\documentclass{amsart}

\usepackage{amsmath,amsthm}
\usepackage {latexsym}
\usepackage{amssymb}

\newcommand{\les}{\lesssim}

\newcommand{\const}{\mbox{\rm const}}
\newcommand{\beeq}{\begin{equation}}
\newcommand{\eneq}{\end{equation}}

\newcommand{\wvec}{{\underline w}}

\newcommand{\bear}{\begin{eqnarray}}
\newcommand{\eear}{\end{eqnarray}}
\newcommand{\beq}{\begin{equation}}
\newcommand{\eeq}{\end{equation}}

\newcommand{\spec}{{\rm spec}}

\newcommand{\half}{\frac{1}{2}}

\newcommand{\eps}{{\varepsilon}}
\newcommand{\R}{{\mathbb R}}
\newcommand{\C}{{\mathbb C}}

\newcommand{\Z}{{\mathbb Z}}

\newcommand{\Compl}{{\mathbb C}}

\newcommand{\calg}{\,{\mathfrak g}}
\newcommand{\calG}{{\mathcal G}}
\newcommand{\calF}{{\mathcal F}}

\newcommand{\calS}{{\mathcal S}}
\newcommand{\calB}{{\mathcal B}}
\newcommand{\calJ}{{\mathcal J}}

\newcommand{\calM}{{\mathcal M}}

\newcommand{\calN}{{\mathcal N}}

\newcommand{\sign}{\mbox{sign}}

\newcommand{\Laplace}{\partial_x^2}

\newcommand{\trip}{|\!|\!|}
\newcommand{\la}{\langle}
\newcommand{\ra}{\rangle}

\newcommand{\vvec}{{\underline v}}
\newcommand{\vece}{{\underline e}}

\def\nn{\nonumber}
\def\calge1{\calg_{\vec{e_1}}}

\def\bm{\left( \begin{array}{cc}}
\def\endm{\end{array}\right)}
\def\ker{{\rm ker}}
\def\Ran{{\rm Ran}}
\def\Hil{{\mathcal H}}
\def\Dom{{\rm Dom}}

\newtheorem{theorem}{Theorem}[section]
\newtheorem{lemma}[theorem]{Lemma}
\newtheorem{defi}[theorem]{Definition}
\newtheorem{cor}[theorem]{Corollary}
\newtheorem{prop}[theorem]{Proposition}
\newtheorem{proposition}[theorem]{Proposition}

\theoremstyle{remark}
\newtheorem{remark}[theorem]{Remark}

\def\gatil{\tilde{\gamma}}

\def\y{{y}}

\def\norm[#1][#2]{\Vert #1 \Vert_{#2}}

\def\vecv{\vvec}

\def\pizer{\pi^{(0)}}

\def\Zzer{Z^{(0)}}

\def\Pim{P_{\rm im}}

\def\trip{|\!|\!|}

\def\Util3{\tilde{U}^{(3)}}
\def\pitil{\tilde{\pi}}
\def\xitil{\tilde{\xi}}
\def\etatil{\tilde{\eta}}

\def\calW{{\mathcal W}}

\def\vecv{\vvec}
\def\vecw{\wvec}
\def\calE{{\mathcal E}}
\numberwithin{equation}{section}

\begin{document}

\title {Stable manifolds for all monic supercritical NLS in one dimension}
\author{J.\ Krieger, W. Schlag}
\thanks{The first author was partially supported by NSF grant DMS-0401177 and the second author by the NSF grant
DMS-0300081 and a Sloan fellowship. The authors wish to thank
Fritz Gesztesy for pointing out reference~\cite{Flug} and helpful
comments on  $\cosh^{-2}(x)$ potentials.}
\address{Harvard University, Dept. of Mathematics, Science Center, 1 Oxford Street, Cambridge, MA 02138, USA}
\email{jkrieger@math.harvard.edu}
\address{253-37 Caltech, Pasadena, CA 91125, USA}
\email{schlag@caltech.edu}

\date{}
\maketitle

\section{Introduction}
\label{sec:intro}

We consider the NLS \beeq \label{eq:NLS}
 i\partial_t \psi + \partial_{x}^2 \psi = -|\psi|^{2\sigma} \psi
\eneq
 on the line $\R$ with $\sigma>2$. This is exactly the $L^2$-supercritical
case and these equations are locally well-posed in
$H^1(\R)=W^{1,2}(\R)$. Let $\phi=\phi(\cdot,\alpha)$ be the ground
state of \beeq \label{eq:ground} -
\phi''+\alpha^2\phi=\phi^{2\sigma+1}. \eneq By this we mean that
$\phi>0$ and $\phi\in C^2(\R)$. It is a classical fact  that such
solutions exist and are unique. In fact, for the case $\alpha=1$
the solutions are
\[ \phi(x,1)=\frac{(\sigma+1)^{\frac{1}{2\sigma}}}{\cosh^{\frac{1}{\sigma}}(\sigma x)},\]
whereas for general $\alpha>0$ they are obtained from this
solution by rescaling:
$\phi(x,\alpha)=\alpha^{\frac{1}{\sigma}}\phi(\alpha x,1)$.

Clearly, the standing wave $\psi=e^{it\alpha^2}\phi$
solves~\eqref{eq:NLS}. We seek an $H^1$-solution $\psi$ of the
form $\psi= W+R$  where
\begin{align}
W(t,x) &= e^{i\theta(t,x)} \phi(x-y(t),\alpha(t)) \label{eq:W}\\
\theta(t,x) &= v(t)x-\int_0^t (v(s)^2-\alpha^2(s))\, ds + \gamma(t) \label{eq:theta} \\
y(t) & = 2\int_0^t v(s)\,ds + D(t) \label{eq:y}
\end{align}
is the usual standing wave with a moving set of parameters
$\pi(t):=(\gamma(t),v(t),D(t),\alpha(t))$, and $R$ is a small
perturbation. Performing a Galilei transform, we may assume that
$W(0,x)=\phi(x,\alpha)$.

\begin{theorem}
\label{thm:main} Fix any $\sigma>2$ in~\eqref{eq:NLS} and any
$\alpha_0>0$.  Let $\Sigma:=\{f\in L^2(\R)\:|\: \trip f\trip <\infty\}$
where\footnote{The weight $\la x\ra$ in the
definition of $\trip\cdot\trip$ can be relaxed to $\la
x\ra^{\half+\eps}$, but we keep it in this form for aesthetic
purposes}
\[ \trip f \trip:= \|f\|_{H^1}+\|\la
x\ra f\|_{L^1\cap L^2}+\|\la x\ra \partial_x f\|_{L^1}
\]
and set
\beeq \label{eq:calB} \calB:=
\Big\{R_0\in L^2(\R)\:|\: \trip R_0 \trip
<\delta\Big\} \eneq
Then there exist a real-linear subspace
$\calS\subset \Sigma$ of co-dimension five and a small $\delta>0$
with the following properties:
there exists a  map $\Phi:\calB\cap \calS\to \Sigma$  such that
\begin{align}
 \trip \Phi(R_0)\trip &\les \trip R_0\trip^2 \qquad \forall R_0\in
 \calB\cap \calS\label{eq:r0til} \\
\trip \Phi(R_0)-\Phi(\tilde{R}_0) \trip &\les  \delta\trip R_0-
\tilde{R}_0\trip \qquad\forall R_0,\tilde{R}_0\in  \calB\cap \calS
\label{eq:lip}
\end{align}
and so that for any $R_0\in~\calB\cap \calS$ the
NLS~\eqref{eq:NLS} has a global $H^1$ solution $\psi(t)$ for
$t\ge0$ with initial condition
$\psi(0)=\phi(\cdot,\alpha_0)+R_0+\Phi(R_0)$. Moreover,
\[ \psi(t)= W(t,\cdot)+R(t) \]
where $W$ as in~\eqref{eq:W} is governed by a path $\pi(t)$ of
parameters which converges to some terminal vector $\pi(\infty)$
such that $\sup_{t\ge0}|\pi(t)-\pi(\infty)|\les \delta^2$ and so
that \beeq \label{eq:Rdec}
 \|R(t)\|_{H^1} \les \delta,\quad \|R(t)\|_\infty \les \delta \la t\ra^{-\frac12},
 \quad \|\la x-y(t)\ra^{-\half-\eps}R(t)\|_\infty
\les \delta \la t\ra^{-1-\eps} \eneq for all $t>0$ and some
$\eps>0$. The solution $\psi(t)$ is unique amongst all solutions
with these initial data and satisfying the above decay assumptions
as well as certain orthogonality relations and decay
assumptions on the path (which will be specified later). Finally, there is
scattering:
\[ R(t)= e^{it\Laplace} f_0 + o_{L^2}(1) \text{\ \ as $t\to\infty$}\]
for some $f_0\in L^2(\R)$.
\end{theorem}

It is  well-known that the {\em supercritical
equation}~\eqref{eq:NLS} is
 {\em orbitally unstable},
see Berestycki and Cazenave~\cite{BerCaz}. This is in contrast to
the {\em orbital stability} of the {\em subcritical equations}
that was proved by Cazenave and Lions~\cite{CazLio} and
Weinstein~\cite{Wei1}, \cite{Wei2}. In fact, the instability
result of Berestycki and Cazenave \cite{BerCaz} shows that one can
have finite time blow-up for initial data
$\psi_0=\phi(\cdot,\alpha)+R_0$ where $R_0$ can be made
arbitrarily small in any reasonable norm. Theorem~\ref{thm:main}
states that we do have the same asymptotic stability and
scattering as in the subcritical case, provided we choose our
initial data on a suitable submanifold.

To understand the origin of~$\calS$, we associate with each
$\alpha>0$  the matrix operator \beeq \label{eq:Hilintro}
\Hil = \Hil(\alpha)=\bm -\partial^2_{x} + \alpha^2 -(\sigma+1)\phi^{2\sigma}(\cdot,\alpha) & -\sigma\phi^{2\sigma}(\cdot,\alpha) \\
\sigma\phi^{2\sigma}(\cdot,\alpha) & \partial_x^2 - \alpha^2 +
(\sigma+1)\phi^{2\sigma}(\cdot,\alpha)
\endm.
\eneq This operator arises by linearizing the NLS \eqref{eq:NLS}
around a standing wave. It is closed on the domain
$W^{2,2}(\R)\times W^{2,2}(\R)$ and its spectrum has the following
form: It is located on $\R\cup i\R$, with essential spectrum equal
to $(-\infty,-\alpha^2]\cup [\alpha^2,\infty)$. The discrete
spectrum equals $\{0,\pm i\gamma(\alpha)\}$, where
$\gamma(\alpha)=\alpha^2\gamma(1)>0$. Both $\pm i\gamma(\alpha)$
are simple eigenvalues with exponentially decaying eigenfunctions,
whereas $0$ is an eigenvalue of geometric multiplicity two and
algebraic multiplicity four (the latter fact goes back to
Weinstein~\cite{Wei1}).

\noindent Next, we introduce the Riesz projection
$P^+_{u}(\alpha)$  such that
\[ \spec(\Hil(\alpha) P_u^+(\alpha))= \{0\}\cup\{i\gamma(\alpha)\}.\]
The notation $P_u^+$ is meant to indicate the unstable modes as
$t\to+\infty$. The real-linear, finite-codimensional subspace
$\calS$ above is precisely the set of $R_0\in \Sigma$ so that\footnote{We will show
below that $P_u^+:\tilde\Sigma\to\tilde\Sigma$ where $\tilde\Sigma=\big\{\binom{f}{\bar f}\::\: f\in\Sigma\big\}$.}
\beeq
\label{eq:orth}
 P_u^+(\alpha_0) \binom{R_0}{\bar{R}_0}=0.
\eneq
The codimension of $\calS$ is  the number of unstable (or
non-decaying) modes of the linearization as $t\to\infty$: four in
the root space and one exponentially unstable mode. The stable
manifold $\calM$ is  the surface described by the parameterization
$R_0\mapsto R_0+\Phi(R_0)$ where $R_0$ belongs to a small ball
$\calB\cap\calS$ inside of~$\calS$. The
inequality~\eqref{eq:r0til} means that $\calS$ is the tangent
space to~$\calM$ at zero, whereas~\eqref{eq:lip} expresses
that~$\calM$  is given in terms of a  Lipschitz parameterization.
It is easy to see that it also the graph of a Lipschitz map
$\tilde{\Phi}: \calS\cap\calB\to\Sigma$. Indeed, define
$\tilde{\Phi}$ as
\[ R_0+P_{\calS}\Phi(R_0) \mapsto R_0 + \Phi(R_0), \]
where $P_\calS$ is the projection onto $\calS$ which is induced by
the Riesz-projection $I-P_u^+(\alpha_0)$ (the latter operates on
$L^2\times L^2$, whereas we need only the first coordinate of this
projection, see Remark~\ref{rem:Jinv} below for the details of
this). The left-hand side is clearly in $\calS$. Moreover, to see
that this map is well-defined as well as Lipschitz, note that
\eqref{eq:lip} implies that
\begin{align*}
& (1-C\delta)\trip R_1-R_0\trip
\le \trip R_1-R_0 + \Phi(R_1)-\Phi(R_0) \trip \le (1+C\delta) \trip R_1-R_0\trip \\
& (1-C\delta)\trip R_1-R_0\trip \le \trip R_1-R_0 +
P_{\calS}\Phi(R_1)-P_{\calS}\Phi(R_0) \trip \le (1+C\delta) \trip
R_1-R_0\trip.
\end{align*}

Since the root-space of $\Hil(\alpha)$ at zero does not destroy
asymptotic stability or instability in the subcritical case, one
would expect that the new unstable phenomena in the supercritical
case should only be connected with the imaginary eigenvalues of
$\Hil(\alpha)$. More precisely, since we are considering
$t\to+\infty$, they should result exclusively from the eigenvalue $i\gamma$, where $\gamma>0$.
Hence,
the true codimension of our stable manifold should be one. In the
following theorem we obtain such a stable hypersurface, by
applying the three-parameter family of Galilei transforms together
with scaling to the manifold $\calM$ from Theorem~\ref{thm:main}.
Since this family acts transversally to $\calM$, we recover four
of the missing dimensions this way.

\begin{theorem}
\label{thm:main2} Fix any $\alpha_0>0$. Then there exist a small
$\delta>0$ and a Lipschitz manifold $\calN$ inside the space
$\Sigma$ of size\footnote{This means that $\calN$ is the graph of
a Lipschitz map $\Psi$ with domain $\calB\cap\tilde\calS$ where
$\tilde\calS$ is a subspace of codimension one and with $\calB$ as
in \eqref{eq:calB}.}  $\delta$ and codimension one so
that $\phi(\cdot,\alpha_0)\in\calN$ with the following property:
for any choice of initial data $\psi(0)\in\calN$ the
NLS~\eqref{eq:NLS} has a global $H^1$ solution $\psi(t)$ for
$t\ge0$.  Moreover,
\[ \psi(t)= W(t,\cdot)+R(t) \]
where $W$ as in~\eqref{eq:W} is governed by a path $\pi(t)$ of
parameters so that $|\pi(0)-(0,0,0,\alpha_0)|\les\delta$ and which
converges to some terminal vector $\pi(\infty)$ such that
$\sup_{t\ge0}|\pi(t)-\pi(\infty)|\les \delta^2$. The solution is
unique under the same conditions as in the preceding theorem.
Finally, \eqref{eq:Rdec} holds and there is scattering:
\[ R(t)= e^{it\Laplace} f_0 + o_{L^2}(1) \text{\ \ as $t\to\infty$ }\]
for some $f_0\in L^2(\R)$.
\end{theorem}

\noindent   This result raises the interesting question of
deciding the behavior of solutions with initial data $\phi(0)\in
\calB\setminus \calN$. It is known that at least for some choices
of such initial data the solution blows up, but the authors do not
know what to expect in general. For related results, see the book
by Li and Wiggins~\cite{LW},  the paper by Tsai and Yau~\cite{TY},
as well as the references there.

Theorems~\ref{thm:main} and~\ref{thm:main2} (but with unweighted norms) were first proved in
three dimensions in~\cite{Sch} for the cubic focusing NLS, albeit
under the assumption that there are no imbedded eigenvalues in the
essential spectrum of the linearized operators. This paper is
closely related to~\cite{Sch}, although it does differ in several important
aspects. First, all spectral properties are proved here analytically, whereas \cite{Sch}
required verifying some spectral properties of the well-known pair of Schr\"odinger
operators $L_-$ and $L_+$ numerically,  and the
absence of imbedded eigenvalues for the systems remained as an assumption.
Second, the (free) dispersive decay in one dimension is $t^{-\frac12}$, which is
not integrable at infinity. This forces us to use an improved decay estimate
which takes the form
\[ \|\la x\ra^{-1}e^{it(-\partial_x^2+V)}P_c\, f\|_{L^\infty(\R)}\le Ct^{-\frac32}\|\la x\ra f\|_{L^1(\R)} \]
for scalar operators $-\partial_x^2+V$ {\em with no resonance at zero}, see~\cite{Sch2}.
Weaker forms of this result on $L^2(\R)$ and with higher weights were known, see Murata~\cite{Mur}
and \cite{BP1}.   However, due to the fact that
we cover here the {\em full supercritical range}, we need to rely on this estimate (or interpolates of it)
rather than any weaker version, albeit for non-selfadjoint matrix operators rather than
scalar operators.
The transition from the scalar case (as in~\cite{Sch2}, say) to the matrix case
requires a more sophisticated functional framework,
which involves developing the scattering theory (i.e., Jost solutions)
of these matrix Schr\"odinger operators as in~\cite{BP1}.
Since the estimates we require here are considerably sharper than the ones in
Buslaev-Perelman, we carefully develop this framework, together
with the necessary spectral theory in the second half of the paper.
This part is of independent interest.
As far as the nonlinear argument is concerned,
the weights inside of our norms will force us to depart
from the contraction procedure employed in ~\cite{Sch}, which takes
place in a fixed Banach space, and use a method of iteration
which adjusts the Banach space depending on the iterate. This issue here is that the weights need to be
centered around a path $y_j(t)$. It is impossible to compare these norms for large $t$
because these paths will then be separated by a distance exceeding one.
This forces us to truncate at a time $T_j$ which grows with~$j$.
This method, which is more involved then
the contraction from~\cite{Sch}, is of independent interest.
\newline {\bf{Notation:}} whenever we use the symbol $\lesssim$, certain
universal multiplicative constants are implied which do not depend on any varying
parameters appearing in the proof.
\section{The linearization, Galilei transforms, and $\calJ$-invariance}
\label{sec:ansatz}

We begin by linearizing around the standing wave.

\begin{lemma}
\label{lem:Zeq}
 Let $0<T\le \infty$.
Assume that
$\pi(t)=(\gamma(t),{v}(t),D(t),\alpha(t)):[0,T)\to\R^4$ belongs to
$C^1([0,T),\R^4)$, and let $W=W(t,x)$ be as in~\eqref{eq:W}. Then
$\psi\in C([0,T),H^1(\R))\cap C^1([0,T),H^{-1}(\R))$
solves~\eqref{eq:NLS} with $\psi=W+R$ iff $Z=\binom{R}{\bar{R}}$
solves the equation
\begin{align}
& i\partial_t Z + \bm \Laplace + (\sigma+1)|W|^{2\sigma} & \sigma|W|^{2(\sigma-1)}W^2 \\
-\sigma|W|^{2(\sigma-1)}\bar{W}^2  & -\Laplace - (\sigma+1)|W|^{2\sigma} \endm Z \nn \\
 & = \dot{v} \binom{-xe^{i\theta} \phi(\cdot-y,\alpha)}{xe^{-i\theta} \phi(\cdot-y,\alpha)} +
\dot{\gamma} \binom{-e^{i\theta}
\phi(\cdot-y,\alpha)}{e^{-i\theta} \phi(\cdot-y,\alpha)}
 + i\dot{\alpha}
\binom{e^{i\theta}\partial_\alpha
\phi(\cdot-y,\alpha)}{e^{-i\theta}\partial_\alpha
\phi(\cdot-y,\alpha)}
 + i\dot{D}
\binom{-e^{i\theta}\partial_x\phi(\cdot-y,\alpha)}{-e^{-i\theta}\partial_x \phi(\cdot-y,\alpha)} \nn \\
&\quad +
\binom{-|R+W|^{2\sigma}(R+W)+|W|^{2\sigma}W+(\sigma+1)|W|^{2\sigma}R
+
\sigma|W|^{2(\sigma-1)}W^2\bar{R}}{|R+W|^{2\sigma}(\bar{R}+\bar{W})-|W|^{2\sigma}\bar{W}-(\sigma+1)|W|^{2\sigma}\bar{R}
- \sigma|W|^{2(\sigma-1)}\bar{W}^2{R}}
\label{eq:Zsys}\\
&:= i\dot{\pi}\partial_\pi W+N(Z,W) \nn
\end{align}
in the sense of $C([0,T),H^1(\R)\times H^1(\R))\cap
C^1([0,T),H^{-1}(\R)\times H^{-1}(\R))$. Here $y$ and $\theta$ are
the functions from~\eqref{eq:y} and~\eqref{eq:theta}, and
$\alpha=\alpha(t)$. For future reference, we denote the matrix
operator on the left-hand side of~\eqref{eq:Zsys} by
$-\Hil(\pi(t))$, i.e., \beeq \label{eq:Hpidef}
\Hil(\pi(t)) := \bm - \Laplace - (\sigma+1)|W|^{2\sigma} & -\sigma|W|^{2(\sigma-1)}W^2 \\
\sigma|W|^{2(\sigma-1)}\bar{W}^2 & \Laplace +
(\sigma+1)|W|^{2\sigma} \endm. \eneq The nonlinear term
in~\eqref{eq:Zsys}, which we denote by
$N(Z,W)=\binom{N_1(Z,W)}{N_2(Z,W)}$, is quadratic in $R,\bar{R}$.
This means that
\[ N(0,W)=\partial_{R}N(0,W)=\partial_{\bar{R}}N(0,W)=0.\]
\end{lemma}
\begin{proof}
It will be convenient to consider the more general NLS \beeq
\label{eq:NLS'}
 i\partial_t \psi + \partial_{x}^2 \psi = -\beta(|\psi|^2)\psi.
\eneq Let $\phi=\phi(\cdot,\alpha(t))$ for ease of notation.
Direct differentiation shows that $W(t,x)$ satisfies
\[ i\partial_tW+\Laplace W = -\beta(|W|^2)W -W(\dot{v}x+\dot{\gamma})-ie^{i\theta}\partial_x \phi\cdot \dot{D}+
ie^{i\theta}\dot{\alpha}\partial_\alpha \phi.\] Hence $W+R$ is a
solution of \eqref{eq:NLS'} iff
\[ i\partial_t R + \Laplace R = -\beta(|W+R|^2)(W+R)+\beta(|W|^2)W
-e^{i\theta}\phi(\dot{v}x+\dot{\gamma}) - ie^{i\theta}\partial_x
\phi\cdot \dot{D} + ie^{i\theta}\dot{\alpha}\partial_\alpha \phi.
\] Joining this equation with its conjugate leads to the system
\begin{align}
& i\partial_t Z + \bm \Laplace + \beta'(|W|^2)|W|^2+\beta(|W|^2) & \beta'(|W|^2)W^2 \\
-\beta'(|W|^2)W^2  & -\Laplace -\beta'(|W|^2)|W|^2-\beta(|W|^2)  \endm Z \nn \\
 & = \dot{v} \binom{-xe^{i\theta} \phi(\cdot-y,\alpha)}{xe^{-i\theta} \phi(\cdot-y,\alpha)} +
\dot{\gamma} \binom{-e^{i\theta}
\phi(\cdot-y,\alpha)}{e^{-i\theta} \phi(\cdot-y,\alpha)}
 + i\dot{\alpha}
\binom{e^{i\theta}\partial_\alpha
\phi(\cdot-y,\alpha)}{e^{-i\theta}\partial_\alpha
\phi(\cdot-y,\alpha)}
 + i\dot{D}
\binom{-e^{i\theta}\partial_x \phi(\cdot-y,\alpha)}{-e^{-i\theta}\partial_x \phi(\cdot-y,\alpha)} \nn \\
&\quad +
\binom{-\beta(|R+W|^2)(R+W)+\beta(|W|^2)W+[\beta'(|W|^2)|W|^2+\beta(|W|^2)]R
+
\beta'(|W|^2)W^2\bar{R}}{\beta(|R+W|^2)(\bar{R}+\bar{W})-\beta(|W|^2)\bar{W}-[\beta'(|W|^2)|W|^2+\beta(|W|^2)]\bar{R}
- \beta'(|W|^2)\bar{W}^2{R}}.\label{eq:Zsys'}
\end{align}
Conversely, if $Z(0)$ is of the form
\[ Z(0) = \binom{Z_1(0)}{\overline{Z_1(0)}}, \]
and $Z(t)$ solves~\eqref{eq:Zsys'}, then $Z(t)$ remains of this
form for all times. This is simply the statement that the
system~\eqref{eq:Zsys'} is invariant under the transformation
\beeq \label{eq:Jdef} \calJ:f\mapsto \overline{Jf}, \qquad J=\bm
0&1 \\1&0 \endm, \quad f=\binom{f_1}{\overline{f_1}}, \eneq which
can be checked by direct verification. This fact allows us to go
back from the system to the scalar equation~\eqref{eq:NLS'}.
Finally, it is easy to see that the nonlinear term
in~\eqref{eq:Zsys'} is quadratic, and that~\eqref{eq:Zsys'}
reduces to~\eqref{eq:Zsys} if  $\beta(u)= u^\sigma$ for all $u>0$.
\end{proof}
\begin{remark} Notice that we have the equality
\begin{equation}\nonumber\begin{split}
& \dot{v} \binom{-xe^{i\theta} \phi(\cdot-y,\alpha)}{xe^{-i\theta}
\phi(\cdot-y,\alpha)} + \dot{\gamma} \binom{-e^{i\theta}
\phi(\cdot-y,\alpha)}{e^{-i\theta} \phi(\cdot-y,\alpha)}
 + i\dot{\alpha}
\binom{e^{i\theta}\partial_\alpha
\phi(\cdot-y,\alpha)}{e^{-i\theta}\partial_\alpha
\phi(\cdot-y,\alpha)}
 + i\dot{D}
\binom{-e^{i\theta}\partial_x\phi(\cdot-y,\alpha)}{-e^{-i\theta}\partial_x \phi(\cdot-y,\alpha)} \nn \\
&=\dot{v} \binom{-(x-y)e^{i\theta}
\phi(\cdot-y,\alpha)}{(x-y)e^{-i\theta} \phi(\cdot-y,\alpha)} +
[\dot{\gamma}+\dot{v}y] \binom{-e^{i\theta}
\phi(\cdot-y,\alpha)}{e^{-i\theta} \phi(\cdot-y,\alpha)}
 + i\dot{\alpha}
\binom{e^{i\theta}\partial_\alpha
\phi(\cdot-y,\alpha)}{e^{-i\theta}\partial_\alpha
\phi(\cdot-y,\alpha)}
 \\&\hspace{11cm}+ i\dot{D}
\binom{-e^{i\theta}\partial_x\phi(\cdot-y,\alpha)}{-e^{-i\theta}\partial_x \phi(\cdot-y,\alpha)} \nn \\
\end{split}\end{equation}
Abusing notation, we shall later refer to this expression as
$i\dot{\tilde{\pi}}\partial_{\pi}W$, where
\begin{equation}\nonumber
\dot{\tilde{\pi}}:=(\dot{\gamma}+\dot{v}y,\dot{v},\dot{D},\dot{\alpha})
\end{equation}
The quantity $\dot{\tilde{\pi}}$ shall play an important role in
the argument to follow.
\end{remark}

The $\calJ$-invariant vectors in $L^2(\R)\times L^2(\R)$ form a
real-linear subspace, namely
\[ \Big\{ \binom{f}{\bar{f}}\:\Big|\: f\in L^2(\R)\Big\}. \]
Writing $f=f_1+if_2$ it can be seen to be isomorphic to the
subspace
\[ \Big\{ \binom{f_1}{f_2}\:\Big|\: f_1,f_2\in L^2(\R), \; f_1,f_2 \text{\ are real-valued}\Big\}, \]
which is clearly linear, but only over $\R$. Throughout the paper,
we need to insure that all vectorial solutions we construct belong
to this subspace. Only then is it possible to revert back to the
scalar NLS~\eqref{eq:NLS}.

To perform estimates, one needs to transform~\eqref{eq:Zsys} to a
resting frame. This requires certain properties of the path
$\pi(t)$. In the following definition, $\eps>0$ and $\delta>0$ are
fixed small constants.

\begin{defi}
\label{def:adm} Let $T>0$, to be thought of as a large
number\footnote{The logic here is that we are working on a
large time interval $[0,T]$, and are not concerned about what happens for
times larger than $T$. In the iteration, we will adjust this
parameter to the stage of the iteration we are at. Of course, $T\to\infty$ in this process} . We say that a
path $\pi:[0,T]\to\R^4$ with
$\pi(t):=(\gamma(t),{v}(t),D(t),\alpha(t))$ is {\em admissible}
provided it belongs to $C^1([0,T],\R^4)$, and the estimate
$|\dot{\pi}(t)|\le \delta^2\la t\ra^{-2-\eps}$ holds. Define a
constant parameter vector $\pi_T=(\gamma_T,{v}_T,D_T,\alpha_T)$ as
\begin{align}
\gamma_T &:= \gamma(T) + 2\int_0^T \int_t^T ({v}(s)\dot{{v}}(s)
-\alpha(s)\dot{\alpha}(s))\,ds\,dt \label{eq:gainf}\\
{v}_T &:= {v}(T)\label{eq:vinf}\\
D_T &:= D(T) -2\int_0^T \int_t^T \dot{{v}}(s)\, ds\,dt \label{eq:Dinf}\\
\alpha_T &:= \alpha(T) \label{eq:alinf}
\end{align}
\end{defi}

\noindent The logic behind these asymptotic parameters can be seen
in the following lemma.

\begin{lemma}
\label{lem:rho_T} Suppose $\pi$ is an admissible path and let
$\theta,y$ and $\theta_T,y_T$ be as in \eqref{eq:theta},
\eqref{eq:y}. Furthermore, define \beeq \label{eq:inf_path} y_T(t)
:= 2t{v}_T + D_T,\; \theta_T(t,x) := {v}_T x -
t({v}_T^2-\alpha_T^2)+\gamma_T \eneq and \beeq \label{eq:rhoinf}
 \rho_T(t,x) := \theta(t,x+y_T)-\theta_T(t,x+y_T).
\eneq Then
\[ |\rho_T(t,x)| \le C_\eps\, \delta^2(1+|x|)\la t\ra^{-\eps}, \quad |y(t)-y_T(t)|\le C_\eps\, \delta^2 \la t\ra^{-\eps},\]
as well as \beeq \label{eq:dott}
 |\dot{\rho}_T(t,x)| \le C_\eps\, \delta^2 (1+|x|)\la t\ra^{-1-\eps}, \quad |\dot{y}(t)-\dot{y}_T(t)|\le C_\eps\, \delta^2\la t
\ra^{-1-\eps} \eneq for all $t$ with $T\geq t\ge0$. The constants
here only depend on $\eps$.
\end{lemma}
\begin{proof}
First, \beeq \label{eq:ti} \theta_T(t,x+y_T)=v_T(x+2tv_T+D_T) +
t(-v_T^2+\alpha_T^2)+\gamma_T =t(-v_T^2+\alpha_T^2)+
v_T(x+D_T)+\gamma_T. \eneq In view of the definition of $\pi_T$,
\begin{align}
\theta(t,x+y_T)-\theta_T(t,x+y_T) &= v(t)(x+2tv_T+D_T) -
\int_0^t (v(s)^2-\alpha^2(s))\,ds + \gamma(t) \nn \\
&\quad  - v_T(x+2tv_T+D_T) + t(v_T^2-\alpha_T^2)-\gamma_T \nn \\
& = ({v}(t)-{v}_T)(x+2t{v}_T+D_T) +2\int_0^T \int_s^T ({v}\dot{{v}}-\alpha\dot{\alpha})(s')\,ds'\, ds \nn \\
&\quad -\gamma_T + \gamma(t) -2\int_t^T\int_s^T ({v}\dot{{v}}-\alpha\dot{\alpha})(s')\,ds'\, ds \nn\\
& =  ({v}(t)-{v}_T)(x+2t{v}_T+D_T) -2\int_t^T\int_s^T ({v}\dot{{v}}-\alpha\dot{\alpha})(s')\,ds'\, ds  \nn \\
&\quad -\gamma(T) + \gamma(t). \label{eq:rhorep}
\end{align}
Using Definition~\ref{def:adm}  implies the desired bound on
$\rho_T$. As for $y(t)-y_T(t)$, the definition of $D_T$ implies
that \beeq \label{eq:yrep} y_T(t)-y(t) = 2tv_T+D_T-2\int_0^t
v(s)\,ds - D(t)=D(T)-D(t)-2\int_t^T\int_s^T \dot{v}(s')\,ds'\,ds,
\eneq which implies the stated estimate for $0\leq t\leq T$.
\end{proof}
\begin{remark} One has similar definitions and properties for
$\pi_{\infty}=(\gamma_\infty,v_\infty,D_\infty,\alpha_\infty)$
etc. One then replaces $[0,T]$ by $[0,\infty)$.
\end{remark}

With $\pi_T=(\gamma_T,{v}_T,D_T,\alpha_T)$ a constant vector,
define the usual {\em Galilei transform} to be \beeq
\label{eq:gal}
 \calg_T(t) = e^{i(\gamma_T+{v}_T x-t|{v}_T|^2)}
\, e^{-i(2t{v}_T+D_T)p}, \eneq where $p:=-i\frac{d}{dx}$. The
action of $\calg_T(t)$ on functions is
\[ (\calg_T(t)f)(x)= e^{i(\gamma_T+{v}_T x-t{v}_T^2)} f(x-2t{v}_T-D_T).\]
It is unitary on $L^2$, and isometric on all $L^p$, and the
commutation property
$e^{it\partial_x^2}\calg_T(0)=\calg_T(t)e^{it\partial_x^2}$ holds.
The inverse of $\calg_T(t)$ is \beeq \label{eq:gal_inv}
\calg_T(t)^{-1} = e^{i(2t{v}_T+D_T)p}\,e^{-i(\gamma_T+{v}_T
x-t{v}_T^2)}= e^{-i(\gamma_T+{v}_T D_T+{v}_T x+t{v}_T^2)}\,
e^{i(2t{v}_T+D_T)p}. \eneq Moreover, the Galilei
transform~\eqref{eq:gal} generates a four-parameter family of
standing waves: Let $\phi(\cdot,\alpha_T)$ be the ground state
of~\eqref{eq:ground} with $\alpha=\alpha_T$. Then \beeq
\label{eq:Winfty}
 W_T(t,\cdot) = \calg_T(t) [e^{it\alpha_T^2}\phi(\cdot,\alpha_T)]
\eneq solves \eqref{eq:NLS}, where $W_T$ is a standing wave as
introduced in~\eqref{eq:W} but with the constant parameter path
$\pi_T$. This can also be written as
\[ W_T(t,x) = e^{i\theta_T(t,x)}\phi(x-y_T(t),\alpha_T), \]
where $y_T, \theta_T$ are as in Lemma~\ref{lem:rho_T}. As usual,
we transform~\eqref{eq:Zsys} to a stationary frame by means of
Galilei transforms. In addition, a modulation will be performed.
The details are as follows.

\begin{lemma} \label{lem:UPDE}
Let $\pi(t)$ and $\pi_T$ be as in Definition~\ref{def:adm}. Given
$Z=\binom{Z_1}{Z_2}$, introduce $U$, as well as $M_T(t), \calG_T
(t)$ as \beeq \label{eq:Udef}
 U(t)=\bm e^{i\omega_T(t)} & 0 \\ 0 & e^{-i\omega_T(t)} \endm \binom{\calg_T(t)^{-1} Z_1(t)}
{\overline{\calg_T(t)^{-1} \overline{Z_2(t)}}} = M_T(t) \calG_T
(t) Z(t), \eneq where $\omega_T(t)=-t\alpha_T^2$. Then $Z(t)$
solves \eqref{eq:Zsys} in the $H^1$ sense iff
 $U=\binom{U_1}{U_2}$ as in~\eqref{eq:Udef} satisfies the following PDE in the $H^1$ sense
(with $\phi_T=\phi(\cdot,\alpha_T)$): \beeq \label{eq:UPDE}
i\dot{U}(t) + \bm \Laplace +(\sigma+1)\phi_T^{2\sigma}-\alpha_T^2
& \sigma\phi_T^{2\sigma} \\ -\sigma\phi_T^{2\sigma} & -\Laplace
-(\sigma+1)\phi_T^{2\sigma} +\alpha_T^2 \endm U =
i\dot{\pi}\partial_\pi \tilde W_T(\pi) +N_T(U,\pi) + V_T U \eneq
where
\begin{align}
 V_T = V_T(t,x) &:= \bm (\sigma+1)(\phi_T^{2\sigma}(x)-\phi^{2\sigma}(x+y_T-y)) & \sigma(\phi_T^{2\sigma}(x)-e^{2i\rho_T} \phi^{2\sigma}(x+y_T-y)) \\
-\sigma(\phi_T^{2\sigma}(x)- e^{-2i\rho_T} \phi^{2\sigma}(x+y_T-y)) & -(\sigma+1)(\phi_T^{2\sigma}(x)-\phi^{2\sigma}(x+y_T-y)) \endm \label{eq:V}\\
 i\dot{\pi}\partial_\pi \tilde W_T(\pi) &:= \dot{{v}} \binom{-(x+y_T)e^{i\rho_T} \phi(x+y_T-y)}{(x+y_T) e^{-i\rho_T} \phi(x+y_T-y)} + \dot{\gamma} \binom{-e^{i\rho_T}\phi(x+y_T-y)}{e^{-i\rho_T}\phi(x+y_T-y)} \label{eq:vdotetc} \\
 & \quad +i\dot{\alpha} \binom{e^{i\rho_T}\partial_\alpha\phi(x+y_T-y)}{e^{-i\rho_T}\partial_\alpha\phi(x+y_T-y)}+i\dot{D}
\binom{-e^{i\rho_T}\partial_x\phi(x+y_T-y)}{-e^{-i\rho_T}\partial_x\phi(x+y_T-y)} \nn \\
N_T(U,\pi) &:= \binom{N_{1T}(U,\pi)}{N_{2T}(U,\pi)} =
\binom{N_{1T}(U,\pi)}{-\overline{N_{1T}(U,\pi)}} \label{eq:NUpi}
\end{align}
and
\begin{align*}
N_{1T}(U,\pi) &=
-|U_1+e^{i\rho_T}\phi(x+y_T-y)|^{2\sigma}(U_1+e^{i\rho_T}\phi(x+y_T-y))+\phi(x+y_T-y)^{2\sigma+1}e^{i\rho_T}\\
& +(\sigma+1)\phi(x+y_T-y)^{2\sigma}U_1 +
\sigma\phi(x+y_T-y)^{2\sigma} e^{2i\rho_T} U_2.
\end{align*}
Here $\rho_T=\rho_T(t,x)$ is as in Lemma~\ref{lem:rho_T} and
$\phi(x+y_T-y)=\phi(x+y_T(t)-y(t),\alpha(t))$. Finally, $Z$ is
$\calJ$-invariant iff $U$ is $\calJ$-invariant, and $U$ is
$\calJ$-invariant iff $U(0)$ is $\calJ$-invariant.
\end{lemma}
\begin{proof}
Throughout this proof we will adhere to the convention that
$\phi=\phi(\cdot,\alpha(t))$ whereas
$\phi_T=\phi(\cdot,\alpha_T)$). Write the equation~\eqref{eq:Zsys}
for $Z$ in the form \beeq \label{eq:Zeq}
 i\partial_t Z -\Hil_T Z = F + (\Hil(\pi(t))-\Hil_T)Z
\eneq where \beeq \label{eq:Hinfty}
 \Hil_T = \bm -\Laplace - (\sigma+1)|W_T|^{2\sigma} & -\sigma |W|^{2(\sigma-1)}W_T^2 \\ \sigma|W|^{2(\sigma-1)}\bar{W}_T^2 & \Laplace +(\sigma+1)|W_T|^{2\sigma} \endm,
\eneq see \eqref{eq:Winfty} and~\eqref{eq:inf_path}. With
$\calG_T(t)$ defined as in~\eqref{eq:Udef}, and with
$p=-i\partial_x$,
\begin{align}
i\frac{d}{dt} \calG_T(t)f &= \binom{i\dot{\calg}_T(t)^{-1}
f_1}{i\overline{\dot{\calg}_T(t)^{-1}
\overline{f_2}}} = \binom{-(2{v}_T p + |{v}_T|^2)\calg_T(t)^{-1} f_1}{-(2{v}_T p - |{v}_T|^2)\overline{\calg_T(t)^{-1} \overline{f_2(t)}}} \nn \\
& = \bm -(2{v}_T p + |{v}_T|^2) & 0 \\
0 & -(2{v}_T p - |{v}_T|^2) \endm \calG_T(t) f \label{eq:gal_der}
\end{align}
for any $f=\binom{f_1}{f_2}$. Furthermore,
\begin{align}
& M_T(t)\calG_T(t) \Hil_T \binom{f_1}{f_2} \label{eq:comm_GH} \\
 &= \bm e^{i\omega_T(t)} & 0 \\ 0 & e^{-i\omega_T(t)} \endm
\binom{-(\Laplace +(\sigma+1)\phi_T^{2\sigma})\calg_T(t)^{-1}f_1 -
\sigma\phi_T^{2\sigma}e^{2i\theta_T(t,x+y_T)}\calg_T(t)^{-1}
f_2}{\phi_T^{2\sigma}e^{-2i\theta_T(t,x+y_T)}
\overline{\calg_T(t)^{-1}\overline{f_1}} + (\Laplace
+(\sigma+1)\phi_T^{2\sigma})\overline{\calg_T(t)^{-1}\overline{f_2}}}
\nn \\
& \quad - \bm e^{i\omega_T(t)} & 0 \\ 0 & e^{-i\omega_T(t)} \endm
\bm -|{v}_T|^2+2i{v}_T \partial_x & 0 \\ 0 & |{v}_T|^2+2i{v}_T \partial_x \endm \calG_T(t) \binom{f_1}{f_2} \nn \\
&= \bm e^{i\omega_T(t)} & 0 \\ 0 & e^{-i\omega_T(t)} \endm
\binom{-(\Laplace +(\sigma+1)\phi_T^{2\sigma})\calg_T(t)^{-1}f_1 -
\sigma\phi_T^{2\sigma} e^{2i[\theta_T(t,x+y_T)-(t|{v}_T|^2+{v}_T(
x+D_T)+\gamma_T)]}\overline{\calg_T(t)^{-1}\overline{f_2}}}{\sigma\phi_T^{2\sigma}e^{-2i[\theta_T(t,x+y_T)-
(t|{v}_T|^2+{v}_T( x+D_T)+\gamma_T)]}\calg_T(t)^{-1}f_1 +
(\Laplace
+(\sigma+1)\phi_T^{2\sigma})\overline{\calg_T(t)^{-1}\overline{f_2}}}
\nn \\
& \quad - \bm e^{i\omega_T(t)} & 0 \\ 0 & e^{-i\omega_T(t)} \endm
\bm -|{v}_T|^2+2i{v}_T\partial_x & 0 \\ 0 &
|{v}_T|^2+2i{v}_T\partial_x
\endm \calG_T(t) \binom{f_1}{f_2}. \nn
\end{align}
Now
\[ \theta_T(t,x+y_T)-(t|v_T|^2+v_T( x+D_T)+\gamma_T) =
t\alpha_T^2,\] see \eqref{eq:ti}. Hence, by the definition of
$\omega(t)$ (and dropping the argument $t$ from $M_T$ and
$\calG_T$ for simplicity),
\begin{align}
\eqref{eq:comm_GH} &= \bm e^{i\omega_T(t)} & 0 \\ 0 &
e^{-i\omega_T(t)} \endm \bm -(\Laplace
+(\sigma+1)\phi_T^{2\sigma}) & -
\sigma\phi_T^{2\sigma} e^{2it\alpha_T^2} \\
\sigma\phi_T^{2\sigma} e^{-2it\alpha_T^2} & \Laplace
+(\sigma+1)\phi_T^{2\sigma} \endm \calG_T
f \nn \\
& \quad - \bm e^{i\omega_T(t)} & 0 \\ 0 & e^{-i\omega_T(t)} \endm
\bm -|{v}_T|^2-2{v}_T p & 0 \\ 0 & |{v}_T|^2-2{v}_T p \endm \calG_T f  \nn \\
&  =  \bm -\Laplace - (\sigma+1)\phi_T^{2\sigma} & -\sigma\phi_T^{2\sigma} \\
\sigma\phi_T^{2\sigma}  & \Laplace +(\sigma+1)\phi_T^{2\sigma}
\endm M_{T}\calG_T f - \bm -|{v}_T|^2-2{v}_T p & 0 \\ 0
& |{v}_T|^2-2{v}_T p \endm M_{T}\calG_T f. \label{eq:H_phi}
\end{align}
Denote the first matrix operator in~\eqref{eq:H_phi} by
$\Hil_\phi$. Hence, in combination with \eqref{eq:gal_der} one
concludes from~\eqref{eq:Zeq} that
\begin{align*}
i\dot{U} &= i\dot{M}_T\calG_T Z + iM_T\dot{\calG}_T Z + M_T\calG_T
\Hil_T +
M_T\calG_T(F+(\Hil(\pi(t))-\Hil_T)Z) \nn \\
&= \bm -\dot{\omega}_T & 0 \\ 0 & \dot{\omega}_T \endm M_T\calG_T
Z +
\bm -(2{v}_T p + |{v}_T|^2) & 0 \\
0 & -(2{v}_T p - |{v}_T|^2) \endm M_T\calG_T Z + \Hil_\phi M_T\calG_T Z \nn\\
&\quad + \bm |{v}_T|^2+2{v}_T p & 0 \\ 0 & -|{v}_T|^2+2{v}_T p \endm M_T\calG_T Z + M_T\calG_T(F+(\Hil(\pi(t))-\Hil_T)Z)  \nn \\
& = \bm -\Laplace +\alpha_T^2 - (\sigma+1)\phi_T^{2\sigma} & -\sigma\phi_T^{2\sigma} \\
\sigma\phi_T^{2\sigma}  & \Laplace -\alpha_T^2
+(\sigma+1)\phi_T^{2\sigma} \endm U(t) +
M_T\calG_T(F+(\Hil(\pi(t))-\Hil_T)\calG_T^{-1}M_T^{-1} U).
\end{align*}
It remains to compute the terms
\begin{align}
 i\dot{\pi}\partial_\pi \tilde W_T(\pi)+N_T(U,\pi) &= M_T(t)\calG_T(t) F(t) \label{eq:Fdef} \\
 V_T &= M_T(t)\calG_T(t)(\Hil(\pi(t))-\Hil_T)\calG_T(t)^{-1}M_T(t)^{-1} \label{eq:Vdef}
\end{align}
In view of \eqref{eq:Zsys}, one has
\begin{align*}
 F &= \dot{v} \binom{-xe^{i\theta} \phi(x-y)}{xe^{-i\theta} \phi(x-y)} +
\dot{\gamma} \binom{-e^{i\theta} \phi(x-y)}{e^{-i\theta} \phi(x-y)} \\
&\quad + i\dot{\alpha} \binom{e^{i\theta}\partial_\alpha
\phi(x-y)}{e^{-i\theta}\partial_\alpha \phi(x-y)}
 + i\dot{D} \binom{-e^{i\theta}\partial_x \phi(x-y)}{-e^{-i\theta}\partial_x \phi(x-y)}
 + \binom{N_1(Z,W)}{N_2(Z,W)}.
\end{align*}
Now
\[ \theta(t,x+y_T)-(\alpha_T^2t+v_T(x+D_T)+t|v_T|^2+\gamma_T)
=\theta(t,x+y_T)-\theta_T(t,x+y_T)=\rho_T(t,x),\] see
\eqref{eq:ti} and Lemma~\ref{lem:rho_T}. Thus, the first term of
$M_T\calG_T F$ is
\begin{align*}
 & \dot{v} \bm e^{i\omega_T} & 0\\ 0 & e^{-i\omega_T} \endm \binom{-(x+y_T)e^{i\theta(t,x+y_T)}\, e^{-i(t|v_T|^2+v_T (x+D_T)+\gamma_T)} \phi(x+y_T-y)}{(x+y_T)
e^{-i\theta(t,x+y_T)}\, e^{i(t|v_T|^2+v_T(x+D_T)+\gamma_T)} \phi(x+y_T-y)} \\
&= \dot{v}  \binom{-(x+y_T)e^{i\rho_T(t,x)}
\phi(x+y_T-y)}{(x+y_T)e^{-i\rho_T(t,x)} \phi(x+y_T-y)}.
\end{align*}
This gives the $\dot{v}$ term in~\eqref{eq:vdotetc}. The other
terms involving $\dot{\alpha}, \dot{\gamma}$, and $\dot{D}$ are
treated similarly, and we skip the details. The cubic term
in~\eqref{eq:Zsys} is also easily transformed, and it leads to the
nonlinear term $N_T(U,\pi)$in~\eqref{eq:NUpi}. We skip that
calculation as well. Finally, it remains to transform
$\Hil(\pi(t))-\Hil_T$. One has
\[
 \Hil(\pi(t))-\Hil_T =
 \bm (\sigma+1)(\phi_T^{2\sigma}(\cdot-y_T)-\phi^{2\sigma}(\cdot-y))
& \sigma (e^{2i\theta_T}\phi_T^{2\sigma}(\cdot-y_T)  - e^{2i\theta} \phi^{2\sigma}(\cdot-y))  \\
-\sigma(e^{-2i\theta} \phi^{2\sigma}(\cdot-y) - e^{-2i\theta_T}
\phi_T^{2\sigma}(\cdot-y_T)) & -
(\sigma+1)(\phi_T^{2\sigma}(\cdot-y_T)-\phi^{2\sigma}(\cdot-y))
\endm
\]
where $\phi_T=\phi(\cdot-y_T(t),\alpha_T)$,
$\phi=\phi(\cdot-y(t),\alpha(t))$ for simplicity. It is easy to
check that
\[
\calG_T(t) (\Hil(\pi(t))-\Hil_T) = \bm (\sigma+1)(\phi_T^{2\sigma}(x)-\phi^{2\sigma}(x+y_T-y)) & * \\
-\sigma e^{-2it\alpha_T^2} (\phi_T^{2\sigma}(x)-e^{-2i\rho_T}
\phi^{2\sigma}(\cdot+y_T-y)) & * \endm \calG_T(t).
\]
After conjugation by the matrix $M_T(t)$ this takes the desired
form~\eqref{eq:V} and we are done. For the final statements
concerning $\calJ$-invariance, observe first that the
transformation~\eqref{eq:Udef} from $Z$ to $U$ preserves
$\calJ$-invariance. Second,  the equation~\eqref{eq:UPDE} is
$\calJ$-invariant, which shows that it suffices to assume the
$\calJ$-invariance of $U(0)$ to guarantee it for all $t\ge0$.  To
check the $\calJ$-invariance of~\eqref{eq:UPDE}, note that the
right-hand side of~\eqref{eq:UPDE} transforms like
\[ \calJ[i\dot{\pi}\partial_\pi \tilde W_T(\pi) + N_T(U,\pi) + V_T U] =
- [i\dot{\pi}\partial_\pi \tilde W_T(\pi) +N_T(\calJ U,\pi) +
V_T\mathcal J U],
\]
while the left-hand side transforms as follows:
\begin{align*}
& \calJ [ i\dot{U}(t) + \bm \Laplace +(\sigma+1)\phi_T^{2\sigma}-\alpha_T^2 & \sigma\phi_T^{2\sigma} \\ -\sigma\phi_T^{2\sigma} & -\Laplace -(\sigma+1)\phi_T^{2\sigma} +\alpha_T^2 \endm U ] \\
& = -i\dot{\calJ U}(t) - \bm \Laplace
+(\sigma+1)\phi_T^{2\sigma}-\alpha_T^2 & \sigma\phi_T^{2\sigma}
\\ -\sigma\phi_T^{2\sigma} & -\Laplace
-(\sigma+1)\phi_T^{2\sigma} +\alpha_T^2 \endm \calJ U
\end{align*}
Combining these statements yields the desired $\calJ$-invariance
of~\eqref{eq:UPDE}.
\end{proof}

Next, we state a standard bound on the nonlinearity $N_T(U,\pi)$.

\begin{lemma}
\label{lem:Nest} The nonlinearity $N_T(U,\pi)$ from
\eqref{eq:NUpi} satisfies
\[ |N_T(U,\pi)(x,t)|\le C\big(|U(x,t)|^{2\sigma+1}+|U(x,t)|^2\phi^{2\sigma-1}(x+y_T-y)\big) \]
for all $x,t\in\R$. Furthermore,
\begin{align*}
 |\partial_x N_T(U,\pi)(x,t)| &\le C\big(|U(x,t)|^{2\sigma}+|U(x,t)|\phi^{2\sigma-1}(x+y_T-y)\big)|\partial_x U(x,t)| \\
& \quad +
C\big(|U(x,t)|^{2\sigma+1}+|U(x,t)|^2\phi^{2\sigma-1}(x+y_T-y)\big).
\end{align*}
$C$ here only depends on $\sigma$.
\end{lemma}
\begin{proof}
Let
\[ F(z,w):= -|z+w|^{2\sigma}(z+w)+|w|^{2\sigma}w + (\sigma+1)|w|^{2\sigma}z + \sigma|w|^{2(\sigma-1)}w^2 \bar{z}\]
for all $z,w\in\Compl$. Then $F(0,w)=0$ as well as $\partial_z
F(0,w)=\partial_{\bar{z}} F(0,w)=0$. Hence,
\[ \sup_{|w|=1} |F(z,w)|\le C(|z|^{2\sigma+1}+|z|^2).\]
Rescaling implies that
\[ |F(z,w)|\le C(|z|^{2\sigma+1}+|z|^2 |w|^{2\sigma-1}).\]
As for the second statement, we use the bound
\[ \sup_{|w|=1} (|\partial_z F(z,w)|+|\partial_{\bar z} F(z,w)|)\le C(|z|^{2\sigma}+|z|),\]
as well as
\[ \sup_{|w|=1} (|\partial_w F(z,w)|+|\partial_{\bar w} F(z,w)|)\le C(|z|^{2\sigma+1}+|z|^2).\]
Via the homogeneity, for all $w\ne0$,
\begin{align*}
 |\partial_z F(z,w)|+|\partial_{\bar z} F(z,w)| &\le C(|z|^{2\sigma}+|z||w|^{2\sigma-1}), \\
 |\partial_w F(z,w)|+|\partial_{\bar w} F(z,w)| &\le C(|z|^{2\sigma+1}|w|^{-1}+|z|^2|w|^{2(\sigma-1)}).
\end{align*}
Thus, if $z=z(x)$ and $w=w(x)$, then
\[ |\partial_x F(z(x),w(x))| \le C(|z|^{2\sigma}+|z||w|^{2\sigma-1})|z_x|+C(|z|^{2\sigma+1}|w|^{-1}+|z|^2|w|^{2(\sigma-1)})|w_x|,\]
which implies that
\begin{align*}
 |\partial_x N_T(U,\pi)(x,t)| &\le C\big(|U(x,t)|^{2\sigma}+|U(x,t)|\phi^{2\sigma-1}(x+y_T-y)\big)|\partial_x U(x,t)| \\
& \quad +
C\big(|U(x,t)|^{2\sigma+1}+|U(x,t)|^2\phi^{2\sigma-1}(x+y_T-y)\big)\frac{|\partial_x\phi(x+y_T-y)|}{\phi(x+y_T-y)}.
\end{align*}
Since $\|\partial_x\phi/\phi\|_\infty\le C$, the lemma follows.
\end{proof}

\section{The linearized problem and the root spaces at zero}
\label{sec:linroot}

Recall that  $\phi=\phi(\cdot,\alpha)$ is the ground state of
$-\Laplace \phi+\alpha^2\phi=\phi^{2\sigma+1}$. Define \beeq
\label{eq:Hilal} \Hil(\alpha) :=   \bm -\Laplace -
(\sigma+1)\phi^{2\sigma}+\alpha^2 & -\sigma\phi^{2\sigma} \\
\sigma\phi^{2\sigma} & \Laplace +(\sigma+1)\phi^{2\sigma}
-\alpha^2 \endm. \eneq Hence the matrix operator on the left-hand
side of~\eqref{eq:UPDE} is equal to~$-\Hil(\alpha_T)$, i.e.,
\eqref{eq:UPDE} can be rewritten as
\[ i\partial_tU -\Hil(\alpha_T) U = i\dot{\pi}\partial_\pi \tilde W_T(\pi)+N_T(U,\pi)+V_TU \text{\ \ or\ \ } i\partial_tU -\Hil_T(t) U
= i\dot{\pi}\partial_\pi \tilde W_T(\pi)+N_T(U,\pi), \] where
$\Hil_T(t):=\Hil(\alpha_T)+V_T(t)$.

In order to prove estimates on \eqref{eq:UPDE}, we will need to
have precise control on the evolution $e^{it\Hil(\alpha)}$.
Sections~\ref{sec:BP}--\ref{sec:spectrum} deal with this issue. In
particular, in Proposition~\ref{prop:spectrum} it is shown that
the the essential spectrum of $\Hil(\alpha)$ equals
$(-\infty,-\alpha^2]\cup[\alpha^2,\infty)$, and that the discrete
spectrum equals $\{0,\pm i\gamma(\alpha)\}$ with
$\gamma=\gamma(\alpha)>0$. Here $0$ is an eigenvalue of geometric
multiplicity two and algebraic multiplicity four, whereas both
$\pm i\gamma$ are simple eigenvalues (we are dealing with the
supercritical case $\sigma>2$). In fact,
$\Hil(\alpha)f^\pm(\alpha)=\pm i\gamma f^\pm(\alpha)$ where
$f^\pm(\alpha)$ are exponentially decaying by an adaptation of
Agmon's argument, see~\cite{Sch}, and similarly $\Hil(\alpha)^*
\tilde{f}^{\pm}(\alpha)=\mp i\gamma \tilde{f}^{\pm}(\alpha)$
(cf.~also Corollary~\ref{cor:imag_proj} below). In \cite{Wei1},
Weinstein showed that the root spaces \beeq \label{eq:wsp}
 \calN = \bigcup_{n=1}^\infty \ker(\Hil(\alpha)^n), \quad
\calN^* = \bigcup_{n=1}^\infty \ker((\Hil(\alpha)^*)^n) \eneq of
$\Hil(\alpha)$ and $\Hil^*(\alpha)$, respectively, are (with
$\phi=\phi(\cdot,\alpha)$)
\begin{align}
 \calN=\calN(\alpha) &= {\rm span} \Big\{ \binom{i\phi}{-i\phi}, \binom{\partial_\alpha\phi}{\partial_\alpha\phi}, \binom{\partial_x\phi}{\partial_x\phi}, \binom{ix\phi}{-ix\phi}\Big\} \label{eq:N}\\
 \calN^*=\calN(\alpha)^* &= {\rm span} \Big\{ \binom{\phi}{\phi}, \binom{i\partial_\alpha\phi}{-i\partial_\alpha\phi}, \binom{i\partial_x\phi}{-i\partial_x\phi}, \binom{x\phi}{x\phi}\Big\}. \label{eq:N*}
\end{align}
In particular, in \eqref{eq:wsp} the kernels are the same starting
with $n=2$. Let $P_d$ be the Riesz projection onto the discrete
spectrum, i.e.,
\[
 P_d = \frac{1}{2\pi i} \oint_\gamma (zI-\Hil)^{-1}\, dz
\] where $\gamma$ is a simple closed curve that encloses the
entire discrete spectrum of~$\Hil$ and lies within the resolvent
set. Moreover, define $P_s=I-P_d$ (``s'' here stands for
``stable''). In Lemma~\ref{lem:L2split} below we show that there
is the direct but not necessarily orthogonal splitting \beeq
\label{eq:L2teilung} L^2(\R)\times L^2(\R) = \calN+{\rm
span}\{f^\pm(\alpha)\}+ \Bigl(\calN^*+  {\rm span}\{\tilde
f^\pm(\alpha)\}  \Bigr)^\perp. \eneq
 Moreover, $P_s$ is exactly the projection onto the orthogonal complement on the right-hand side with kernel equal to the sum of the first two terms,
i.e., it is the projection onto the orthogonal complement which is
induced by this splitting.

\noindent The main estimates on $e^{it\Hil}P_s$ are as follows,
see Sections~\ref{sec:scat}-\ref{sec:weigh}:
\begin{itemize}
\item  $\sup_{t}\|e^{it\Hil}P_sf\|_2 \le C \|f\|_2\quad$  and
$\quad\sup_{t}\|e^{it\Hil}P_sf\|_{H^1} \le C \|f\|_{H^1}$ \item
$\|\la x\ra e^{it\Hil}P_sf\|_2 \le C (\la t\ra \|f\|_{H^1}+\|\la
x\ra f\|_2)$ \item  $\|\la x\ra^{-\theta} e^{it\Hil}P_sf\|_\infty
\le C|t|^{-\half-\theta} \|\la x\ra^{\theta}f\|_1$ for all
$0\le\theta\le1$.
\end{itemize}
In order to apply these estimates to~\eqref{eq:UPDE}, we need to
project $U$ onto $\Ran(P_s)$. Following common practice, see
Soffer, Weinstein~\cite{SofWei1}, \cite{SofWei2}, and Buslaev,
Perelman~\cite{BP1}, we will make an appropriate choice of the
path~$\pi(t)$ in order to insure that $U(t)$ is perpendicular to
$\calN^*$. However, for technical reasons it is necessary to
impose an orthogonality condition onto a {\em time-dependent}
family of functions rather than $\calN^*$ itself. We introduce
this family in the following definition. In view of
Lemma~\ref{lem:rho_T}, it approaches $\calN^*$ in the limit $t\to
T$.

\begin{defi}
\label{def:rootspace} Assume that $\pi$ is an admissible path and
let $y,\theta$ be as in~\eqref{eq:y}, \eqref{eq:theta},
$y_T,\theta_T$ as in~\eqref{eq:inf_path}, and $\rho_T$ as
in~\eqref{eq:rhoinf}. With these functions, define\footnote{This
notation is a bit inaccurate, since $\tilde{\xi}_{\ell}$ depends on
both $T$ and the path $\pi$ chosen. We will later explicitly
denote the path dependence by $\xi_{\ell}(\pi)$, and imply a time $T$
explicitly chosen for each path $\pi$.}
\begin{align*}
\xitil_1(t) &:=
\binom{e^{i\rho_T}\phi(\cdot+(y_T-y)(t),\alpha(t))}{e^{-i\rho_T}\phi(\cdot+(y_T-y)(t),\alpha(t))},\qquad
\xitil_2(t) := \binom{ie^{i\rho_T}\partial_\alpha\phi(\cdot+(y_T-y)(t),\alpha(t))}{-ie^{-i\rho_T}\partial_\alpha\phi(\cdot+(y_T-y)(t),\alpha(t))} \\
\xitil_{3}(t) &:=
\binom{e^{i\rho_T}(x+(y_T-y)(t))\phi(\cdot+(y_T-y)(t),\alpha(t))}{e^{-i\rho_T}(x+(y_T-y)(t))\phi(\cdot+(y_T-y)(t),\alpha(t))},\qquad
\xitil_{4}(t) :=
\binom{ie^{i\rho_T}\partial_x\phi(\cdot+(y_T-y)(t),\alpha(t))}{-ie^{-i\rho_T}
\partial_x\phi(\cdot+(y_T-y)(t),\alpha(t))}.
\end{align*}
We similarly introduce the notation
\begin{align*}
\xi_1(t) &:=
\binom{e^{i\theta(t,x)}\phi(\cdot-y(t),\alpha(t))}{e^{-i\theta(t,x)}\phi(\cdot-y(t),\alpha(t))},\qquad
\xi_2(t) := \binom{ie^{i\theta(t,x)}\partial_\alpha\phi(\cdot-y(t),\alpha(t))}{-ie^{-i\theta(t,x)}\partial_\alpha\phi(\cdot-y(t),\alpha(t))} \\
\xi_{3}(t) &:=
\binom{e^{i\theta(t,x)}(x-y(t))\phi(\cdot-y(t),\alpha(t))}{e^{-i\theta(t,x)}(x-y(t))\phi(\cdot-y(t),\alpha(t))},\qquad
\xi_{4}(t) :=
\binom{ie^{i\theta(t,x)}\partial_x\phi(\cdot-y(t),\alpha(t))}{-ie^{-i\theta(t,x)}
\partial_x\phi(\cdot-y(t),\alpha(t))}.
\end{align*}
We also introduce other families $\{\etatil_j\}_{j=1}^4$ and $\{\eta_j\}_{j=1}^4$ by
\beeq
\label{eq:etadef}
 \etatil_j = \bm -i & 0\\ 0 & i \endm \xitil_j,\, \eta_j = \bm -i & 0\\ 0 & i \endm \xi_j,\text{\ \ for any\ \ }1\le j\le 4.
\eneq
\end{defi}

\noindent By inspection, $\calJ \xitil_j=\xitil_j$, $\calJ
\xi_j=\xi_j$ for $1\le j\le 4$ and we chose $\etatil_j$, $\eta_j$
in such a way that $\calJ \etatil_j=\etatil_j$, $\calJ
\eta_j=\eta_j$  for each $j$. Clearly, while the $\xitil_j$,
$\xi_j$ correspond to $\Hil^*$, the $\etatil_j$, $\eta_j$
correspond to~$\Hil$, see~\eqref{eq:N} and~\eqref{eq:N*}. Next, we
modify the $\gamma$ parameter.

\begin{lemma}
\label{lem:sigma_tild} Let $\pi(t)$ be an admissible path as in
Definition~\ref{def:adm}. Set \beeq \label{eq:dotga}
\dot{\gatil}(t):=\dot{\gamma}(t)+\dot{v}(t) y(t) \eneq
 and $\gatil(0):=0$, i.e.,
\[  \gatil(t):= \int_0^t \Big[\dot{\gamma}(s)+\dot{v}(s) y(s)\Big]\,ds. \]
Then the function $\dot{\pi}\partial_\pi \tilde W_{T}(\pi)$ on the
right-hand side of~\eqref{eq:UPDE} satisfies
\[
\dot{\pi}\partial_\pi \tilde W_{T}(\pi)=
-\dot{D}\etatil_{4}-\dot{v}\etatil_{3} + \dot{\alpha}\etatil_2 -
\dot{\gatil}\etatil_1
\]
where the functions $\{\etatil_j\}_{j=1}^4$ are as
in~\eqref{eq:etadef}. Similarly, we have (see \eqref{eq:Zsys})
\begin{equation}\nonumber
\dot{\pi}\partial_{\pi}W=-\dot{D}\eta_{4}-\dot{v}\eta_{3} +
\dot{\alpha}\eta_2 - \dot{\gatil}\eta_1
\end{equation}
\end{lemma}
\begin{proof} By inspection. \end{proof}

The following lemma records some useful facts about the two
families in Definition~\ref{def:rootspace}.

\begin{lemma}
\label{lem:orth}
Let $\phi=\phi(\cdot,\alpha(t))$ be the ground
state of~\eqref{eq:ground} and let $\{\xitil_j\}_{j=1}^4$ and
$\{\etatil_j\}_{j=1}^4$ be as in Definition~\ref{def:rootspace}.
Then
\begin{align*}
\la\xitil_1,\etatil_j \ra &= 2 \la\partial_\alpha\phi,\phi\ra \text{\ \ if $j=2$ and $=0$ else}, \\
\la\xitil_2,\etatil_j \ra &= -2 \la\partial_\alpha\phi,\phi\ra \text{\ \ if $j=1$ and $=0$ else}, \\
\la\xitil_{3},\etatil_j \ra &= -\la\phi,\phi \ra \text{\ \ if $j=4$ and $=0$ else}, \\
\la \xitil_{4},\etatil_j \ra &= \la\phi,\phi \ra \text{\ \ if
$j=3$ and $=0$ else}.
\end{align*}
Here $\partial_\alpha \la\phi,\phi\ra =
2\la\partial_\alpha\phi,\phi\ra=(2\sigma^{-1}-1)\alpha^{-1}\|\phi\|_2^2$.
\end{lemma}
\begin{proof} A simple calculation.\end{proof}
We can now derive the usual modulation equations for the
admissible path $\pi$ under the orthogonality condition \beeq
\label{eq:OC}
 \la U(t),\xitil_j(t) \ra =0
\eneq for all $t\ge 0$, $1\le j\le 4$.

\begin{lemma}
\label{lem:modul} Assume that $\pi$ is an admissible path and that
$U$ is an $H^1$ solution of~\eqref{eq:UPDE} with an initial
condition $U(0)$ which satisfies the  orthogonality
condition~\eqref{eq:OC} at time $t=0$. Then $U$ satisfies the
orthogonality assumptions~\eqref{eq:OC} for all times iff $\pi$
satisfies the following modulation equations (with
$\phi=\phi(\cdot,\alpha(t))$)
\begin{align*}
\dot{\alpha}(2\sigma^{-1}-1)\alpha^{-1}\|\phi\|_2^2 &= \la
U,\dot{\pitil}\tilde{\calS}_{1}(.+y_{T}-y)\ra
+  \la i N_T(U,\pi),\xitil_1 \ra \\
\dot{\gatil}(2\sigma^{-1}-1)\alpha^{-1}\|\phi\|_2^2 &= \la
U,\dot{\pitil}\tilde{\calS}_{2}(.+y_{T}-y)\ra + \la i N_T(U,\pi),\xitil_2 \ra \\
\dot{D} \|\phi\|_2^2 &= \la
U,\dot{\pitil}\tilde{\calS}_{3}(.+y_{T}-y)\ra
+ \la i N_T(U,\pi),\xitil_{3} \ra \\
-\dot{v} \|\phi\|_2^2 &= \la
U,\dot{\pitil}\tilde{\calS}_{4}(.+y_{T}-y) \ra + \la
i N_T(U,\pi),\xitil_{4} \ra.
\end{align*}
In these formulae we denote by
$\dot{\pitil}\tilde{\calS}_{\ell}(.+y_{T}-y)$ a linear combination of
four rapidly decaying smooth functions with coefficients
$\dot{\pitil}$ and centered at $(y_{T}-y)(t)$. The right-hand side
is real-valued, consisting of scalar products of $\calJ$-invariant
vectors. We denote (as before)
$\tilde{\pi}(t)=(\tilde{\gamma}(t),v(t),D(t),\alpha(t))$.
\end{lemma}
\begin{proof} The orthogonality conditions \eqref{eq:OC} are
equivalent to the conditions $\la Z,\xi_{\ell}\ra=0$, $\ell=1,2,3,4$.
Observe that assuming this, we have an equality of the form
\begin{equation}\nonumber
\la[i\partial_{t}-\Hil(\pi(t))]Z,\xi_{\ell}\ra = \la
Z,\dot{\tilde{\pi}}\calS_{\ell}(.-y(t))\ra
\end{equation}
for suitable rapidly decaying functions $\calS_{\ell}(\cdot)$. More
precisely, without assuming the orthogonality, we have for
suitable $\lambda_{\ell}\in\C$
\begin{equation}\nonumber
\la[i\partial_{t}-\Hil(\pi(t))]Z,\xi_{\ell}\ra = \la
Z,\dot{\tilde{\pi}}\calS_{\ell}(.-y(t))\ra+\lambda_{\ell}\la
Z,\xi_{\ell}\ra+i\partial_{t}\la Z,\xi_{\ell}\ra
\end{equation}
The statement of the lemma follows from this, the fact that
solutions to first order linear ODE vanish identically if they
vanish at one point, and the preceding lemmata.
\end{proof}
\begin{remark} The proof reveals that one may cast this system in
the $Z$-picture schematically as follows:
\begin{equation}\nonumber
\la i\dot{\tilde{\pi}}\partial_{\pi}W(\pi),\xi_{\ell}(\pi)\ra =\la
Z(t),i\dot{\tilde{\pi}}\calS_{\ell}(\pi)(t)\ra -\la
N(Z,\pi),\xi_{\ell}(t)\ra
\end{equation}
Equivalently, in the $U$-picture, this can be written
\begin{equation}\nonumber
\la i\dot{\tilde{\pi}}\partial_{\pi}W(\pi),\xi_{\ell}(\pi)\ra =\la
U(t),i\dot{\tilde{\pi}}\tilde{\calS_{\ell}}(.+y_{T}-y)(t)\ra -\la
N_T(U,\pi),\tilde{\xi}_{\ell}(t)\ra
\end{equation}

\end{remark}

\section{Constructing the solution: The iteration scheme}
\label{sec:iterate}

According to Lemma~\ref{lem:Zeq}, in order to solve the
NLS~\eqref{eq:NLS} with $\psi(t)=W(t)+R(t)$, we need to find an
admissible path $\pi(t)$ as well as a function
\[ Z\in C([0,\infty),H^1(\R)\times H^1(\R))\cap C^1([0,\infty),H^{-1}(\R)\times H^{-1}(\R))\]
so that $Z(t)$ is $\calJ$-invariant and such that $(\pi(t),Z(t))$
together satisfy~\eqref{eq:Zsys}. This will be accomplished by
means of an iteration argument. To explain it, we will need to
deal with several paths simultaneously. Therefore, our notations
will need to indicate relative to which paths Galilei transforms,
root spaces, etc.~are defined. For example, $\calG_T(\pi)(t)$ will
mean the (vector) Galilei transform from~\eqref{eq:Udef} defined
in terms of~$\pi$, and $\{\xitil_j(\pi)(t)\}_{j=1}^4$ will be the
set of functions from Definition~\ref{def:rootspace} which are
obtained from~$\pi$ together with a time $T$ specified explicitly
in conjunction with $\pi$. The iteration scheme is based on the
linearized equation~\eqref{eq:Zsys}. In principle, we want a
suitable bounded subset of some Banach space $\tilde{X}_*$,
containing $(Z,\pi)$ (defined globally in time), such that given
some $(\pizer,\Zzer)\in \tilde{X}_*$ with
$\Zzer=\binom{R^{(0)}}{\bar{R}^{(0)}}$, we can solve for
\begin{align}
 i\partial_t Z(t) + \bm \Laplace + (\sigma+1)|W(\pizer)|^2 & \sigma W^2(\pizer) \\ -\sigma\bar{W}^2(\pizer) & -\Laplace -(\sigma+1)|W(\pizer)|^2 \endm Z(t)
&= i\dot{{\pi}}\partial_\pi W(\pizer) + N(\Zzer,\pizer).
\nn 
\end{align}
The vector $\dot{\tilde{\pi}}$ is determined by means of the
modulation equations. The initial datum $Z(0)$ should be chosen in
such a fashion that the solution $Z(t)$ does not grow in time, i. e.,
such that the exponentially growing mode remain controlled, see below. The
initial datum $\pi(0)$ is a fixed constant. Then we need to show
that $(Z,\pi)$ again lies in the same subset of $\tilde{X}_*$.
Unfortunately, this straightforward approach runs into severe
difficulties due to the fact that our norms (defining the
underlying Banach space) contain weights centered around the usual $y$-curve determined by
the path $\pi^{(0)}$, and $\pi$ diverges from $\pi^{(0)}$ at infinity,
leading to an incompatible norm. In particular, there is no
absolute $\tilde{X}_*$ we can work with, but only relative
versions of the form $\tilde{X}_{*}(\pi^{(0)})$, depending on the
given path. Our way out of this consists in adapting the time
intervals on which the iterates are constructed. More precisely,
we shall solve the above equation on progressively longer
intervals, which are chosen in such a fashion that the paths differ
very little on them. In particular, the weighted norms are all
compatible on such intervals. The details of the iterative
construction are rather involved and we now present them:

\begin{defi} We let $X_*$ be the subset of the function
space:
\begin{align*} X&:=\Big\{ (\pi,U)\in  {\rm Lip}([0,\infty),\R^4)\times
\big[L^\infty((0,\infty),H^1(\R)\times H^1(\R))\cap L^\infty_{\rm
loc}((0,\infty),Y\times Y)\big] \:\big|\:\\&
\pi(0)=(\alpha_0,0,0,0)\Big\}
\end{align*}
where $Y=\{f\in H^1(\R)\:|\: \la x\ra f\in L^2,\,\partial_x f\in
L^q(\R)\}$ for a very large number\footnote{It will be seen that
$q\rightarrow\infty$ as $\sigma\rightarrow 2$.} $q$, for which the
following norm is $<\infty$:
\begin{align}
\|(\pi,U)\|_{X_*}:=&\sup_{0\leq t<\infty}\Big\{\la
t\ra^{2+\epsilon}[|\dot{\alpha}(t)|+|\dot{v}(t)|+|\dot{\gatil}(t)|+|\dot{D}(t)|]
 +\|U(t)\|_2 +  \|\partial_x U(t)\|_2 \label{eq:X*norm} \\
&+\sup_{0\leq\theta\leq 1}\la t\ra^{-\theta}\|\la x\ra^{\theta}
U(t)\|_2+ \sup_{0\leq\theta\leq\frac{1}{2}+\epsilon}\la
t\ra^{\frac12+\theta}\|\la x\ra^{-\theta} U(t)\|_\infty +\la
t\ra^{1+\eps}\|\la x\ra^{-\half-2\eps} \partial_x U(t)\|_q \Big\} \nn
\end{align}
We define
$X_{*}([0,T])$ analogously, replacing $\infty$ by $T$. We also
introduce time-localized versions of this norm, as follows:
\begin{equation}
\nonumber
||(\pi,U)||_{X_*([0,T])}:=\inf_{(\tilde{\pi},\tilde{U})|_{[0,T]}=(\pi,U)}||(\tilde{\pi},\tilde{U})||_{X_*}
\end{equation}
In the last line, $(\tilde{\pi},\tilde{U})$ are to be in $X$.
Finally, we introduce $||(\pi,Z)||_{\tilde{X}_*(\pi^{(0)})}$. The
definition of the latter is the same as for $||.||_{X_{*}}$,
except that we replace $\la x\ra$ be $\la x-y^{(0)}(t)\ra$. Here
the quantity $y^{(0)}(t)$ is given by \eqref{eq:y} with respect to
the path $\pi^{(0)}$. Time-localized versions of this norm are
defined as before.
\end{defi}

With the norm \eqref{eq:X*norm} the space $X_*$ becomes a Banach space.
With these tools we can now detail the iterative step and the a
priori estimates: we shall again use the notation
\begin{equation}\nonumber
\dot{\tilde{\pi}}=(\dot{\tilde{\gamma}},\dot{v},\dot{D},\dot{\alpha}),\,\dot{\tilde{\gamma}}=\dot{\gamma}+\dot{v}y
\end{equation}
Also, we shall use the notation
\begin{equation}\nonumber
\tilde{\pi}(t)=\int_{0}^{t}\dot{\tilde{\pi}}(s)ds + (0,0,0,\alpha_0).
\end{equation}

\begin{theorem}\label{iterate: a priori est}
Let $T_{i}=i+\delta^{-1}$, where $\delta>0$ is as
in the preceding definition, $i\in{\mathbf{Z}}_{\geq 0}$. There
exists $\delta_{0}>0$ such that if $0\leq\delta<\delta_{0}$, there
exist positive numbers $A,C$ with the following properties: Assume
that for $i\geq 1$, we are given $(\pi^{(j)},Z^{(j)})\in
X_{*}([0,T_i])$, $0\leq j\leq i-1$, all $Z^{(j)}$
$\calJ$-invariant, with the properties
\begin{equation}\nonumber
||(\pi^{j},Z^{(j)})||_{\tilde{X}_*(\pi^{(j-1)})([0,T_{j}])}<
C\delta,\,\max_{1\leq j\leq i-1}
T_{j}\sup_{t\in[0,T_{j}]}|\tilde{\pi}^{(j)}-\tilde{\pi}^{(j-1)}|(t)<A
\end{equation}
Also, assume that $Z^{(j)}$ is constant past time $t=T_{j}$, and
that $\pi^{(j)}$ is a straight line past $t=T_{j}$. Then given
$R_{0}$ as in Theorem~\ref{thm:main}, there exists a canonical
procedure for determining $\calJ$-invariant initial
data\footnote{of course, these initial data are not just close to
$ \bm R_{0}\\
\overline{R_{0}} \endm$, but uniquely defined perturbation thereof} $Z^{(i)}(0)$ satisfying $\la
Z^{(i)}(0),\xi_{\ell}(\pi^{(i-1)})\ra=0$, such that the following
conclusion applies: the combined system
\begin{align}
 &i\partial_t Z^{(i)}(t) + \bm \Laplace + (\sigma+1)|W(\pi^{(i-1)})|^{2\sigma} & \sigma |W|^{2\sigma-2}(\pi^{(i-1)})W^{2}(\pi^{(i-1)}) \\ -\sigma |W|^{2\sigma-2}(\pi^{(i-1)})W^{2}(\pi^{(i-1)}) & -\Laplace -(\sigma+1)|W(\pi^{(i-1)})|^{2\sigma} \endm Z^{(i)}(t)
\nn \\
&\hspace{8cm}= i\dot{{\pi}}^{(i)}\partial_\pi W(\pi^{(i-1)}) +
N(Z^{(i-1)},\pi^{(i-1)}), \nn\\
&Z^{(i)}(t)|_{t=0}=Z^{(i)}(0)
\label{eq:Zsyszer}\\
&\la
i\dot{{\pi}}^{(i)}\partial_{\pi}W(\pi^{(i-1)}),\xi_{\ell}(\pi^{(i-1)})\ra
=\la
Z^{(i)}(t),i\dot{\tilde{\pi}}^{(i-1)}\calS_{\ell}(\pi^{(i-1)})(t)\ra
-\la
N(Z^{(i-1)},\pi^{(i-1)}),\xi_{\ell}(\pi^{(i-1)}(t))\ra,\nn \\
&\pi^{(i)}(0)=(0,0,0,\alpha_{0}) \nn
\end{align}
has a solution on $[0,T_{i}]$ satisfying the inequalities
\begin{equation}\label{apriori estimates}
||(\pi^{(i)},Z^{(i)})||_{\tilde{X}_{*}(\pi^{(i-1)})([0,T_{i}])}<
C\delta,\,
T_{i}\sup_{t\in[0,T_{i}]}|\tilde{\pi}^{(i)}-\tilde{\pi}^{(i-1)}|(t)<A
\end{equation}
Also, $Z^{(i)}$ is $\calJ$-invariant, and $\delta$ can be made
small independently of $A,C$.
\end{theorem}
In light of the theorem, we can make the following
\begin{defi}\label{Iteration} Let $\delta_{0}>0$, $A,C$ be as in the preceding theorem. The iterates
$Z^{(i)}$ are defined as follows:
put $Z^{(0)}:=\bm R_{0}\\ \overline{R_{0}}\endm$,
$\pi^{(0)}:=(\alpha_{0},0,0,0)$. Then determine
$(Z^{(i)},\pi^{(i)})$, $i\geq 1$ from the preceding theorem: given
$(\pi^{(0)},Z^{(0)}),\ldots,(\pi^{(i-1)},Z^{(i-1)})$, one
constructs $(\pi^{(i)},Z^{(i)})$ on $[0,T_{i}]$ and extends
$Z^{(i)}$ beyond $T_{i}$ as a constant, and $\pi^{(i)}$ as a
straight line.
\end{defi}
We now prove the theorem.
\begin{proof}
In order to avoid confusion, we shall stick to the following
conventions: we let $\xitil_{\ell}(\pi^{(i)})$ etc.~denote the
functions defined  in Definition~\ref{def:rootspace}, with the
translations $(y_{T}-y)(t)$ replaced by
$(y^{(i)}_{T_{i}}-y^{(i)})(t)$ etc.~and time $T_{i}$ as above.
Also, we denote the basis of the root space of the operator
\begin{align}\nonumber
\bm
\triangle+(\sigma+1)\phi(x,\alpha^{(i-1)}_{T_{i}})^{2\sigma}-(\alpha^{(i-1)}_{T_{i}})^{2}&\sigma\phi(x,\alpha^{(i-1)}_{T_{i}})^{2\sigma}\\
-\sigma\phi(x,\alpha^{(i-1)}_{T_{i}})^{2\sigma}&-\triangle-(\sigma+1)\phi(x,\alpha^{(i-1)}_{T_{i}})^{2\sigma}+(\alpha^{(i-1)}_{T_{i}})^{2}\\
\endm
\end{align}
by $\tilde{\eta}_{\ell}(\alpha^{(i-1)}_{T_{i}})$, and let
eigenvectors corresponding to the imaginary eigenvalues $\pm
i\gamma(\alpha^{(i-1)}_{T_{i}})$ be
$f^{\pm}(\alpha^{(i-1)}_{T_{i}})$, both of which are chosen to be
$\calJ$-invariant. We shall make the following ansatz for
$Z^{(i)}(0)$:
\begin{align}\label{data}
Z^{(i)}(0)=\bm R_{0}\\
\overline{R_{0}}\endm+h^{(i)}f^{+}(\alpha^{(i-1)}_{T_{i}})+\sum_{j=1}^{4}a_{j}^{(i)}\tilde{\eta}_{j}(\alpha^{(i-1)}_{T_{i}}),\quad h^{(i)},\,
a^{(i)}_{j}\in\R
\end{align}
As stated in the theorem, we want this to satisfy the
orthogonality relations
\[ \la
Z^{(i)}(0),\xi_{\ell}(\pi^{(i-1)}(0))\ra=0.\]
The assumptions in the theorem imply in particular that
$|\tilde{\eta}_{j}(\alpha^{(i-1)}_{T_{i}})-\tilde{\eta}(\alpha_{0})|\lesssim
\delta$. This entails that for $\delta$ sufficiently small, we may
uniquely solve the system
\begin{align}\nonumber
\Big\la\bm R_{0}\\
\overline{R_{0}}\endm+h^{(i)}f^{+}(\alpha^{(i-1)}_{T_{i}})+\sum_{j=1}^{4}a_{j}^{(i)}\tilde{\eta}_{j}(\alpha^{(i-1)}_{T_{i}}),\xi_{\ell}(\pi^{(i-1)}(0))\Big\ra
=0
\end{align}
for the $a_{j}^{(i)}$, for given $h^{(i)}\in \R$ and fixed $R_{0}$
as in Theorem~\ref{thm:main}. Moreover, the fact that both
$\tilde{\eta}_{j}(\alpha^{(i-1)}_{T_{i}})$,
$\xi_{\ell}(\pi^{(i-1)}(0))$ are $\calJ$ invariant implies that
$a_{j}^{(i)}(h^{(i)},\alpha^{(i-1)}_{T_{i}})\in\R$ for
$h^{(i)}\in\R$. We shall then determine $h^{(i)}$ in such a fashion
that the corresponding solution $Z^{(i)}(t)$ does not grow with
time. In order to carry out estimates, we apply a Gauge
transformation to \eqref{eq:Zsyszer}. Specifically, we put
\begin{equation}\nonumber
U^{(i)}=M_{T_{i}}(\pi^{(i-1)})\calG_{T_{i}}(\pi^{(i-1)})Z^{(i)},
\end{equation}
where we have
\begin{align*}\nonumber
M_{T_{i}}(\pi^{(i-1)})(t)=\bm
e^{i\alpha^{(i-1)}(T_{i})^{2}t}&0\\0&
e^{-i\alpha^{(i-1)}(T_{i})^{2}t}\\\endm
\end{align*}
\begin{align*}\nonumber
\calG_{T_{i}}(\pi^{(i-1)})(t,x)=\bm
e^{-i(\gamma^{(i-1)}_{T_{i}}+v^{(i-1)}_{T_{i}}D^{(i-1)}_{T_{i}}+v^{(i-1)}_{T_{i}}x+t(v^{(i-1)}_{T_{i}})^{2})}e^{i(2tv^{(i-1)}_{T_{i}}+D^{(i-1)}_{T_{i}})p}\\
\overline{e^{-i(\gamma^{(i-1)}_{T_{i}}+v^{(i-1)}_{T_{i}}D^{(i-1)}_{T_{i}}+v^{(i-1)}_{T_{i}}x+t(v^{(i-1)}_{T_{i}})^{2})}e^{i(2tv^{(i-1)}_{T_{i}}+D^{(i-1)}_{T_{i}})p}}\\\endm
\end{align*}
The notation here is the same as in Section~\ref{sec:ansatz}. We
apply the same procedure to the $Z^{(j)}$, thereby introducing
quantities $U^{(j)}$, $j\leq i-1$. Proceeding as in
Section~\ref{sec:ansatz}, we derive the following equation in the
gauged picture:
\begin{equation}\nonumber\begin{split}
&i\partial_{t}U^{(i)}-\Hil(\alpha^{(i-1)}_{T_{i}})U^{(i)}\\&=
M_{T_{i}}(\pi^{(i-1)})\calG_{T_{i}}(\pi^{(i-1)})[i\dot{{\pi}}^{(i)}\partial_{\pi}W(\pi^{(i-1)})+N(Z^{(i-1)},\pi^{(i-1)})+(\Hil(\pi^{(i-1)}(t))-\Hil(\pi^{(i-1)}_{T_{i}})(t))Z^{(i)}]\\
\end{split}\end{equation}
Here  $\Hil(\pi^{(i-1)}_{T_{i}}(t))$ is defined as
$\Hil(\pi^{(i-1)}(t))$ with the path $y^{(i-1)}(t)$ replaced by
the straight line path given by
$2tv^{(i-1)}_{T_{i}}+D^{(i-1)}_{T_{i}}$. Written out, the
following equation results:
\begin{align}
&i\dot{U}^{(i)}(t) + \bm \Laplace
+(\sigma+1)\phi_{T_{i}}^{2\sigma}-(\alpha^{(i-1)}_{T_{i}})^2 &
\sigma\phi_{T_{i}}^{2\sigma} \\ -\sigma\phi_{T_{i}}^{2\sigma} &
-\Laplace -(\sigma+1)\phi_{T_{i}}^{2\sigma}
+(\alpha^{(i-1)}_{T_{i}})^2
\endm U^{(i)} \nn\\
& = i\dot{{\pi}}^{(i)}\partial_\pi \tilde
W_{T_{i}}(\pi^{(i-1)})
+N_{T_{i}}(M_{T_{i-1}}(\pi^{(i-1)})\calG_{T_{i-1}}(\pi^{(i-1)})\calG_{T_{i-2}}(\pi^{(i-2)})^{-1}M_{T_{i-2}}(\pi^{(i-2)})^{-1}U^{(i-1)},\pi^{(i-1)})
\label{eq:UPDE1} \\
&+ V^{(i-1)}U^{(i)} \nn
\end{align}
We use the abbreviation
$\phi_{T_{i}}=\phi(.,\alpha^{(i-1)}_{T_{i}})$, as well as the following:
\begin{equation}
\nn
\begin{split}
 &V^{(i-1)} = V_{T_{i}}(\pi^{(i-1)}) =\\
& \bm (\sigma+1)(\phi_{T_{i}}^{2\sigma}-\phi^{2\sigma}(\cdot+y^{(i-1)}_{T_{i}}-y^{(i-1)})) & \sigma(\phi_{T_{i}}^{2\sigma}-e^{2i\rho^{(i-1)}_{T_{i}}} \phi^{2\sigma}(\cdot+y^{(i-1)}_{T_{i}}-y^{(i-1)})) \\
-\sigma(\phi_{T_{i}}^{2\sigma}- e^{-2i\rho^{(i-1)}_{T_{i}}} \phi^{2\sigma}(\cdot+y^{(i-1)}_{T_{i}}-y^{(i-1)})) & -(\sigma+1)(\phi_{T_{i}}^{2\sigma}-\phi^{2\sigma}(\cdot+y^{(i-1)}_{T_{i}}-y^{(i-1)})) \endm \label{eq:V1}\\
\end{split}\end{equation}
\begin{align}
 &i\dot{{\pi}}^{(i)}\partial_\pi \tilde W_{T_{i}}(\pi^{(i-1)}) := \nn \\
&\dot{{v}}^{(i)} \binom{-(\cdot+y^{(i-1)}_{T_{i}}-y^{(i-1)})e^{i\rho^{(i-1)}_{T_{i}}} \phi(\cdot+y^{(i-1)}_{T_{i}}-y^{(i-1)})}{(\cdot+y^{(i-1)}_{T_{i}}-y^{(i-1)}) e^{-i\rho^{(i-1)}_{T_{i}}} \phi(\cdot+y^{(i-1)}_{T_{i}}-y^{(i-1)})} \\&+ (\dot{\gamma}^{(i)}+\dot{v}^{(i)}y^{(i-1)}) \binom{-e^{i\rho^{(i-1)}_{T_{i}}}\phi(\cdot+y^{(i-1)}_{T_{i}}-y^{(i-1)})}{e^{-i\rho^{(i-1)}_{T_{i}}}\phi(\cdot+y^{(i-1)}_{T_{i}}-y^{(i-1)})} \label{eq:vdotetc1} \\
 & \quad +i\dot{\alpha}^{(i)}
 \binom{e^{i\rho^{(i-1)}_{T_{i}}}\partial_\alpha\phi(\cdot+y^{(i-1)}_{T_{i}}-y^{(i-1)})}{e^{-i\rho^{(i-1)}_{T_{i}}}\partial_\alpha\phi(\cdot+y^{(i-1)}_{T_{i}}-y^{(i-1)})}+i\dot{D}^{(i)}
\binom{-e^{i\rho^{(i-1)}_{T_{i}}}\partial_x\phi(\cdot+y^{(i-1)}_{T_{i}}-y^{(i-1)})}{-e^{-i\rho^{(i-1)}_{T_{i}}}\partial_x\phi(\cdot+y^{(i-1)}_{T_{i}}-y^{(i-1)})} \nn \\
&N_{T_{i}}(U,\pi^{(i-1)}) :=
\binom{N_{1T_{i}}(U,\pi^{(i-1)})}{N_{2T_{i}}(U,\pi^{(i-1)})} =
\binom{N_{1T_{i}}(U,\pi^{(i-1)})}{-\overline{N_{1T_{i}}(U,\pi^{(i-1)})}}
\label{eq:NUpi1}
\end{align}
and we have
\begin{align*}
&N_{1T_{i}}(U,\pi^{(i-1)}) =\\&
-|U_{1}+e^{i\rho^{(i-1)}_{T_{i}}}\phi(\cdot+y^{(i-1)}_{T_{i}}-y^{(i-1)})|^{2\sigma}(U_{1}+e^{i\rho^{(i-1)}_{T_{i}}}\phi(\cdot+y^{(i-1)}_{T_{i}}-y^{(i-1)}))\\&+\phi(\cdot+y^{(i-1)}_{T_{i}}-y^{(i-1)})^{2\sigma+1}e^{i\rho^{(i-1)}_{T_{i}}}
+(\sigma+1)\phi(\cdot+y^{(i-1)}_{T_{i}}-y^{(i-1)})^{2\sigma}U_{1}\\& +
\sigma\phi(\cdot+y^{(i-1)}_{T_{i}}-y^{(i-1)})^{2\sigma}
e^{2i\rho^{(i-1)}_{T_{i}}} U_2.
\end{align*}
As can be seen from \eqref{eq:vdotetc1}, we need to define
\beeq
\label{eq:yi-1} \dot{\tilde\gamma}^{(i)} = \dot{\gamma}^{(i)}+\dot{v}^{(i)}y^{(i-1)}
\eneq
and then, as usual, $\dot{\tilde\pi}^{(i)}=(\dot{\tilde\gamma}^{(i)}, \dot{v}^{(i)}, \dot{D}^{(i)}, \dot{\alpha}^{(i)} )$
and $\tilde{\pi}^{(i)}(0)=(0,0,0,\alpha_0)$.   As above, we will sometimes write
$\dot{\tilde{\pi}}^{(i)}\partial_\pi \tilde
W_{T_{i}}(\pi^{(i-1)})$ instead of $\dot{{\pi}}^{(i)}\partial_\pi \tilde
W_{T_{i}}(\pi^{(i-1)})$ to emphasize that we are working with $\tilde{\pi}^{(i)}$ rather than $\pi^{(i)}$ itself.

We shall analyze \eqref{eq:UPDE1} on the interval $[0,T_{i}]$, and
establish control over $||U^{(i)}||_{X_{*}([0,T_{i}])}$, upon
defining $h$ suitably. More precisely, we shall establish the
following

\begin{prop}\label{(U,pi)} There exists a canonically determined value
$h(R_{0},Z^{(i-1)},\pi^{(i-1)},\pi^{(i-2)})\in\R$ such that with
initial data \eqref{data} and the assumptions of
Theorem~\ref{iterate: a priori est}, we have for some universal
constant $C_{0}$ independent of $A,C$
\begin{equation}\nn \label{crux1}
||(\pi^{(i)},U^{(i)})||_{X_*([0,T_{i}])}\lesssim
C^{2}[(A+1)\delta]^{2}+C^{2\sigma+1}\delta^{2\sigma+1}+C_{0}\delta
\end{equation}
\end{prop}
This proposition allows us to retrieve the a priori bound on
$(\pi^{(i)},U^{(i)})$, respectively $(\pi^{(i)},Z^{(i)})$; in
order to establish Theorem~\ref{iterate: a priori est}, we still
need to retrieve control over
\[ T_{i}\sup_{0\leq t\leq
T_{i}}|\pi^{(i)}(t)-\pi^{(i-1)}(t)|.\]
This follows from the next

\begin{proposition}\label{pi-pi}
Assume $h$ is chosen as in the preceding theorem, and moreover
$A,C>1$ as above. Introduce the norm
\begin{equation}\nonumber\begin{split}
||(\pi,Z)||_{Y^{(i)}([0,T])}:=&\sup_{t\in
[0,T]}\la t\ra^{1+\frac{\epsilon}{2}}|\dot{\pi}(t)|+\sup_{t\in[0,T]}[\la t\ra^{-1}||Z(t)||_{L_{x}^{2}}+||\la x-y^{(i-1)}(t)\ra^{-\theta}Z(t)||_{L_{x}^{\infty}}]\\
\end{split}\end{equation}
Then we have the inequality\footnote{The statement here is far
from optimal, but it is all that is needed to close the
iteration.}
\begin{equation}\nn 
\begin{split}
&||(\tilde{\pi}^{(i)}-\tilde{\pi}^{(i-1)},Z^{(i)}-Z^{(i-1)})||_{Y^{(i)}([0,T_{i}])}\\&\lesssim
[A^{2}C\delta+A^{2}(C\delta)^{2\sigma}]||(\tilde{\pi}^{(i-1)}-\tilde{\pi}^{(i-2)},Z^{(i-1)}-Z^{(i-2)})||_{Y^{(i-1)}([0,T_{i-1}])}+[1+(A+1)^{2}(C\delta)^{2}+(C\delta)^{2\sigma}]
T_{i}^{-1}
\end{split}
\end{equation}
\end{proposition}
Let us now assume that these two propositions hold. Then we can finish the proof
Theorem~\ref{iterate: a priori est}.
Observe that if we iterate the inequality of the second proposition, we get
\begin{equation}\nonumber\begin{split}
&\sup_{t\in[0,T_{i}]}\la
t\ra^{1+\frac{\epsilon}{2}}|\dot{\tilde{\pi}}^{(i)}-\dot{\tilde{\pi}}^{(i-1)}|(t)\lesssim
[1+(A+1)^{2}(C\delta)^{2}+(C\delta)^{2\sigma}]T_{i}^{-1}\\&+A^{2}[C\delta+(C\delta)^{2\sigma}]
[1+(A+1)^{2}(C\delta)^{2}+(C\delta)^{2\sigma}]T_{i-1}^{-1}+\ldots\\&\hspace{7cm}+(A^{2}[C\delta+(C\delta)^{2\sigma}])^{i}[1+(A+1)^{2}(C\delta)^{2}+(C\delta)^{2\sigma}]T_{0}^{-1}
\end{split}\end{equation}
We conclude that if we choose $A,C_{0}$ large enough and then
$\delta$ small enough, we can bound
\begin{equation}\nonumber
||(\pi^{(i)},U^{(i)})||_{X_*([0,T_{i}])}<C\delta,\,T_{i}\sup_{t\in[0,T_{i}]}\la
t\ra^{1+\frac{\epsilon}{2}}|\dot{\tilde{\pi}}^{(i)}-\dot{\tilde{\pi}}^{(i-1)}|(t)<A,
\end{equation}
which suffices to close the iteration.
\end{proof}

Thus we are left with proving the two propositions. This will be
established by means of a sequence of estimates.

\begin{proof}[Proof of Proposition~\ref{(U,pi)}] We shall estimate the
various parts constituting the norm
$||(\pi^{(i)},U^{(i)})||_{X_*([0,T_{i}])}$. We commence with the
parts concerning $U^{(i)}$. We shall decompose $U^{(i)}$ into its
dispersive, root, and hyperbolic part with respect to the operator
$\Hil(\alpha^{(i-1)}_{T_{i}})$. Thus we split
\begin{equation}\nonumber
U^{(i)}=U^{(i)}_{dis}+U^{(i)}_{root}+U^{(i)}_{hyp}
\end{equation}
and write
\begin{equation}\nonumber
U^{(i)}_{root}(t)=\sum_{j=1}^{4}\tilde{a}_{j}^{(i)}(t)\tilde{\eta}_{j}(\alpha^{(i-1)}_{T_{j}}),\quad
U^{(i)}_{hyp}=b^{(i)+}(t)f^{+}(\alpha^{(i-1)}_{T_{i}})+b^{(i)-}(\alpha^{(i-1)}_{T_{i}})
\end{equation}
We mean here that $U^{(i)}_{root}(t)=P_0(\alpha^{(i-1)}_{T_{j}}) U^{(i)}$, and
$U^{(i)}_{hyp}=P_{\pm i\gamma}(\alpha^{(i-1)}_{T_{j}}) U^{(i)}$.
The orthogonality condition $\la Z^{(i)},\xi_{\ell}(\pi^{(i-1)})\ra =
0$, $\ell=1,2,3,4$, which is equivalent to $\la
U^{(i)},\tilde{\xi}_{\ell}(\pi^{(i-1)})\ra = 0$, implies an equation
of the form
\begin{equation}\nonumber
\la U^{(i)}_{dis},\tilde{\xi}_{\ell}(\pi^{(i-1)})\ra +\la
U^{(i)}_{hyp},\tilde{\xi}_{\ell}(\pi^{(i-1)})\ra
+\sum_{j=1}^{4}\tilde{a}_{j}^{(i)}(t)\la\tilde{\eta}_{j}(\alpha^{(i-1)}_{T_{j}}),\tilde{\xi}_{\ell}(\pi^{(i-1)})\ra
= 0
\end{equation}
Our a priori assumptions on $\pi^{(i-1)}$ imply, upon choosing
$\delta$ small enough, that the $4\times 4$ matrix with entries
$\la\tilde{\eta}_{j}(\alpha^{(i-1)}_{T_{j}}),\tilde{\xi}_{\ell}(\pi^{(i-1)})\ra$
is nonsingular. In particular, we can deduce formulae
$a_{j}=a_{j}(U^{(i)}_{root},U^{(i)}_{hyp},\pi^{(i-1)})$. We
feed this information back into \eqref{eq:UPDE1}. Then we
specialize this equation to the dispersive and hyperbolic parts:
for the dispersive part, we get
\beeq
\nonumber\begin{split}
 &i\dot{U^{(i)}}_{dis}(t) + \bm \Laplace
+(\sigma+1)\phi_{T_{i}}^{2\sigma}-(\alpha^{(i-1)}_{T_{i}})^2 &
\sigma\phi_{T_{i}}^{2\sigma} \\ -\sigma\phi_{T_{i}}^{2\sigma} &
-\Laplace -(\sigma+1)\phi_{T_{i}}^{2\sigma}
+(\alpha^{(i-1)}_{T_{i}})^2
\endm U^{(i)}_{dis}\\
& = P_{s}(\alpha^{(i-1)}_{T_{i}})\Big\{i\dot{\tilde{\pi}}^{(i)}\partial_\pi \tilde
W_{T_{i}}(\pi^{(i-1)}) +
V^{(i-1)}[U^{(i)}_{dis}+U^{(i)}_{hyp}+U^{(i)}_{root}(U^{(i)}_{dis},U^{(i)}_{hyp})]\\
&+N_{T_{i}}(M_{T_{i-1}}(\pi^{(i-1)})\calG_{T_{i-1}}(\pi^{(i-1)})\calG_{T_{i-2}}(\pi^{(i-2)})^{-1}M_{T_{i-2}}(\pi^{(i-2)})^{-1}U^{(i-1)},\pi^{(i-1)})\Big\}
\end{split}
\eneq
For the hyperbolic part, we need to get a condition on
$b^{(i)\pm}$: denoting the hyperbolic projection of the right-hand
side of \eqref{eq:UPDE1} as
\begin{equation}\nonumber
\bm g^{+}(U^{(i)},U^{(i-1)},\pi^{(i)},\pi^{(i-1)},\pi^{(i-2)})\\
g^{-}(U^{(i)},U^{(i-1)},\pi^{(i)},\pi^{(i-1)},\pi^{(i-2)})\endm,
\end{equation}
we can formulate the system
\[ \frac{d}{dt}\binom{b^{(i)+}}{b^{(i)-}} + \bm -\gamma(\alpha^{(i-1)}_{T_{i}}) & 0 \\ 0 & \gamma(\alpha^{(i-1)}_{T_{i}}) \endm
\binom{b^{(i)+}}{b^{(i)-}} = \binom{g^+}{g^-},
\]
In order to control the growth of $b^{(i)+}$, we use the following
lemma already used in \cite{Sch}:
\begin{lemma}
\label{lem:ODE_stable} Consider the two-dimensional ODE
\[ \dot{x}(t)-A_0 x(t) = f(t),\qquad x(0)=\binom{x_1(0)}{x_2(0)}\]
where $f=\binom{f_1}{f_2}\in L^\infty([0,\infty),\Compl^2)$ and $
A_0 = \bm \gamma & 0 \\ 0 & -\gamma \endm $ where $\gamma>0$. Then
$x(t)=\binom{x_1(t)}{x_2(t)}$ remains bounded for all times iff
\beeq \label{eq:stable} 0 = x_1(0)+\int_0^\infty e^{-\gamma t}
f_1(t)\, dt. \eneq Moreover, in that case \beeq
 x_1(t) = - \int_t^\infty e^{-(s-t)\gamma} f_1(s)\, ds, \quad \label{eq:Duh_stab}
 x_2(t) = e^{-t\gamma}x_2(0)+ \int_0^t e^{-(t-s)\gamma} f_2(s)\, ds.
\eneq
for all $t\ge0$.
\end{lemma}
\begin{proof}
Clearly, $x_1(t) = e^{t\gamma}x_1(0) + \int_0^t e^{(t-s)\gamma}
f_1(s)\, ds $ and $ x_2(t) = e^{-t\gamma}x_2(0) + \int_0^t
e^{-(t-s)\gamma} f_2(s)\, ds$. If \[\lim_{t\to\infty}
e^{-t\gamma}x_1(t)=0,\] then $ 0 = x_1(0) + \int_0^\infty
e^{-s\gamma} f_1(s)\, ds,$ which is \eqref{eq:stable}. Conversely,
if this holds, then $ x_1(t) = -e^{t\gamma} \int_t^\infty
e^{-s\gamma } f_1(s)\, ds,$ and the lemma is proved.
 \end{proof}
Of course we are working on a finite time interval, but we use
this lemma to motivate our choice of $b^{(i)+}(0)$, namely
\[
b^{(i)+}(0)=
-\int_{0}^{T_i}e^{-\gamma(\alpha^{(i-1)}_{T_{i}})t}g^{+}(U^{(i)},U^{(i-1)},\pi^{(i)},\pi^{(i-1)},\pi^{(i-2)})(t)dt.
\]
We now claim that there is a unique choice of $h^{(i)}$ for which
we have
\begin{align}\nonumber
P_{Im}^{+}(\alpha^{(i-1)}_{T_{i}})\Big\{\calG_{T_{i}}(\pi^{(i-1)}(0))\Big[\bm R_{0}\\
\overline{R_{0}}\endm +h^{(i)}
f^{+}(\alpha^{(i-1)}_{T_{i}})+\sum_{j=1}^{4}a_{j}^{(i)}(h^{(i)},\alpha^{(i-1)}_{T_{i}})\tilde{\eta}_{j}(\alpha^{(i-1)}_{T_{i}})\Big]
\Big\}
=b^{(i)+}(0)
\end{align}
This follows from the fact that due to our
assumptions\footnote{provided we choose $\delta>0$ small enough},
we have
\begin{align}\nonumber
&P_{Im}(\alpha^{(i-1)}_{T_{i}})\Big[(\calG_{T_{i}}(\pi^{(i-1)}(0))-I)\bm R_{0}\\
\overline{R_{0}}\endm\Big]=O(\delta),\\&P_{Im}(\alpha^{(i-1)}_{T_{i}})\Big[(\calG_{T_{i}}(\pi^{(i-1)}(0))-I)
\sum_{j=1}^{4}a_{j}^{(i)}(h^{(i)},\alpha^{(i-1)}_{T_{i}})\tilde{\eta}_{j}(\alpha^{(i-1)}_{T_{i}})\Big]=O\Big(\delta
h^{(i)}\Big)
\end{align}
Indeed, the dependence of $h^{(i)}$ on $b^{(i)+}(0)$ is linear.
One also sees from the $\calJ$-invariance of
$f^{+}(\alpha^{(i-1)}_{T_{i}})$ as well as the root space
representatives and the $\calJ$-invariance of the equations that
$b^{(i)+}(t)$ is always real-valued, whence so is $h^{(i)}$. With
this $h^{(i)}$, we can then define
\begin{equation}\nonumber
U^{(i)}_{dis}(0):=P_{s}(\alpha^{(i-1)}_{T_{i}})\Big\{\calG_{T_{i}}(\pi^{(i-1)}(0))\Big[\bm R_{0}\\
\overline{R_{0}}\endm +h^{(i)}
f^{+}(\alpha^{(i-1)}_{T_{i}})+\sum_{j=1}^{4}a_{j}^{(i)}(h^{(i)},\alpha^{(i-1)}_{T_{i}})\tilde{\eta}_{j}(\alpha^{(i-1)}_{T_{i}})\Big]
\Big\}
\end{equation}
Similarly, we can uniquely specify $b^{(i)-}(0)$. Of course,
solving for $U^{(i)}_{dis}$, $U^{(i)}_{hyp}$ is complicated by the
fact that the unknowns are on both sides of the equations (and
indeed also implicitly determine the initial data). Thus we need
to run a contraction argument to solve for them. Specifically,
denoting
\[ h^{(i)}=h^{(i)}(U^{(i)}_{dis},U^{(i)}_{hyp},U^{(i-1)},\pi^{(i-1)},\pi^{(i-2)})\]
in the sense just established, we introduce a map
$F_{T_{i}}(U^{(i-1)},\pi^{(i-1)},\pi^{(i-2)})$, which sends a
given pair $(\pi^{*},U^{*})$ satisfying the orthogonality
relations \[\la U^{*},\tilde{\xi}_{\ell}(\pi^{(i-1)})\ra =0,\quad \ell=1,2,3,4\]
into another one (satisfying the same orthogonality relations)
$(\pi,U)$ as follows:
\beeq\label{Ueqn}\begin{split}
 &i\dot{U}_{dis}(t) + \bm \Laplace
+(\sigma+1)\phi_{T_{i}}^{2\sigma}-(\alpha^{(i-1)}_{T_{i}})^2 &
\sigma\phi_{T_{i}}^{2\sigma} \\ -\sigma\phi_{T_{i}}^{2\sigma} &
-\Laplace -(\sigma+1)\phi_{T_{i}}^{2\sigma}
+(\alpha^{(i-1)}_{T_{i}})^2
\endm U_{dis}\\& = P_{s}(\alpha^{(i-1)}_{T_{i}})\big[i\dot{\tilde{\pi}}^{*}\partial_\pi \tilde
W_{T_{i}}(\pi^{(i-1)}) +
V^{(i-1)}(U^{*}_{dis}+U^{*}_{hyp}+U^{*}_{root}(U^{*}_{dis},U^{*}_{hyp}))\\
&+N_{T_{i}}(M_{T_{i-1}}(\pi^{(i-1)})\calG_{T_{i-1}}(\pi^{(i-1)})\calG_{T_{i-2}}(\pi^{(i-2)})^{-1}M_{T_{i-2}}(\pi^{(i-2)})^{-1}U^{(i-1)},\pi^{(i-1)})\big]\\
\end{split}
\eneq
\begin{equation}\label{eqn}\begin{split}
&U_{dis}(0)=P_{s}(\alpha^{(i-1)}_{T_{i}})\Big\{\calG_{T_{i}}(\pi^{(i-1)}(0))\Big[\bm R_{0}\\
\overline{R_{0}}\endm
+h^{(i)}(U^{*}_{dis},U^{*}_{hyp},U^{(i-1)},\pi^{(i-1)},\pi^{(i-2)})
f^{+}(\alpha^{(i-1)}_{T_{i}})\\&\hspace{4cm}+\sum_{j=1}^{4}a_{j}^{(i)}(h^{(i)}(U^{*}_{dis},U^{*}_{hyp},U^{(i-1)},\pi^{(i-1)},\pi^{(i-2)}),\alpha^{(i-1)}_{T_{i}})\tilde{\eta}_{j}(\alpha^{(i-1)}_{T_{i}})\Big]\Big\}
\end{split}\end{equation}
\begin{equation}\label{Uhypeqn}
 \frac{d}{dt}\binom{b^{+}}{b^{-}} + \bm
-\gamma(\alpha^{(i-1)}_{T_{i}}) & 0 \\ 0 &
\gamma(\alpha^{(i-1)}_{T_{i}}) \endm \binom{b^{+}}{b^{-}} = \bm
g^+(U^{*},U^{(i-1)},\pi^{*},\pi^{(i-1)},\pi^{(i-2)})\\g^-(U^{*},U^{(i-1)},\pi^{*},\pi^{(i-1)},\pi^{(i-2)})\endm
\end{equation}
\begin{equation}\nonumber
b^{+}(0)=b^{+}(U^{*},U^{(i-1)},\pi^{*},\pi^{(i-1)},\pi^{(i-2)}),\,b^{-}(0)=b^{-}(U^{*},U^{(i-1)},\pi^{*},\pi^{(i-1)},\pi^{(i-2)})
\end{equation}
\begin{equation}\label{path}
\begin{split}
&\la
i\dot{\tilde{\pi}}\partial_\pi \tilde
W_{T_{i}}(\pi^{(i-1)}),\tilde \xi_{\ell}(\pi^{(i-1)})\ra
=\la
U^{*}(t),i\dot{\tilde{\pi}}^{(i-1)}\tilde{\calS_{\ell}}(.+y^{(i-1)}_{T}-y^{(i-1)})(t)\ra \\
&-\la
N_{T_{i}}(M_{T_{i-1}}(\pi^{(i-1)})\calG_{T_{i-1}}(\pi^{(i-1)})\calG_{T_{i-2}}(\pi^{(i-2)})^{-1}
M_{T_{i-2}}(\pi^{(i-2)})^{-1}U^{(i-1)},\pi^{(i-1)}),\tilde{\xi}_{\ell}(\pi^{(i-1)}(t))\ra
\end{split}
\end{equation}
Here $\tilde{\calS_{\ell}}$ is from Lemma~\ref{lem:modul}.
Of course we have written
$U_{hyp}=b^{+}f^{+}(\alpha^{(i-1)}_{T_{i}})+b^{-}f^{-}(\alpha^{(i-1)}_{T_{i}})$,
and this in addition to $U_{dis}$ uniquely determines the root
part $U_{root}$, on account of the orthogonality relations $\la
U,\tilde{\xi}_{\ell}(\pi^{(i-1)})\ra=0$, which we assume to hold. Our
task is to find a fixed point for the affine map
\[ F_{T_{i}}(U^{(i-1)},\pi^{(i-1)},\pi^{(i-2)})\;:\; (\pi^*,U^*)\mapsto (\pi,U) \]
In order to do this,
we need to show that it is a contraction with respect to the norm
$||.||_{X_{*}([0, T_{i}])}$. This will then (finally!) define the
iterate $(\pi^{(i)},U^{(i)})$ and thereby (undoing the Gauge)
$(\pi^{(i)}, Z^{(i)})$. Note the equations~\eqref{Ueqn}--\eqref{path} can all be solved
by integration in terms of the initial data and the right-hand sides.

We shall show here that
$F_{T_{i}}(U^{(i-1)},\pi^{(i-1)},\pi^{(i-2)})$ sends the ball
\begin{equation}\nonumber
||(\pi,U)||_{X_{*}([0,T_{i}])}<C_{0}
C^{2}[(A+1)\delta]^{2}+C^{2\sigma+1}\delta^{2\sigma+1}+C_{0}\delta=:M(A,C,C_{0},\delta)
\end{equation}
into itself, provided $\delta$ is small enough in relation to
$A,C,C_{0}$, and provided the latter quantities are large enough (relative to
some absolute constant, and with $C_{0}$ small enough in relation to
$C$). The same estimates, upon considering a suitable difference
equation, will establish the contraction property, as
well as the inequality in Proposition~\ref{(U,pi)}.\\
{\bf{(A):}} Estimating
$||U_{dis}||_{L_{t}^{\infty}L_{x}^{2}([0,T_{i}])}$: for this we
use the fact that
\begin{equation}\nonumber
||U_{dis}||_{L_{t}^{\infty}L_{x}^{2}([0,T_{i}])}\lesssim
\sup_{t\in
[0,T_{i}]}||e^{it\Hil(\alpha^{(i-1)}_{T_{i}})}U_{dis}(0)||_{L_{x}^{2}}+||\text{right-hand
side of \eqref{Ueqn}}||_{L_{t}^{1}L_{x}^{2}([0,T_{i}])}
\end{equation}
We commence by estimating the various terms on the right-hand side
of \eqref{Ueqn}: first, observe that
\begin{equation}\nonumber
P_{s}(\alpha^{(i-1)}_{T_{i}})[i\dot{\tilde{\pi}}^{*}\partial_\pi
\tilde W_{T_{i}}(\pi^{(i-1)})]=
P_{s}(\alpha^{(i-1)}_{T_{i}})[i\dot{\tilde{\pi}}^{*}\partial_\pi
\tilde W_{T_{i}}(\pi^{(i-1)})-i\dot{\tilde{\pi}}^{*}\partial_\pi
\tilde W_{T_{i}}(\pi^{(i-1)}_{T_{i}})],
\end{equation}
where $\partial_\pi \tilde W_{T_{i}}(\pi^{(i-1)}_{T_{i}})$ is the
same as $\partial_\pi \tilde W_{T_{i}}(\pi^{(i-1)})$ with the path
$\pi^{(i-1)}$ replaced by the straight line path
$\pi^{(i-1)}_{T_{i}}$, see Definition~\ref{def:rootspace}, Lemma~\ref{lem:sigma_tild}, and Lemma~\ref{lem:rho_T}. Therefore, by
our assumptions on $\pi^{(i-1)}$, we get
\begin{equation}
\nonumber\begin{split}
||P_{s}(\alpha^{(i-1)}_{T_{i}})[i\dot{\tilde{\pi}}^{*}\partial_\pi
\tilde
W_{T_{i}}(\pi^{(i-1)})]||_{L_{t}^{1}L_{x}^{2}([0,T_{i}])}&\lesssim
\delta C [\sup_{t\in [0,T_{i}]}\la
t\ra^{2+\epsilon}|\dot{\tilde{\pi}}^{*}|]\int_{0}^{T_{i}}\la
t\ra^{-2-\epsilon}dt \\
&\lesssim C\delta
M(A,C,C_{0},\delta)
\end{split}\end{equation}
We have used the following simple
\begin{lemma} Under the assumptions of Proposition~\ref{(U,pi)},
we have
\begin{equation}\nonumber
|\partial_\pi \tilde W_{T_{i}}(\pi^{(i-1)})-\partial_\pi \tilde
W_{T_{i}}(\pi^{(i-1)}_{T_{i}})|(t)\lesssim C_{\epsilon}\,\delta\la t\ra
^{-\epsilon}
\end{equation}
\end{lemma}
\begin{proof} This follows from the definition of $\partial_\pi \tilde
W_{T_{i}}(\pi^{(i-1)})$ (with $y,\,y_{T}$ replaced by
$y^{(i-1)},\,y^{(i-1)}_{T}$), and the fact that
\begin{equation}\nonumber
|y^{(i-1)}(t)-y^{(i-1)}_{T}(t)|\lesssim C\delta \la
t\ra^{-\epsilon},\,|\rho_T(t,x)| \le C_\eps\, \delta^2(1+|x|)\la
t\ra^{-\eps}
\end{equation}
as follows from Lemma~\ref{lem:rho_T}.
\end{proof}
Next, using the definition of $V^{(i-1)}$ in \eqref{eq:V1}, as
well as the linear dependence of $U_{root}$ on $U_{dis},U_{hyp}$
we get
\begin{equation}
\nonumber\begin{split}
&||V^{(i-1)}(U^{*}_{dis}+U^{*}_{hyp}+U^{*}_{root}(U^{*}_{dis},U^{*}_{hyp}))||_{L_{t}^{1}L_{x}^{2}([0,T_{i}])}\\
&\lesssim \|V^{(i-1)}\la
x\ra^{\theta}\|_{L_{t}^{\infty}L_{x}^{2}([0,T_{i}])}[\|\la
x\ra^{-\theta}U^{*}_{dis}\|_{L_{t}^{1}L_{x}^{\infty}([0,T_{i}])}+
\|\la x\ra^{-\theta}U^{*}_{hyp}\|_{L_{t}^{1}L_{x}^{\infty}([0,T_{i}])}+\|\la x\ra^{-\theta}U^{*}_{root}\|_{L_{t}^{1}L_{x}^{\infty}([0,T_{i}])}]\\
&\lesssim C\delta \|\la t \ra ^{1+\epsilon}\la
x\ra^{-\theta}U^{*}\|_{L_{x}^{\infty}}\int_{0}^{T_{i}}<t>^{-1-\epsilon}dt
\lesssim C\delta M(A,C,C_{0},\delta)\\
\end{split}\end{equation}
We proceed to the last and most complicated term of the
nonlinearity. To simplify notation, denote
\begin{equation}\nonumber
\tilde{U}^{(i-1)}:=M_{T_{i-1}}(\pi^{(i-1)})\calG_{T_{i-1}}(\pi^{(i-1)})\calG_{T_{i-2}}(\pi^{(i-2)})^{-1}M_{T_{i-2}}(\pi^{(i-2)})^{-1}U^{(i-1)}
\end{equation}
Thus we need to estimate
$\|N_{T_{i}}(\tilde{U}^{(i-1)},\pi^{(i-1)})\|_{L_{t}^{1}L_{x}^{2}([0,T_{i}])}$.
Using Lemma~\ref{lem:Nest}, we get
\begin{equation}\nonumber
|N_{T_{i}}(\tilde{U}^{(i-1)},\pi^{(i-1)})|(t,x)\lesssim
|\tilde{U}^{(i-1)}(x,t)|^{2\sigma+1}+|\tilde{U}^{(i-1)}(x,t)|^2\phi^{2\sigma-1}(x+(y^{(i-1)}_{T_{i}}-y^{(i-1)})(t))
\end{equation}
We first estimate the contribution of the local term. We have
\begin{equation}\nonumber\begin{split}
&||\tilde{U}^{(i-1)}(x,t)^2\phi^{2\sigma-1}(x+(y^{(i-1)}_{T_{i}}-y^{(i-1)})(t))||_{L_{t}^{1}L_{x}^{2}([0,T_{i}])}
\\&\lesssim
||\la x\ra^{-2\theta}\tilde{U}^{(i-1)}(x,t)^2||_{L_{t}^{1}L_{x}^{\infty}}||\la x\ra^{2\theta}\phi^{2\sigma-1}(x+(y^{(i-1)}_{T_{i}}-y^{(i-1)})(t))||_{L_{t}^{\infty}L_{x}^{2}}\\
&\lesssim (A+1)C\delta ||\la t\ra^{1+\epsilon}\la
x\ra^{-\theta}U^{(i-1)}||_{L_{t}^{\infty}L_{x}^{\infty}}[\int_{0}^{T_{i}}\la
t\ra ^{-1-\epsilon}t^{-\frac{1}{2}}dt]\lesssim (A+1)C^{2}\delta^{2} \\
\end{split}\end{equation}
We have crudely bounded for $t\in[0,T_{i}]$
\begin{equation}\label{comp}
\sup_{x\in \R}\la x\ra^{-\theta}|\tilde{U}^{(i-1)}|(t,x)\lesssim
(A+1) \sup_{x\in\R}\la x\ra^{-\theta}|U^{(i-1)}|(t,x)
\end{equation}
Finally, we need to estimate
$||(\tilde{U}^{(i-1)})^{2\sigma+1}||_{L_{t}^{1}L_{x}^{2}([0,T_{i}])}$.
This is straightforward, we have
\begin{equation}\nonumber
||(\tilde{U}^{(i-1)})^{2\sigma+1}||_{L_{t}^{1}L_{x}^{2}([0,T_{i}])}
\lesssim
||U^{(i-1)}||_{L_{t}^{1}L_{x}^{\infty}}^{2\sigma}||U^{(i-1)}||_{L_{t}^{\infty}L_{x}^{2}}
\end{equation}
Note that the translations and phase functions distinguishing
$\tilde{U}^{(i-1)}$ from $U^{(i-1)}$ are irrelevant here. One
bounds the preceding by
\begin{equation}\nonumber
\lesssim
C\delta\int_{1}^{T_{i}}t^{-\sigma}||t^{\frac{1}{2}}U^{(i-1)}(t)||_{L_{x}^{\infty}}^{2\sigma}dt
+C\delta
\int_{0}^{1}[||\partial_{x}U^{(i-1)}(t)||_{L_{x}^{2}}^{2\sigma}+||U^{(i-1)}(t)||_{L_{x}^{2}}^{2\sigma}]dt,
\end{equation}
where we have used Sobolev's inequality in the last step. The
expression can be bounded by $\lesssim (C\delta)^{2\sigma+1}$. To
finish the estimation of $U_{dis}$, we still need to handle the
free contribution, i.e.\ $\sup_{t\in
[0,T_{i}]}||e^{it\Hil(\alpha^{(i-1)}_{T_{i}})}U_{dis}(0)||_{L_{x}^{2}}$.
For this we need to carefully keep track of the definition of
$U_{dis}(0)$, which was
\begin{equation}\nonumber\begin{split}
&U_{dis}(0)=P_{s}(\alpha^{(i-1)}_{T_{i}})\Big\{\calG_{T_{i}}(\pi^{(i-1)}(0))\Big[\bm R_{0}\\
\overline{R_{0}}\endm
+h^{(i)}(U^{*}_{dis},U^{*}_{hyp},U^{(i-1)},\pi^{(i-1)},\pi^{(i-2)})
f^{+}(\alpha^{(i-1)}_{T_{i}})\\&\hspace{4cm}+\sum_{j=1}^{4}a_{j}^{(i)}(h^{(i)}(U^{*}_{dis},U^{*}_{hyp},U^{(i-1)},\pi^{(i-1)},\pi^{(i-2)}),\alpha^{(i-1)}_{T_{i}})\tilde{\eta}_{j}(\alpha^{(i-1)}_{T_{i}})\Big]\Big\}
\end{split}\end{equation}
Recall that
$h^{(i)}(U^{*}_{dis},U^{*}_{hyp},U^{(i-1)},\pi^{(i-1)},\pi^{(i-2)})$
depended linearly on
$b^{(i)+}(U^{*}_{dis},U^{*}_{hyp},U^{(i-1)},\pi^{(i-1)},\pi^{(i-2)})$,
which in turn was given by the expression
\begin{equation}\nonumber
-\int_{0}^{T_i}e^{-\gamma(\alpha^{(i-1)}_{T_{i}})t}g^{+}(U^{*},U^{(i-1)},\pi^{*},\pi^{(i-1)},\pi^{(i-2)})(t)\, dt,
\end{equation}
where $g^{+}(...)$ is as in \eqref{Uhypeqn}. But
this part is estimated exactly like above (indeed, we have an extra
exponentially decaying weight), and one winds up with the estimate
\begin{equation}\nonumber
\left|\int_{0}^{T_i}e^{-\gamma(\alpha^{(i-1)}_{T_{i}})t}g^{+}(U^{*},U^{(i-1)},\pi^{*},\pi^{(i-1)},\pi^{(i-2)})(t)\,dt\right|
\lesssim (A+1)(C\delta)^{2}+(C\delta)^{5}+C\delta
M(A,C,C_{0},\delta)
\end{equation}
The fifth power here comes from $2\sigma+1>5$.
Keeping in mind the linear dependence of the $a_{j}^{(i)}$ on
$h^{(i)}$, and choosing $C_{0}$ such that
\begin{equation}\nonumber
\Big\|P_{s}(\alpha^{(i-1)}_{T_{i}})\calG_{T_{i}}(\pi^{(i-1)}(0))\bm R_{0}\\
\overline{R_{0}}\endm \Big\|_{L_{x}^{2}}<\frac{C_{0}}{\Lambda}\delta,
\end{equation}
we get the estimate (for large $\Lambda$)
\begin{equation}\nonumber
||U_{dis}(0)||_{L_{x}^{2}}<\frac{C_{0}}{\Lambda}\delta+
C_{1}[(A+1)(C\delta)^{2}+(C\delta)^{5}+C\delta
M(A,C,C_{0},\delta)]
\end{equation}
for suitable $C_{1}$. Using the approximate unitarity of the
evolution $e^{it\Hil(\alpha^{(i-1)}_{T_{i}})}$ acting
$P_{s}(\alpha^{(i-1)}_{T_{i}})(L^{2})$, we finally obtain
\begin{equation}\nonumber
||U_{dis}||_{L_{t}^{\infty}L_{x}^{2}([0,T_{i}])}\leq
\frac{C_{0}}{\tilde{\Lambda}}\delta+C_{2}[(A+1)(C\delta)^{2}+(C\delta)^{5}+C\delta
M(A,C,C_{0},\delta)]
\end{equation}
for suitable $C_{2}$ (independent of all other constants) and large
$\tilde{\Lambda}$, provided $\Lambda$ was chosen large enough.
\\{\bf{(B):}} Estimating
$||U_{root}||_{L_{t}^{\infty}L_{x}^{2}([0,T_{i}])}+||U_{hyp}||_{L_{t}^{\infty}L_{x}^{2}([0,T_{i}])}$.
To complete the estimate for
$||U||_{L_{t}^{\infty}L_{x}^{2}([0,T_{i}])}$, we still need to
estimate the contributions from the hyperbolic and root part. For
the former, we use \eqref{Uhypeqn}, as well as the condition for
$b^{+}(0)$ and Lemma~\ref{lem:ODE_stable}, which results in
\begin{equation}\nonumber
b^{+}(t)=-\int_{t}^{T_{i}}e^{-\gamma(\alpha^{(i-1)}_{T_{i}})(s-t)}g^{+}(U^{*},U^{(i-1)},\pi^{*},\pi^{(i-1)},\pi^{(i-2)})(t)\,dt
\end{equation}
\begin{equation}\nonumber
b^{-}(t)=e^{-t\gamma(\alpha^{(i-1)}_{T_{i}})}b^{-}(0)+
\int_{0}^{t}e^{\gamma(\alpha^{(i-1)}_{T_{i}})(s-t)}g^{-}(U^{*},U^{(i-1)},\pi^{*},\pi^{(i-1)},\pi^{(i-2)})(t)\,dt
\end{equation}
for all $t\in[0,T_{i}]$. One can then bound this by the same kind
of expression as
$\|U_{dis}\|_{L_{t}^{\infty}L_{x}^{2}([0,T_{i}])}$. Finally, the
fact that $U_{root}$ depends linearly on $U_{dis}$, $U_{hyp}$
implies the same kind of bound for it. This completes estimating
the $\|.\|_{L_{t}^{\infty}L_{x}^{2}}$-contribution.\\
{\bf{(C):}} The contribution of
$\|\partial_{x}U\|_{L_{t}^{\infty}L_{x}^{2}}$. Using
Corollary~\ref{cor:weightp}, as well as the Duhamel parametrix, we
see that we need to estimate
\begin{equation}\nonumber
\|\partial_{x}[\text{right-hand side of \eqref{Ueqn} without the
$P_{s}$}]\|_{L_{t}^{1}L_{x}^{2}([0,T_{i}])}
\end{equation}
However, this follows from almost identical estimates. One simply
substitutes $\|\la t\ra ^{1+\epsilon}\la
x\ra^{-\frac{1}{2}-2\epsilon}\partial_{x}U\|_{L_{t}^{\infty}L_{x}^{q}}$
where before we used $\|\la t\ra
^{1+\epsilon}\la x\ra^{-\theta}U\|_{L_{t}^{\infty}L_{x}^{\infty}[0,T_{i}]}$.\\
{\bf{(D):}} The contribution of $\sup_{0\leq\theta\leq 1}\la t\ra
^{-\theta}\|\la x\ra^{\theta}U\|_{L_{x}^{2}}$. This is treated
just like the preceding cases, using the Duhamel formula in
addition to Lemma~\ref{lem:moments} (more precisely, an
interpolate of this and the approximate $L^{2}$ conservation). The
details are very similar to previous calculations, and we skip
them.\\
{\bf{(E):}} Estimating the weighted norm $||\la t\ra^{1+\epsilon}
\la
x\ra^{-\theta}U_{dis}||_{L_{t}^{\infty}L_{x}^{\infty}([0,T_{i}])}$.
Using Duhamel's formula as well as Lemma~\ref{lem:theta}, we see
that we have\footnote{The $s$ in $P_{s}(\alpha^{(i-1)}_{T_{i}})$
stands for 'stable' and has nothing to do with the
integration variable $s$.}
\begin{equation}\nonumber\begin{split}
&\|\la x\ra^{-\theta}U_{dis}(t)\|_{L_{x}^{\infty}} \lesssim \|\la
x\ra^{-\theta}e^{i t
\Hil(\alpha^{(i-1)}_{T_{i}})}U_{dis}(0)\|_{L_{x}^{\infty}}\\&+
\int_{0}^{t-1}\la t-s\ra^{-1-\epsilon}\big\|\la
x\ra^{\theta}P_{s}(\alpha^{(i-1)}_{T_{i}})[i\dot{\tilde{\pi}}^{*}\partial_\pi
\tilde W_{T_{i}}(\pi^{(i-1)}) +
V^{(i-1)}(U^{*}_{dis}+U^{*}_{hyp}+U^{*}_{root}(U^{*}_{dis},U^{*}_{hyp}))\\
&+N_{T_{i}}(M_{T_{i-1}}(\pi^{(i-1)})\calG_{T_{i-1}}(\pi^{(i-1)})\calG_{T_{i-2}}(\pi^{(i-2)})^{-1}M_{T_{i-2}}(\pi^{(i-2)})^{-1}U^{(i-1)},\pi^{(i-1)})](s)\big\|_{L_{x}^{1}}\,ds\\
&+\int_{t-1}^{t}(t-s)^{-\frac{1}{2}}\|P_{s}(\alpha^{(i-1)}_{T_{i}})[i\dot{\tilde{\pi}}^{*}\partial_\pi
\tilde W_{T_{i}}(\pi^{(i-1)}) +
V^{(i-1)}(U^{*}_{dis}+U^{*}_{hyp}+U^{*}_{root}(U^{*}_{dis},U^{*}_{hyp}))\\
&+N_{T_{i}}(M_{T_{i-1}}(\pi^{(i-1)})\calG_{T_{i-1}}(\pi^{(i-1)})\calG_{T_{i-2}}(\pi^{(i-2)})^{-1}M_{T_{i-2}}(\pi^{(i-2)})^{-1}U^{(i-1)},\pi^{(i-1)})](s)\|_{L_{x}^{1}}\,ds\\
\end{split}
\end{equation}
For the last integral expression, we have rather crudely discarded
the weight $\la x\ra^{-\theta}$, since this term turns out to be
very small. We proceed as in case {\bf{(A)}} by treating all the
different terms. First, consider the first integral expression on
the right. We estimate for $t\in [0,T_{i}]$
\begin{equation}\nonumber\begin{split}
&\int_{0}^{t-1}\la t-s\ra^{-1-\epsilon}\|\la
x\ra^{\theta}P_{s}(\alpha^{(i-1)}_{T_{i}})[i\dot{\tilde{\pi}}^{*}\partial_\pi
\tilde W_{T_{i}}(\pi^{(i-1)})](s)\|_{L_{x}^{1}}ds\\
&=\int_{0}^{t-1}\la t-s\ra^{-1-\epsilon}\|\la
x\ra^{\theta}P_{s}(\alpha^{(i-1)}_{T_{i}})[i\dot{\tilde{\pi}}^{*}\partial_\pi
\tilde W_{T_{i}}(\pi^{(i-1)})-i\dot{\tilde{\pi}}^{*}\partial_\pi
\tilde W_{T_{i}}(\pi^{(i-1)}_{T_{i}})](s)\|_{L_{x}^{1}}ds\\
&\lesssim C\delta \int_{0}^{t-1}\la t-s\ra^{-1-\epsilon}\la
s\ra^{-2-\epsilon}ds\lesssim C\delta M(A,C,C_{0},\delta)\la
t\ra^{-1-\epsilon}\\
\end{split}\end{equation}
This is acceptable in light of the definition of
$\|.\|_{X_{*}([0,T_{i}])}$. We have used the fact that we may
safely move the weight $\la x\ra^{\theta}$ past
$P_{s}(\alpha^{(i-1)}_{T_{i}})$ and then discard it, due to the
local nature of the expression. Next, we have
\begin{equation}\nonumber\begin{split}
&\int_{0}^{t-1}\la t-s\ra^{-1-\epsilon}\|\la
x\ra^{\theta}P_{s}(\alpha^{(i-1)}_{T_{i}})[V^{(i-1)}(U^{*}_{dis}+U^{*}_{hyp}+U^{*}_{root}(U^{*}_{dis},U^{*}_{hyp}))](s)\|_{L_{x}^{1}}\,ds\\
&\lesssim \int_{0}^{t-1}\la t-s\ra^{-1-\epsilon}\, \|\la
x\ra^{2\theta}V^{(i-1)}(s)\|_{L_{x}^{1}}\|\la
x\ra^{-\theta}(U^{*}_{dis}+U^{*}_{hyp}+U^{*}_{root}(U^{*}_{dis},U^{*}_{hyp})(s)\|_{L_{x}^{\infty}}\,ds\\
&\lesssim C\delta M(A,C,C_{0},\delta) \int_{0}^{t-1}\la
t-s\ra^{-1-\epsilon}\la s\ra^{-1-\epsilon}ds\lesssim \la
t\ra^{-1-\epsilon}C\delta M(A,C,C_{0},\delta),\\
\end{split}\end{equation}
as desired. Next, we consider
\begin{equation}\nonumber
\int_{0}^{t-1}\la t-s\ra^{-1-\epsilon}\|\la
x\ra^{\theta}P_{s}(\alpha^{(i-1)}_{T_{i}})[N_{T_{i}}(\tilde{U}^{(i-1)},\pi^{(i-1)})](s)\|_{L_{x}^{1}}\,ds
\end{equation}
As before, we move the weight past the operator
$P_{s}(\alpha^{(i-1)}_{T_{i}})$ and discard the latter, obtaining
the expression
\begin{equation}\nonumber
\int_{0}^{t-1}\la t-s\ra^{-1-\epsilon}\|\la
x\ra^{\theta}[N_{T_{i}}(\tilde{U}^{(i-1)},\pi^{(i-1)})](s)\|_{L_{x}^{1}}\,ds
\end{equation}
We use Lemma~\ref{lem:Nest}, which implies
\begin{equation}\nonumber\begin{split}
&\int_{0}^{t-1}\la t-s\ra^{-1-\epsilon}\|\la
x\ra^{\theta}N_{T_{i}}(\tilde{U}^{(i-1)},\pi^{(i-1)})(s)\|_{L_{x}^{1}}ds\\
&\lesssim
\int_{0}^{t-1}\la
t-s\ra^{-1-\epsilon}\|\la x\ra^{\theta}[|\tilde{U}^{(i-1)}(x,s)|^{2\sigma+1}+|\tilde{U}^{(i-1)}(x,s)|^2\phi^{2\sigma-1}(x+(y^{(i-1)}_{T_{i}}-y^{(i-1)})(s))]\|_{L_{x}^{1}}\, ds
\end{split}\end{equation}
We treat each of the two summands separately. First, consider the
local term. The weight is simply absorbed here, i.e., we can estimate
\begin{equation}\nonumber\begin{split}
&\int_{0}^{t-1}\la
t-s\ra^{-1-\epsilon}\|\la x\ra^{\theta}|\tilde{U}^{(i-1)}(x,s)|^2\phi^{2\sigma-1}(x+(y^{(i-1)}_{T_{i}}-y^{(i-1)})(s)\|_{L_{x}^{1}}\,ds\\
&\lesssim\int_{0}^{t-1}\la t-s\ra^{-1-\epsilon}\|\la
x\ra^{-\theta}\tilde{U}^{(i-1)}(x,s)\|_{L_{x}^{\infty}}\|\tilde{U}^{(i-1)}(x,s)\|_{L_{x}^{\infty}}\|\la
x\ra^{2\theta}\phi^{2\sigma-1}(x+(y^{(i-1)}_{T_{i}}-y^{(i-1)})(s))\|_{L_{x}^{1}}\,ds
\end{split}\end{equation}
Again exploiting \eqref{comp} (actually, we only need to pay $A$
here), we can bound this integral by
\begin{equation}\nonumber
\lesssim A(C\delta)^{2}\int_{0}^{t-1}\la t-s\ra^{-1-\epsilon}\la
s\ra ^{-\frac{3}{2}}ds\lesssim A(C\delta)^{2}\la
t\ra^{-1-\epsilon}
\end{equation}
Next, consider the nonlocal term. Here, we can no longer absorb
the weights, and therefore need to carefully keep track of the
exact powers. For this purpose,  we assume (as we may) that
$0<\epsilon<\sigma-1-2\theta$. Using Lemma~\ref{lem:moments}, we get
\begin{equation}\nonumber\begin{split}
&\int_{0}^{t-1}\la t-s\ra^{-1-\epsilon}\|\la
x\ra^{\theta}\tilde{U}^{(i-1)}(x,s)|^{2\sigma+1}\|_{L_{x}^{1}}\,ds \\
&\lesssim\int_{0}^{t-1}\la t-s\ra^{-1-\epsilon}\|\la
x\ra^{\theta}
\tilde{U}^{(i-1)}(x,s)\|_{L_{x}^{2}}\|\tilde{U}^{(i-1)}(x,s)\|_{L_{x}^{2}}\|\tilde{U}^{(i-1)}(x,s)\|_{L_{x}^{\infty}}^{2\sigma-1}\,ds\\
&\lesssim A(C\delta)^{2\sigma+1}\int_{0}^{t-1}\la
t-s\ra^{-1-\epsilon}s^{-\sigma+\frac{1}{2}+\theta}\,ds
\lesssim A(C\delta)^{2\sigma+1}\la t\ra^{-1-\epsilon}
\end{split}\end{equation}
We have crudely bounded
\begin{equation}\nonumber
\sup_{x\in\R}|\la x\ra^{\theta} \tilde{U}^{(i-1)}(x,s)|\lesssim A
\sup_{x\in\R}|\la x\ra U^{(i-1)}(x,s)|
\end{equation}
We proceed to estimating the fringe integral over the interval
$[t-1,t]$. The estimate here is even simpler: we get
\begin{equation}\nonumber
\int_{t-1}^{t}(t-s)^{-\frac{1}{2}}\big\|P_{s}(\alpha^{(i-1)}_{T_{i}})\big[i\dot{\tilde{\pi}}^{*}\partial_\pi
\tilde W_{T_{i}}(\pi^{(i-1)})(s)\big]\big\|_{L_{x}^{1}}\, ds\lesssim C\delta
\la t\ra^{-2-\epsilon}M(A,C,C_{0},\delta)
\end{equation}
Next, we have
\begin{equation}\nonumber
\int_{t-1}^{t}(t-s)^{-\frac{1}{2}}\|P_{s}(\alpha^{(i-1)}_{T_{i}})[
V^{(i-1)}(U^{*}_{dis}+U^{*}_{hyp}+U^{*}_{root}(U^{*}_{dis},U^{*}_{hyp}))(s)]\|_{L_{x}^{1}}ds
\lesssim C\delta M(A,C,C_{0},\delta) \la t\ra^{-1-\epsilon}
\end{equation}
Of course we exploit here that
\begin{equation}\nonumber
\|\la x\ra^{-\theta}U^{*}\|_{L_{x}^{\infty}}\leq \la
t\ra^{-1-\epsilon}M(A,C,C_{0},\delta)
\end{equation}
The remaining terms are more of the same and omitted\footnote{One
has to be  careful here, since the paths $\pi^{(i-1)},\pi^{(i-2)}$
diverge a bit more past time $T_{i-1}=T_{i}-1$. Thus one needs to
replace $A$ by $A+1$.}. We still need to estimate the free
contribution. First, note that for $0\leq t\leq 1$ we have by
Sobolev's inequality as well as the approximate unitarity of the free
evolution, see Lemma~\ref{cor:inter}
\begin{equation}\nonumber
\|\la
x\ra^{-\theta}e^{it\Hil(\alpha^{(i-1)}_{T_{i}})}U_{dis}(0)\|_{L_{x}^{\infty}}\lesssim
\|e^{it\Hil(\alpha^{(i-1)}_{T_{i}})}U_{dis}(0)\|_{H^{1}}\lesssim
\|U(0)\|_{H^{1}}
\end{equation}
The latter is estimated as in parts {\bf{(A)}}, {\bf{(C)}}. Next,
assuming $t\geq 1$, we have (using Lemma~\ref{lem:theta})
\begin{equation}\nonumber
\|\la
x\ra^{-\theta}e^{it\Hil(\alpha^{(i-1)}_{T_{i}})}U_{dis}(0)\|_{L_{x}^{\infty}}\lesssim
\la t\ra^{-1-\epsilon}\|\la x\ra^{\theta}U_{dis}(0)\|_{L_{x}^{1}}
\end{equation}
Using that
\begin{equation}\nonumber\begin{split}
&U_{dis}(0)=P_{s}(\alpha^{(i-1)}_{T_{i}})\Big\{\calG_{T_{i}}(\pi^{(i-1)}(0))\Big[\bm R_{0}\\
\overline{R_{0}}\endm
+h^{(i)}(U^{*}_{dis},U^{*}_{hyp},U^{(i-1)},\pi^{(i-1)},\pi^{(i-2)})
f^{+}(\alpha^{(i-1)}_{T_{i}})\\&\hspace{4cm}+\sum_{j=1}^{4}a_{j}^{(i)}(h^{(i)}(U^{*}_{dis},U^{*}_{hyp},U^{(i-1)},\pi^{(i-1)},\pi^{(i-2)}),\alpha^{(i-1)}_{T_{i}})\tilde{\eta}_{j}(\alpha^{(i-1)}_{T_{i}})\Big] \Big\}
\end{split}\end{equation}
This shows that first, the following term needs to be estimated:
\begin{equation}\nonumber
\Big\|\la x\ra^{\theta}P_{s}(\alpha^{(i-1)}_{T_{i}})\calG_{T_{i}}(\pi^{(i-1)}(0))\bm R_{0}\\
\overline{R_{0}}\endm \Big\|_{L_{x}^{1}} \lesssim
\Big\|\la x\ra^{\theta}\calG_{T_{i}}(\pi^{(i-1)}(0))\bm R_{0}\\
\overline{R_{0}}\endm\Big\|_{L_{x}^{1}}\lesssim \Big\|\la x\ra^\theta \bm R_{0}\\
\overline{R_{0}}\endm\Big\|_{L_{x}^{1}}
\end{equation}
which is majorized by $\frac{C_{0}}{\Lambda}\delta$ for $C_{0}$
chosen sufficiently large and large $\Lambda$. The remaining terms
constituting $U_{dis}(0)$, being local, are estimated similarly,
in light of the earlier comments on $h^{(i)}(...)$ etc. More
precisely, one obtains
\begin{equation}\nonumber
\|\la x\ra^\theta U_{dis}(0)\|_{L_{x}^{1}}\leq
\frac{C_{0}}{\Lambda}\delta +C\delta M(A,C,C_{0},\delta)+
C_{3}[(A+1)(C\delta)^{2}+(C\delta)^{2\sigma+1}]
\end{equation}
for an absolute constant $C_{3}$. We are done with {\bf{(E)}}. \\
{\bf{(F):}} The estimate for $\sup_{0\leq t\leq T_{i}}\la
t\ra^{-1-\epsilon}\|\la x\ra^{-\theta} U_{hyp}\|_{L_{x}^{\infty}}$
and $\sup_{0\leq t\leq T_{i}}\la t\ra^{-1-\epsilon}\|\la
x\ra^{-\theta} U_{root}\|_{L_{x}^{\infty}}$. The estimate follows
once we establish it for $U_{hyp}$, on account of the linear
dependence of $U_{root}$ on $U_{dis}$, $U_{hyp}$. To see it for
$U_{hyp}$, use (recall the terminology from the beginning of the
proof of Proposition~\ref{(U,pi)})
\begin{equation}\nonumber
b^{+}(t)=-\int_{t}^{T_{i}}e^{-\gamma(\alpha^{(i-1)}_{T_{i}})(s-t)}g^{+}(U^{(*)},U^{(i-1)},\pi^{(i)},\pi^{(i-1)},\pi^{(i-2)})(t)\,dt
\end{equation}
\begin{equation}\nonumber
b^{-}(t)=e^{-t\gamma(\alpha^{(i-1)}_{T_{i}})}b^{-}(0)+
\int_{0}^{t}e^{\gamma(\alpha^{(i-1)}_{T_{i}})(s-t)}g^{-}(U^{*},U^{(i-1)},\pi^{*},\pi^{(i-1)},\pi^{(i-2)})(t)dt
\end{equation}
The exponentially decaying weight accounts for the integrability
of the integrands. But then the desired decay rate of $\la
t\ra^{-1-\epsilon}$ follows easily from the preceding estimates
for $g^{\pm}(...)$.\\
{\bf{(G):}} The estimate for $\sup_{0\leq t\leq T_{i}}\la
t\ra^{-1-\epsilon}\|\la
x\ra^{-\frac{1}{2}-2\epsilon}\partial_{x}U_{dis}(t)\|_{L_{x}^{q}}$.
This is handled by using Lemma~\ref{cor:weightp}, in addition to
the Duhamel's formula. One reduces to estimating the
differentiated nonlinearity, which is handled just as before,
using $\la t\ra^{-1-\epsilon}\|\la
x\ra^{-\frac{1}{2}-2\epsilon}\partial_{x}U(t)\|_{L_{x}^{q}}$
instead of $\la t\ra^{-1-\epsilon}\|\la
x\ra^{-\theta}\partial_{x}U(t)\|_{L_{x}^{\infty}}$ in some
places.\\
{\bf{(H):}} Estimating $\sup_{0\leq t\leq  T_{i}}\la
t\ra^{2+\epsilon}|\dot{\tilde{\pi}}|(t)$. For this we of course use
\eqref{path}. The following terms need to be estimated:
\begin{equation}\nonumber
|\la
U^{*}(t),i\dot{\tilde{\pi}}^{(i-1)}\tilde{\calS_{\ell}}(.+y^{(i-1)}_{T_{i}}-y^{(i-1)})(t)\ra|
\lesssim
\|U^{*}(t)\|_{L_{x}^{2}}\|\tilde{\calS_{\ell}}(.+y^{(i-1)}_{T_{i}}-y^{(i-1)})(t)\|_{L_{x}^{2}}
|\dot{\tilde{\pi}}^{(i-1)}(t)|\lesssim C\delta \la t\ra
^{-2-\epsilon}
\end{equation}
Of course this estimate is rather crude, but it suffices for our
purposes. Next, consider the expression $\la
N(\tilde{U}^{(i-1)}(t),\pi^{(i-1)}(t)),\tilde{\xi}_{\ell}(\pi^{(i-1)}(t))\ra$.
Using Lemma~\ref{lem:Nest}, we reduce this to the following two
estimates. It suffices to consider $T_{i}\geq t\geq 1$.
\begin{align*}
&\Big|\la
|\tilde{U}^{(i-1)}(x,t)|^2\phi^{2\sigma-1}(x+(y^{(i-1)}_{T_{i}}-y^{(i-1)})(t)),\tilde{\xi}_{\ell}(\pi^{(i-1)}(t))\ra \Big|\\
&\lesssim (A+1)^{2}(\la t\ra^{1+\epsilon}\|\la
x\ra^{-\theta}U^{(i-1)}\|_{L_{x}^{\infty}}^{2}\la t\ra
^{-2-2\epsilon}\lesssim ((A+1)C\delta)^{2}\la
t\ra^{-2-\epsilon}\\
\Big|\la |\tilde{U}^{(i-1)}(x,t)|^{2\sigma+1},\tilde{\xi}_{\ell}(\pi^{(i-1)}(t))\ra\Big| &\lesssim
t^{-\sigma-\frac{1}{2}}[t^{\frac{1}{2}}\|\tilde{U}^{(i-1)}(x,t)\|_{L_{x}^{\infty}}]^{2\sigma+1}
\|\tilde{\xi}_{\ell}(\pi^{(i-1)}(t)\|_{L_{x}^{1}}\\
& \lesssim \la t\ra^{-\sigma-\frac{1}{2}}(C\delta)^{2\sigma+1}
\end{align*}
Finally we have controlled all the components constituting
$\|(\pi,U)\|_{X_*([0,T_{i}])}$. Gathering the preceding estimates
from {\bf{(A)}}-{\bf{(H)}}, we get
\begin{equation}\nonumber
 \|(\pi,U)\|_{X_*([0,T_{i}])}\leq \frac{C_{0}}{\tilde{\Lambda}}\delta +C\delta M(A,C,C_{0},\delta)+
\tilde{C}_{3}[(A+1)^{2}(C\delta)^{2}+(C\delta)^{2\sigma+1}]
\end{equation}
for suitable large $\tilde{\Lambda}$ (if we choose $C_{0}$ large
enough) and $\tilde{C}_{3}$ an absolute constant. We conclude that
if $C$ is sufficiently large in relation to $C_{0},\tilde{C}_{3}$,
and $\delta$ sufficiently small in relation to $C$, as well as
$C_{0}$ sufficiently large in relation to $\tilde{C}_{3}$, we get
\begin{equation}\nonumber
\|(\pi,U)\|_{X_*([0,T_{i}])}\leq M(A,C,C_{0},\delta),
\end{equation}
which establishes the a priori inequality we need. Since the
contraction step follows along the same lines, the map
$F_{T_{i}}(U^{(i-1)},\pi^{(i-1)},\pi^{(i-2)})$ has a fixed point.
This is the next iterate $(\pi^{(i)},U^{(i)})$. This establishes
Proposition~\ref{(U,pi)}.
\end{proof}

\begin{proof}[Proof of Proposition~\ref{pi-pi}] We need to analyze the
difference equation at the level of the $Z$, which will be
accomplished by transforming $Z^{(i)}-Z^{(i-1)}$ into a suitable
gauge, similarly to the preceding (note, however, that we need to
consider the difference of the $Z^{(i)}$, since it is easy to see that
the difference of the $U^{(i)}$ cannot be controlled in a reasonable way).
Starting from \eqref{eq:Zsyszer} etc.\ we
arrive at the following equations, valid on $[0,T_{i-1}]$:
\begin{equation}\label{diff2}\begin{split}
&i\partial_{t}(Z^{(i)}-Z^{(i-1)})
-\Hil(\pi^{(i-1)})(Z^{(i)}-Z^{(i-1)})=[i\dot{\tilde{\pi}}^{(i)}\partial_{\pi}W(\pi^{(i-1)})-i\dot{\tilde{\pi}}^{(i-1)}\partial_{\pi}W(\pi^{(i-2)})]
\\&\hspace{2cm}+[N(Z^{(i-1)},\pi^{(i-1)})-N(Z^{(i-2)},\pi^{(i-2)})] +
[\Hil(\pi^{(i-2)})-\Hil(\pi^{(i-1)})]Z^{(i-1)}\\
\end{split}\end{equation}
\begin{equation}\nn
\begin{split}
&i[\dot{\tilde{\pi}}^{(i)}-\dot{\tilde{\pi}}^{(i-1)}]=\la
Z^{(i)}(t),i\dot{\tilde{\pi}}^{(i-1)}\tilde{\calS}_{\ell}(\pi^{(i-1)})(t)\ra
-\la N(Z^{(i-1)},\pi^{(i-1)}),\tilde{\xi}_{\ell}(\pi^{(i-1)})(t)\ra-
\\&(\la
Z^{(i-1)}(t),i\dot{\tilde{\pi}}^{(i-2)}\tilde{\calS}_{\ell}(\pi^{(i-2)})(t)\ra
-\la
N(Z^{(i-2)},\pi^{(i-2)}),\tilde{\xi}_{\ell}(\pi^{(i-2)})(t)\ra)\\
\end{split}\end{equation}
We commence by analyzing the first equation. Introduce the gauged
quantity
\begin{equation}\nonumber
U^{(i,i-1)}:=M_{T_{i}}(\pi^{(i-1)})\calG_{T_{i}}(\pi^{(i-1)})(Z^{(i)}-Z^{(i-1)}),
\end{equation}
It satisfies the following equation:
\begin{equation}\label{gauged difference eqn.}\begin{split}
&[i\partial_{t}-\Hil(\alpha_{T_{i}}^{(i-1)})]U^{(i,i-1)}=
M_{T_{i}}(\pi^{(i-1)})\calG_{T_{i}}(\pi^{(i-1)})\Big\{[i\dot{\tilde{\pi}}^{(i)}\partial_{\pi}W(\pi^{(i-1)})-i\dot{\tilde{\pi}}^{(i-1)}\partial_{\pi}W(\pi^{(i-2)})]
\\&+[N(Z^{(i-1)},\pi^{(i-1)})-N(Z^{(i-2)},\pi^{(i-2)})]+
[\Hil(\pi^{(i-2)})-\Hil(\pi^{(i-1)})]Z^{(i-1)}\\&+[\Hil(\pi^{(i-1)}(t)-\Hil(\pi^{(i-1)}_{T_{i}}(t))](Z^{(i)}-Z^{(i-1)}) \Big\}
\end{split}\end{equation}
Introduce the norm $\|(\cdot,\cdot)\|_{\tilde{Y}([0,T])}$ like
$\|(\cdot,\cdot)\|_{Y^{(i)}([0,T])}$ but with the weight $\la
x-y^{(i-1)}(t)\ra$ replaced by $\la x\ra$. We estimate
the various constituents of this norm. \\
{\bf{(A):}} The contribution from $\|\la
t\ra^{-1}U^{(i,i-1)}_{dis}\|_{L_{t}^{\infty}L_{x}^{2}([0,T_{i-1}])}$.
For $t\in [0,T_{i-1}]$, we shall use
\begin{equation}\nonumber
\|\la t\ra^{-1}U^{(i,i-1)}_{dis}(t)\|_{L_{x}^{2}}\lesssim \la
t\ra^{-1}\|U^{(i,i-1)}_{dis}(0)\|_{L_{x}^{2}}+\la
t\ra^{-1}\|P_{s}(\alpha^{(i-1)}_{T_{i}}[\text{left-hand side of
\eqref{diff2}}]\|_{L_{t}^{1}L_{x}^{2}([0,t])}
\end{equation}
we need to estimate the following terms:
\begin{equation}\nonumber\begin{split}
&i\dot{\tilde{\pi}}^{(i)}\partial_{\pi}W(\pi^{(i-1)})-i\dot{\tilde{\pi}}^{(i-1)}\partial_{\pi}W(\pi^{(i-2)})
\\&=i(\dot{\tilde{\pi}}^{(i)}-\dot{\tilde{\pi}}^{(i-1)})\partial_{\pi}W(\pi^{(i-1)})+\dot{\tilde{\pi}}^{(i-1)}(\partial_{\pi}W(\pi^{(i-1)})-
\partial_{\pi}W(\pi^{(i-2)}))\\
\end{split}\end{equation}
As for the a priori estimates, we can estimate
\begin{equation}\nonumber\begin{split}
&\|i(\dot{\tilde{\pi}}^{(i)}-\dot{\tilde{\pi}}^{(i-1)})P_{s}[M_{T_{i}}(\pi^{(i-1)})\calG_{T_{i}}(\pi^{(i-1)})\partial_{\pi}W(\pi^{(i-1)})]\|_{L_{t}^{1}L_{x}^{2}([0,t])}
\\&=\|i(\dot{\tilde{\pi}}^{(i)}-\dot{\tilde{\pi}}^{(i-1)})P_{s}[M_{T_{i}}(\pi^{(i-1)})\calG_{T_{i}}(\pi^{(i-1)})(\partial_{\pi}W(\pi^{(i-1)})-\partial_{\pi}W(\pi^{(i-1)}_{T_{i}}))]\|_{L_{t}^{1}L_{x}^{2}([0,t])}\\
&\lesssim C\delta \sup_{t\in [0,T_{i-1}]}\la
t\ra^{1+\frac{\epsilon}{2}} |(\dot{\tilde{\pi}}^{(i)}-\dot{\tilde{\pi}}^{(i-1)})(t)|\;\int_{0}^{T_{i-1}}\la
t\ra^{-1-\frac{\epsilon}{2}}\,dt\\
\end{split}\end{equation}
Moreover, we have for $t\in[0,T_{i-1}]$
\begin{equation}\nonumber\begin{split}
&\|P_{s}[M_{T_{i}}(\pi^{(i-1)})\calG_{T_{i}}(\pi^{(i-1)})\dot{\tilde{\pi}}^{(i-1)}(\partial_{\pi}W(\pi^{(i-1)})-
\partial_{\pi}W(\pi^{(i-2)}))]\|_{L_{t}^{1}L_{x}^{2}[0,t]}
\\&\lesssim
(\sup_{s\in[0,T_{i-1}]}\la s^{2+\epsilon}\ra|\dot{\pi}^{(i-1)}|(s))\int_{0}^{t}\la s\ra^{-1-\frac{\epsilon}{2}}ds \sup_{s\in[0,T_{i-1}]}|\tilde{\pi}^{(i-1)}-\tilde{\pi}^{(i-2)}|(s)\\
&\lesssim C\delta
\sup_{s\in[0,T_{i-1}]}\la s\ra^{1+\frac{\epsilon}{2}}|\dot{\tilde{\pi}}^{(i-1)}-\dot{\tilde{\pi}}^{(i-2)}|(s)\\
\end{split}\end{equation}
We have used the following simple
\begin{lemma}\label{tedious} The following inequality holds
\begin{equation}\nonumber
\sup_{t\in[0,T_{i-1}]}|\partial_{\pi}W(\pi^{(i-1)})(t)-\partial_{\pi}W(\pi^{(i-2)})(t)|\lesssim
\sup_{t\in[0,T_{i-1}]}\la t\ra^{1+\frac{\epsilon}{2}}
|\tilde{\pi}^{(i-1)}(t)-\tilde{\pi}^{(i-2)}(t)|
\end{equation}
\end{lemma}
\begin{proof} Note that $\partial_{\pi}W(\pi^{(i-1)})$ consists
precisely of the generalized root functions, translated by
$y^{(i-1)}$ and twisted by a phase $e^{i\theta^{(i-1)}}$, and
similarly for $\partial_{\pi}W(\pi^{(i-2)})$. But we have for
$t\in[0,T_{i-1}]$
\begin{equation}\nonumber\begin{split}
&|y^{(i-1)}(t)-y^{(i-2)}(t)|\leq
\int_{0}^{t}|v^{(i-1)}(s)-v^{(i-2)}(s)|ds+
|D^{(i-1)}(t)-D^{(i-2)}(t)|\\&\lesssim \sup_{s\in[0,T_{i-1}]}\la
s\ra^{1+\frac{\epsilon}{2}}|v^{(i-1)}(s)-v^{(i-2)}(s)| +
|D^{(i-1)}(t)-D^{(i-2)}(t)|\\
 \end{split}\end{equation}
and similarly, we have
\begin{equation}\nonumber\begin{split}
|\theta^{(i-1)}(t)-\theta^{(i-2)}(t)|\leq
&|v^{(i-1)}(t)-v^{(i-2)}(t)|\,|x|+\left|\int_{0}^{t}((v^{(i-1)})^{2}(s)-(v^{(i-2)}(s))^{2})\,ds\right|\\&
+\left|\int_{0}^{t}((\alpha^{(i-1)})^{2}(s)-(\alpha^{(i-2)}(s))^{2})\,ds\right|+|\gamma^{(i-1)}(t)-\gamma^{(i-2)}(t)|\\
&\lesssim
|v^{(i-1)}(t)-v^{(i-2)}(t)|\,|x|+|\gamma^{(i-1)}(t)-\gamma^{(i-2)}(t)|\\&+
 \sup_{s\in[0,T_{i-1}]}\la
s\ra^{1+\frac{\epsilon}{2}}[|v^{(i-1)}(s)-v^{(i-2)}(s)|+|\alpha^{(i-1)}(s)-\alpha^{(i-2)}(s)|]
\\&+|\gamma^{(i-1)}(t)-\gamma^{(i-2)}(t)|\\
\end{split}\end{equation}
Of course the factor $x$ gets absorbed because of the local nature
of the term. Finally note that
\begin{equation}\nonumber
|\gamma^{(i-1)}(t)-\gamma^{(i-2)}(t)|\lesssim \sup_{s\in[0,t]}\la
s\ra^{1+\frac{\epsilon}{2}}|\tilde{\pi}^{(i-1)}-\tilde{\pi}^{(i-2)}|(s),
\end{equation}
and the lemma follows.
\end{proof}
In order to estimate the second term in the nonlinearity, we use that
\begin{equation}\nonumber\begin{split}
[N(Z^{(i-1)},\pi^{(i-1)})-N(Z^{(i-2)},\pi^{(i-2)})]&=[Z^{(i-1)}-Z^{(i-2)}]\int_{0}^{1}\partial_{Z}N(s[Z^{(i-1)}-Z^{(i-2)}]+Z^{(i-2)},\pi^{(i-1)})\,ds
\\&+[\pi^{(i-1)}-\pi^{(i-2)}]\int_{0}^{1}\partial_{\pi}N(Z^{(i-2)},\pi^{(i-2)}+s(\pi^{(i-1)}-\pi^{(i-2)}))\,ds\\
\end{split}\end{equation}
Now we have for $t\in [0,T_{i-1}]$
\begin{equation}\nonumber\begin{split}
&\la t\ra^{-1}\Big\|P_{s}[M_{T_{i}}(\pi^{(i-1)})\calG_{T_{i}}(\pi^{(i-1)})[Z^{(i-1)}-Z^{(i-2)}]\\&\hspace{5cm}\int_{0}^{1}\partial_{Z}N(s[Z^{(i-1)}-Z^{(i-2)}]+Z^{(i-2)},\pi^{(i-1)})ds]\Big\|_{L_{t}^{1}L_{x}^{2}[0,t]}\\
&\lesssim
[AC\delta+(C\delta)^{2\sigma}]\sup_{s\in[0,T_{i-1}]}\la s\ra^{-1}\|Z^{(i-1)}-Z^{(i-2)}\|_{L_{x}^{2}}(s)
\end{split}\end{equation}
We have exploited that
\begin{equation}\nonumber\begin{split}
&|\partial_{Z}N(s[Z^{(i-1)}-Z^{(i-2)}]+Z^{(i-2)},\pi^{(i-1)})|\\&\hspace{2cm}\lesssim
|s[Z^{(i-1)}-Z^{(i-2)}]+Z^{(i-2)}|^{2\sigma}+\phi(.-y^{(i-1)})|s[Z^{(i-1)}-Z^{(i-2)}]+Z^{(i-2)}|\\
\end{split}\end{equation}
which follows from Lemma~\ref{lem:Nest}. Next, using reasoning as
in the proof of Lemma~\ref{tedious} in combination with
Lemma~\ref{lem:Nest}, we obtain
\begin{equation}\nonumber\begin{split}
&\Big|[\pi^{(i-1)}-\pi^{(i-2)}](t)\int_{0}^{1}\partial_{\pi}N(Z^{(i-2)},\pi^{(i-2)}+s(\pi^{(i-1)}-\pi^{(i-2)}))(t)\,ds\Big|
\\&\lesssim
[\sup_{s\in[0,t]}\la
s\ra^{1+\frac{\epsilon}{2}}|\tilde{\pi}^{(i-1)}-\tilde{\pi}^{(i-2)}|(s)]|Z^{(i-2)}(t)|^{2}\sup_{0\leq
s\leq 1}\phi(.-y^{(i-1)}+s(y^{(i-1)}-y^{(i-2)}))(t)
\end{split}\end{equation}
We can further estimate this by 
\begin{equation}\nonumber\begin{split}
&\Big|[\sup_{s\in[0,t]}\la
s\ra^{1+\frac{\epsilon}{2}}|\tilde{\pi}^{(i-1)}-\tilde{\pi}^{(i-2)}|(s)]|Z^{(i-2)}(t)|^{2}\sup_{0\leq
s\leq 1}\phi(.-y^{(i-1)}+s(y^{(i-1)}-y^{(i-2)}))(t)\Big|\\
&\lesssim \la t\ra^{1+\frac{\epsilon}{2}}[\la
x-y^{(i-1)}\ra^{-\theta}|Z^{(i-2)}|(t)]^{2}\\&\sup_{0\leq s\leq
1}\la
x-y^{(i-1)}\ra^{2\theta}\phi(.-y^{(i-1)}+s(y^{(i-1)}-y^{(i-2)}))(t)\,
\sup_{0\leq s\leq
T_{i-1}}|\tilde{\pi}^{(i-1)}-\tilde{\pi}^{(i-2)}|(s)\\
\end{split}\end{equation}
and we bound this in turn by
\begin{equation}\nonumber
\lesssim (AC\delta)^{2}\la t\ra^{-1-\frac{\epsilon}{2}}\sup_{0\leq
s\leq T_{i-1}}\la
s\ra^{1+\frac{\epsilon}{2}}|\dot{\tilde{\pi}}^{(i-1)}-\dot{\tilde{\pi}}^{(i-2)}|(s),
\end{equation}
which can be comfortably integrated against $t$. Now we estimate
the third term in the nonlinearity. Note that it is here that we
are forced to use the extra weight $\la t\ra ^{-1}$:
\begin{equation}\nonumber
\la
t\ra^{-1}\Big\|P_{s}\Big\{M_{T_{i}}(\pi^{(i-1)})\calG_{T_{i}}(\pi^{(i-1)})[\Hil(\pi^{(i-2)})-\Hil(\pi^{(i-1)})]Z^{(i-1)}\Big\}\Big\|_{L_{t}^{1}L_{x}^{2}[0,t]}
\lesssim AC\delta
\sup_{s\in[0,T_{i-1}]}|\tilde{\pi}^{(i-1)}-\tilde{\pi}^{(i-2)}|(s)
\end{equation}
One again uses a version of Lemma~\ref{tedious}. Furthermore, we obtain
\begin{equation}\nonumber\begin{split}
&\sup_{0\leq t\leq T_{i-1}}\la t\ra
^{-1}\big\|[\Hil(\pi^{(i-1)}(\cdot))-\Hil(\pi^{(i-1)}_{T_{i}}(\cdot))](Z^{(i)}-Z^{(i-1)})\big\|_{L_{t}^{1}L_{x}^{2}[0,t]}\\
&\lesssim C\delta \sup_{0\leq s\leq T_{i-1}}\|\la
\cdot-y^{(i-1)}(s)\ra^{-\theta}(Z^{(i)}-Z^{(i-1)})(s)\|_{L_{x}^{\infty}}\\
\end{split}\end{equation}
Next, consider $\la
t\ra^{-1}\|U^{(i,i-1)}_{dis}(0)\|_{L_{x}^{2}}$. Note that
\begin{equation}\nonumber
U^{(i,i-1)}(0)=U^{i}(0)-M_{T_{i}}(\pi^{(i-1)})\calG_{T_{i}}(\pi^{(i-1)})(0)\calG_{T_{i-1}}^{-1}(\pi^{(i-2)})(0)M_{T_{i}}^{-1}(\pi^{(i-2)})U^{(i-1)}(0)
\end{equation}
Moreover, we have
\begin{align*}\nonumber
&U^{(i)}(0)=M_{T_{i}}(\pi^{(i-1)})\calG_{T_{i}}(\pi^{(i-1)})(0)\Big[\bm R_{0}\\
\bar{R}_{0}\endm
+h^{(i-1)}f^{+}(\alpha_{T_{i}}^{(i-1)})+\sum_{k=1}^{4}a^{(i-1)}_{k}\tilde{\eta}_{k}(\alpha_{T_{i}}^{(i-1)})\Big]
\end{align*}
\begin{align*}\nonumber
&U^{(i-1)}(0)=M_{T_{i-1}}(\pi^{(i-2)})\calG_{T_{i-1}}(\pi^{(i-2)})(0)\Big[\bm R_{0}\\
\bar{R}_{0}\endm
+h^{(i-2)}f^{+}(\alpha_{T_{i-1}}^{(i-2)})+\sum_{k=1}^{4}a^{(i-2)}_{k}\tilde{\eta}_{k}(\alpha_{T_{i-1}}^{(i-2)})\Big]
\end{align*}
The conclusion is that
\begin{equation}\nonumber
U^{(i,i-1)}(0)=O(\delta[\tilde{\pi}^{(i-1)}_{T_{i}}-\tilde{\pi}^{(i-2)}_{T_{i-1}}])+O(h^{(i-1)}-h^{(i-2)})
\end{equation}
Thus we need  to estimate the difference $h^{(i-1)}-h^{(i-2)}$.
For this, we decompose $U^{(i,i-1)}$ with respect to the operator
$\Hil(\alpha^{(i-1)}_{T_{i}})$:
\begin{equation}\nonumber
U^{(i,i-1)}=\sum_{j=1}^{4}\tilde{a}_{j}^{(i,i-1)}(t)\tilde{\eta}^{(i-1)}_{j}(\alpha^{(i-1)}_{T_{i}})
+b^{+(i,i-1)}(t)f^{+}(\alpha^{(i-1)}_{T_{i}})+b^{-(i,i-1)}(t)f^{-}(\alpha^{(i-1)}_{T_{i}})+U^{(i,i-1)}_{dis}
\end{equation}
In other words, we have
\begin{equation}\nonumber
U^{(i,i-1)}_{root}=\sum_{j=1}^{4}\tilde{a}_{j}^{(i,i-1)}(t)\tilde{\eta}^{(i-1)}_{j}(\alpha^{(i-1)}_{T_{i}})
,\,U^{(i,i-1)}_{hyp}=b^{+(i,i-1)}(t)f^{+}(\alpha^{(i-1)}_{T_{i}})+b^{-(i,i-1)}(t)f^{-}(\alpha^{(i-1)}_{T_{i}})
\end{equation}
In particular, we have the identity
\begin{equation}\nonumber\begin{split}
&P^{+}_{Im}(\alpha_{T_{i}}^{(i-1)})\Big\{M(\pi^{(i-1)}_{T_{i}})(0)\calG_{T_{i}}(\pi^{(i-1)})(0)\Big[h^{(i-1)}f^{+}(\alpha_{T_{i}}^{(i-1)})-h^{(i-2)}f^{+}(\alpha_{T_{i}}^{(i-2)})\\&+\sum_{k=1}^{4}(a_{k}^{(i-1)}-a_{k}^{(j-1)})\tilde{\eta}_{k}(\alpha_{T_{i}}^{(i-1)})
+\sum_{k}a_{k}^{(i-2)}[\tilde{\eta}_{k}(\alpha_{T_{i}}^{(i-1)})-\tilde{\eta}_{k}(\alpha_{T_{i}}^{(j-1)})]\Big]\Big\}=b^{(i,i-1)(+)}(0)f^{+}(\alpha^{(i-1)}_{T_{i}})\\
\end{split}\end{equation}
In order to control the difference $h^{(i-1)}-h^{(i-2)}$, we
need to control $b^{(i,i-1)(+)}(0)$. Indeed, assume we control the
latter. Then we have
\begin{equation}\nonumber
h^{(i-1)}f^{+}(\alpha_{T_{i}}^{(i-1)})-h^{(i-2)}f^{+}(\alpha_{T_{i-1}}^{(i-2)})=[h^{(i-1)}-h^{(i-2)}]f^{+}(\alpha_{T_{i}}^{(i-1)})
+h^{(i-2)}[f^{+}(\alpha_{T_{i}}^{(i-1)})-f^{+}(\alpha_{T_{i-1}}^{(i-2)})]
\end{equation}
The conclusion is that
\begin{equation}\nonumber\begin{split}
&|h^{(i-1)}-h^{(i-2)}|\\&\lesssim
\big|P^{+}_{Im}(\alpha_{T_{i}}^{(i-1)})\big\{M_{T_{i}}(\pi^{(i-1)})(0)\calG_{T_{i}}(\pi^{(i-1)})(0)[h^{(i-1)}f^{+}(\alpha_{T_{i}}^{(i-1)})-h^{(i-2)}f^{+}(\alpha_{T_{i-1}}^{(i-2)})]\big\}\big|+\delta|\alpha^{(i-1)}_{T_{i}}-\alpha^{(i-2)}_{T_{i-1}}|\\
&+\big|P^{+}_{Im}(\alpha_{T_{i}}^{(i-1)})\big\{M_{T_{i}}(\pi^{(i-1)})(0)(\calG_{T_{i}}(\pi^{(i-1)})-I)(0)[h^{(i-1)}f^{+}(\alpha_{T_{i}}^{(i-1)})-h^{(i-2)}f^{+}(\alpha_{T_{i-1}}^{(i-2)})]\big\}\big|\\
\end{split}\end{equation}
On account of the a priori estimates established thus far, we
conclude that
\begin{equation}\nonumber
|(\calG_{T_{i}}(\pi^{(i-1)})(0)-I)f^{\pm}(\alpha_{T_{i}}^{(i-1)})|\lesssim
C\delta
\end{equation}
Furthermore, we have
\begin{equation}\nonumber\begin{split}
&\Big|P^{+}_{Im}(\alpha_{T_{i}}^{(i-1)})\Big[M_{T_{i}}(\pi^{(i-1)})\calG_{T_{i}}(\pi^{(i-1)})\sum_{k=1}^{4}(a_{k}^{(i-1)}-a_{k}^{(i-2)})\tilde{\eta}_{k}(\alpha_{T_{i}}^{(i-1)})\Big]\Big|(0)\\
&=\Big|P^{+}_{Im}(\alpha_{T_{i}}^{(i-1)})\Big[(M_{T_{i}}(\pi^{(i-1)})\calG_{T_{i}}(\pi^{(i-1)})-I)\sum_{k=1}^{4}(a_{k}^{(i-1)}-a_{k}^{(i-2)})\tilde{\eta}_{k}(\alpha_{T_{i}}^{(i-1)})\Big]\Big|(0)\\
&\lesssim \delta [|h^{(i-1)}-h^{(i-2)}|+|\alpha^{(i-1)}_{T_{i}}-\alpha^{(i-2)}_{T_{i-1}}|]\\
\end{split}\end{equation}
We use here that
$a_{k}^{(i-1)}=a_{k}^{(i-1)}(h^{(i-1)},\alpha^{(i-1)}_{T_{i}})$
etc.  Also, we observe that\footnote{Use that
$a_{k}^{(i-1)}=O(\delta)$.}
\begin{equation}\nonumber\begin{split}
&\Big|P^{+}_{Im}(\alpha_{T_{i}}^{(i-1)})\Big[M_{T_{i}}(\pi^{(i-1)})\calG_{T_{i}}(\pi^{(i-1)})\sum_{k}a_{k}^{(i-1)}
\big(\tilde{\eta}_{k}(\alpha_{T_{i}}^{(i-1)})-\tilde{\eta}_{k}(\alpha_{T_{i-1}}^{(i-2)})\big)\Big]\Big|
\lesssim
\delta|\tilde{\pi}_{T_{i}}^{(i-1)}-\tilde{\pi}_{T_{i-1}}^{(i-2)}|
\end{split}\end{equation}
Finally, we observe that thanks to the a priori estimates thus far
established
\begin{equation}\nonumber
|\tilde{\pi}_{T_{i}}^{(i-1)}-\tilde{\pi}_{T_{i-1}}^{(i-2)}|\leq
|\tilde{\pi}_{T_{i-1}}^{(i-1)}-\tilde{\pi}_{T_{i-1}}^{(i-2)}|+|\dot{\tilde{\pi}}^{(i-1)}|\lesssim
\sup_{t\in[0,T_{i-1}]}|\tilde{\pi}^{(i-1)}-\tilde{\pi}^{(i-2)}|(t)+T_{i-1}^{-2-\epsilon}
\end{equation}
Putting all of these observations together, we get that
\begin{equation}\nonumber
|h^{(i-1)}-h^{(i-2)}|\lesssim
b^{+(i,i-1)}(0)+\delta\|(\tilde{\pi}^{(i-1)}-\tilde{\pi}^{(i-2)},Z^{(i-1)}-Z^{(i-2)})\|_{Y^{(i-1)}([0,T_{i-1}])}+T_{i-1}^{-2-\epsilon}
\end{equation}
and thus it suffices to estimate $b^{+(i,i-1)}(0)$, as claimed
earlier. To do so, we observe the ODE system satisfied by the
$b^{+(i,i-1)}$ as follows:
\begin{align*}\nonumber
&\frac{d}{dt}\bm b^{(i,i-1)+}(t)\\
b^{(i,i-1)-}(t)\endm - \bm \sigma(\alpha^{(i-1)}_{T_{i}})&
0\\0&-\sigma(\alpha^{(i-1)}_{T_{i}})\endm\bm b^{(i,i-1)+}(t)\\
b^{(i,i-1)-}(t)\endm \\&=\bm
g^{+}(U^{(i,i-1)},U^{(i-1)},U^{(i-2)},\pi^{(i-1)},\pi^{(i-2)},\dot{\pi}^{(i)}-\dot{\pi}^{(i-1)})\\g^{-}(U^{(i,i-1)}U^{(i-1)},U^{(i-2)},\pi^{(i-1)},\pi^{(i-2)},\pi^{(i)}-\pi^{(i-1)})\endm
\end{align*}
where
\begin{equation}\nonumber\begin{split}
&g^{\pm}(U^{(i,i-1)},U^{(i-1)},\pi^{(i-1)},\pi^{(i-2)},\dot{\pi}^{(i)}-\dot{\pi}^{(i-1)})\\
&=P^{\pm}_{Im}(\alpha^{(i-1)}_{T_{i}})\Big\{M_{T_{i}}(\pi^{(i-1)})\calG_{T_{i}}(\pi^{(i-1)})[\text{right-hand
side of \eqref{diff2}}\\
&+V^{(i-1)}(Z^{(i)}-Z^{(i-1)})]\Big\}
\end{split}\end{equation}
On the other hand, by construction of $U^{(i,i-1)}$ as difference
of two $L^{2}$-bounded functions, we know that $b^{(i,i-1)}(t)$ is
controlled by an a priori bound $C\delta$ on $[0,T_{i}]$. The
solution $b^{+(i,i-1)}(t)$ of the ODE is given by
\begin{equation}\nonumber
b^{+(i,i-1)}(t)=e^{t\gamma(\alpha_{T_{i}}^{(i-1)})}b^{+(i,i-1)}(0)+\int_{0}^{t}e^{(t-s)\gamma(\alpha_{T_{i}}^{(i-1)})}g^{+}(,)ds
\end{equation}
We conclude that
\begin{equation}\nonumber
b^{+(i,i-1)}(0)=-\int_{0}^{T_{i-1}}e^{-s\gamma(\alpha_{T_{i}}^{(i-1)})}g^{+}(,)ds+O(e^{-\gamma(\alpha_{T_{i}}^{(i-1)})T_{i-1}})
\end{equation}
Thus we need to estimate $g^{+}(.)$. However, this can be done as in our previous work concerning the bound on
$\|\la t\ra^{-1}(Z^{(i)}-Z^{(i-1)})\|_{L_{x}^{2}}$. This completes
estimating the dispersive part $U^{(i,i-1)}_{dis}$.\\
{\bf{(B):}} We also need to estimate the root as well as the
hyperbolic part for $t\in[0,T_{i-1}]$. For the latter, one again
invokes the ODE system, as before. However, there is a technical
difficulty since we no longer work on an infinite interval and can
no longer conclude that there is just one possible value for
$b^{+(i,i-1)}(0)$ given in terms of an integral representation. We
do have the representation (for $t\in[0,T_{i-1}]$)
\begin{equation}\nonumber
b^{+(i,i-1)}(t)=e^{t\gamma(\alpha_{T_{i}}^{(i-1)})}b^{+(i,i-1)}(0)+\int_{0}^{t}e^{(t-s)\gamma(\alpha_{T_{i}}^{(i-1)})}g^{+}(,)ds.
\end{equation}
As before, this allows us to conclude that
\begin{equation}\nonumber
b^{+(i,i-1)}(t)=O(e^{(t-T_{i-1})\gamma(\alpha_{T_{i}}^{(i-1)})})-\int_{t}^{T_{i-1}}e^{(t-s)\gamma(\alpha_{T_{i}}^{(i-1)})}g^{+}(,)ds
\end{equation}
The first term on the left is only $O(1)$, though, as
$t\rightarrow T_{i-1}$. Note that this problem doesn't occur for
the establishment of the a priori estimates since we define
$b^{+(i)}(0)$ etc to be given by the integral expression. To
remedy this, observe that we also have
\begin{equation}\nonumber
U^{(i,i-1)}(t)=U^{(i)}(t)-M_{T_{i}}(\pi^{(i-1)})(t)\calG_{T_{i}}(\pi^{(i-1)})(t)\calG_{T_{i-1}}(\pi^{(i-2)})(t)^{-1}M_{T_{i-1}}(\pi^{(i-2)})(t)^{-1}U^{(i-1)}(t)
\end{equation}
Observe that as long as $T_{i-1}-t \geq C\log|T_{i-1}|$, one
obtains the desired estimate for $b^{(i,i-1)}(t)$. If $T_{i-1}-t<
C\log|T_{i-1}|$, one invokes the a priori estimates, which produce
a bound of the form
\begin{equation}\nonumber
\la t\ra^{-1}\|U^{(i,i-1)}_{hyp}(t)\|_{L_{x}^{2}}\lesssim \la
t\ra^{-1}[\|U^{(i)}(t)\|_{L_{x}^{\infty}}+\|U^{(i-1)}(t)\|_{L_{x}^{\infty}}]\lesssim
C\delta T_{i-1}^{-\frac{3}{2}},
\end{equation}
which is more than what we need. We proceed to the root part.
Using the notation from above, we have
\begin{equation}\nonumber
U^{(i,i-1)}_{root}(t)=\sum_{j=1}^{4}\tilde{a}_{j}^{(i,i-1)}(t)\tilde{\eta}^{(i-1)}_{j}(\alpha^{(i-1)}_{T_{i}})
\end{equation}
From the iterative construction, we have
\begin{equation}\nonumber
\la U^{(i,i-1)},\tilde{\xi}_{\ell}(\pi^{(i-1)}) \ra = \la
Z^{(i,i-1)}, \xi_{\ell}(\pi^{(i-1)})\ra =\la Z^{(i-1)},
-\xi_{\ell}(\pi^{(i-2)})+\xi_{\ell}(\pi^{(i-1)})\ra
\end{equation}
Hence
\begin{equation}\nonumber
\la U^{(i,i-1)}_{root},\tilde{\xi}_{\ell}(\pi^{(i-1)}) \ra =-\la
U^{(i,i-1)}_{dis},\tilde{\xi}_{\ell}(\pi^{(i-1)}) \ra-\la
U^{(i,i-1)}_{hyp},\tilde{\xi}_{\ell}(\pi^{(i-1)}) \ra+\la Z^{(i-1)},
-\xi_{\ell}(\pi^{(i-2)})+\xi_{\ell}(\pi^{(i-1)})\ra
\end{equation}
From here one obtains the bound
\begin{equation}\nonumber
\sup_{t\in [0,T_{i-1}]}\la
t\ra^{-1}\|U^{(i,i-1)}_{root}(t)\|_{L_{x}^{2}}\lesssim
AC\delta\|(Z^{(i-1)}-Z^{(i-2)},\pi^{(i-1)}-\pi^{(i-2)})\|_{Y^{(i-1)}[0,T_{i-1}]}+C\delta
T_{i}^{-1-\epsilon}
\end{equation}
{\bf{(C):}} Finally, we consider $\sup_{t\in [T_{i-1},T_{i}]}\la
t\ra^{-1}\|U^{(i,i-1)}\|_{L_{x}^{2}}$: we use the a priori
estimates to obtain
\begin{equation}\nonumber
\la t\ra^{-1}\|U^{(i,i-1)}\|\leq \la t\ra^{-1}
[\|Z^{(i)}\|_{L_{x}^{2}}+\|Z^{(i-1)}\|_{L_{x}^{2}}]\lesssim
C\delta T_{i}^{-1}
\end{equation}
\\
We proceed to the weighted $L^{\infty}$-norm
$\sup_{t\in[0,T_{i}]}\|\la
x-\y^{(i-1)}(t)\ra^{-\theta}(Z^{(i)}-Z^{(i-1)})(t)\|_{L_{x}^{\infty}}$.
We follow the same strategy as before by changing the gauge,
working with
\begin{equation}\nonumber
U^{(i,i-1)}=M_{T_{i}}(\pi^{(i-1)})(t)\calG_{T_{i}}(\pi^{(i-1)})(t)(Z^{(i)}-Z^{(i-1)})(t)
\end{equation}
{\bf{(D):}} the estimate for $\|\la
x\ra^{-\theta}U^{(i,i-1)}_{dis}\|_{L_{x}^{\infty}([0,T_{i-1}])}$.
We use \eqref{gauged difference eqn.}, Lemma~\ref{lem:theta} as
well as the Duhamel parametrix, to obtain
\begin{equation}\label{weighted difference}\begin{split}
&\|\la
x\ra^{-\theta}U^{(i,i-1)}_{dis}(t)\|_{L_{x}^{\infty}}\lesssim
\|\la x\ra^{-\theta}e^{i
t\Hil(\alpha^{(i-1)}_{T_{i}})}U^{(i,i-1)}_{dis}(0)\|_{L_{x}^{\infty}}\\&\hspace{4cm}+\int_{0}^{t-1}(t-s)^{-1-\epsilon}\|\la
x\ra^{\theta}P_{s}(\alpha^{(i-1)}_{T_{i}})[\text{left-hand side of
\eqref{gauged difference
eqn.}}(s)]\|_{L_{x}^{1}}ds\\
&\hspace{4cm}+\int_{t-1}^{t}(t-s)^{-\frac{1}{2}}\|P_{s}(\alpha^{(i-1)}_{T_{i}})[\text{left-hand side of \eqref{gauged difference
eqn.}}(s)]\|_{L_{x}^{1}}ds
\end{split}\end{equation}
We start with the first integral on the right, which can be
decomposed into the following terms:
\begin{equation}\nonumber\begin{split}
&\int_{0}^{t-1}(t-s)^{-1-\epsilon}\|\la x\ra^{\theta}P_{s}(\alpha^{(i-1)}_{T_{i}})[M_{T_{i}}(\pi^{(i-1)})\calG_{T_{i}}(\pi^{(i-1)})[i\dot{\tilde{\pi}}^{(i)}\partial_{\pi}W(\pi^{(i-1)})\\&\hspace{9cm}-i\dot{\tilde{\pi}}^{(i-1)}\partial_{\pi}W(\pi^{(i-2)})]](s)\|_{L_{x}^{1}}\,ds\\
\end{split}\end{equation}
We split this into two contributions, the first of which is
\begin{equation}\nonumber\begin{split}
&\int_{0}^{t-1}(t-s)^{-1-\epsilon}\big\|\la x\ra^{\theta}P_{s}(\alpha^{(i-1)}_{T_{i}})\big\{M_{T_{i}}(\pi^{(i-1)})\calG_{T_{i}}(\pi^{(i-1)})[i(\dot{\tilde{\pi}}^{(i)}-\dot{\tilde{\pi}}^{(i-1)})\partial_{\pi}W(\pi^{(i-1)})](s)\big\}\big\|_{L_{x}^{1}}\,ds\\
\end{split}\end{equation}
Observe that one may move the $\la x\ra^{\theta}$ inside,
replacing it by $\la x-y^{(i-1)}(t)\ra^{\theta}$. Then one may
replace the expression
\[P_{s}(\alpha^{(i-1)}_{T_{i}})[M_{T_{i}}(\pi^{(i-1)})\calG_{T_{i}}(\pi^{(i-1)})\partial_{\pi}W(\pi^{(i-1)})]\]
by
\begin{equation}\nonumber
P_{s}(\alpha^{(i-1)}_{T_{i}})\big\{M_{T_{i}}(\pi^{(i-1)})\calG_{T_{i}}(\pi^{(i-1)})
[\partial_{\pi}W(\pi^{(i-1)})-\partial_{\pi}W(\pi^{(i-1)}_{T_{i}})]\big\},
\end{equation}
see Lemma~\ref{lem:rho_T}. Thus, we need to estimate
\begin{equation}\nonumber
\delta\int_{0}^{t-1}(t-s)^{-1-\epsilon}\|\la
x-y^{(i-1)}(s)\ra^{\theta}i(\dot{\tilde{\pi}}^{(i)}-\dot{\tilde{\pi}}^{(i-1)})[\partial_{\pi}W(\pi^{(i-1)})(s)-\partial_{\pi}W(\pi^{(i-1)}_{T_{i}})]\|_{L_{x}^{1}}\,ds
\end{equation}
Due to the local nature of this expression, one can suppress the
factor $\la x-y^{(i-1)}(s)\ra^{\theta}$, and we can estimate the
above by
\begin{equation}\nonumber
\delta\sup_{t\in[0,T_{i-1}]}\la
t\ra^{1+\frac{\epsilon}{2}}|\dot{\tilde{\pi}}^{(i)}-\dot{\tilde{\pi}}^{(i-1)}|(t)\int_{0}^{t-1}(t-s)^{-1-\epsilon}\la
s\ra^{-1-\frac{\epsilon}{2}}\,ds\lesssim
C\delta\sup_{t\in[0,T_{i-1}]}\la
t\ra^{1+\frac{\epsilon}{2}}|\dot{\tilde{\pi}}^{(i)}-\dot{\tilde{\pi}}^{(i-1)}|(t)
\end{equation}
The second contribution is the following
\begin{equation}\nonumber\begin{split}
&\int_{0}^{t-1}(t-s)^{-1-\epsilon}\Big\|\la
x\ra^{\theta}P_{s}(\alpha^{(i-1)}_{T_{i}})
\Big\{M_{T_{i}}(\pi^{(i-1)})\calG_{T_{i}}(\pi^{(i-1)})[\dot{\tilde{\pi}}^{(i-1)}(\partial_{\pi}W(\pi^{(i-1)})-
\partial_{\pi}W(\pi^{(i-2)}))(s)]\Big\}\Big\|_{L_{x}^{1}}\,ds\\
&\lesssim \int_{0}^{t-1}(t-s)^{-1-\epsilon}\|\la
x-y^{(i-1)}(s)\ra^{\theta}\dot{\tilde{\pi}}^{(i-1)}(\partial_{\pi}W(\pi^{(i-1)})-
\partial_{\pi}W(\pi^{(i-2)}))(s)\|_{L_{x}^{1}}\,ds\\
&\lesssim[\sup_{t\in
[0,T_{i-1}]}|\tilde{\pi}^{(i-1)}-\tilde{\pi}^{(i-2)}|(t)]
AC\delta\int_{0}^{t-1}(t-s)^{-1-\epsilon}s^{-1-\frac{\epsilon}{2}}ds\lesssim AC\delta\sup_{t\in [0,T_{i-1}]}<t>^{1+\frac{\epsilon}{2}}|\dot{\tilde{\pi}}^{(i-1)}-\dot{\tilde{\pi}}^{(i-2)}|(t)\\
\end{split}\end{equation}
We have used Lemma~\ref{tedious}.\\
Next, we treat the difference
$N(Z^{(i-1)},\pi^{(i-1)})-N(Z^{(i-2)},\pi^{(i-2)})$ in the
nonlinearity. We split
\begin{equation}\nonumber\begin{split}
&N(Z^{(i-1)},\pi^{(i-1)})-N(Z^{(i-2)},\pi^{(i-2)})=N(Z^{(i-1)},\pi^{(i-1)})-N(Z^{(i-2)},\pi^{(i-1)})
\\&\hspace{8cm}+N(Z^{(i-2)},\pi^{(i-1)})-N(Z^{(i-2)},\pi^{(i-2)})\\
\end{split}\end{equation}
Proceeding as in the estimate of
$\|U^{(i,i-1)}_{dis}\|_{L_{x}^{2}}$, we write
\begin{equation}\nonumber\begin{split}
&\int_{0}^{t-1}(t-s)^{-1-\epsilon}\Big\|\la
x\ra^{\theta}P_{s}(\alpha^{(i-1)}_{T_{i}})\Big\{M_{T_{i}}(\pi^{(i-1)})\calG_{T_{i}}(\pi^{(i-1)})[
N(Z^{(i-1)},\pi^{(i-1)})-N(Z^{(i-2)},\pi^{(i-1)})](s)\Big\}\Big\|_{L_{x}^{1}}\,ds\\
&=\int_{0}^{1}\int_{0}^{t-1}(t-s)^{-1-\epsilon}\Big\|\la
x\ra^{\theta}P_{s}(\alpha^{(i-1)}_{T_{i}})\Big \{M_{T_{i}}(\pi^{(i-1)})\calG_{T_{i}}(\pi^{(i-1)})[
(Z^{(i-1)}-Z^{(i-2)})\\&\hspace{8cm}\partial_{Z}N(Z^{(i-2)}+\mu(Z^{(i-1)}-Z^{(i-2)}),\pi^{(i-1)})](s)\Big]\Big\|_{L_{x}^{1}}\,ds d\mu\\
\end{split}\end{equation}
We use Lemma~\ref{lem:Nest} to reduce this to two estimates, a
local and nonlocal one. First, we get (by moving the multiplier $\la
x\ra^{\theta}$ inside and then removing
$P_{s}(\alpha^{(i-1)}_{T_{i}})$ etc.)
\begin{equation}\nonumber\begin{split}
&\int_{0}^{t-1}(t-s)^{-1-\epsilon}\|\la x-y^{(i-1)}(s)\ra^{\theta}[Z^{(i-1)}-Z^{(i-2)}](s)\\
&\hspace{7cm}|Z^{(i-2)}+\mu(Z^{(i-1)}-Z^{(i-2)})|(s)\phi(.-y^{(i-1)}(s))\|_{L_{x}^{1}}\,ds\\
&\lesssim
A^{2}C\delta[\sup_{t\in[0,T_{i-1}]}\|\la x-y^{(i-2)}(t)\ra^{-\theta}(Z^{(i-1)}-Z^{(i-2)})(t)\|_{L_{x}^{\infty}}]\int_{0}^{t-1}(t-s)^{-1-\epsilon}\la s\ra^{-1-\epsilon}\, ds
\end{split}\end{equation}
For the nonlocal term, we estimate if we assume (as we may) that
$0<\epsilon<\sigma-1-2\theta$
\begin{equation}\nonumber\begin{split}
&\int_{0}^{t-1}(t-s)^{-1-\epsilon}\|\la x-y^{(i-1)}(s)\ra^{\theta}[Z^{(i-1)}-Z^{(i-2)}](s)\\&\hspace{7cm}|Z^{(i-2)}+\mu(Z^{(i-1)}-Z^{(i-2)})|^{2\sigma}(s)\|_{L_{x}^{1}}\,ds\\
&\lesssim
A^{\theta}[\sup_{t\in[0,T_{i-1}]}\|\la x-y^{(i-2)}(t)\ra^{-\theta}(Z^{(i-1)}-Z^{(i-2)})(t)\|_{L_{x}^{\infty}}]\\&\hspace{2cm}\int_{0}^{t-1}(t-s)^{-1-\epsilon}\|\la x-y^{i-1}(s)\ra^{2\theta}|Z^{(i-2)}+\mu(Z^{(i-1)}-Z^{(i-2)})|^{2\sigma}(s)\|_{L_{x}^{1}}\,ds\\
&\lesssim A^{2}(C\delta)^{4}[\sup_{t\in[0,T_{i-1}]}\|\la x-y^{(i-2)}(t)\ra^{-\theta}(Z^{(i-1)}-Z^{(i-2)})(t)\|_{L_{x}^{\infty}}]\\
\end{split}\end{equation}
We have used Lemma~\ref{lem:moments}. Further, using
Lemma~\ref{tedious} as well as Lemma~\ref{lem:Nest}, we get
\begin{equation}\begin{split}
&|N(Z^{(i-2)},\pi^{(i-1)})-N(Z^{(i-2)},\pi^{(i-2)})|(s)\\&\lesssim
\sup_{0\leq \mu\leq 1}\la
s\ra|\tilde{\pi}^{(i-1)}-\tilde{\pi}^{(i-2)}|(s)
[|Z^{(i-2)}|^{2}(s)\phi(.-\mu
y^{(i-2)}(s)-(1-\mu)y^{(i-1)}(s))+|Z^{(i-2)}|(s)^{2\sigma+1}]\\
\end{split}\end{equation}
Plugging this back in, we get for the contribution of the local
term
\begin{equation}\nonumber\begin{split}
&\int_{0}^{t-1}(t-s)^{-1-\epsilon}\la s\ra \sup_{0\leq\mu\leq
1}\|\la
x-y^{(i-1)}(s)\ra^{\theta}|\tilde{\pi}^{(i-1)}-\tilde{\pi}^{(i-2)}|(s)
[|Z^{(i-2)}|^{2}(s)\phi(.-\mu y^{(i-2)}(s)\|_{L_{x}^{1}}ds\\
&\lesssim
A^{2}C^{2}\delta[\sup_{s\in[0,T_{i-1}]}\la s\ra^{1+\frac{\epsilon}{2}}|\dot{\tilde{\pi}}^{(i-1)}-\dot{\tilde{\pi}}^{(i-2)}|(s)]\int_{0}^{t-1}(t-s)^{-1-\epsilon}\la s\ra^{-1-\epsilon}ds,\\
\end{split}\end{equation}
while for the contribution of the nonlocal term, we get
\begin{equation}\nonumber\begin{split}
&\int_{0}^{t-1}(t-s)^{-1-\epsilon}\la s\ra \sup_{0\leq\mu\leq
1}\|\la
x-y^{(i-1)}(s)\ra^{\theta}|\tilde{\pi}^{(i-1)}-\tilde{\pi}^{(i-2)}|(s)
\,|Z^{(i-2)}|(s)^{2\sigma+1}\|_{L_{x}^{1}}\,ds\\
&\lesssim A(C\delta)^{2\sigma}[\sup_{s\in[0,T_{i-1}]}\la
s\ra^{1+\frac{\epsilon}{2}}|\dot{\tilde{\pi}}^{(i-1)}-\dot{\tilde{\pi}}^{(i-2)}|(s)]\int_{0}^{t-1}(t-s)^{-1-\epsilon}\la
s\ra^{\frac{3}{2}+\theta-\sigma}ds\\
\end{split}\end{equation}
Finally, we need to control the following two local terms:
\begin{equation}\nonumber\begin{split}
&\int_{0}^{t-1}(t-s)^{-1-\epsilon}\|\la x-y^{(i-1)}(s)\ra^{\theta}[\Hil(\pi^{(i-2)})-\Hil(\pi^{(i-1)})](s)Z^{(i-1)}\|_{L_{x}^{1}}ds\\
&\lesssim \sup_{s\in[0,T_{i-1}]}\la s\ra^{1+\frac{\epsilon}{2}}|\dot{\tilde{\pi}}^{(i-1)}-\dot{\tilde{\pi}}^{(i-2)}|(s)A^{2}C\delta\int_{0}^{t-1}(t-s)^{-1-\epsilon}s^{-\frac{\epsilon}{2}}ds\\
\end{split}\end{equation}
where we have used Lemma~\ref{tedious}. Similarly, by Lemma~\ref{lem:rho_T},
\begin{equation}\nonumber\begin{split}
&\int_{0}^{t-1}(t-s)^{-1-\epsilon}\|\la x-y^{(i-1)}(s)\ra
^{\theta}[\Hil(\pi^{(i-1)})(t)-\Hil(\pi^{(i-1)}_{T_{i}}(t))](Z^{(i)}-Z^{(i-1)})\|_{L_{x}^{1}}ds
\\&\lesssim
C\delta[\sup_{s\in[0,T_{i-1}]}\|\la x-y^{(i-2)}(s)\ra^{-\theta}(Z^{(i)}-Z^{(i-1)})\|_{L_{x}^{\infty}}]\int_{0}^{t-1}(t-s)^{-1-\epsilon}s^{-\epsilon}ds\\
\end{split}\end{equation}
This completes estimating the
first integral expression in \eqref{weighted difference}. The second
is handled similarly, and is hence omitted. We still need to estimate
the contribution from the free term. The fact that the difference
$(Z^{(i)}-Z^{(i-1)})(0)$ is a local term ensures that this can be
accomplished exactly along the same lines as for
$\|(Z^{(i)}-Z^{(i-1)})(0)\|_{L_{x}^{2}}$, which we dealt with
earlier. This completes part {\bf{(D)}}.\\
{\bf{(E):}} The estimate for $\|\la x\ra
^{-\theta}U^{(i,i-1)}_{hyp}\|_{L_{x}^{\infty}}$. This is handled
exactly like the corresponding contribution of $\|\la
t\ra^{-1}U^{(i,i-1)}_{hyp}\|_{L_{x}^{2}}$, on account of the fact
that this term is local.\\
{\bf{(F):}} The estimate for $\|\la x\ra
^{-\theta}U^{(i,i-1)}_{root}\|_{L_{x}^{\infty}}$. Use again that
\begin{equation}\nonumber
\la U^{(i,i-1)}_{root},\tilde{\xi}_{\ell}(\pi^{(i-1)}) \ra =-\la
U^{(i,i-1)}_{dis},\tilde{\xi}_{\ell}(\pi^{(i-1)}) \ra-\la
U^{(i,i-1)}_{hyp},\tilde{\xi}_{\ell}(\pi^{(i-1)}) \ra+\la Z^{(i-1)},
\xi_{\ell}(\pi^{(i-2)})-\xi_{\ell}(\pi^{(i-1)})\ra
\end{equation}
One combines the previous estimates for $\|\la
x\ra^{-\theta}U^{(i,i-1)}_{dis}\|_{L_{x}^{\infty}}$, $\|\la
x\ra^{-\theta}U^{(i,i-1)}_{hyp}\|_{L_{x}^{\infty}}$ with the
following:
\begin{equation}\nonumber
|\la Z^{(i-1)}, \xi_{\ell}(\pi^{(i-2)})-\xi_{\ell}(\pi^{(i-1)})\ra|(t)
\lesssim AC\delta \sup_{s\in[0,T_{i-1}]}\la
s\ra^{1+\frac{\epsilon}{2}}|\dot{\tilde{\pi}}^{(i-1)}-\dot{\tilde{\pi}}^{(i-2)}|(s)
\end{equation}
{\bf{(G):}} Estimating $\sup_{t\in[T_{i-1},T_{i}]}\|\la
x\ra^{-\theta}U^{(i,i-1)}(t)\|_{L_{x}^{\infty}}$. This is again
accomplished by using the a priori estimates. We get for
$t\in[T_{i-1},T_{i}]$
\begin{equation}\nonumber
\|\la x\ra^{-\theta}U^{(i,i-1)}(t)\|_{L_{x}^{\infty}}\lesssim
(A+1)C\delta T_{i}^{-1-\epsilon},
\end{equation}
taking into account the translation effect of
\begin{equation}\nonumber
M_{T_{i-1}}(\pi^{(i-1)})(t)\calG_{T_{i-1}}(\pi^{(i-1)})(t)\calG_{T_{i-2}}(\pi^{(i-2)})(t)^{-1}M_{T_{i-2}}(\pi^{(i-2)})(t)^{-1}
\end{equation}
Putting these estimates together completes the estimation for the
norms involving $Z^{(i)}-Z^{(i-1)}$.\\
{\bf{(H):}} We proceed to estimating $\sup_{t\in[0,T_{i}]}\la
t\ra^{1+\frac{\epsilon}{2}}|\dot{\tilde{\pi}}^{(i)}-\dot{\tilde{\pi}}^{(i-1)}|(t)$.
First, consider the case $t\in[0,T_{i-1}]$. Using the equation
\begin{equation}\label{diff1}\begin{split}
&i[\dot{\tilde{\pi}}^{(i)}-\dot{\tilde{\pi}}^{(i-1)}]=\la
Z^{(i)}(t),i\dot{\tilde{\pi}}^{(i-1)}\tilde{\calS}_{\ell}(\pi^{(i-1)})(t)\ra
-\la N(Z^{(i-1)},\pi^{(i-1)}),\tilde{\xi}_{\ell}(\pi^{(i-1)})(t)\ra-
\\&\hspace{3cm}(\la
Z^{(i-1)}(t),i\dot{\tilde{\pi}}^{(i-2)}\tilde{\calS}_{\ell}(\pi^{(i-2)})(t)\ra
-\la
N(Z^{(i-2)},\pi^{(i-2)}),\tilde{\xi}_{\ell}(\pi^{(i-2)})(t)\ra)\\
\end{split}\end{equation}
First we have
\begin{equation}\nonumber
|\la
Z^{(i)}-Z^{(i-1)},i\dot{\tilde{\pi}}^{(i-1)}\calS_{\ell}(\pi^{(i-1)})\ra(t)|\lesssim
C\delta\la t\ra^{-1-\epsilon}\sup_{t\in [0,T_{i-1}]}\la
t\ra^{-1}\|(Z^{(i)}-Z^{(i-1)})(t)\|_{L_{x}^{2}}
\end{equation}
Of course this term will then be moved to the left. Next, we have
\begin{equation}\nonumber
|\la
Z^{(i-1)},[\dot{\tilde{\pi}}^{(i-1)}-\dot{\tilde{\pi}}^{(i-2)}]\calS_{\ell}(\pi^{(i-2)})\ra
|\lesssim \la
t\ra^{-1-\frac{\epsilon}{2}}\|Z^{(i-1)}\|_{L_{x}^{2}}\sup_{s\in[0,T_{i-1}]}<s>^{1+\frac{\epsilon}{2}}|\dot{\tilde{\pi}}^{(i-1)}-\dot{\tilde{\pi}}^{(i-2)}|(s)
\end{equation}
Further, using a variant of lemma~\ref{tedious}, we get
\begin{equation}\nonumber
|\la Z^{(i-1)},
\dot{\tilde{\pi}}^{(i-1)}(\calS_{\ell}(\pi^{(i-1)})-\calS_{\ell}(\pi^{(i-2)})(t)\ra|\lesssim
\|Z^{(i-1)}\|_{L_{x}^{2}}\la t\ra^{1+\frac{\epsilon}{2}}
|\dot{\tilde{\pi}}^{(i-1)}(t)| |(\pi^{(i-1)}-\pi^{(i-2)})(t)|,
\end{equation}
which is bounded by (always keeping in mind that
$t\in[0,T_{i-1}]$)
\begin{equation}\nonumber
C\delta\la t\ra^{-1-\frac{\epsilon}{2}}\sup_{s\in[0,T_{i-1}]}\la
s\ra^{1+\frac{\epsilon}{2}}|\dot{\tilde{\pi}}^{(i-1)}-\dot{\tilde{\pi}}^{(i-2)}|(s)
\end{equation}
The difference
\begin{equation}\nonumber
\la
N(Z^{(i-1)},\pi^{(i-1)}),\tilde{\xi}_{\ell}(\pi^{(i-1)})(t)\ra-\la
N(Z^{(i-2)},\pi^{(i-2)}),\tilde{\xi}_{\ell}(\pi^{(i-2)})(t)\ra
\end{equation}
is handled as before. Indeed,  one obtains
\begin{equation}\nonumber\begin{split}
&|\la
N(Z^{(i-1)},\pi^{(i-1)}),\tilde{\xi}_{\ell}(\pi^{(i-1)})(t)\ra-\la
N(Z^{(i-1)},\pi^{(i-2)}),\tilde{\xi}_{\ell}(\pi^{(i-1)})(t)\ra|\\
&\lesssim
\la t\ra^{1+\frac{\epsilon}{2}}|\tilde{\pi}^{(i-1)}-\tilde{\pi}^{(i-2)}|(t)|\partial_{\pi}N(Z^{(i-1)},\mu\pi^{(i-1)}+(1-\mu)\pi^{(i-2)}(t))|\\
&\lesssim A^{2}[C^{2}\delta^{2}+(C\delta)^{2\sigma}]
\la t\ra^{-1-\frac{\epsilon}{2}}\sup_{t\in[0,T_{i-1}]}\la t\ra^{1+\frac{\epsilon}{2}}|(\dot{\tilde{\pi}}^{(i-1)}-\dot{\tilde{\pi}}^{(i-1)})(t)|\\
\end{split}\end{equation}
\begin{equation}\nonumber\begin{split}
&|\la
N(Z^{(i-1)},\pi^{(i-2)}),\tilde{\xi}_{\ell}(\pi^{(i-1)})(t)\ra-\la
N(Z^{(i-2)},\pi^{(i-2)}),\tilde{\xi}_{\ell}(\pi^{(i-1)})(t)\ra|\\
&\lesssim |Z^{(i-1)}-Z^{(i-2)}|(t)|\partial_{Z}N(\mu Z^{(i-1)}
+(1-\mu)Z^{(i-2)},\pi^{(i-2)})(t)|\\&\lesssim
A^{2}[C\delta+(C\delta)^{2\sigma}]\la t\ra^{-1-\epsilon}
\sup_{s\in
[0,T_{i-1}]}\|\la x-y^{(i-1)}(s)\ra^{-\theta}(Z^{(i-1)}-Z^{(i-2)})\|_{L_{x}^{\infty}}\\
\end{split}\end{equation}
{\bf{(I)}}: Finally, consider the case $t\in[T_{i-1},T_{i}]$. Then
we estimate 
\begin{equation}\nonumber\begin{split}
&|\la
Z^{(i)}(t),i\dot{\tilde{\pi}}^{(i-1)}\tilde{\calS}_{\ell}(\pi^{(i-1)})(t)\ra
-\la N(Z^{(i-1)},\pi^{(i-1)}),\tilde{\xi}_{\ell}(\pi^{(i-1)})(t)\ra-
\\&\hspace{3cm}(\la
Z^{(i-1)}(t),i\dot{\tilde{\pi}}^{(i-2)}\tilde{\calS}_{\ell}(\pi^{(i-2)})(t)\ra
-\la
N(Z^{(i-2)},\pi^{(i-2)}),\tilde{\xi}_{\ell}(\pi^{(i-2)})(t)\ra)|\\
&\lesssim [(A+1)^{2}(C\delta)^{2}+(C\delta)^{2\sigma}] T_{i}^{-2-\epsilon}\\
\end{split}\end{equation}
Putting all of the preceding estimates {\bf{(A)}}-{\bf{(I)}}
together, we obtain the claim of Proposition~\ref{pi-pi}.
\end{proof}

We have now shown that Definition~\ref{Iteration} results in
iterates satisfying the a priori estimates \eqref{apriori
estimates}. We now need to show that the $(\pi^{(i)},Z^{(i)})$
converge in a suitable sense. We have the following
\begin{theorem}\label{thm:convergence} For $i\geq 1$, $j\geq i$, the following inequality
holds\footnote{This convergence is rather slow, of course, and can
be significantly improved by choosing $T_{i}$ less
conservatively.} :
\begin{equation}\nonumber
\|(\tilde{\pi}^{(i)}-\tilde{\pi}^{(j)},Z^{(i)}-Z^{(j)})\|_{Y^{(i)}([0,T_{i}])}\lesssim
i^{-1}
\end{equation}
There exists $(\pi,Z)\in X_{*}$ solving \eqref{eq:Zsys} in the
$H^{1}$-sense in addition to satisfying the orthogonality
relations
\begin{equation}\nonumber
\la Z,\xi_{\ell}(\pi(t))\ra=0\quad \forall t\in[0,\infty)
\end{equation}
with the following property: for every $T>0$, we have
\begin{equation}\nonumber
\lim_{i\rightarrow
\infty}\|(\pi-\pi^{(i)},Z-Z^{(i)})\|_{\tilde{X}_{*}(\pi)([0,T])}\rightarrow
0
\end{equation}
In particular, we get
\begin{equation}\nonumber
\|(\pi,Z)\|_{\tilde{X}_{*}(\pi)}\leq C\delta
\end{equation}
Also, we have
\begin{equation}\nonumber
Z(0)=\bm R_{0}\\
\overline{R_{0}}\endm+h
f^{+}(\alpha_{\infty})+\sum_{j=1}^{4}a_{j}\tilde{\eta}_{j}(\alpha_{\infty})
\end{equation}
for suitable $h,\,a_{j}(h)\in\R$, where $h=h(R_{0})$ depends in a
Lipschitz continuous fashion on $R_{0}$. Finally, $(\pi,Z)$ is the
unique solution with these initial data and satisfying the above
bounds and orthogonality relations.
\end{theorem}
\begin{proof} We need to show that the $(\pi^{(i)}, Z^{(i)})$ form
a Cauchy sequence in a suitable sense. This follows from the
following pair of inequalities. Let $j\geq i$.
\begin{equation}\nonumber\begin{split}
&\|(\tilde{\pi}^{(i)}-\tilde{\pi}^{(j)},Z^{(i)}-Z^{(j)})\|_{Y^{(i)}([0,T_{i}])}\\&\lesssim
[A^{2}C\delta+A^{2}(C\delta)^{2\sigma}]\|(\tilde{\pi}^{(i-1)}-\tilde{\pi}^{(j-1)},Z^{(i-1)}-Z^{(j-1)})\|_{Y^{(i-1)}([0,T_{i-1}])}+[1+(A+1)^{2}(C\delta)^{2}+(C\delta)^{2\sigma}]
T_{i}^{-1}\\
\end{split}\end{equation}
\begin{equation}\nonumber
T_{i}\sup_{j\geq
i}\sup_{t\in[0,T_{i}]}|\tilde{\pi}^{(i)}-\tilde{\pi}^{(j)}|(t)\lesssim
A
\end{equation}
The proof of these follows along the exact same lines as the proof
of Proposition~\ref{pi-pi}. One uses the fact that on account of
the a priori estimates already established, we have ($j\geq 0$)
\begin{equation}\nonumber
T_{0}\sup_{j\geq
i}\sup_{t\in[0,T_{0}]}|\pi^{(0)}-\pi^{(j)}|(t)\lesssim A
\end{equation}
If one iterates the first of the above two inequalities, one
obtains
\begin{equation}\nonumber
\|(\tilde{\pi}^{(i)}-\tilde{\pi}^{(j)},Z^{(i)}-Z^{(j)})\|_{Y^{(i)}([0,T_{i}])}\lesssim
i^{-1}
\end{equation}
Now choose $T>0$ and $i_{0}$ such that $T_{i_{0}}>T$. Then we see
that for $i\geq i_{0}$, $j\geq i$, we have
\begin{equation}\nonumber
\|(\tilde{\pi}^{(i)}-\tilde{\pi}^{(j)},Z^{(i)}-Z^{(j)})\|_{Y^{(i)}([0,T])}\lesssim
i^{-1}
\end{equation}
In particular, the numbers $\la
t\ra^{1+\frac{\epsilon}{2}}\dot{\tilde{\pi}}^{(i)}(t)$ converge
uniformly on $[0,T]$, whence also $\pi^{(i)}|_{[0,T]}$ converge
toward some Lipschitz continuous path $\pi_{T}$ on $[0,T]$.
Actually, this path is $C^{1}$, since $\dot{\pi}$ is locally the
uniform limit of continuous functions. Note that if we define
$y_{T}(t)=\int_{0}^{t}v_{T}(s)ds+D_{T}(s)$, we have
\begin{equation}\nonumber
y^{(i)}|_{[0,T]}(t)\rightarrow  y_{T}(t),
\end{equation}
and consequently we have
\begin{equation}\nonumber
\|.\|_{Y^{(i)}[0,T]}\rightarrow \|.\|_{Y[0,T]},
\end{equation}
where $\|.\|_{Y[0,T]}$ is defined like $\|.\|_{Y^{(i)}[0,T]}$ with
$y^{(i-1)}(t)$ replaced by $y(t)$. In summary, we get
\begin{equation}\nonumber
\lim_{i,j\rightarrow\infty}\|(\tilde{\pi}^{(i)}-\tilde{\pi}^{(j)},Z^{(i)}-Z^{(j)})\|_{Y[0,T]}=0.
\end{equation}
This in conjunction to the a priori estimates implies that the
$Z^{(i)}|_{[0,T]}$ converge point-wise toward some function
$Z_{T}\in X_{*}([0,T])$, which satisfies
\begin{equation}\label{local convergence}
\lim_{i\rightarrow\infty}\|(\pi_{T}-\pi^{(i)},Z_{T}-Z^{(i)}|_{[0,T]})\|_{X_{*}[0,T]}=0
\end{equation}
Clearly $(\pi_{T},Z_{T})$ weakly solves \eqref{eq:Zsys} and
satisfies the orthogonality relations on $[0,T]$. Indeed, we can
improve this statement by observing that the norm
$\|.\|_{Y^{(i)}([0,T_{i}])}$ in Proposition~\ref{pi-pi} may be
strengthened to also include
\begin{equation}\nonumber
\sup_{t\in [0,T_{i}]}[\la t\ra^{-1}\|Z(t)\|_{H^{1}}+\|\la
x-y^{(i-1)}\ra^{-\frac{1}{2}-2\epsilon}\partial_{x}Z(t)\|_{L_{x}^{q}}],\,\,\,q\,\text{as
in}\,\|.\|_{X_{*}}
\end{equation}
The justification for this is as in the proof of
Proposition~\ref{(U,pi)}. This in particular entails that
\begin{equation}\nonumber
Z_{T}\in C([0,T],H^{1}(\R)\times H^{1}(\R))\cap
C^{1}([0,T),H^{-1}(\R)\times H^{-1}(\R)),
\end{equation}
and $Z$ solves the equation in the $H^{1}$-sense. Replacing $T$ by
a larger $\tilde{T}$, we can compatibly extend $(\pi_{T}, Z_{T})$
to a larger time interval, whence all the way to $[0,\infty)$.
Then the a priori estimates imply that the $(\pi,Z)$ thus
constructed lies in $X_{*}$, as well as
\begin{equation}\nonumber
C([0,\infty),H^{1}(\R)\times H^{1}(\R))\cap
C^{1}([0,\infty),H^{-1}(\R)\times H^{-1}(\R))
\end{equation}
In particular, the limit
$\alpha_{\infty}:=\lim_{t\rightarrow\infty}\alpha(t)$ exists. The
estimate \eqref{local convergence} implies that
\begin{equation}\nonumber
Z^{(i)}(0)\rightarrow Z(0)
\end{equation}
in the $L^{2}$-sense. In particular, recalling
\begin{equation}\nonumber
Z^{(i)}(0)=\bm R_{0}\\
\bar{R}_{0}\endm
+h^{(i-1)}f^{+}(\alpha_{T_{i}}^{(i-1)})+\sum_{k=1}^{4}a^{(i-1)}_{k}\tilde{\eta}_{k}(\alpha_{T_{i}}^{(i-1)}),
\end{equation}
we get $h^{(i-1)}\rightarrow h^+$ for suitable $h^{+}\in\R$, and
similarly  $a^{(i-1)}_{k}\rightarrow a_{k}$ for suitable
$a_{k}\in\R$. Moreover, on account of the a priori estimates, we
have
\begin{equation}\nonumber
|\alpha_{T_{i}}^{(i-1)}-\alpha^{(i-1)}(T)|<T^{-1-\epsilon},\,|\alpha_{\infty}-\alpha(T)|<T^{-1-\epsilon}
\end{equation}
for $i$ suitably large. Hence
\begin{equation}\nonumber
|\alpha_{T_{i}}^{(i-1)}-\alpha_{\infty}|<2T^{-1-\epsilon}+|\alpha(T)-\alpha^{(i-1)}(T)|
\end{equation}
and therefore
\begin{equation}\nonumber
\lim_{i\rightarrow\infty}|\alpha^{(i-1)}_{T_{i}}-\alpha_{\infty}|<2T^{-1-\epsilon}
\end{equation}
Letting $T\rightarrow\infty$, we get
$\lim_{i\rightarrow\infty}\alpha^{(i-1)}_{T_{i}}=\alpha_{\infty}$,
whence indeed
\begin{equation}\label{initial data}
Z(0)=\bm R_{0}\\
\bar{R}_{0}\endm
+h^{}f^{+}(\alpha_{\infty})+\sum_{k=1}^{4}a_{k}\tilde{\eta}_{k}(\alpha_{\infty}).
\end{equation}
We now verify that  $h^{+},\,a_{k}$ depend in a Lipschitz
continuous manner on $R_{0}$. We claim the following:
\begin{lemma}\label{Lipschitz} Let $(\pi,Z)$ be the solution associated with $\bm
R_{0}\\ \overline{R_{0}}\endm$, and let $(\pi^{*},Z^{*})$
be the solution associated with $\bm R_{0}^{*}\\
\overline{R_{0}^{*}}\endm$. Then the following inequality
holds\footnote{We shall assume in the following that $\trip
R_{0}-R_{0}^{*}\trip<1$.}
\begin{equation}\nonumber
\|\la t\ra^{-1} [Z-Z^{*}]\|_{L_{t}^{\infty}L_{x}^2([0,\trip
R_{0}-R_{0}^{*}\trip^{-1}])}+\|(\tilde{\pi}-\tilde{\pi^{*}})(t)\|_{L_{t}^{\infty}([0,\trip
R_{0}-R_{0}^{*}\trip^{-1}])}\lesssim \trip R_{0}-R_{0}^{*}\trip
\end{equation}
\end{lemma}

Assuming this lemma for now, we introduce
$U(t)=M_{\infty}(\pi)(t)\calG_{\infty}(\pi)(t)Z(t)$, and
analogously for $U^{*}$. The preceding considerations imply that
we may write
\begin{equation}\nonumber
U=U_{dis}+U_{root}+U_{hyp}\,\,\text{\ with respect to
$\Hil(\alpha_{\infty})$},
\end{equation}
and we can write
\begin{equation}\nonumber
U_{hyp}=b^{+}(t)f^{+}(\alpha_{\infty})+b^{-}(t)f^{-}(\alpha_{\infty}),
\end{equation}
where
\begin{equation}\nonumber
b^{+}(0)=-\int_{0}^{\infty}e^{-\gamma(\alpha_{\infty})s}g^{+}(U,\pi)(s)ds
\end{equation}
and we have
\begin{equation}\nonumber
g^{+}(U,\pi)=P^{+}_{Im}(\alpha_{\infty})[i\dot{\tilde{\pi}}\partial_\pi
\tilde W_{\infty}(\pi) +N(U,\pi)+VU]
\end{equation}
\begin{align}
 V = V(t) &:= \bm (\sigma+1)(\phi_\infty^{2\sigma}(x)-\phi^{2\sigma}(x+y_\infty-y)) & \sigma(\phi_\infty^{2\sigma}(x)-e^{2i\rho_\infty} \phi^{2\sigma}(x+y_\infty-y)) \\
-\sigma(\phi_\infty^{2\sigma}(x)- e^{-2i\rho_\infty} \phi^{2\sigma}(x+y_\infty-y)) & -(\sigma+1)(\phi_\infty^{2\sigma}(x)-\phi^{2\sigma}(x+y_\infty-y)) \endm \label{eq:V*}\\
 i\dot{\pi}\partial_\pi \tilde W(\pi) &:= \dot{{v}} \binom{-(x+y_\infty)e^{i\rho_\infty} \phi(x+y_\infty-y)}{(x+y_\infty) e^{-i\rho_\infty} \phi(x+y_\infty-y)} + \dot{\gamma} \binom{-e^{i\rho_\infty}\phi(x+y_\infty-y)}{e^{-i\rho_\infty}\phi(x+y_\infty-y)} \label{eq:vdotetc*} \\
 & \quad +i\dot{\alpha} \binom{e^{i\rho_\infty}\partial_\alpha\phi(x+y_\infty-y)}{e^{-i\rho_\infty}\partial_\alpha\phi(x+y_\infty-y)}+i\dot{D}
\binom{-e^{i\rho_\infty}\partial_x\phi(x+y_\infty-y)}{-e^{-i\rho_\infty}\partial_x\phi(x+y_\infty-y)} \nn \\
\end{align}
and $N(U,\pi)$ is defined as in \eqref{eq:NUpi} with $T$ replaced
by $\infty$. Plugging in the above estimate, one
easily\footnote{Repeating estimates as in the proofs of
Proposition~\ref{(U,pi)} etc.\ and breaking the integrals into two
parts, one over the interval $[0,\trip
R_{0}-R_{0}^{*}\trip^{-1}]$, the other over its complement, where
the a priori estimates are used.} obtains the bound
\begin{equation}\nonumber
\left|\int_{0}^{\infty}e^{-\gamma(\alpha_{\infty})s}g^{+}(U,\pi)(s)ds-\int_{0}^{\infty}e^{-\gamma(\alpha^{*}_{\infty})s}g^{+}(U^{*},\pi^{*})(s)ds\right|\lesssim \delta\trip R_{0}-R^{*}_{0}\trip
\end{equation}
Now one uses that (with an analogous equation determining
$b^{*+}(0)$)
\begin{equation}\nonumber
b^{+}(0)=P^{+}_{Im}(\alpha_{\infty})M_{\infty}(\pi)(0)\calG_{\infty}(\pi)(0)\Big[\bm R_{0}\\
\bar{R}_{0}\endm
+hf^{+}(\alpha_{\infty})+\sum_{k=1}^{4}a_{k}\tilde{\eta}_{k}(\alpha_{\infty})\Big]
\end{equation}
Observing that
\begin{equation}\nonumber
P^{+}_{Im}(\alpha_{\infty})\bm R_{0}\\
\bar{R}_{0}\endm=[P^{+}_{Im}(\alpha_{\infty})-P^{+}_{Im}(\alpha_{0})]\bm R_{0}\\
\bar{R}_{0}\endm,
\end{equation}
one infers from the preceding that
\begin{equation}\nonumber
|h^{}-h^{*}|\lesssim \delta \trip R_{0}-R^{*}_{0}\trip
\end{equation}
Furthermore, exploiting the orthogonality relations determining
$U_{root}$, $U^{*}_{root}$, one obtains a similar estimate for
$a_{k}-a^{*}_{k}$ for all $k$. The argument just given also easily
implies the bound
\begin{equation}\nonumber
|h| +\sum_{k=1}^{4}|a_{k}|\lesssim \trip R_{0}\trip^{2}
\end{equation}
We now turn to the proof of lemma~\ref{Lipschitz}, which is based
on a recursive inequality:\\
\begin{proof}(Lemma~\ref{Lipschitz}) We recycle the notation from
the proof of Proposition~\ref{(U,pi)} etc.\ We claim the following
pair of inequalities hold true\footnote{Possibly shrinking
$\delta$ and growing $A,C$ a bit.} for $j\geq i$:
\begin{equation}\nn
\begin{split}
&\|(\tilde{\pi}^{(i)}-\tilde{\pi}^{*(i)},Z^{(i)}-Z^{*(i)})\|_{Y^{(i)}([0,\min\{T_{i},\trip
R_{0}-R^{*}_{0}\trip^{-1}\}])}\\&\lesssim
[A^{2}C\delta+A^{2}(C\delta)^{2\sigma}]\|(\tilde{\pi}^{(i-1)}-\tilde{\pi}^{*(i-1)},Z^{(i-1)}-Z^{*(i-1)})\|_{Y^{(i-1)}([0,\min\{T_{i-1},\trip
R_{0}-R^{*}_{0}\trip^{-1}\}])}\\&+[1+(A+1)^{2}(C\delta)^{2}+(C\delta)^{2\sigma}]
T_{i}^{-1}+\trip R_{0}-R^{*}_{0}\trip\\
\end{split}\end{equation}
\begin{equation}\nonumber
\min\{T_{i},\trip
R_{0}-R^{*}_{0}\trip^{-1}\}\sup_{t\in[0,\min\{T_{i},\trip
R_{0}-R^{*}_{0}\trip^{-1}\}]}|\tilde{\pi}^{(i)}-\tilde{\pi}^{*(i)}|(t)<A
\end{equation}
The proof of this proceeds inductively, assuming
\begin{equation}\nonumber
\sup_{0\leq j\leq i-1}\min\{T_{j},\trip
R_{0}-R^{*}_{0}\trip^{-1}\}\sup_{t\in[0,\min\{T_{j},\trip
R_{0}-R^{*}_{0}\trip^{-1}\}]}|\tilde{\pi}^{(j)}-\tilde{\pi}^{*(j)}|(t)<A
\end{equation}
Playing the usual iteration game, it is then easy to see that one
can retrieve the latter inequality for $j=i$ from the foregoing
inequality. To prove the first inequality, one writes the
difference equation satisfied by $Z-Z^{*}$ in the following
fashion:
\begin{equation}\label{diff2*}\begin{split}
&i\partial_{t}(Z^{(i)}-Z^{*(i)})
-\Hil(\pi^{(i-1)})(Z^{(i)}-Z^{*(i)})=[i\dot{\tilde{\pi}}^{(i)}\partial_{\pi}W(\pi^{(i-1)})-i\dot{\tilde{\pi}}^{*(i)}\partial_{\pi}W(\pi^{*(i-1)})]
\\&\hspace{2cm}+[N(Z^{(i-1)},\pi^{(i-1)})-N(Z^{*(i-1)},\pi^{*(i-1)})] +
[\Hil(\pi^{*(i-1)})-\Hil(\pi^{(i-1)})]Z^{*(i)}\\
\end{split}\end{equation}
\begin{equation}\label{diff1*}\begin{split}
&i[\dot{\tilde{\pi}}^{(i)}-\dot{\tilde{\pi}}^{*(i)}]=\la
Z^{(i)}(t),i\dot{\tilde{\pi}}^{(i-1)}\tilde{\calS}_{\ell}(\pi^{(i-1)})(t)\ra
-\la N(Z^{(i-1)},\pi^{(i-1)}),\tilde{\xi}_{\ell}(\pi^{(i-1)})(t)\ra-
\\&(\la
Z^{*(i)}(t),i\dot{\tilde{\pi}}^{*(i-1)}\tilde{\calS}_{\ell}(\pi^{*(i-1)})(t)\ra
-\la
N(Z^{*(i-1)},\pi^{*(i-1)}),\tilde{\xi}_{\ell}(\pi^{*(i-1)})(t)\ra)\\
\end{split}\end{equation}
The estimation on the interval $[0,\min\{T_{i},\trip
R_{0}-R^{*}_{0}\trip^{-1}\}]$ follows then almost verbatim the
proof of Proposition~\ref{pi-pi}. Note that the paths
$\pi,\pi^{*}$ lead to compatible norms on this interval.
\end{proof}
We have almost completed the proof of Theorem~\ref{thm:convergence}.
All that is left is the uniqueness part. For this, consider the
solutions $(\pi,Z)\in\tilde{X}_{*}(\pi)$, $(\pi^{*},Z^{*})\in
\tilde{X}_{*}(\pi^{*})$ with identical initial data:
\begin{equation}\nonumber
Z(0)=\bm R_{0}\\
\bar{R}_{0}\endm
+h^{}f^{+}(\alpha_{\infty})+\sum_{k=1}^{4}a_{k}\tilde{\eta}_{k}(\alpha_{\infty})
=Z^{*}(0)=\bm R_{0}\\
\bar{R}_{0}\endm
+h^{*}f^{+}(\alpha^{*}_{\infty})+\sum_{k=1}^{4}a^{*}_{k}\tilde{\eta}_{k}(\alpha^{*}_{\infty})]
\end{equation}
We study the difference equation for $Z-Z^{*}$, $\pi-\pi^{*}$. On
account of the a priori bounds, the two paths
$y(t)=\int_{0}^{t}v(s)ds+D(t)$,
$y^{*}(t)=\int_{0}^{t}v^{*}(s)ds+D^{*}(t)$ differ by $A$ (say) on
$[0,T]$ where $T=\delta^{-1}$. We estimate (the norm $\|.\|_{Y}$
here is either with respect to $\pi$ or $\pi^{*}$)
\begin{equation}\nonumber
\|(\tilde{\pi}-\tilde{\pi}^{*},Z-Z^{*})\|_{Y([0,T+10])}\lesssim
\delta \|(\tilde{\pi}-\tilde{\pi}^{*},Z-Z^{*})\|_{Y([0,T])}+\delta
T^{-1}
\end{equation}
Implicit in the equation are the constants $A,C$ from before,
which we assume to be chosen once and for all. From the above, we
obtain in particular that
\begin{equation}\nonumber
\la T+10\ra \sup_{t\in[0,T+10]}|\tilde{\pi}-\tilde{\pi}^{*}|(t)<A
\end{equation}
Now one repeats the same argument with $T+10$ instead $T$ etc. The
conclusion is that
\begin{equation}\nonumber
\|(\tilde{\pi}-\tilde{\pi}^{*},Z-Z^{*})\|_{Y([0,T])}=0 \text{\ \ for all\ \ } T>0,
\end{equation}
whence the two solutions agree.
\end{proof}

\begin{proof}(Theorem~\ref{thm:main}) In light of
theorem~\ref{thm:convergence}, we are almost done, we only need to
verify the scattering statement. Given a solution $(\pi,Z)$ on
$[0,\infty)$ constructed as in the preceding proof and  with
$Z(0)$ as given by \eqref{initial data}, we define
\begin{equation}\nonumber
\Phi(R_0)=\text{upper entry
of}\,\,[h^{}f^{+}(\alpha_{\infty})+\sum_{k=1}^{4}a_{k}\tilde{\eta}_{k}(\alpha_{\infty})]
\end{equation}
Define as usual
$U(t)=M_{\infty}(\pi)(t)\calG_{\infty}(\pi)(t)Z(t)$. We first seek
a representation of the form
\begin{equation}\nonumber
U(t)=e^{-it\Hil(\alpha_{\infty})}U_{1}+o_{L^{2}}(1)
\end{equation}
for a suitable $U_{1}\in P_{s}(L^{2}(\R))$. Define
\begin{equation}\nonumber
U_{1}:=
U_{dis}(0)-i\int_{0}^{\infty}e^{ir\Hil(\alpha_{\infty})}P_{s}(\alpha_{\infty})[F(r)]\,dr,
\end{equation}
where
\begin{equation}\nonumber
F(t):=i\partial_{t}U(t)-\Hil(\alpha_{\infty})U(t)
\end{equation}
Clearly, we then have
\begin{equation}\nonumber
U(t)-e^{-it\Hil(\alpha_{\infty})}U_{1}=U_{root}+U_{hyp}-i\int_{t}^{\infty}e^{ir\Hil(\alpha_{\infty})}P_{s}(\alpha_{\infty})[F(r)]\,dr
\end{equation}
On the other hand, the preceding estimates as well as the fact
that both $U_{root},U_{hyp}$ are local and satisfy suitable
$L^{\infty}$ decay estimates imply that
\begin{equation}\nonumber
U_{root}+U_{hyp}-i\int_{t}^{\infty}e^{ir\Hil(\alpha_{\infty})}P_{s}(\alpha_{\infty})[F(r)]\,dr=o_{L^{2}}(1)
\end{equation}
It remains to show that one has scattering for the evolution of
$\Hil(\alpha_\infty)$. This is a standard Cook's method argument.
Indeed, write
\[ \Hil(\alpha_\infty) = \bm -\Laplace+\alpha_\infty^2 & 0 \\ 0 & \Laplace -\alpha_\infty^2 \endm
+ \bm -(\sigma+1)\phi_\infty^2 & -\sigma\phi_\infty^2 \\
\sigma\phi_\infty^2 & (\sigma+1)\phi_\infty^2 \endm =:
\Hil_0(\alpha_\infty) + V,
\]
where $\phi_\infty:=\phi(\cdot,\alpha_\infty)$. Then
\[
e^{-it\Hil(\alpha_\infty)} U_1 = e^{-it\Hil_0(\alpha_\infty)}
U_1 -i \int_0^t e^{-i(t-s)\Hil_0(\alpha_\infty)} V e^{-is
\Hil(\alpha_\infty)} U_1\, ds\] and thus \beeq \label{eq:scat_2}
 U_2:=\lim_{t\to\infty} e^{it\Hil_0(\alpha_\infty)} e^{-it\Hil(\alpha_\infty)} U_1
\eneq exists as a strong $L^2$ limit. Indeed, this follows from
\[ \int_0^\infty \|e^{is\Hil_0(\alpha_\infty)} V e^{-is \Hil(\alpha_\infty)} U_1\|_2\, ds
\les \int_0^\infty \|\la x\ra^{-\theta}e^{-is \Hil(\alpha_\infty)}
U_1\|_{\infty}\, ds
\]
with the latter integral being controlled as follows:
\begin{equation}\nonumber\begin{split}
& \|\la x\ra^{-\theta}e^{-is \Hil(\alpha_\infty)} U_1\|_{\infty}
=\big\|\la x\ra^{-\theta}\big[e^{-is
\Hil(\alpha_\infty)}U_{dis}(0)-i\int_{0}^{\infty}e^{i(r-s)
\Hil(\alpha_{\infty})}P_{s}(\alpha_{\infty})F(r)\,dr\big]\big\|_{\infty}\\
\end{split}\end{equation}
The first term on the right is controlled by
Corollary~\ref{cor:weightp} as well as Sobolev's inequality:
\begin{equation}\nonumber
\|\la x\ra^{-\theta}e^{-is
\Hil(\alpha_\infty)}U_{dis}(0)\|_{L_{x}^{\infty}}\lesssim \la
s\ra^{-1-\epsilon}
\end{equation}
For the integral term we rewrite it as
\begin{equation}\nonumber\begin{split}
&\Big\|\la
x\ra^{-\theta}\int_{0}^{\infty}e^{i(r-s)\Hil(\alpha_{\infty})}P_{s}(\alpha_{\infty})[F(r)]\,dr\Big\|_{\infty}\\
&=\Big\|\la
x\ra^{-\theta}\int_{0}^{\infty}e^{i(r-s)\Hil(\alpha_{\infty})}P_{s}(\alpha_{\infty})[i\dot{\pi}\partial_\pi
\tilde W_\infty(\pi)(r) +N_\infty(U,\pi)(r) + V_\infty U(r)]\,dr\Big\|_{\infty}\\
\end{split}\end{equation}
Proceeding as in the proof of Proposition~\ref{(U,pi)} etc., and
keeping the a priori estimates in mind, we easily bound this
integral by
\begin{equation}\nonumber
\lesssim\int_{0}^{\infty}\min\{|r-s|^{-1-\epsilon},|r-s|^{-\frac{1}{2}}\}\la
r\ra^{-1-\epsilon}\,dr
\end{equation}
Putting the preceding observations together, we obtain that
\begin{equation}\nonumber
\int_0^\infty \|\la x\ra^{-\theta}e^{-is \Hil(\alpha_\infty)}
U_1\|_{\infty}\, ds \lesssim \int_{0}^{\infty}[\la
s\ra^{-1-\epsilon}+\int_{0}^{\infty}\min\{|r-s|^{-1-\epsilon},|r-s|^{-\frac{1}{2}}\}\la
r\ra^{-1-\epsilon}\,dr]\,ds<\infty,
\end{equation}
as desired. It follows that
\[ U(t) = e^{-it\Hil_0(\alpha_\infty)} U_2 + o_{L^2}(1). \]
Finally,
\[ Z(t) = \calG_\infty(t)^{-1} M(t)^{-1}U(t) = e^{-it\Hil_0} \calG_\infty^{-1}(0) U_2 + o_{L^2}(1),\]
where $\Hil_0=\bm -\Laplace & 0\\ 0 & \Laplace\endm$. Setting
$\calG_\infty^{-1}(0) U_2 = \binom{f_0}{\bar{f}_0}$ and $Z(t)=
\binom{R(t)}{\bar{R}(t)}$, we obtain
\[ R(t) = e^{it\Laplace} f_0 +  o_{L^2}(1),\]
and the theorem is proved.
\end{proof}

\begin{proof}(Proof of Theorem~\ref{thm:main2})
The idea is as follows: Given $\alpha_0$, consider the
NLS~\eqref{eq:NLS} with initial data $\phi(\cdot,\alpha_0)+R_0$.
Applying the usual four-parameter family of symmetries (Galilei
giving three parameters, scaling one --- scaling here is the same
as the parameter $\alpha$), we transform this to $W(0,\cdot)+R_1$
where $W(0,x)$ is a soliton with a general parameter vector
$\pi_0$ which is close to $(0,0,0,\alpha_0)$. Hence, we can apply
Theorem~\ref{thm:main} to conclude that these initial data will
give rise to global solutions with the desired properties as long
as $W(0,x)+R_1$ lies on the stable manifold associated
with~$W(0,x)$. To prove that we obtain four dimensions back in
this fashion requires checking that the derivatives of $W(0,x)$ in
its parameters are transverse to the linear space $\calS$ of
Theorem~\ref{thm:main}. However, these derivatives are basically
the elements of the root space $\calN$ of $\Hil(\alpha_0)$,
whereas we know that $\calS$ is perpendicular to the root space
$\calN^*$ of $\Hil(\alpha_0)^*$. More precisely, it is easily
verified that these derivatives are
\begin{align}\nonumber
\bm\partial_{\alpha}\phi\\\partial_{\alpha}\phi\endm,\,\bm
i\phi\\i\phi\endm,\,\bm
ix\phi\\-ix\phi\endm,\,\bm\partial_{x}\phi\\\partial_{x}\phi\endm
\end{align}
But Lemma~\ref{lem:orth} implies that no nonzero vector in $\calN$
is perpendicular to $\calN^*$, which proves that $\calN$ is
transverse to $\calS$, as desired.
\end{proof}

\section{The scattering theory for Schr\"odinger systems}
\label{sec:BP}

This section presents the scattering theory for matrix
Hamiltonians on the line which was developed by Buslaev and
Perelman~\cite{BP1}. Since the presentation in~\cite{BP1} is
somewhat sketchy, and since we need to refine some of the
estimates in~\cite{BP1}, we give full details.

\begin{defi}
\label{def:matrix} In what follows,
\[
\Hil_0 := \bm -\partial_{xx}+1 & 0 \\ 0 & \partial_{xx}-1 \endm,
\; V(x) = \bm V_{1}(x) & V_{2}(x) \\ -V_{2}(x) & -V_{1}(x) \endm,
\; \Hil:=\Hil_0+V.
\]
We will assume that $V$ as well as all its derivatives are
exponentially decaying: \beeq \label{eq:Vdec} \|V^{(k)}(x)\|\le
C_k\,e^{-\gamma |x|} \quad \forall \,k\ge0 \eneq with some
$0<\gamma<1$. Moreover, all entries of $V$ are real-valued, and we
will also assume that $V$ is even: $V(x)=V(-x)$.
\end{defi}

The decay and regularity assumptions can be relaxed to polynomial
decay and a finite number of derivatives, but we do not dwell on
this issue. Let the usual Pauli matrices be given by
\[
\sigma_1=\bm 0&1\\1&0\endm, \;\sigma_2=\bm 0&-i\\i&0\endm,\;
\sigma_3=\bm 1&0\\0&-1\endm.
\]
Note that any $V$ as in Definition~\ref{def:matrix} satisfies
\beeq \label{eq:pauli} \sigma_3 \Hil^* \sigma_3 = \Hil,\;
\sigma_1\Hil\sigma_1=-\Hil. \eneq The following three lemmas
construct a basis of the solution space to $ \Hil
f=(\lambda^2+1)f$ with prescribed asymptotics at infinity. These
are of course analogues of the Jost solutions in the scalar case.
Throughout, $\mu=\sqrt{\lambda^2+2}$.

\begin{lemma}
\label{lem:f3}
 For every $\lambda\in\R$ there exists a solution
$f_3(x,\lambda)$ of the equation
\[ \Hil f_3(\cdot,\lambda) = (\lambda^2+1) f_3(\cdot,\lambda) \]
with the property that $f_3(x,\lambda)\sim e^{-\mu x}\binom{0}{1}$
as $x\to\infty$. Moreover, $f_3$ is smooth in both variables and
satisfies the estimates \beeq \label{eq:f3est} \Bigl|
\partial_\lambda^\ell\partial_x^{k} \Big[e^{\mu x}f_3(x,\lambda) -
\binom{0}{1}\Big]\Bigr| \leq C_k\,\mu^{-1-\ell} e^{-\gamma x}
\eneq for all $x\ge0$ and $k,\ell\ge0$. Finally,
$\sup_{\lambda\in\R}\sup_{x\in\R} |e^{\mu x}f_3(x,\lambda)|\le
C(V).$
\end{lemma}
\begin{proof}
We set \beeq \label{eq:f3def}
f_3(x,\lambda):= e^{-\mu x}\binom{0}{1}+ \int_x^\infty \bm \frac{\sin(\lambda(y-x))}{\lambda} & 0 \\
0 & -\frac{\sinh(\mu(y-x))}{\mu} \endm V(y) f_3(y,\lambda)\, dy.
\eneq Equivalently, with \beeq \label{eq:f3K}
 K(x,y;\lambda) := \bm \frac{\sin(\lambda(y-x))}{\lambda} & 0 \\
0 & -\frac{\sinh(\mu(y-x))}{\mu} \endm e^{\mu(x-y)}, \eneq we have
\beeq \label{eq:volf3} e^{\mu x}f_3(x,\lambda)= \binom{0}{1}+
\int_x^\infty K(x,y;\lambda) V(y) e^{\mu y} f_3(y,\lambda)\, dy.
\eneq Since for all $\lambda\in\R$
\[ \sup_{y\ge x} | K(x,y;\lambda)| \le \sup_{y\ge x} \Big\{
(y-x) e^{\mu(x-y)} + \frac{\sinh(\mu(x-y))}{\mu} e^{\mu(x-y)}
\Big\} \le C\mu^{-1},
\]
with a universal constant $C$, we conclude that $e^{\mu
x}f_3(\cdot,\lambda)$ solves a Volterra integral equation and thus
\beeq \label{eq:f3bd}
 \sup_{\lambda \in\R}\sup_{x\in\R} |e^{\mu
x}f_3(x,\lambda)| \le C(V). \eneq Thus, we obtain that
\[ \Big|e^{\mu x}f_3(x,\lambda) - \binom{0}{1}\Big| \le
\int_x^\infty C\mu^{-1}\, e^{-\gamma y}\, dy \le C\mu^{-1}
e^{-\gamma x}\] for all $x\ge0$. The estimate~\eqref{eq:f3est}
follows by differentiating the Volterra equation~\eqref{eq:volf3}.
Indeed, since $K(x,x,\lambda)=0$ and $\partial_x
K(x,y;\lambda)=-\partial_y K(x,y;\lambda)$, integration by parts
yields
\[ \partial_x (e^{\mu x}f_3(x,\lambda)) = \int_x^\infty
K(x,y;\lambda) V'(y) e^{\mu y}f_3(y,\lambda)\, dy + \int_x^\infty
K(x,y;\lambda)V(y)\partial_y(e^{\mu y}f_3(y,\lambda))\, dy.\] By
the usual estimates for Volterra equations as well
as~\eqref{eq:Vdec} and~\eqref{eq:f3bd},
\[ |\partial_x (e^{\mu x}f_3(x,\lambda))| \le C\mu^{-1} \int_x^\infty
e^{-\gamma y}\, dy \le C\,\mu^{-1} e^{-\gamma x},
\]
for $x\ge0$, as claimed. The higher derivatives in~$x$ follow in a
similar fashion. Indeed, integrating by parts one verifies
inductively that \beeq \label{eq:f3xkder}
\partial_x^k  (e^{\mu x}f_3(x,\lambda)) = \sum_{j=0}^k \binom{k}{j} \int_x^\infty V^{(k-j)}(y)\partial_y^j[e^{\mu y}f_3(y,\lambda)]\, dy,
\eneq which implies the bounds
\[ |\partial_x^k  (e^{\mu x}f_3(x,\lambda))| \le C_k\, \mu^{-1}e^{-\gamma x}\]
for all $x\ge0$. As far as the derivatives in $\lambda$ are
concerned, differentiating \eqref{eq:f3K} in~$\lambda$ reveals
that
\[ \sup_{y\ge x}|\partial_\lambda^\ell K(x,y;\lambda)| \le C_\ell\, \mu^{-\ell-1}
\]
for all $\ell\ge0$. Apply $\partial_\lambda^\ell$
to~\eqref{eq:f3xkder}. Induction in $\ell$ implies the
estimate~\eqref{eq:f3est}.
\end{proof}

Next, we find a pair of oscillatory solutions.

\begin{lemma}
\label{lem:f1} For all $\lambda\in\R$ there exist solutions
$f_1(\cdot,\lambda), f_2(\cdot,\lambda)$ of
\[ \Hil f_j(\cdot,\lambda) = (1+\lambda^2) f_j(\cdot,\lambda)\]
($j=1,2$) and with the property that
$f_2(\cdot,\lambda)=\overline{f_1(\cdot,\lambda)}$ and
\[
f_1(x,\lambda) = \binom{e^{ix\lambda}}{0} + O(\mu^{-1}\,e^{-\gamma
x})
\]
as $x\to\infty$. The constant in the $O$-term is uniform in
$\lambda\in \R$. Moreover, $f_1$ is smooth in both variables and
there exists a constant\footnote{This constant becomes large as
$V$ becomes large.} $x_0\ge0$ only depending on $V$ such that
\beeq \label{eq:f1der} \Big|\partial_{\lambda}^\ell\partial_x^k
\Big[e^{-i\lambda x}f_1(x,\lambda)-\binom{1}{0}\Big]\Big| \le
C_{k,\ell}\,\mu^{-1+k}x^\ell\,e^{-\gamma x} \eneq for all
$k,\ell\ge0$, $\lambda\in\R$,  and $x\ge x_0$. Finally, the same
bound holds for all $x\ge0$ provided $|\lambda|\ge\lambda_0$ where
$\lambda_0\ge0$ is some  constant${}^1$ that only depends on $V$.
\end{lemma}
\begin{proof} Since $V$ has real entries, any solution $f_1(\cdot,\lambda)$ gives rise to another
solution $\overline{f_2(\cdot,\lambda)}$. Hence it will suffice to
find $f_1$. We seek a solution of the form \beeq \label{eq:f1ans}
 f_1(x,\lambda) = \binom{1}{0}v(x,\lambda) + f_3(x,\lambda) u(x,\lambda)
\eneq where $v(x,\lambda)\sim e^{ix\lambda}$,  and
$f_3(x,\lambda)=\binom{f_3^{(1)}(x,\lambda)}{f_3^{(2)}(x,\lambda)}$
is as in Lemma~\ref{lem:f3}. Clearly,
\begin{align*}
\Hil( f_3(\cdot,\lambda) u) &= u \Hil(f_3(\cdot,\lambda)) +
\binom{-\partial_{xx}(uf_3^{(1)}(\cdot,\lambda))+u\partial_{xx}
f_3^{(1)}(\cdot,\lambda)}{\partial_{xx}(uf_3^{(1)}(\cdot,\lambda)-u\partial_{xx} f_3^{(1)}(\cdot,\lambda)}  \\
&= u (1+\lambda^2)f_3(\cdot,\lambda)) +
\binom{-\partial_{xx}(uf_3^{(1)}(\cdot,\lambda))+u\partial_{xx}
f_3^{(1)}(\cdot,\lambda)}{\partial_{xx}(uf_3^{(1)}(\cdot,\lambda)-u\partial_{xx}
f_3^{(1)}(\cdot,\lambda)}.
\end{align*}
In order to have $\Hil
f_1(\cdot,\lambda)=(1+\lambda^2)f_1(\cdot,\lambda)$, we therefore
need
\begin{align}
 0 &= (\Hil-(\lambda^2+1)) f_1(\cdot,\lambda) \nonumber \\
 & =
\binom{(-\partial_{xx}-\lambda^2+V_{11}) v(\cdot,\lambda)}{V_{21}
v(\cdot,\lambda)} +\binom{-\partial_{xx} u(\cdot,\lambda)
f_3^{(1)}(\cdot,\lambda)-2\partial_x u(\cdot,\lambda)
\partial_x f_3^{(1)}(\cdot,\lambda)}{\partial_{xx} u(\cdot,\lambda)
f_3^{(2)}(\cdot,\lambda)+2\partial_x u(\cdot,\lambda) \partial_x
f_3^{(2)}(\cdot,\lambda)} \label{eq:f1one}
\end{align}
The homogeneous equation
\[ y''\, f_3^{(2)}(\cdot,\lambda) + 2 y'\,
\partial_x f_3^{(2)}(\cdot,\lambda) =0\]
has the solution
\[ y'(x) = C\, (f_3^{(2)}(x,\lambda))^{-2}\]
which is well-defined provided $x\ge x_0$ by Lemma~\ref{lem:f3}.
Here $x_0\ge0$ is a large constant independent of~$\lambda$ (for
large $|\lambda|$, we can take $x_0=0$). From the second
coordinate in~\eqref{eq:f1one} and the usual "variation of
constants" method we obtain that
\[ C'(x)= -f_3^{(2)}(x,\lambda)V_{21}(x)v(x,\lambda) \]
which implies, together with the boundary condition $C(\infty)=0$,
that \beeq \label{eq:u'eq} u'(x,\lambda) =
[f_3^{(2)}(x,\lambda)]^{-2}\int_x^\infty f_3^{(2)}(y,\lambda)\,
V_{21}(y) v(y,\lambda)\,dy. \eneq From the first coordinate
in~\eqref{eq:f1one} we conclude that (dropping $\lambda$ from
$f_3,v,u$ for simplicity)
\begin{align*}
 v(x) &= e^{ix\lambda} - \int_x^\infty
\frac{\sin(\lambda(y-x))}{\lambda}
\big(u''(y)f_3^{(1)}(y)+2u'(y)f_3^{(1)}(y)'-V_{11}(y)v(y)\big)\,
dy \\
&= e^{ix\lambda} + \int_x^\infty
\frac{\sin(\lambda(y-x))}{\lambda}
\big(V_{21}(y)\frac{f_3^{(1)}}{f_3^{(2)}}(y)+V_{11}(y)\big)v(y)\,
dy \\
& - 2 \int_x^\infty \frac{\sin(\lambda(y-x))}{\lambda}
\Big(-\frac{{f_3^{(2)}}(y)'}{f_3^{(2)}(y)}\,f_3^{(1)}(y)+
f_3^{(1)}(y)'\Big) [f_3^{(2)}(y)]^{-2}\int_y^\infty f_3^{(2)}(z)
V_{21}(z) v(z)\,dz\,dy.
\end{align*}
This can be written as \beeq
 \label{eq:veq}
 v(x,\lambda) = e^{ix\lambda} +
\int_x^\infty K(x,y;\lambda)v(y,\lambda)\, dy \eneq where
$K(x,y;\lambda):= K_1(x,y;\lambda)+K_2(x,y;\lambda)$ according to
the splitting
\begin{align}
& K_1(x,y;\lambda) =
\frac{\sin(\lambda(y-x))}{\lambda}\Big(V_{21}(y)\frac{f_3^{(1)}}{f_3^{(2)}}(y;\lambda)+V_{11}(y)\Big) \label{eq:K1}\\
& K_2(x,y;\lambda) = \label{eq:K2}\\
& -2 \int_x^y \frac{\sin(\lambda(z-x))}{\lambda}
\Big(-\frac{{f_3^{(2)}}(z;\lambda)'}{f_3^{(2)}(z;\lambda)}\,f_3^{(1)}(z;\lambda)+
f_3^{(1)}(z;\lambda)'\Big) [f_3^{(2)}(z;\lambda)]^{-2}\, dz\,
f_3^{(2)}(y;\lambda) V_{21}(y). \nonumber
\end{align}
Since $y\ge x\ge x_0\ge0$, Lemma~\ref{lem:f3} implies that
\begin{align*}
 |K_1(x,y;\lambda)| &\le C\,(y-x)(1+|\lambda|(y-x))^{-1} e^{-\gamma y} \\
 |K_2(x,y;\lambda)| &\le \int_x^y (z-x)(1+|\lambda|(z-x))^{-1}e^{(\mu-\gamma)z}\, dz\, e^{-(\mu+\gamma)z}
 \le C \frac{y-x}{1+|\lambda|(y-x)} \mu^{-1}e^{-2\gamma y}.
\end{align*}
In the final estimate we used that $\frac{a}{1+|\lambda|a}$ is
increasing in $a$. Hence, \eqref{eq:veq} is a Volterra equation
with a solution $v(x,\lambda)$ on the interval $[x_0,\infty)$
satisfying
\begin{align*}
 |v(x,\lambda)-e^{ix\lambda}| &\le C\int_x^\infty \frac{y-x}{1+|\lambda|(y-x)}\, e^{-\gamma y}\, dy = Ce^{-\gamma x}\int_0^\infty \frac{u}{1+|\lambda|u}\,e^{-\gamma u}\, du\\
&\le C(1+|\lambda|)^{-1}e^{-\gamma x}
\end{align*}
for all $x\ge x_0$. Thus, in view of \eqref{eq:u'eq},
\[ |u'(x,\lambda)| \le C\, e^{2\mu x} \int_x^\infty e^{-(\mu+\gamma)
y} \, dy \le C\,\mu^{-1}\, e^{(\mu-\gamma)x}
\] for all $x\ge x_0$ and $\lambda\in\R$. Hence, assuming that
$u(x_0,\lambda)=0$, we obtain that
\[ |u(x,\lambda)| \le C \mu^{-2}\, e^{(\mu-\gamma)x}. \]
By the preceding,
\[ \Big|f_1(x,\lambda) - \binom{e^{ix\lambda}}{0}\Big| \le
C\mu^{-1}\,e^{-\gamma x}
\]
with a constant $C$ which is uniform in $x\ge x_0$ and
$\lambda\in\R$. Now continue $f_1(\cdot,\lambda)$ to the left of
$x_0$ by means of the existence and uniqueness theorem.

As far as the derivatives are concerned, set
\[
\tilde{K}(x,y;\lambda):=e^{-i\lambda(x-y)}K(x,y;\lambda), \quad
\tilde{v}(x,\lambda):=e^{-i\lambda x} v(x,\lambda),
\]
and similarly with $\tilde K_1, \tilde K_2$, see
\eqref{eq:K1},\eqref{eq:K2}. Then \beeq \label{eq:volt2}
 \partial_\lambda \tilde v(x,\lambda) = \int_x^\infty \partial_\lambda \tilde K(x,y;\lambda)\tilde v(y,\lambda)\, dy
+\int_x^\infty \tilde K(x,y;\lambda)\partial_\lambda\tilde
v(y,\lambda)\, dy. \eneq In view of \eqref{eq:f3est},
\[ |\partial_\lambda f_3(x,\lambda)| \le C\,x e^{-\mu x} \]
for large $x$. Consequently, for $y\ge x\ge x_0$,
\begin{align*}
|\partial_\lambda \tilde K_1(x,y;\lambda)| &=
\Big|\partial_\lambda \Big[e^{-i\lambda(x-y)}
\frac{\sin(\lambda(y-x))}{\lambda}\Big(V_{21}(y)\frac{f_3^{(1)}}{f_3^{(2)}}(y;\lambda)+V_{11}(y)\Big)\Big]\Big|  \\
& \le C\frac{(y-x)^2}{1+|\lambda|(y-x)}e^{-\gamma y} +
C\mu^{-1}\,\frac{y-x}{1+|\lambda|(y-x)}ye^{-2\gamma y}.
\end{align*}
To obtain this bound, it is helpful to introduce
$\phi(u):=\frac{\sin (u)}{u}$. Then $|\phi^{(k)}(u)|\le C_k\,
(1+|u|)^{-k-1}$ for all $k\ge0$ and
\[ \frac{\sin(\lambda(y-x))}{\lambda} = (y-x)\phi(\lambda(y-x)). \]
Similarly,
\begin{align*}
&|\partial_\lambda \tilde K_2(x,y;\lambda)| \\
&=\Big|\partial_\lambda \int_x^y
e^{-i\lambda(x-y)}\frac{\sin(\lambda(z-x))}{\lambda}
\Big(-\frac{{f_3^{(2)}}(z;\lambda)'}{f_3^{(2)}(z;\lambda)}\,f_3^{(1)}(z;\lambda)+
f_3^{(1)}(z;\lambda)'\Big) [f_3^{(2)}(z;\lambda)]^{-2}\, dz\,
f_3^{(2)}(y;\lambda)
V_{21}(y) \Big| \\
&\le C\int_x^y \frac{(y-x)^2}{1+|\lambda|(y-x)} e^{-(\mu+\gamma)z} e^{2\mu z}\, dz e^{-(\mu+\gamma)y} \\
& + C\int_x^y \frac{y-x}{1+|\lambda|(y-x)} (1+\mu z)\mu^{-1}e^{-(\mu+\gamma)z} e^{2\mu z}\, dz e^{-(\mu+\gamma)y} \\
&\le C \frac{(y-x)^2}{1+|\lambda|(y-x)} \mu^{-1}e^{-2\gamma y} + C
\frac{y-x}{1+|\lambda|(y-x)} \mu^{-1}y e^{-2\gamma y}.
\end{align*}
The conclusion is that
\[
\Big| \int_x^\infty \partial_\lambda \tilde K(x,y;\lambda)\tilde
v(y,\lambda)\, dy\Big| \le C\, \mu^{-1}e^{-\gamma x}
\]
and therefore also
\[ |\partial_\lambda\tilde v(x,\lambda)|\le C\mu^{-1}e^{-\gamma x},\]
see \eqref{eq:volt2}. Inserting this into \eqref{eq:f1ans} yields
\begin{align*}
|\partial_\lambda[e^{-ix\lambda} f_1(x,\lambda)]| \le
C\mu^{-1}xe^{-\gamma x}
\end{align*}
for all $x\ge x_0$. The case of higher derivatives in $\lambda$ is
similar. Indeed, \beeq \label{eq:volt_ell}
 \partial_\lambda^\ell \tilde v(x,\lambda) = \sum_{j=0}^\ell \binom{\ell}{j} \int_x^\infty \partial_\lambda^j \tilde K(x,y;\lambda)
\partial_\lambda^{\ell-j}\tilde v(y,\lambda)\, dy.
\eneq As before, Lemma~\ref{lem:f3} leads to the bounds
\[
|\partial_\lambda^\ell \tilde K(x,y;\lambda)| \le
C_\ell\,\frac{(y-x)^{\ell+1}+y-x}{1+|\lambda|(y-x)}e^{-\gamma y}.
\]
In combination with \eqref{eq:volt_ell}, this estimate inductively
yields
\[ |\partial_\lambda^\ell \tilde v(x,\lambda)|\le C_\ell\, \mu^{-1}e^{-\gamma x} \text{\ \ or\ \ }
|\partial_\lambda^\ell v(x,\lambda)|\le C_\ell\, x^\ell
\mu^{-1}e^{-\gamma x}\] for large $x$. Inserting this bound into
the defining equation \eqref{eq:u'eq} for $u'$ implies that
\[
|\partial_\lambda^\ell u'(x,\lambda)| \le C_\ell\, x^\ell
\mu^{-1}e^{(\mu-\gamma)x}.
\]
Since $\partial_\lambda^\ell u(x_0,\lambda)=0$ for all $\lambda$,
\[
|\partial_\lambda^\ell u(x,\lambda)| \le C_\ell\, x^\ell
\mu^{-2}e^{(\mu-\gamma)x}.
\]
Finally, inserting this estimate into \eqref{eq:f1ans} we obtain
\[
|\partial_\lambda^\ell [e^{-ix\lambda} f_1(x,\lambda)]| \le
C_\ell\, x^\ell\mu^{-1}\, e^{-\gamma x},
\]
which is \eqref{eq:f1der} with $k=0$. The case of $x$-derivatives,
i.e., $k\ge1$ in~\eqref{eq:f1der} follows by similar
considerations. We skip the details.
\end{proof}

\begin{remark}
\label{rem:x0} As already noted, we can take $x_0=0$ in the
previous proof for large $|\lambda|$. This allows us to state that
\[
\sup_{x\ge 0} \Big|
\partial_x f_1(x,\lambda)-i\lambda\binom{e^{ix\lambda}}{0}
\Big| \le C,
\]
for large $|\lambda|$, which will be useful later. Another
important (but simple) observation concerns the point $\lambda=0$.
There $f_1(\cdot,\lambda), f_2(\cdot,\lambda)$ are identical.
However, it is simple to obtain a pair of linearly independent
solutions at $\lambda=0$. Indeed, just take $f_1(\cdot,0)$ and
$\partial_\lambda f_1(\cdot,0)$. Note that the asymptotic behavior
of $\partial_\lambda f_1(\cdot,0)$ is $ix$ as $x\to\infty$.
Alternatively, one can also work with the pair
\[ f_1(\cdot,\lambda),\qquad \frac{f_1(\cdot,\lambda)-f_2(\cdot,-\lambda)}{\lambda} \]
which is independent for all $\lambda\in\R$.
\end{remark}

Next, we construct an exponentially growing solution at $+\infty$.
We will later modify $\tilde f_4$ to obtain $f_4$, hence the
notation.

\begin{lemma}
\label{lem:f4} There exists a solution $\tilde f_4(\cdot,\lambda)$
of $\Hil \tilde f_4(\cdot,\lambda)=(1+\lambda^2)\tilde
f_4(\cdot,\lambda)$ with the property that
\[ \tilde f_4(\cdot,\lambda) = e^{\mu x}\binom{0}{1} +
O\big((1+|\lambda|)^{-1}e^{(\mu-\gamma)x}\big)
\]
as $x\to\infty$. The constant in the $O$-term is uniform in
$\lambda$. Moreover, $f_4(\cdot,\lambda)$ is a smooth function of
its arguments and there exists a constant\footnote{This constant
becomes large as $V$ becomes large.}  $x_1\ge 0$ only depending on
$V$ such that \beeq \label{eq:tilf4der}
\Big|\partial_{\lambda}^\ell\partial_x^k \Big[e^{-\mu x}\tilde
f_4(x,\lambda)-\binom{0}{1}\Big]\Big| \le
C_{k,\ell}\,\mu^{-1+k}x^\ell\,e^{-\gamma x} \eneq for all
$k,\ell\ge0$, $\lambda\in\R$,  and $x\ge x_1$. Finally, the same
bound holds for all $x\ge0$ provided $|\lambda|\ge\lambda_1$ where
$\lambda_1\ge0$ is some  constant${}^2$ that only depends on $V$.
\end{lemma}
\begin{proof}
Make the ansatz
\begin{align}
\tilde f_4(x,\lambda) &= e^{\mu x} \binom{0}{1} + \int_x^\infty
\bm 0 & 0
\\ 0 & -\frac{1}{2\mu}e^{\mu(x-y)} \endm V(y) \tilde f_4(y,\lambda)\, dy \nonumber \\
& + \int_{x_1}^x \bm \frac{\sin(\lambda(x-y))}{\lambda} & 0  \\
0 & -\frac{1}{2\mu} e^{-\mu(x-y)} \endm V(y)\tilde
f_4(y,\lambda)\, dy, \label{eq:f4tildef}
\end{align}
where $x_1\ge0$ is some constant that will be chosen large.
Clearly, if $\tilde f_4(x,\lambda)$ grows at most like $e^{\mu x}$
as $x\to\infty$, then this integral equation is well-defined.
Moreover, it is easy to check that a solution of this integral
equation satisfies $\Hil \tilde
f_4(\cdot,\lambda)=(1+\lambda^2)\tilde f_4(\cdot,\lambda)$. To
find a solution, we solve
\begin{align}
f(x,\lambda) &=  \binom{0}{1} + \int_x^\infty \bm 0 & 0
\\ 0 & -\frac{1}{2\mu}\endm V(y) f(y,\lambda)\, dy
\nn \\
& + \int_{x_1}^x \bm \frac{\sin(\lambda(x-y))}{\lambda}e^{-\mu(x-y)} & 0 \\
0 & -\frac{1}{2\mu} e^{-2\mu(x-y)} \endm V(y)f(y,\lambda)\, dy
\label{eq:Tdef}
\end{align}
by the contraction principle. Thus, $\tilde f_4(x,\lambda):=
e^{\mu x}f(x,\lambda)$ will be the desired solution. Denote the
right-hand side by $T$. Then
\begin{align*}
(Tf-Tg)(x) &= \int_x^\infty \bm 0 & 0
\\ 0 & -\frac{1}{2\mu}\endm V(y) [f(y)-g(y)]\, dy
\\
& + \int_{x_1}^x \bm \frac{\sin(\lambda(x-y))}{\lambda}e^{-\mu(x-y)} & 0 \\
0 & -\frac{1}{2\mu} e^{-2\mu(x-y)} \endm V(y)[f(y)-g(y)]\, dy,
\end{align*}
which implies that
\begin{align*}
|(Tf-Tg)(x)| &\le C(1+|\lambda|)^{-1}\int_x^\infty e^{-\gamma y}
|f(y)-g(y)|\, dy
\\
& + C(1+|\lambda|)^{-1}\int_{x_1}^x [\mu (x-y)e^{-\mu(x-y)} +
e^{-2\mu(x-y)}] e^{-\gamma
y} |f(y)-g(y)|\, dy \\
&\le C(1+|\lambda|)^{-1} e^{-\gamma x_1} \sup_{y\ge x_1}
|f(y)-g(y)|
\end{align*}
for all $x\ge x_1$. Hence, $T$ is a contraction in
\[ \big\{ f\in C([x_1,\infty),\Compl^2)\:|\: \sup_{x\ge x_1} |f(x)|\le 2 \big\} \]
provided $x_1$ is large (if $|\lambda|$ is large, we can take
$x_1=0$). If $f(\cdot,\lambda)$ is the fixed-point of $T$, then
\[ \Big|f(x,\lambda)-\binom{0}{1}\Big| \le C(1+|\lambda|)^{-1}e^{-\gamma x},\]
for all $x\ge x_1$. The estimate on the derivatives follows by
differentiating the equation~\eqref{eq:Tdef}.
\end{proof}

\begin{remark}
\label{rem:x02} As noted, for large $|\lambda|$, we can take
$x_1=0$ in the previous proof. This allows us to state that
\[
\Big| \tilde{f}_4(0,\lambda)-\binom{0}{1} \Big| \le
C\,|\lambda|^{-1},
\]
as well as
\[
\Big| \partial_x \tilde{f}_4(0,\lambda)-\binom{0}{\mu} \Big| \le
C,
\]
 for large $|\lambda|$.
\end{remark}

Here we record a useful property of these solutions.

\begin{cor}
\label{cor:signs} The solutions $f_1,f_2,f_3,\tilde f_4$ from the
previous lemmas satisfy
\begin{align*}
f_1(\cdot,-\lambda) &= \overline{f_1(\cdot,\lambda)} =
f_2(\cdot,\lambda),\qquad f_2(\cdot,-\lambda) =
\overline{f_2(\cdot,\lambda)} =
f_1(\cdot,\lambda) \\
f_3(\cdot,-\lambda) &= \overline{f_3(\cdot,\lambda)} =
f_3(\cdot,\lambda),\qquad \tilde f_4(\cdot,-\lambda) =
\overline{\tilde f_4(\cdot,\lambda)} = \tilde f_4(\cdot,\lambda)
\end{align*}
for all $\lambda\in\R$. Moreover, $f_1(\cdot,0)=f_2(\cdot,0)$.
\end{cor}
\begin{proof}
This can be seen by inspecting the various integral equations
defining these solutions. Indeed, since $V$ has real entries,
\eqref{eq:f3def} and \eqref{eq:f4tildef} are invariant under both
conjugation and the substitution $\lambda\to -\lambda$. Finally,
\eqref{eq:f1ans}, \eqref{eq:u'eq}, and \eqref{eq:veq} imply that
conjugation of $f_1(\cdot,\lambda)$ is equivalent to
$\lambda\to-\lambda$.
\end{proof}

The following two lemmas introduce the Wronskian in the matrix
context. We will need to use the property that $\sigma_3 V^*
\sigma_3=V$, see Definition~\ref{def:matrix} and~\eqref{eq:pauli}.

\begin{lemma}
\label{lem:wronski} For two differentiable functions $f,g$ taking
values in $\Compl^2$ let
\[ W[f,g](x) := \la f'(x),g(x)\ra - \la f(x),g'(x) \ra \]
where $\la \cdot,\cdot\ra$ is the {\em real} scalar product.
Suppose that $(\Hil-z)f=0$ and $(\Hil-z)g=0$.
Then
\[ W[f,g] = \const. \]
Moreover, \begin{align*}
W[f_1(\cdot,\lambda),f_2(\cdot,\lambda)]&=2i\lambda,\;
W[f_3(\cdot,\lambda),\tilde f_4(\cdot,\lambda)] = -2\mu,\\
W[f_1(\cdot,\lambda),f_3(\cdot,\lambda)]&=W[f_2(\cdot,\lambda),f_3(\cdot,\lambda)]=0,
\end{align*}
where $f_1,f_2,f_3,\tilde f_4$ are as in
Lemmas~\ref{lem:f3}-\ref{lem:f4}. There exists a unique choice of
$c_1,c_2\in \Compl$ so that \beeq \label{eq:f4def}
 f_4(\cdot,\lambda):=\tilde
f_4(\cdot,\lambda)-c_1(\lambda) f_1(\cdot,\lambda) - c_2(\lambda)
f_2(\cdot,\lambda)
 \eneq satisfies
\[ W[f_1(\cdot,\lambda),f_4(\cdot,\lambda)]= W[f_2(\cdot,\lambda),f_4(\cdot,\lambda)]=0.\]
 Furthermore,
$\bar{f}_4(\cdot,\lambda)=f_4(\cdot,-\lambda)=f_4(\cdot,\lambda)$.
Finally, we also have $W[f_3(\cdot,\lambda),f_4(\cdot,\lambda)] =
-2\mu$ and $c_j(\lambda)=O(\lambda^{-1})$ for $j=1,2$ as
$|\lambda|\to\infty$.
\end{lemma}
\begin{proof}
Compute
\begin{align*}
 \frac{d}{dx} W[f,g](x) &= \la \sigma_3 f'', \sigma_3 g\ra - \la
 \sigma_3 f,\sigma_3 g'' \ra \\
 &= \la (1-z)\sigma_3 f + Vf, \sigma_3 g\ra - \la \sigma_3
 f,(1-z)\sigma_3 g + Vg\ra \\
 &= \la \sigma_3 V f,g\ra - \la V^*\sigma_3 f,g \ra =0.
\end{align*}
The statements about the Wronskians follow from the asymptotics in
the previous lemmas. Now define
\[f_4(\cdot,\lambda)=\tilde f_4(\cdot,\lambda)-c_1(\lambda) f_1(\cdot,\lambda) -
c_2(\lambda) f_2(\cdot,\lambda) \] so that \[
W[f_1(\cdot,\lambda),f_4(\cdot,\lambda)]=
W[f_2(\cdot,\lambda),f_4(\cdot,\lambda)]=0. \] Thus, we need that
\begin{align*}
 0 & = W[f_1(\cdot,\lambda),\tilde
f_4(\cdot,\lambda)]-c_2(\lambda)W[f_1(\cdot,\lambda),f_2(\cdot,\lambda)],\\
0 &= W[f_2(\cdot,\lambda),\tilde
f_4(\cdot,\lambda)]-c_1(\lambda)W[f_2(\cdot,\lambda),f_1(\cdot,\lambda)],
\end{align*}
 and therefore
\[ c_2(\lambda) = (2i\lambda)^{-1}W[f_1(\cdot,\lambda),\tilde f_4(\cdot,\lambda)], \quad c_1(\lambda)=
-(2i\lambda)^{-1} W[f_2(\cdot,\lambda),\tilde f_4(\cdot,\lambda)].
\]
It follows from \eqref{eq:f1der} and \eqref{eq:tilf4der} that
$c_1(\lambda)=O(\lambda^{-1})$ and $c_2(\lambda)=O(\lambda^{-1})$.
By inspection,
\[ \bar{c}_1(\lambda) = (2i\lambda)^{-1} W[\tilde{f}_4(\cdot,\lambda),f_1(\cdot,\lambda)] = c_2(\lambda),\]
which implies that $\bar{f}_4(\cdot,\lambda)= f_4(\cdot,\lambda)$
by Corollary~\ref{cor:signs}. Also,
\[c_1(-\lambda) =
(2i\lambda)^{-1} W[f_2(\cdot,-\lambda),\tilde
f_4(\cdot,-\lambda)]= (2i\lambda)^{-1} W[f_1(\cdot,\lambda),\tilde
f_4(\cdot,\lambda)] = c_2(\lambda),\] and by
Corollary~\ref{cor:signs} again
$f_4(\cdot,-\lambda)=f_4(\cdot,\lambda)$. Finally,
$W[f_3,f_4]=W[f_3,\tilde f_4]=-2\mu$ as claimed.
\end{proof}

We will need the following analogue of \eqref{eq:tilf4der} for
$f_4$.

\begin{cor}
\label{cor:f4der} Let $f_4$ be as in \eqref{eq:f4def}. Then
\[
\Big|\partial_{\lambda}^\ell\partial_x^k \Big[e^{-\mu x}
f_4(x,\lambda)-\binom{0}{1}\Big]\Big| \le
C_{k,\ell}\,\mu^{-1+k}x^\ell\,e^{-\gamma x}
\]
for all $k,\ell\ge0$.
\end{cor}
\begin{proof} This is an immediate consequence of Lemma~\ref{lem:wronski} and \eqref{eq:tilf4der}.
\end{proof}

Recall  that we are assuming that $V$ is {\em even}. In that case,
set
\[ g_j(x,\lambda) = f_j(-x,\lambda) \text{\ \ for\ \ }1\le j\le
4.\] Since $V(x)=V(-x)$, these functions are again solutions of
\[\Hil g_j(\cdot,\lambda)=(1+\lambda^2)g_j(\cdot,\lambda)\]
which have the same asymptotic behavior as $x\to -\infty$ as the
$f_j$ when $x\to\infty$.

\begin{lemma}
\label{lem:matwron}
 Suppose $F,G$ are $2\times 2$ matrix solutions
of $\Hil F= zF$, $\Hil G=zG$, with some $z\in \Compl$. Then the
matrix Wronskian
\[ \calW[F,G](x) := {F'}^t(x) G(x) - F^t(x) G'(x) \]
is independent of $x$. Now suppose that $W[F,F]=0$ or $W[G,G]=0$.
Then
\[ \det \calW[F,G] =0 \]
iff there exist vectors $a,b\in \Compl^2$, not both zero, such
that
\[ F(x)a+G(x)b=0\]
for all $x\in\R$.
\end{lemma}
\begin{proof}
By assumption, $F=[\phi_1\;\phi_2]$, $G=[\psi_1\; \psi_2]$ where
$(\Hil-z)\phi_j=0$ and $(\Hil-z)\psi_j=0$ for $j=1,2$. Hence,
\[
\calW[F,G](x) = \bm W[\phi_1,\psi_1](x) & W[\phi_1,\psi_2](x) \\
W[\phi_2,\psi_1](x) & W[\phi_2,\psi_2](x) \endm.
\]
By Lemma~\ref{lem:wronski}, each of the entries is independent
of~$x$ and thus $\calW[F,G]$. For the second statement, compute
\begin{align}
\bm \calW[F,G] & \calW[F,F] \\ \calW[G,G] & \calW[G,F] \endm &=
\bm {F'}^t & - F^t
\\ {G'}^t & - G^t \endm \bm G & F\\ G' & F' \endm \nonumber\\
&= \bm 0 & I \\ I & 0 \endm \bm G & F\\ G' & F' \endm^t \bm 0 & -I
\\ I & 0 \endm \bm G & F\\ G' & F' \endm. \label{eq:fehs}
\end{align}
Since $\calW[G,F]=-\calW[F,G]^t$ and $\calW[F,F]=0$ (or
$\calW[G,G]=0$) by assumption, we conclude that \beeq
\label{eq:det_sq}
 [\det(\calW[F,G])]^2 = \Big[\det \bm G & F\\ G'
& F'
\endm\Big]^2.
\eneq Now suppose that there exist vectors $a,b\in \Compl^2$, not
both zero, such that
\[ F(x)a+G(x)b=0\]
for all $x\in\R$. Then, clearly,
\[ \bm G(x) & F(x)\\ G'(x) & F'(x) \endm \binom{b}{a} = 0 \]
for every $x\in\R$, and $\det\calW[F,G]=0$. Conversely, if
$\det\calW[F,G]=0$, then for any given $x\in \R$ there exists
$\vecv=\vecv(x)\in \Compl^4$ such that
\[
\bm G(x) & F(x)\\ G'(x) & F'(x) \endm \vecv = 0.
\]
Fix $x=0$, say, and let $\vecv=\binom{b}{a}$ with
$a,b\in\Compl^2$. Then the column vector $y(x):=G(x)b+F(x)a$ is a
solution of $(\Hil-z) y=0$ with $y(0)=0$ and $y'(0)=0$. By the
uniqueness theorem, $y(x)=0$ for all $x\in\R$, and we are done.
\end{proof}

\begin{defi}
\label{def:FG} With $f_j(\cdot,\lambda),g_j(\cdot,\lambda)$ as
above, set for each $\lambda\in\R$,
\begin{align*}
F_1(\cdot,\lambda) &:= (f_1(\cdot,\lambda),f_3(\cdot,\lambda)),
\quad F_2(\cdot,\lambda) = (f_2(\cdot,\lambda),f_4(\cdot,\lambda))
\\
G_1(\cdot,\lambda) &:= (g_2(\cdot,\lambda),g_4(\cdot,\lambda)),
\quad G_2(\cdot,\lambda) =
(g_1(\cdot,\lambda),g_3(\cdot,\lambda)).
\end{align*}
\end{defi}

\begin{remark}
\label{rem:free} At this point it may be helpful to consider the
case $V=0$. Then
\begin{align*}
 F_1(x,\lambda) &=\bm e^{ix\lambda} & 0 \\ 0 & e^{-\mu x} \endm,\;
 F_2(x,\lambda) =\bm e^{-ix\lambda} & 0 \\ 0 & e^{\mu x} \endm, \\
 G_1(x,\lambda) &=\bm e^{ix\lambda} & 0 \\ 0 & e^{-\mu x} \endm,\;
 G_2(x,\lambda) =\bm e^{-ix\lambda} & 0 \\ 0 & e^{\mu x} \endm.
\end{align*}
Hence, in that case $F_1=G_1$ and $F_2=G_2$.
\end{remark}

We record some simple but useful symmetry properties of these
matrix solutions.

\begin{lemma}
\label{lem:FGsymm} For all $\lambda\in\R$,
\begin{align*}
G_1(x,\lambda) &= F_2(-x,\lambda), \quad G_2(x,\lambda) =
F_1(-x,\lambda) \\
\overline{F_1(\cdot,\lambda)} &= F_1(\cdot,-\lambda), \quad
\overline{F_2(\cdot,\lambda)} = F_2(\cdot,-\lambda) \\
\overline{G_1(\cdot,\lambda)} &= G_1(\cdot,-\lambda), \quad
\overline{G_2(\cdot,\lambda)} = G_2(\cdot,-\lambda).
\end{align*}
\end{lemma}
\begin{proof}
This is an immediate consequence of the definitions,
Corollary~\ref{cor:signs}, and Lemma~\ref{lem:wronski}.
\end{proof}

\begin{lemma}
\label{lem:AB} For every $\lambda\in\R\setminus\{0\}$ there exist
unique constant $2\times 2$ matrices $A=A(\lambda), B=B(\lambda)$
with complex entries so that \beeq \label{eq:ABdef}
F_1(\cdot,\lambda) = G_1(\cdot,\lambda) A(\lambda) +
G_2(\cdot,\lambda) B(\lambda). \eneq Then
$A(-\lambda)=\overline{A(\lambda)}$,
$B(-\lambda)=\overline{B(\lambda)}$ and
\begin{align}
G_2(\cdot,\lambda) &= F_2(\cdot,\lambda) A(\lambda) + F_1(\cdot,\lambda) B(\lambda) \label{eq:ABagain}\\
\calW[F_1(\cdot,\lambda),G_2(\cdot,\lambda)] &= A(\lambda)^t(2i\lambda p -2\mu q) \label{eq:wron1}\\
\calW[F_1(\cdot,\lambda),G_1(\cdot,\lambda)] &=
-B(\lambda)^t(2i\lambda p -2\mu q) \label{eq:wron2},
\end{align}
where $p=\bm 1&0 \\ 0&0\endm$, and $q=\bm 0&0\\0&1\endm$.
\end{lemma}
\begin{proof}
For each $\lambda\ne0$, the columns of $G_1(\cdot,\lambda),
G_2(\cdot,\lambda)$ form a basis of the kernel of
$\Hil-(\lambda^2+1)$. Hence, the columns of $F_1(\cdot,\lambda)$
are linear combinations of these solutions, hence the existence of
$A(\lambda),B(\lambda)$. Replacing $x$ with $-x$
in~\eqref{eq:ABdef} implies~\eqref{eq:ABagain}. To
check~\eqref{eq:wron1}, compute
\begin{align*}
\calW[F_1,G_2] &= (A^t {G_1'}^t + B^t {G_2'}^t)G_2 - (A^tG_1^t+B^t
G_2^t) G_2' \\
 &= A^t({G_1'}^tG_2-G_1^tG_2')+ B^t({G_2'}^tG_2-G_2^tG_2') \\
&= -A^t({F_2'}^t F_1-F_2^t F_1')- B^t({F_1'}^t F_1-F_1^t F_1') \\
& = -A^t \bm W[f_2,f_1] & W[f_2,f_3] \\ W[f_4,f_1] & W[f_4,f_3]
\endm - B^t \bm 0 & W[f_1,f_3] \\ W[f_3,f_1] & 0
\endm\\
&= A^t (2i\lambda p -2\mu q),
\end{align*}
where the last line follows from Lemma~\ref{lem:wronski}.
For~\eqref{eq:wron2}, we compute
\begin{align*}
\calW[F_1,G_1] &= (A^t {G_1'}^t + B^t {G_2'}^t)G_1 - (A^tG_1^t+B^t
G_2^t) G_1' \\
 &= A^t({G_1'}^tG_1-G_1^tG_1')+ B^t({G_2'}^tG_1-G_2^tG_1') \\
&= -A^t({F_2'}^t F_2-F_2^t F_2')- B^t({F_1'}^t F_2-F_1^t F_2') \\
& = -A^t \bm 0 & W[f_2,f_4] \\ W[f_4,f_2] & 0
\endm - B^t \bm W[f_1,f_2] & W[f_1,f_4] \\ W[f_3,f_2] & W[f_3,f_4]
\endm\\
&= - B^t (2i\lambda p -2\mu q),
\end{align*}
where the last line again follows from Lemma~\ref{lem:wronski}.
Finally, by Lemma~\ref{lem:FGsymm},
\[
F_1(\cdot,-\lambda)= G_1(\cdot,-\lambda)A(-\lambda) +
G_2(\cdot,-\lambda) B(-\lambda)
\]
is the same as
\[
F_1(\cdot,\lambda)= G_1(\cdot,\lambda)\overline{A(-\lambda)} +
G_2(\cdot,\lambda) \overline{B(-\lambda)}
\]
so that $A(\lambda)=\overline{A(-\lambda)}$,
$B(\lambda)=\overline{B(-\lambda)}$ for all $\lambda\in\R$.
\end{proof}

The following corollary is natural in view of
Remark~\ref{rem:free}. Indeed, the limit $|\lambda|\to\infty$
should correspond to $V\simeq 0$.

\begin{cor}
\label{cor:ABasymp} $A(\lambda)$ and $B(\lambda)$ are smooth for
$\lambda\ne0$. Furthermore, $\lambda pA(\lambda)$, $qA(\lambda)$,
$\lambda p B(\lambda)$ and $q B(\lambda)$ are smooth functions of
$\lambda\in\R$. As $|\lambda|\to\infty$,
\[ A(\lambda)=I+O(\lambda^{-1}),\quad
B(\lambda)=O(\lambda^{-1}).\]
\end{cor}
\begin{proof}
The regularity statements are immediate from \eqref{eq:wron1}
and~\eqref{eq:wron2}. By Remark \ref{rem:x0},
\begin{align*}
&\calW[F_1(\cdot,\lambda),G_2(\cdot,\lambda)] \\
&= [f_1'(0,\lambda)\; f_3'(0,\lambda)]^t [g_1(0,\lambda)\;
g_3(0,\lambda)] - [f_1(0,\lambda)\;
f_3(0,\lambda)]^t [g_1'(0,\lambda)\; g_3'(0,\lambda)] \\
&= [f_1'(0,\lambda)\; f_3'(0,\lambda)]^t [f_1(0,\lambda)\;
f_3(0,\lambda)] + [f_1(0,\lambda)\;
f_3(0,\lambda)]^t [f_1'(0,\lambda)\; f_3'(0,\lambda)] \\
&= \bm i\lambda +O(1) & O(1) \\
O(1) & -\mu + O(1)\endm \Big(I+O(\lambda^{-1}) \Big) +
\Big(I+O(\lambda^{-1}) \Big)
\bm i\lambda +O(1) & O(1) \\
O(1) & -\mu + O(1)\endm \\
&= \bm 2i\lambda & 0 \\ 0 & -2\mu\endm + O(1) = 2i\lambda p -2\mu
q + O(1) .
\end{align*}
It now follows from \eqref{eq:wron1} that $A^t = I +
O(\lambda^{-1})$ as $|\lambda|\to\infty$. Similarly, by
Remark~\ref{rem:x02} and the property that
$c_j(\lambda)=O(\lambda^{-1})$,
\begin{align*}
&\calW[F_1(\cdot,\lambda),G_1(\cdot,\lambda)] \\
&= [f_1'(0,\lambda)\; f_3'(0,\lambda)]^t [g_2(0,\lambda)\;
g_4(0,\lambda)] - [f_1(0,\lambda)\;
f_3(0,\lambda)]^t [g_2'(0,\lambda)\; g_4'(0,\lambda)] \\
&= [f_1'(0,\lambda)\; f_3'(0,\lambda)]^t [f_2(0,\lambda)\;
f_4(0,\lambda)] + [f_1(0,\lambda)\;
f_3(0,\lambda)]^t [f_2'(0,\lambda)\; f_4'(0,\lambda)] \\
&= \bm i\lambda +O(1) & O(1) \\
O(1) & -\mu + O(1)\endm \Big(I+O(\lambda^{-1}) \Big) +
\Big(I+O(\lambda^{-1}) \Big)
\bm -i\lambda +O(1) & O(1) \\
O(1) & \mu + O(1)\endm \\
&=  O(1),
\end{align*}
and the desired bound follows from \eqref{eq:wron2}.
\end{proof}

The following lemma establishes relations between $A$ and $B$
which are analogous to those satisfied by the reflection and
transmission coefficients in scalar scattering theory.

\begin{lemma}
\label{lem:quad} For each $\lambda\ne0$, the matrices $A(\lambda),
B(\lambda)$ satisfy the following  relations:
\begin{align}
-2i\lambda p &= -2i\lambda A(\lambda)^* p A(\lambda) -2\mu
A(\lambda)^* q B(\lambda) + 2\mu B(\lambda)^* qA(\lambda) +
2i\lambda B(\lambda)^*p B(\lambda) \label{eq:quad1} \\
0 &=  A^t(2i\lambda p - 2\mu q) B - B^t (2i\lambda p - 2\mu q)A
\label{eq:quad2} \\
0 &= 2\mu A^*(\lambda)q-2i\lambda B^*(\lambda)p-2\mu
qA(\lambda)-2i\lambda p B(\lambda) \label{eq:linrel}
\end{align}
\end{lemma}
\begin{proof}
Compute $\calW[F_1(\cdot,\lambda),\bar{F}_1(\cdot,\lambda)]$ in
two different ways. For the most part, we will suppress $\lambda$
in our notation for the sake of simplicity. Then, on the one hand,
\[ \calW[\bar{F}_1,{F}_1] = \bm W[\bar{f}_1,{f}_1] & W[\bar{f}_1,{f}_3]
\\ W[\bar{f}_3,{f}_1] & W[\bar{f}_3,{f}_3] \endm = \bm -2i\lambda & 0
\\ 0 & 0 \endm = -2i\lambda p.\]
And on the other hand,
\begin{align}
\calW[\bar{F}_1,{F}_1] &= \calW[\bar{G}_1 \bar{A} +
\bar{G}_2 \bar B, G_1 A + G_2 B] \nonumber \\
&= A^* \calW[\bar{G}_1,G_1]A + A^* \calW[\bar{G}_1,G_2]B +
B^*\calW[\bar{G}_2,G_1]A + B^*\calW[\bar{G}_2,G_2]B. \label{eq:zw}
\end{align}
Next, we compute each of the matrix Wronskians on the right-hand
side of~\eqref{eq:zw}. Before doing so, we calculate
(see~\eqref{eq:f4def})
\[
W[\bar{f}_4,f_4] = W[f_4,f_4]
 =0.
\]
Then
\begin{align*}
\calW[\bar{G}_1,G_1] &= -\calW[\bar{F}_2,F_2] \\
&= -\bm W[\bar{f}_2,f_2] & W[\bar{f}_2,f_4] \\ W[\bar{f}_4,f_2] &
W[\bar{f}_4,f_4] \endm = -\bm W[f_1,f_2] & W[f_1,f_4] \\
\overline{W[{f}_4,f_1]} &
W[\bar{f}_4,f_4] \endm \\
& = - \bm 2i\lambda & 0 \\ 0 & 0 \endm = -2i\lambda p,
\end{align*}
and
\begin{align*}
\calW[\bar{G}_1,G_2] &= -\calW[\bar{F}_2,F_1] \\
&= -\bm W[\bar{f}_2,f_1] & W[\bar{f}_2,f_3] \\ W[\bar{f}_4,f_1] &
W[\bar{f}_4,f_3] \endm = -\bm W[f_1,f_1] & W[f_1,f_3] \\
\overline{W[{f}_4,f_2]} &
\overline{W[f_4,f_3]} \endm \\
& = -\bm 0 & 0 \\ 0 & 2\mu\endm = -2\mu q.
\end{align*}
Thus,
\[ \calW[\bar{G}_2,G_1] = \overline{\calW[G_2,\bar{G}_1]} = -\overline{\calW[\bar{G}_1,G_2]^t}
= 2\mu q. \] Finally,
\begin{align*}
\calW[\bar{G}_2,G_2] &= -\calW[\bar{F}_1,F_1] \\
&= -\bm W[\bar{f}_1,f_1] & W[\bar{f}_1,f_3] \\ W[\bar{f}_3,f_1] &
W[\bar{f}_3,f_3] \endm = -\bm W[f_2,f_1] & W[f_2,f_3] \\
W[{f}_3,f_1] &
{W[f_3,f_3]} \endm \\
& = \bm 2i\lambda & 0 \\ 0 & 0 \endm = 2i\lambda p.
\end{align*}
Inserting this into \eqref{eq:zw} yields
\[
-2i\lambda p= -2i\lambda A^* p A -2\mu A^* q B + 2\mu B^* qA +
2i\lambda B^*p B,
\]
as claimed. For the second quadratic relation, we compute
$\calW[F_1(\cdot,\lambda),{F}_1(\cdot,\lambda)]$ in two different
ways. On the one hand, \beeq \label{eq:F1F1} \calW[{F}_1,{F}_1] =
\bm W[{f}_1,{f}_1] & W[{f}_1,{f}_3]
\\ W[{f}_3,{f}_1] & W[{f}_3,{f}_3] \endm = 0.\eneq
And on the other hand,
\begin{align*}
\calW[{F}_1,{F}_1] &= \calW[{G}_1 {A} +
{G}_2  B, G_1 A + G_2 B] \nonumber \\
&= A^t \calW[{G}_1,G_1]A + A^t \calW[{G}_1,G_2]B +
B^t\calW[{G}_2,G_1]A + B^t\calW[{G}_2,G_2]B \\
& = -A^t \calW[{F}_1,F_1]A - A^t \calW[{F}_1,F_2]B -
B^t\calW[{F}_2,F_1]A - B^t\calW[{F}_2,F_2]B.
\end{align*}
By Lemma~\ref{lem:wronski},
\[ \calW[{F}_2,{F}_2] = \bm W[{f}_2,{f}_2] & W[{f}_2,{f}_4]
\\ W[{f}_4,{f}_2] & W[{f}_4,{f}_4] \endm = 0.\]
By the same lemma,
\[ \calW[{F}_1,{F}_2] = \bm W[{f}_1,{f}_2] & W[{f}_1,{f}_4]
\\ W[{f}_3,{f}_2] & W[{f}_3,{f}_4] \endm = 2i\lambda p - 2\mu q,\]
and therefore,
\[ \calW[{F}_2,F_1] = -\calW[{F}_1,{F}_2]^t = -2i\lambda p + 2\mu
q. \] The conclusion is that
\[ 0 = - A^t(2i\lambda p - 2\mu q) B +
B^t (2i\lambda p - 2\mu q)A,\] which is~\eqref{eq:quad2}.

Finally, to obtain \eqref{eq:linrel}, we compute
$\calW[\bar{F}_1,G_2]$ in two different ways: On the one hand,
\begin{align*}
\calW[\bar{F}_1,G_2] &= \calW[\bar{G}_1
\bar{A}+\bar{G}_2\bar{B},G_2] \\
&= A^* \calW[\bar{G}_1,G_2] + B^* \calW[\bar{G}_2,G_2] \\
&= -A^* \calW[\bar{F}_2,F_1] - B^* \calW[\bar{F}_1,F_1] \\
& = - 2\mu A^*q + 2i\lambda B^* p,
\end{align*}
and on the other hand,
\begin{align*}
\calW[\bar{F}_1,G_2] &= \calW[\bar{F}_1,F_2 A + F_1 B] \\
&= \calW[\bar{F}_1,F_2]A + \calW[\bar{F}_1,F_1] B \\
&= -2\mu q A - 2i\lambda p B,
\end{align*}
as claimed.
\end{proof}

Next we turn to the important question of invertibility of
$A(\lambda)$.

\begin{lemma}
\label{lem:Ainv} If $\lambda\ne0$, then the following are
equivalent:
\begin{itemize}
\item $\det A(\lambda)=0$ \item $E=\lambda^2+1$ is an eigenvalue
of $\Hil$ \item $\det \calW[F_1(\cdot,\lambda),G_2(\cdot,\lambda)]
=0$.
\end{itemize}
\end{lemma}
\begin{proof}
Since $\lambda\ne0$, the first and third properties are equivalent
by~\eqref{eq:wron1}. Now  suppose that $\lambda^2+1$ is an
eigenvalue. Then $F_1(\cdot,\lambda)\binom{0}{1}$ and
$G_2(\cdot,\lambda)\binom{0}{1}$ are linearly dependent (note: the
sign of $\lambda$ is irrelevant here). Since
\[ F_1(\cdot,\lambda)\binom{0}{1} = G_1(\cdot,\lambda)A(\lambda)\binom{0}{1}
+ G_2(\cdot,\lambda)B(\lambda)\binom{0}{1},\] we conclude that
\[ A(\lambda)\binom{0}{1} =0, \quad B(\lambda)\binom{0}{1} =
\binom{0}{\alpha}\] for some $\alpha\ne0$. In particular,
$A(\lambda)$ is singular, as claimed.

Conversely, let $A(\lambda)v=0$ for some $v\in\Compl^2$, $v\ne0$.
Then \eqref{eq:quad1} implies that
\begin{align*}
-2i\lambda \la pv,v\ra &= -2i\lambda \la A(\lambda)^* p
A(\lambda)v,v\ra -2\mu \la A(\lambda)^* q B(\lambda)v,v\ra \\
& \quad + 2\mu \la B(\lambda)^* qA(\lambda)v,v\ra + 2i\lambda \la
B(\lambda)^*p
B(\lambda)v,v\ra \\
&= 2i\lambda \|pBv\|^2,
\end{align*}
and hence $-2i\lambda \|pv\|^2 = 2i\lambda \|pBv\|^2$. Since
$\lambda\ne0$, this implies that $pv=0$ and $pBv=0$. In other
words, both $v$ and $B(\lambda)v$ are parallel to $\binom{0}{1}$.
By the previous paragraph, this implies that $\lambda^2+1$ is an
eigenvalue.
\end{proof}

Unlike the scalar case, where rapid decay of the potential insures
absence of embedded eigenvalues, this is not the case in the
system case. Indeed, take $V_2=0$ and $V_1<0$. If $|V_1|$ is
sufficiently large, then there is $E>1$ and $f\in L^2(\R)$ so that
\[ (\partial_{xx}-1-V_1)f=Ef.\]
This implies that $\Hil\binom{0}{f}=E\binom{0}{f}$ so that $E$
becomes an embedded eigenvalue of~$\Hil$.

The case $E=1$ requires more care.

\begin{defi}
\label{def:res} We say that $E=\pm 1$ is a {\em resonance} of
$\Hil$  provided $\Hil f=\pm f$ has a solution $f\in
L^\infty\setminus L^2$.
\end{defi}

First, we characterize the solutions $f\in L^\infty\setminus L^2$.

\begin{lemma}
\label{lem:0sol} Any solution $f\in L^\infty\setminus L^2$ of
$\Hil f=f$ is of the form
\[ f(x)= C_{\pm}\binom{1}{0} + O(e^{\mp \gamma x}) \text{\ \ as\ \ }x\to\pm\infty\]
where both $C_+\ne0$ and $C_-\ne0$. Similarly, solutions of $\Hil
f=-f$ are of the form
\[ f(x)= C_{\pm}\binom{0}{1} + O(e^{\mp \gamma x}) \text{\ \ as\ \ }x\to\pm\infty\]
\end{lemma}
\begin{proof}
In view of Remark~\ref{rem:x0}, any such $f$ has to be a linear
combination of
\[ f_1(\cdot,0),\; \partial_\lambda f_1(\cdot,0),\; f_3(\cdot,0),\; f_4(\cdot,0). \]
Clearly, only $f_1(\cdot,0), f_3(\cdot,\lambda)$ can occur in this
linear combination when $x\to\infty$. Similarly when $x\to
-\infty$.
\end{proof}

There is a characterization similar to Lemma~\ref{lem:Ainv} for
the endpoint $E=1$ ($E=-1$ is analogous).

\begin{lemma}
\label{lem:D} $E=1$ is a resonance or an eigenvalue of $\Hil$ iff
$\det\calW[F_1(\cdot,0),G_2(\cdot,0)] =0$.
\end{lemma}
\begin{proof}
By definition and the proof of Lemma~\ref{lem:0sol}, $E=1$ is a
resonance or an eigenvalue of $\Hil$ iff there exist $a,b\in
\Compl^2$ not zero such that $F_1(\cdot,0)a$ and $G_2(\cdot,0)b$
are linearly dependent. In view of~\eqref{eq:F1F1} and
Lemma~\ref{lem:matwron} this is in turn equivalent to
$\det\calW[F_1(\cdot,0),G_2(\cdot,0)] =0$, as claimed.
\end{proof}

\noindent Define \beeq \label{eq:Ddef}
 D(\lambda):=
\calW[F_1(\cdot,\lambda),G_2(\cdot,\lambda)] \eneq for all
$\lambda\in\R$.
 By \eqref{eq:wron1},
\begin{align}
(2i\lambda p-2\mu
q)A(\lambda)&=\calW[F_1(\cdot,\lambda),G_2(\cdot,\lambda)]^t  =
-\calW[G_2(\cdot,\lambda),F_1(\cdot,\lambda)] \nonumber\\
&= \calW[F_1(\cdot,\lambda),G_2(\cdot,\lambda)]=D(\lambda), \label{eq:AD}\\
D(\lambda)^t &= -\calW[G_2(\cdot,\lambda),F_1(\cdot,\lambda)]=
\calW[F_1(\cdot,\lambda),G_2(\cdot,\lambda)]=D(\lambda)\label{eq:Dsymm}
\\
D(-\lambda) &= \calW[F_1(\cdot,-\lambda),G_2(\cdot,-\lambda)] =
\overline{\calW[F_1(\cdot,\lambda),G_2(\cdot,\lambda)]}=
\overline{D(\lambda)} \label{eq:Dbar}
\end{align}
for all $\lambda\in\R$. Combining the previous two lemmas
therefore yields

\begin{cor}
\label{cor:imbed} The following properties are equivalent:
\begin{itemize}
\item There are no eigenvalues in $[1,\infty)$ and $E=1$ is not a
resonance \item $D(\lambda)$ is invertible for all $\lambda\in\R$.
\end{itemize}
In that case,
\[ A(\lambda)^{-1} = D(\lambda)^{-1}(2i\lambda p-2\mu q)\]
for all $\lambda\in\R$.
\end{cor}

The point of the final statement is that it should be viewed as a
definition of the left-hand side in case $\lambda=0$. Note
that~\eqref{eq:wron2} allows us to conclude that \beeq
\label{eq:Blim} \lim_{\lambda\to 0} (2i\lambda p B(\lambda) -
2qB(\lambda)) = \calW[F_2(\cdot,0),G_2(\cdot,0)]. \eneq

\section{Scattering solutions, the resolvent, and the distorted Fourier transform.}
\label{sec:scat}

From now on, we shall assume that the conditions of
Corollary~\ref{cor:imbed} hold. We will call such Hamiltonians
{\em admissible}.

\begin{defi}
\label{def:admH} We say that $\Hil$  is {\em admissible} if it
satisfies the requirements of Definition~\ref{def:matrix}, if
there are no eigenvalues in the essential spectrum
$(-\infty,-1]\cup[1,\infty)$, and if the edges $\pm1$ are not
resonances.
\end{defi}

Later we will prove that the linearization of NLS around a ground
state is admissible. It turns out that this class of $\Hil$ admits
the construction of scattering solutions for all energies $|E|>1$,
see Lemma~\ref{lem:scat}. We start with a rather obvious lemma
about the smoothness of $D(\lambda)^{-1}$.

\begin{lemma}
\label{lem:Dsmooth} Let $\Hil$ be admissible. Then both
$D(\lambda)$ and $D^{-1}(\lambda)$ are smooth functions in
$\lambda\in\R$. Moreover, $D(\lambda)^{-1}\lambda = O(1)$ as
$|\lambda|\to\infty$.
\end{lemma}
\begin{proof}
Since both $F_1(x,\lambda)$ and $G_2(x,\lambda)$ are smooth
functions in $\lambda$, it follows from \eqref{eq:Ddef} that
$D(\lambda)$, and therefore also $\det(D(\lambda))$, are smooth.
Since Corollary~\ref{cor:imbed} implies that
$\det(D(\lambda))\ne0$ for all $\lambda\in\R$, we conclude from
Cramer's rule that $D^{-1}(\lambda)$ is smooth for all $\lambda$.
Finally, the asymptotics of $D^{-1}(\lambda)$
 follows from
Corollary~\ref{cor:ABasymp}.
\end{proof}

For the remainder of this section, admissibility of $\Hil$ will be
a standing assumption and we will not mention it further.

\begin{lemma}
\label{lem:scat} Let $\vece=\binom{1}{0}$. Then for all
$\lambda\in\R$
\begin{align}
\calF(x,\lambda) &:= 2i\lambda F_1(x,\lambda) D(\lambda)^{-1}
\vece \label{eq:Fscat}\\
\calG(x,\lambda) &:= 2i\lambda G_2(x,\lambda) D(\lambda)^{-1}
\vece \label{eq:Gscat}
\end{align}
are bounded solutions of $\Hil f=(1+\lambda^2)f$. Moreover, their
asymptotics are given by\footnote{Here
$O(\lambda(1+|\lambda|)^{-2})$  stands for a function which
vanishes linearly as $\lambda\to0$ and decays like $\lambda^{-1}$
as $\lambda\to\infty$.}
\begin{align*}
\calF(x,\lambda) &= s(\lambda)[e^{ix\lambda}\vece +
O((1+|\lambda|)^{-1}e^{-\gamma x})] +
O(\lambda(1+|\lambda|)^{-2}e^{-\mu x})  \text{\ \ \
as\  \ } x\to\infty \\
\calF(x,\lambda) &= [e^{ix\lambda} +
r(\lambda)e^{-ix\lambda}]\vece +
O(\lambda(1+|\lambda|)^{-2}e^{\gamma x}) \text{\ \ \ as\  \ }
x\to-\infty
\\
\calG(x,\lambda) &= s(\lambda)[e^{-ix\lambda}\vece +
O((1+|\lambda|)^{-1}e^{\gamma x})] +
O(\lambda(1+|\lambda|)^{-2}e^{\mu x})  \text{\ \ \
as\  \ } x\to-\infty \\
\calG(x,\lambda) &= [e^{-ix\lambda} +
r(\lambda)e^{ix\lambda}]\vece +
O(\lambda(1+|\lambda|)^{-2}e^{-\gamma x}) \text{\ \ \ as\  \ }
x\to\infty
\end{align*}
where  \[s(\lambda)\vece=2i\lambda p D(\lambda)^{-1}\vece \text{\
\ and\ \ } r(\lambda)\vece = 2i\lambda p B(\lambda)
D(\lambda)^{-1}\vece\] are smooth functions for all $\lambda\in\R$
with $s(0)=0, r(0)=-1$.  The matrix $S(\lambda):=\bm s(\lambda) &
r(\lambda)\\r(\lambda) & s(\lambda) \endm$ is unitary. In fact,
one has
\[ S(\lambda)^*=S(\lambda)^{-1}=S(-\lambda)\]
for all $\lambda\in\R$.
\end{lemma}
\begin{proof}
By Lemmas~\ref{lem:f1} and \ref{lem:f3},
\begin{align*}
\calF(x,\lambda) &= 2i\lambda F_1(x,\lambda)pD(\lambda)^{-1} \vece
+ 2i\lambda F_1(x,\lambda) qD(\lambda)^{-1} \vece \\
&= s(\lambda) f_1(x,\lambda) + O(\lambda(1+|\lambda|)^{-1}e^{-\mu
x}) = s(\lambda)[e^{ix\lambda}\vece +
O((1+|\lambda|)^{-1}e^{-\gamma x})] +
O(\lambda(1+|\lambda|)^{-1}e^{-\mu x})
\end{align*}
as $x\to\infty$. On the other hand, as $x\to-\infty$,
\begin{align*}
\calF(x,\lambda) &= G_1(x,\lambda)\vece + G_2(x,\lambda)
B(\lambda)A(\lambda)^{-1}\vece \\
&= f_2(-x,\lambda) + G_2(x,\lambda) p
B(\lambda)A(\lambda)^{-1}\vece + G_2(x,\lambda) q
B(\lambda)A(\lambda)^{-1}\vece \\
&= f_2(-x,\lambda) + r(\lambda) f_1(-x,\lambda) +
O(\lambda(1+|\lambda|)^{-1} e^{\gamma x}) \\
&= [e^{ix\lambda} + r(\lambda) e^{-ix\lambda}]\vece +
O(\lambda(1+|\lambda|)^{-1} e^{\gamma x}).
\end{align*}
Here we used that $\lambda^{-1}qB(\lambda)A(\lambda)^{-1}\vece=2i
B(\lambda)D^{-1}(\lambda)\vece$ is smooth in $\lambda$. The
asymptotics for $\calG$ now follow since
$\calG(x,\lambda)=\calF(-x,\lambda)$.

As far as the unitarity is concerned, \eqref{eq:wron1} implies
that for $\lambda\ne0$,
 \begin{align*} & -2i\lambda \la
pA^{-1}(\lambda)\vece,A^{-1}(\lambda)\vece\ra
\\ &= -2i\lambda \la A(\lambda)^* p
A(\lambda)A(\lambda)^{-1}\vece,A(\lambda)^{-1}\vece\ra -2\mu \la
A(\lambda)^* q B(\lambda)A(\lambda)^{-1}\vece,A(\lambda)^{-1}\vece\ra \\
& + 2\mu \la B(\lambda)^*
qA(\lambda)A(\lambda)^{-1}\vece,A(\lambda)^{-1}\vece\ra +
2i\lambda \la B(\lambda)^*p
B(\lambda)A(\lambda)^{-1}\vece,A(\lambda)^{-1}\vece\ra.
\end{align*}
Since $s(\lambda)\vece=pA(\lambda)^{-1}\vece$ and $r(\lambda)\vece
= B(\lambda)A(\lambda)^{-1}\vece$ for $\lambda\ne0$, we obtain
from this that
\[ |s(\lambda)|^2+|r(\lambda)|^2=1,\]
which also extends to $\lambda=0$ by continuity. On the other
hand, \eqref{eq:linrel} implies that
\begin{align*}
0 &=  2\mu \la A^*(\lambda)q
A(\lambda)^{-1}\vece,A(\lambda)^{-1}\vece\ra-2i\lambda
\la B^*(\lambda)pA(\lambda)^{-1}\vece,A(\lambda)^{-1}\vece\ra\\
& -2\mu \la
qA(\lambda)A(\lambda)^{-1}\vece,A(\lambda)^{-1}\vece\ra-2i\lambda
\la p B(\lambda)A(\lambda)^{-1}\vece,A(\lambda)^{-1}\vece\ra \\
& = 2\mu \la q A(\lambda)^{-1}\vece,\vece\ra-2i\lambda
\la pA(\lambda)^{-1}\vece,pB(\lambda)A(\lambda)^{-1}\vece\ra\\
& -2\mu \la q\vece,A(\lambda)^{-1}\vece\ra-2i\lambda \la p
B(\lambda)A(\lambda)^{-1}\vece,pA(\lambda)^{-1}\vece\ra \\
&= -2i\lambda s(\lambda)\bar{r}(\lambda) -2i\lambda
r(\lambda)\bar{s}(\lambda)
\end{align*}
so that $0=s(\lambda)\bar{r}(\lambda)+r(\lambda)\bar{s}(\lambda)$,
as claimed. Finally, \begin{align*}
 s(-\lambda)\vece &= -2i\lambda
p D(-\lambda)^{-1}\vece =\overline{ 2i\lambda
D(\lambda)^{-1}\vece} =
\overline{s(\lambda)}\vece, \\
 r(-\lambda)\vece &= -2i\lambda
p B(-\lambda)D(-\lambda)^{-1}\vece =\overline{ 2i\lambda
B(\lambda)D(\lambda)^{-1}\vece} = \overline{r(\lambda)}\vece
\end{align*}
for all $\lambda\in\R$ which proves that
$S(-\lambda)=\overline{S(\lambda)}^t = S(\lambda)^*$.
\end{proof}

The solutions $\calF(\cdot,\lambda),\calG(\cdot,\lambda)$ are
fundamental for several reasons, one being that they form the
(distorted) Fourier basis associated with $\Hil$. This will be
clarified later. First, we show that any globally bounded solution
for $\lambda\ne0$ is a linear combination of these two solutions.

\begin{lemma}
\label{lem:FG} Any  solution $f$ of $\Hil f = (\lambda^2+1)f$ with
$\lambda\ne 0$ and $f\in L^\infty(\R)$ is a linear combination of
$\calF(\cdot,\lambda)$ and $\calG(\cdot,\lambda)$.
\end{lemma}
\begin{proof}
By (the proof of) Lemma~\ref{lem:matwron}, see \eqref{eq:det_sq},
the matrix
\[ \bm F_1(x,\lambda) & G_2(x,\lambda) \\
 F_1'(x,\lambda) & G_2'(x,\lambda) \endm \]
is invertible for all $x\in\R$. As noted in the proof of that
lemma, this means that the four columns of $F_1, G_2$ are linearly
independent. Hence, there exist $\vecv,\vecw\in\Compl^2$ so that
\[ f(x) = F_1(x,\lambda) \vecv + G_2(x,\lambda)\vecw \]
for all $x\in\R$. As $x\to-\infty$,
\[ f(x) = G_1(x,\lambda) A(\lambda)\vecv + G_2(x,\lambda)[
B(\lambda)\vecv + \vecw]
\] remains bounded iff $A(\lambda)\vecv$
is parallel to $\vece$. Similarly, as $x\to\infty$,
\[ f(x) = F_1(x,\lambda) [\vecv + B(\lambda)\vecw] + F_2(x,\lambda)
A(\lambda)\vecw\] remains bounded iff $A(\lambda)\vecw$ is
parallel to $\vece$. Hence,
\[ \vecv = \alpha A^{-1}(\lambda)\vece,\qquad \vecw=\beta
A^{-1}(\lambda) \vecw\] for some constants
$\alpha,\beta\in\Compl$. This implies that
\[ f(x)= \alpha \calF(x,\lambda) + \beta \calG(x,\lambda)\]
for all $x\in\R$, as desired.
\end{proof}

We can now obtain expressions for the resolvent kernel on the
essential spectrum.

\begin{lemma}
\label{lem:resolv} For all $\lambda\ge0$,
\begin{align*}
\big(\Hil-(\lambda^2+1+i0)\big)^{-1}(x,y) &= \left\{
\begin{array}{ll}
-F_1(x,\lambda)D^{-1}(\lambda)
G_2^t(y,\lambda)\sigma_3 & \text{\ \ if\ }x\ge y \\
-G_2(x,\lambda)D^{-1}(\lambda) F_1^t(y,\lambda)\sigma_3 & \text{\
\ if\ }x\le y
\end{array}
\right. \\
\big(\Hil-(\lambda^2+1-i0)\big)^{-1}(x,y) &= \left\{
\begin{array}{ll}
-F_1(x,-\lambda)D^{-1}(-\lambda)
G_2^t(y,-\lambda)\sigma_3 & \text{\ \ if\ }x\ge y \\
-G_2(x,-\lambda)D^{-1}(-\lambda) F_1^t(y,-\lambda)\sigma_3 &
\text{\ \ if\ }x\le y
\end{array}
\right.
\end{align*}
\end{lemma}
\begin{proof}
There exist matrices $\calM_1(y,\lambda)$ and $\calM_2(y,\lambda)$
so that
\[
\big(\Hil-(\lambda^2+1+i0)\big)^{-1}(x,y) = \left\{
\begin{array}{ll}
F_1(x,\lambda)\calM_1(y,\lambda) & \text{\ \ if\ }x\ge y \\
G_2(x,\lambda)\calM_2(y,\lambda) & \text{\ \ if\ }x\le y
\end{array}
\right.
\]
The choice of $F_1(x,\lambda)$ and $G_2(x,\lambda)$ for
$\lambda\ge0$ is due to the fact that these are the only functions
that remain bounded for $\lambda+i\epsilon$ as $x\to\infty$ or
$x\to-\infty$, respectively.  As usual, one needs compatibility
conditions at $x=y$. Here they take the form
\begin{align*}
F_1(x,\lambda)\calM_1(x,\lambda) -
G_2(x,\lambda)\calM_2(x,\lambda) &= 0 \\
F_1'(x,\lambda)\calM_1(x,\lambda) -
G_2'(x,\lambda)\calM_2(x,\lambda) &= -\sigma_3.
\end{align*}
To see why, observe that for any Schwartz function, say, $h(x)=
\binom{h_1(x)}{h_2(x)}$ we need to ensure that
\[ f(x):= G_2(x,\lambda)\int_{x}^\infty  \calM_2(y,\lambda)h(y)\, dy +
F_1(x,\lambda)\int_{-\infty}^x \calM_1(y,\lambda) h(y)\, dy \]
satisfies
\[ (\Hil-(\lambda^2+1))f = (\sigma_3(-\partial_{xx}+1) + V) f(x) =
h(x)\] for all $x\in\R$. Direct differentiation leads to the
conditions above.
 In matrix notation,
\[
\bm F_1(x,\lambda) & G_2(x,\lambda) \\
 F_1'(x,\lambda) & G_2'(x,\lambda) \endm
 \binom{\calM_1(x,\lambda)}{-\calM_2(x,\lambda)} =
 \binom{0}{-\sigma_3}.
\]
By Lemma~\ref{lem:matwron} (or more precisely, its proof), the
$4\times4$ matrix on the left-hand side is invertible for all
$x\in\R$. In fact, in view of \eqref{eq:fehs} and~\eqref{eq:F1F1}
we have the following explicit expression for the inverse:
\begin{align*}
\bm F_1 & G_2 \\ F_1' & G_2'\endm^{-1} &= \bm -D^{-t} & 0 \\ 0 &
D^{-1}\endm \bm 0 & I \\ I & 0 \endm \bm F_1^t & {F_1'}^t \\ G_2^t
& {G_2'}^t \endm \bm 0 & -I \\ I & 0 \endm \\
&= \bm -D^{-t}{G_2'}^t & D^{-t}G_2^t \\ D^{-1}{F_1'}^t &
-D^{-1}F_1^t \endm.
\end{align*}
Consequently,
\[ \binom{\calM_1(x,\lambda)}{\calM_2(x,\lambda)} =
\binom{-D^{-t}(\lambda)G_2^t(x,\lambda)\sigma_3}{-D^{-1}(\lambda)F_1^t(x,\lambda)\sigma_3}
=\binom{-D^{-1}(\lambda)G_2^t(x,\lambda)\sigma_3}{-D^{-1}(\lambda)F_1^t(x,\lambda)\sigma_3},
\]
as claimed. The case of $-i0$ is basically identical, and we are
done.
\end{proof}

\begin{remark}
If $V=0$, then it is easy to check that Lemma~\ref{lem:resolv}
yields
\[ (\Hil-(\lambda^2+1\pm i0))^{-1}(x,y) = \bm \mp \frac{e^{\pm
i\lambda|x-y|}}{2i\lambda} & 0 \\ 0 & -\frac{e^{-\mu|x-y|}}{2\mu}
\endm,\]
which is also an immediate consequence of the standard formulas
for the one-dimensional scalar free resolvent.
\end{remark}

Next, we need to express the jump of the resolvent across the
spectrum $[1,\infty)$.

\begin{lemma}
\label{lem:jump} Let
\[ \calE(x,\lambda):= [\calF(x,\lambda)\;\;\calG(x,\lambda)]\]
for all $\lambda\in\R$. Then \beeq \label{eq:jump}
\big(\Hil-(\lambda^2+1+i0)\big)^{-1}(x,y) -
\big(\Hil-(\lambda^2+1-i0)\big)^{-1}(x,y) = \frac{-1}{2i\lambda}
\calE(x,\lambda)\calE^*(y,\lambda)\sigma_3
 \eneq
for all $\lambda\ge0$.
\end{lemma}
\begin{proof} Let $\lambda\ge0$ and set
\[ \calS(x,y;\lambda):=
\big(\Hil-(\lambda^2+1+i0)\big)^{-1}(x,y) -
\big(\Hil-(\lambda^2+1-i0)\big)^{-1}(x,y).
\]
For fixed $y\in\R$
\[ (\Hil-(\lambda^2+1)) \calS(\cdot,y;\lambda) = 0\]
so that
\[ x\mapsto \calS(x,y;\lambda)\]
is a globally bounded solution. It vanishes identically for
$\lambda=0$. If $\lambda>0$, then Lemma~\ref{lem:FG} implies that
\[ \calS(x,y;\lambda) = \calE(x;\lambda) \calM(y;\lambda) \]
for some matrix $\calM(y;\lambda)$. Lemma~\ref{lem:resolv} implies
that \[ \sigma_3\, \overline{\calS(x,y;\lambda)}^t \sigma_3 =
-\calS(x,y;\lambda). \] Hence,
\[ \sigma_3\, \calM(y;\lambda)^*\calE(x;\lambda)^*\,\sigma_3= - \calE(y;\lambda)
\calM(x;\lambda).\] Since $\calE(x,\lambda)$ is invertible for
$\lambda>0$ and every $x\in\R$, we conclude that
\[ \calE(y;\lambda)^{-1}\,\sigma_3 \calM(y;\lambda)^*= -
\calM(x;\lambda)\sigma_3 \calE(x;\lambda)^{-*} =: C(\lambda)
\]
is a matrix which depends only on $\lambda$. Moreover, we see that
$C(\lambda)^*=-C(\lambda)$. Consequently, \beeq \label{eq:Siden}
 \calS(x,y;\lambda)
= -\calE(x;\lambda) C(\lambda) \calE(y;\lambda)^{*} \sigma_3.
\eneq To determine $C(\lambda)$, we invoke the asymptotics of both
sides of this equation as $x\to\infty$ and $y\to-\infty$. In view
of Lemma~\ref{lem:resolv} the left-hand side satisfies
\begin{align}
\calS(x,y;\lambda) &= -\bm e^{ix\lambda} & 0\\ 0 & 0\endm
D^{-1}(\lambda)\bm e^{-iy\lambda} & 0\\ 0&0\endm \sigma_3 \nonumber\\
& \quad + \bm e^{-ix\lambda} & 0\\ 0 & 0\endm D^{-1}(-\lambda)\bm
e^{iy\lambda} & 0\\ 0&0\endm \sigma_3 + o(1) \nonumber\\
& = - e^{i(x-y)\lambda} \frac{s(\lambda)}{2i\lambda} p -
e^{-i(x-y)\lambda} \frac{s(-\lambda)}{2i\lambda} p + o(1)
\label{eq:links}
\end{align} in this limit. The matrix $C(\lambda)$ can be written
as \[ C(\lambda) = \bm i\alpha & z\\ -\bar{z} & i\beta\endm \]
where $\alpha,\beta\in\R$ and $z\in\Compl$.
 By Lemma~\ref{lem:scat} and~\eqref{eq:Siden},
\begin{align}
&\calS(x,y;\lambda) \nonumber\\
&= -\bm s(\lambda)e^{ix\lambda} & e^{-ix\lambda}+r(\lambda)e^{ix\lambda}\\
0 & 0\endm C(\lambda)\bm e^{-iy\lambda} +
\bar{r}(\lambda)e^{iy\lambda}& 0 \\ \bar{s}(\lambda)e^{iy\lambda}
\endm
\sigma_3 + o(1) \nonumber \nonumber\\
&= -\Big[ e^{i(x+y)\lambda} (i\alpha s(\lambda)\bar{r}(\lambda) +
z|s(\lambda)|^2 + i\beta r(\lambda)\bar{s}(\lambda) -
\bar{z}|r(\lambda)|^2) \nonumber \\
& + e^{i(x-y)\lambda} (i\alpha s(\lambda)-r(\lambda)\bar{z})  +
e^{-i(x-y)\lambda}(-\bar{z}\bar{r}(\lambda)+i\beta\bar{s}(\lambda))
-\bar{z} e^{-i(x+y)\lambda} \Big]p + o(1). \label{eq:rechts}
\end{align}
Comparing \eqref{eq:links} with \eqref{eq:rechts} yields $z=0$ and
\[ i\alpha s(\lambda) = \frac{s(\lambda)}{2i\lambda}, \qquad
 i\beta \bar{s}(\lambda) = \frac{s(-\lambda)}{2i\lambda} = \frac{\bar{s}(\lambda)}{2i\lambda}. \]
 This implies that $\alpha=\beta=-\frac{1}{2\lambda}$ or
 \[ C(\lambda) = \frac{1}{2i\lambda} I,\]
and the lemma follows.
\end{proof}

Let $P_d$ be the Riesz projection onto the discrete spectrum,
i.e., \beeq \label{eq:rp}
 P_d = \frac{1}{2\pi i} \oint_\gamma (zI-\Hil)^{-1}\, dz
\eneq where $\gamma$ is a simple closed curve that encloses the
entire discrete spectrum of~$\Hil$ and lies within the resolvent
set. Moreover, define $P_s=I-P_d$ (``s'' here stands for
``stable'').

We now recall a general representation formula for the expression
$\la e^{it\Hil}P_s\phi,\psi\ra$ from Section~7 of~\cite{Sch}. It
is elementary and was probably known before. Although \cite{Sch}
dealt with dimension three, these particular statements are
independent of the dimension (basically, \eqref{eq:ac} follows
from the Hille-Yoshida theorem). Nevertheless, we present the
proof for the reader's convenience.

\begin{lemma}
\label{lem:rep} Assume that $\Hil$ is admissible. Then there is
the representation \beeq e^{it\Hil} = \frac{1}{2\pi
i}\int_{|\lambda|\ge 1} e^{it\lambda}\,
[(\Hil-(\lambda+i0))^{-1}-(\Hil-(\lambda-i0))^{-1}]\,d\lambda +
\sum_{j} e^{it\Hil} P_{\zeta_j}, \label{eq:ac} \eneq where the sum
runs over the entire discrete spectrum $\{\zeta_j\}_j$  and
 $P_{\zeta_j}$ is the Riesz projection corresponding to the eigenvalue $\zeta_j$.
The formula~\eqref{eq:ac} and the convergence of the integral are
to be understood in the following weak sense: If
$\phi,\psi\in\calS$, then
\[
\la e^{it\Hil}\phi,\psi\ra = \lim_{R\to\infty} \frac{1}{2\pi
i}\int_{R\ge|\lambda|\ge 1}\!\! e^{it\lambda} \big\la
[(\Hil-(\lambda+i0))^{-1}-(\Hil-(\lambda-i0))^{-1}]\phi,\psi\big\ra\,d\lambda
+ \sum_{j} \la e^{it\Hil}P_{\zeta_j}\phi,\psi \ra.
\]
for all $t$. The integrand here is well-defined in view of
Lemma~\ref{lem:jump}.
\end{lemma}
\begin{proof}
We start by checking the following limiting absorption principle
\beeq \label{eq:lim_ap} \sup_{|\lambda|\ge \lambda_0,\,\eps>0}
|\lambda|^{\half}\|(\Hil-(\lambda\pm i\eps))^{-1}\| <\infty, \eneq
where the norm is that of $L^{2,\sigma}(\R)\times L^{2,\sigma}(\R)
\to L^{2,-\sigma}(\R)\times L^{2,-\sigma}(\R)$ where
$\sigma>\half$, say, and $\lambda_0$ is large. In the free case
(i.e., $V=0$ and $\Hil=\Hil_0$) this bound is an immediate
consequence of the explicit form of the resolvent. If $V\ne0$,
then we write
\[ (\Hil-(\lambda\pm i0))^{-1}= (\Hil_0-(\lambda\pm i0))^{-1}\big(I+V(\Hil_0-(\lambda\pm i0))^{-1}\big)^{-1}.\]
This is to be understood as identity between operators
$L^{2,\sigma}(\R)\times L^{2,\sigma}(\R) \to
L^{2,-\sigma}(\R)\times L^{2,-\sigma}(\R)$ where $\sigma>\half$.
Note the the  inverse
\[ \big(I+V(\Hil_0-(\lambda\pm i0))^{-1}\big)^{-1}\]
exists as operator on $L^{2,\sigma}(\R)\times L^{2,\sigma}(\R) \to
L^{2,\sigma}(\R)\times L^{2,\sigma}(\R)$ because
\[ \sup_{|\lambda|\ge\lambda_0} \|V(\Hil_0-(\lambda\pm i0))^{-1}\| < \half \]
for $\lambda_0$ large. Hence, \eqref{eq:lim_ap} holds.

The evolution $e^{it\Hil}$ is defined via the Hille-Yoshida
theorem. Indeed, let $a>0$ be large. Then $i\Hil-a$ satisfies
(with $\rho$ the resolvent set)
\[
 \rho(i\Hil-a)\supset (0,\infty)\text{\ \ and\ \ }\|(i\Hil-a-\lambda)^{-1}\|_{2\to2}\le \lambda^{-1}\text{\ \ for all\ \ }\lambda>0.
\]
The estimate here follows from
\begin{align*}
\|(i\Hil-a-\lambda)^{-1}\| &\le \|(i\Hil_0-a-\lambda)^{-1}\| \big\|\big(I+iV(i\Hil_0-a-\lambda)^{-1}\big)^{-1}\big\| \\
&\le (\lambda+a)^{-1}\frac{1}{1-C(a+\lambda)^{-1}} =
\frac{1}{\lambda+a-C}\le \lambda^{-1}
\end{align*}
provided $a$ is large. Hence $\{e^{t(i\Hil-a)}\}_{t\ge0}$ is a
contractive semigroup, so that $\|e^{it\Hil}\|_{2\to2}\le
e^{|t|a}$ for all $t\in\R$. If $\Re z>a$, then there is the
Laplace transform \beeq \label{eq:Lap}
 (i\Hil-z)^{-1}=-\int_0^\infty e^{-tz}\,e^{it\Hil}\,dt
\eneq as well as its inverse (with $b>a$ and $t>0$) \beeq
\label{eq:invLap} e^{it\Hil} = -\frac{1}{2\pi i}
\int_{b-i\infty}^{b+i\infty} e^{tz}\, (i\Hil-z)^{-1}\, dz. \eneq
While \eqref{eq:Lap} converges in the norm sense,
defining~\eqref{eq:invLap} requires more care. The claim is that
for any $\phi,\psi\in \Dom(\Hil)=W^{2,2}\times W^{2,2}$, \beeq
\label{eq:limR}
 \la e^{it\Hil} \phi,\psi \ra =
- \lim_{R\to\infty} \frac{1}{2\pi i} \int_{b-iR}^{b+iR} e^{tz}\,
\la (i\Hil-z)^{-1}\phi,\psi\ra\, dz. \eneq To verify this, let
$t>0$ and use \eqref{eq:Lap} to conclude that \bear -\frac{1}{2\pi
i} \int_{b-iR}^{b+iR} e^{tz}\, \la (i\Hil-z)^{-1}\phi,\psi\ra\, dz
&=&
 \frac{1}{2\pi i} \int_{b-iR}^{b+iR} e^{tz}\,  \int_0^\infty e^{-sz}\,\la e^{is\Hil}\,\phi,\psi\ra\, dsdz \nn \\
&=& \frac{1}{\pi} \int_0^\infty
e^{(t-s)b}\,\frac{\sin((t-s)R)}{t-s}\,\la
e^{is\Hil}\,\phi,\psi\ra\, ds. \label{eq:dirker} \eear Since
$e^{(t-s)b}\,\la e^{is\Hil}\,\phi,\psi\ra$ is a $C^1$ function
in~$s$ (recall $\phi\in \Dom(\Hil)$)
 as well as  exponentially decaying in $s$ (because of $b>a$),
it follows from standard properties of the Dirichlet kernel that
the limit in~\eqref{eq:dirker} exists and equals $\la
e^{it\Hil}\phi,\psi\ra$, as claimed. Note that if $t<0$, then the
limit is zero. Therefore, it follows that for any $b>a$,
\begin{align*}
 \la e^{it\Hil} \phi,\psi \ra &=
- \lim_{R\to\infty} \Big\{ \frac{1}{2\pi i} \int_{b-iR}^{b+iR}
e^{tz}\, \la (i\Hil-z)^{-1}\phi,\psi\ra\, dz
- \frac{1}{2\pi i} \int_{-b-iR}^{-b+iR} e^{tz}\, \la (i\Hil-z)^{-1}\phi,\psi\ra\, dz \Big\} \\
&= \lim_{R\to\infty} \frac{1}{2\pi i} \int_{-R}^{R}
e^{it\lambda}\, \la [e^{-bt}(\Hil-(\lambda+ib))^{-1} -
e^{bt}(\Hil-(\lambda-ib))^{-1}]\phi,\psi\ra\, d\lambda.
\end{align*}
Next, assume that $\phi,\psi$ are as in the statement of the
theorem, and shift the contour in the previous integrals by
sending $b\to0+$. More precisely, we apply the residue theorem to
the contour integrals over the rectangles with vertices $\pm
R+ib$, $\pm R+i0$ and the reflected one below the real axis.  The
horizontal segments on the real axis need to avoid the poles,
which can be achieved by
 surrounding each of the at most finitely many real poles of the resolvent $(\Hil-z)^{-1}$
by a small semi-circle. Combining each such semi-circle with its
reflection yields a small closed loop and the resulting integral
is precisely the Riesz projection corresponding to that real
eigenvalue. The Riesz projections corresponding to eigenvalues on
the imaginary axis are obtained as residues. On the other hand, we
also need to show that the contribution by the horizontal segments
is zero in the limit $R\to\infty$. This, however, follows from the
limiting absorption principle~\eqref{eq:lim_ap}. The lemma
follows.
\end{proof}

Under our assumptions, $\Hil$ can only have finitely many points
in its discrete spectrum, each of which is an eigenvalue of finite
algebraic multiplicity (however, the geometric and algebraic
multiplicities my differ for each one of them). Let $P_d$ denote
the Riesz projection onto the discrete spectrum. It is given by
the Cauchy integral of the resolvent around a simple closed curve
which surrounds the entire discrete spectrum but avoids the
essential spectrum. We write $P_s=I-P_d$ for the ``stable''
projection.

We can now state the Fourier expansion theorem. So far, our
analysis has been restricted to the right half of the essential
spectrum, i.e.,  $[1,\infty)$. To extend this to the left half, it
will be convenient (but not essential) to use a further property
of $V$, see Definition~\ref{def:matrix} and~\eqref{eq:pauli}.
Namely,
\[
\sigma_1 V\sigma_1 = -V \text{\ \ and\ \ }\sigma_1 \Hil \sigma_1 =
-\Hil \text{\ \ where\ \ }\sigma_1=\bm 0&1\\1&0
\endm.
\]
Therefore, if we denote the scattering solutions of
Lemma~\ref{lem:scat} by $\calF_+(x,\lambda),\calG_+(x,\lambda)$,
then the corresponding ones for the negative essential spectrum
are
\[ \calF_{-}(x,\lambda) := \sigma_1\calF_{+}(x,\lambda),\quad
\calG_{-}(x,\lambda) := \sigma_1\calG_{+}(x,\lambda).\]

\begin{prop}
\label{prop:fourier} Let
\[ e_{\pm}(x,\lambda) = \left\{\begin{array}{ll} \calF_{\pm}(x,\lambda) &
\text{\ \ if\ \ }\lambda\ge 0 \\
\calG_{\pm}(x,-\lambda) & \text{\ \ if\ \ }\lambda\le 0
\end{array} \right.
\]
Then for every $\phi,\psi\in \calS$,
\[ \la P_s \phi,\psi\ra = \frac{1}{2\pi} \int_{-\infty}^\infty \la
\phi, \sigma_3 e_{+}(\cdot,\lambda)\ra \overline{\la \psi,
e_{+}(\cdot,\lambda) \ra} \, d\lambda + \frac{1}{2\pi}
\int_{-\infty}^\infty \la \phi, \sigma_3 e_{-}(\cdot,\lambda)\ra
\overline{\la \psi, e_{-}(\cdot,\lambda) \ra} \, d\lambda.
\]
The integrals on the right-hand side are absolutely convergent. In
fact, the integrand is rapidly decaying.
\end{prop}
\begin{proof}
We start from the representation formula
\begin{align*}
 \la P_s
\phi,\psi\ra &= \frac{1}{2\pi i}
\Big\{\int_{-\infty}^{-1}+\int_1^\infty\Big\} \big \la
\big((\Hil-(u+i0))^{-1}- (\Hil-(u-i0))^{-1}\big)\phi,\psi
\big\ra\, du \\
&= \frac{1}{2\pi i} \int_0^\infty 2\lambda \big \la
\big((\Hil-(\lambda^2+1+i0))^{-1}-
(\Hil-(\lambda^2+1-i0))^{-1}\big)\phi,\psi \big\ra\, d\lambda \\
& + \frac{1}{2\pi i} \int_0^\infty 2\lambda \big \la
\big((\Hil-(-\lambda^2-1+i0))^{-1}-
(\Hil-(-\lambda^2-1-i0))^{-1}\big)\phi,\psi \big\ra\, d\lambda
\end{align*}
which holds in the principal value sense. This was proved
in~\cite{Sch}, see Lemma~44. Since for $\lambda>0$
\begin{align*}
&(\Hil-(\lambda^2+1+i0))^{-1}(x,y)-
(\Hil-(\lambda^2+1-i0))^{-1}(x,y) \\
&= \frac{-1}{2i\lambda}
[e_+(x,\lambda)\;\;e_+(x,-\lambda)][e_+(y,\lambda)\;\;e_+(y,-\lambda)]^*\sigma_3
\\
&(\Hil-(-\lambda^2-1+i0))^{-1}(x,y)-
(\Hil-(-\lambda^2-1-i0))^{-1}(x,y) \\
&= \frac{-1}{2i\lambda}
[e_-(x,\lambda)\;\;e_-(x,-\lambda)][e_-(y,\lambda)\;\;e_-(y,-\lambda)]^*\sigma_3,
\end{align*}
this representation formula can be rewritten as
\begin{align*}
 \la P_s
\phi,\psi\ra &= \frac{1}{2\pi} \int_0^\infty \Big\la
\binom{\overline{e_+(y,\lambda)}^t}{\overline{e_+(y,-\lambda)}^t}\sigma_3\phi(y),
\binom{\overline{e_+(x,\lambda)}^t}{\overline{e_+(x,-\lambda)}^t}\psi(x)
\Big\ra\, d\lambda \\
& \quad +\frac{1}{2\pi} \int_0^\infty \Big\la
\binom{\overline{e_-(y,\lambda)}^t}{\overline{e_-(y,-\lambda)}^t}\sigma_3\phi(y),
\binom{\overline{e_-(x,\lambda)}^t}{\overline{e_-(x,-\lambda)}^t}\psi(x)
\Big\ra\, d\lambda \\
&= \frac{1}{2\pi} \int_{-\infty}^\infty \la \phi, \sigma_3
e_{+}(\cdot,\lambda)\ra \overline{\la \psi, e_{+}(\cdot,\lambda)
\ra} \, d\lambda + \frac{1}{2\pi} \int_{-\infty}^\infty \la \phi,
\sigma_3 e_{-}(\cdot,\lambda)\ra \overline{\la \psi,
e_{-}(\cdot,\lambda) \ra} \, d\lambda,
\end{align*}
as claimed. So far, we need to interpret the right-hand side in
the principal value sense. Since $e_{\pm}$ are solutions, i.e.,
they satisfy $\Hil
e_{\pm}(\cdot,\lambda)=(1+\lambda^2)e_{\pm}(\cdot,\lambda)$, we
obtain
\[ \la \phi,\sigma_3 e_+\ra = \la \Hil^m \phi,\sigma_3 e_+\ra (1+\lambda^2)^{m}\]
and the rapid decay follows.
\end{proof}

In fact, the same proof also yields the following representation
of the time evolution.

\begin{cor}
\label{cor:Evol} With the same notation as in the previous
theorem,
\begin{align*}
 \la e^{it\Hil}P_s \phi,\psi\ra &= \frac{e^{it}}{2\pi}
\int_{-\infty}^\infty e^{it\lambda^2}\la \phi, \sigma_3
e_{+}(\cdot,\lambda)\ra \overline{\la \psi, e_{+}(\cdot,\lambda)
\ra} \, d\lambda \\
& + \frac{e^{-it}}{2\pi}  \int_{-\infty}^\infty
e^{-it\lambda^2}\la \phi, \sigma_3 e_{-}(\cdot,\lambda)\ra
\overline{\la \psi, e_{-}(\cdot,\lambda) \ra} \, d\lambda,
\end{align*}
with absolutely convergent integrals.
\end{cor}

Formally, this can be written as
\[
e^{it\Hil}P_s \phi = \frac{e^{it}}{2\pi} \int_{-\infty}^\infty
e^{it\lambda^2}\la \phi, \sigma_3
e_{+}(\cdot,\lambda)\ra \, e_{+}(\cdot,\lambda) \, d\lambda \\
 + \frac{e^{-it}}{2\pi}  \int_{-\infty}^\infty  e^{-it\lambda^2}\la
\phi, \sigma_3 e_{-}(\cdot,\lambda)\ra \, e_{-}(\cdot,\lambda) \,
d\lambda.
\]

One easy consequence of Proposition~\ref{prop:fourier} is the
following stability bound on the evolution. Not too surprisingly,
it can also be established independently of the scattering theory
from above. In fact, it is a relatively straightforward
consequence of Kato's smoothing theory which does not depend on
the dimension. See~\cite{Sch} for the three-dimensional case. The
argument which is presented there only uses the representation
from Lemma~\ref{lem:rep} and carries over to the one-dimensional
case as well.

\begin{lemma}
\label{lem:L2stable} Let $\Hil$ be admissible. Then the following
{\em stability} bound holds: \beeq \sup_{t\in\R} \|e^{it\Hil}
P_s\|_{2\to2} \le C. \label{eq:L2stable} \eneq
\end{lemma}
\begin{proof}
In view of Corollary~\ref{cor:Evol},
\begin{align*}
 |\la e^{it\Hil}P_s \phi,\psi\ra| &\le
\int_{-\infty}^\infty |\la \phi, \sigma_3
e_{+}(\cdot,\lambda)\ra|\, |\la \psi, e_{+}(\cdot,\lambda) \ra| \,
d\lambda
 + \int_{-\infty}^\infty
|\la \phi, \sigma_3 e_{-}(\cdot,\lambda)\ra|\,
|\la \psi, e_{-}(\cdot,\lambda) \ra| \, d\lambda, \\
&\le 2\max_{\pm}\Big(\int_{-\infty}^\infty |\la \phi, \sigma_3
e_{\pm}(\cdot,\lambda)\ra|^2\,d\lambda\Big)^{\half}
\Big(\int_{-\infty}^\infty|\la \psi, e_{\pm}(\cdot,\lambda) \ra|
\, d\lambda\Big)^{\half} \le C\|\phi\|_2\|\psi\|_2,
\end{align*}
where the final bound follows from the asymptotics of
Lemma~\ref{lem:scat}.
\end{proof}

Next, we state the natural bound on $xe^{it\Hil}P_s$:

\begin{lemma}
\label{lem:moments} Let $\Hil$ be admissible. Then
\[ \| x e^{it\Hil}P_s f\|_2 \le C\la t\ra\|f\|_{H^1}+C\| x f\|_2\]
for all $t\ge0$.
\end{lemma}
\begin{proof}
Clearly,
\begin{align*}
 i\frac{d}{dt} \la |x|^2 e^{it\Hil}P_s f, e^{it\Hil}P_s f\ra &= \la (\Hil^* |x|^2-|x|^2\Hil) e^{it\Hil}P_s f, e^{it\Hil}P_s f \ra \\
&= -\la [\partial_x^2, |x|^2]\sigma_3  e^{it\Hil}P_s f, e^{it\Hil}P_s f \ra + \la (V^*-V)|x|^2 e^{it\Hil}P_s f, e^{it\Hil}P_s f \ra \\
&= -\la (2+2x\partial_x)\sigma_3  e^{it\Hil}P_s f, e^{it\Hil}P_s f
\ra + \la (V^*-V)|x|^2 e^{it\Hil}P_s f, e^{it\Hil}P_s f \ra.
\end{align*}
In particular,
\begin{align*}
 \Big|\;i\frac{d}{dt} \la |x|^2 e^{it\Hil}P_s f, e^{it\Hil}P_s f\ra \;\vert_{t=0}\Big|
&\le \Big|\la (2+2x\partial_x)\sigma_3 P_s f, P_s f \ra\Big| +\Big| \la (V^*-V)|x|^2 P_s f, P_s f \ra \Big| \\
&\le C (\|\la x\ra f\|_2 \|f\|_{H^1} + \|f\|_2^2).
\end{align*}
Thus,
\begin{align*}
 -\frac{d^2}{dt^2} \la |x|^2 e^{it\Hil}P_s f, e^{it\Hil}P_s f\ra
&= -2\la (\Hil^*x\partial_x-x\partial_x\Hil^*)\sigma_3  e^{it\Hil}P_s f, e^{it\Hil}P_s f \ra \\
& \quad + \la (\Hil^*(V^*-V)|x|^2-(V^*-V)|x|^2\Hil) e^{it\Hil}P_s f, e^{it\Hil}P_s f \ra  \\
&= 2\la (\partial_x^2(x\partial_x)-x\partial_x\partial_x^2)
e^{it\Hil}P_s f, e^{it\Hil}P_s f \ra
+2\la x\partial_x(V^*)\sigma_3  e^{it\Hil}P_s f, e^{it\Hil}P_s f \ra \\
&\quad - \la (\partial_x^2 \sigma_3(V^*-V)|x|^2-(V^*-V)|x|^2\sigma_3\partial_x^2) e^{it\Hil}P_s f, e^{it\Hil}P_s f \ra \\
&+ \la (V^*(V^*-V)|x|^2-(V^*-V)|x|^2V) e^{it\Hil}P_s f,
e^{it\Hil}P_s f \ra,
\end{align*}
which implies that
\[
\sup_{t}\Big|\frac{d^2}{dt^2} \la |x|^2 e^{it\Hil}P_s f,
e^{it\Hil}P_s f\ra \Big| \le \sup_{t}C\|e^{it\Hil}P_s
f\|_{H^1}^2\le C\|f\|_{H^1}^2.
\]
A Taylor expansion of degree two therefore yields
\begin{align*}
\| xe^{it\Hil}P_s f\|_2^2 &\le C\,\|\la x\ra f\|_2^2 + C\,t(\|\la x\ra f\|_2 \|f\|_{H^1} + \|f\|_2^2) + C\,t^2 \|f\|_{H^1}^2 \\
&\le C\,(\la t\ra^2 \|f\|_{H^1}^2 + \|\la x\ra f\|_2^2),
\end{align*}
and the lemma follows.
\end{proof}

\section{Dispersive estimates: The unweighted case}
\label{sec:noweigh}

In this section, we prove dispersive estimates on $e^{it\Hil}P_s$
in the $L^1(\R)\to L^\infty(\R)$ sense. For the scalar case, see
for example~\cite{GolSch}.

\begin{prop}
\label{prop:disp1} Let $\Hil$ be admissible, see
Definition~\ref{def:admH}. Then for all $t\ne0$,
\[ \|e^{it\Hil}P_s\|_{1\to\infty}\le C|t|^{-\half}\]
with some $C=C(V)$.
\end{prop}
\begin{proof}
We will follow the proof strategy of the one-dimensional case
of~\cite{GolSch}. To do so, we start from the representation
formula
\begin{align*}
\la e^{it\Hil} P_s \phi,\psi \ra &= \frac{e^{it}}{\pi i}
\int_0^\infty e^{it\lambda^2}\lambda \big \la
\big((\Hil-(\lambda^2+1+i0))^{-1}-
(\Hil-(\lambda^2+1-i0))^{-1}\big)\phi,\psi \big\ra\, d\lambda \\
& + \frac{e^{-it}}{\pi i} \int_0^\infty e^{-it\lambda^2} \lambda
\big \la \big((\Hil-(-\lambda^2-1+i0))^{-1}-
(\Hil-(-\lambda^2-1-i0))^{-1}\big)\phi,\psi \big\ra\, d\lambda,
\end{align*}
which holds in the principal value sense if $\phi,\psi\in\calS$.
It will suffice to deal with energies $E\ge1$ since the second
integral will satisfy the same bounds as the first by symmetry.
Let $\chi$ be a smooth, even, and compactly supported bump
function so that $\chi(\lambda)=1$ for $|\lambda|\le \lambda_1$,
where $\lambda_1=\lambda_1(V)$ will be specified later. On the
support of $1-\chi(\lambda)$ (the ``high energy case'') we will
use a Born series expansion. More precisely, since
\begin{align*}
(\Hil_0-(\lambda^2+1\pm i0))^{-1}(x,y) &= \bm \big(-\partial_{xx}-(\lambda^2\pm i0)\big)^{-1}(x,y)  & 0\\
0 &  (-\partial_{xx}-\mu^2)^{-1}(x,y)  \endm \\
& = \bm \mp \frac{e^{\pm i\lambda|x-y|}}{2i\lambda} & 0 \\ 0 &
-\frac{e^{-\mu|x-y|}}{2\mu}
\endm,
\end{align*}
satisfies
\[ \|(\Hil_0-(\lambda^2+1\pm i0))^{-1}\|_{1\to\infty} \le C|\lambda|^{-1},\]
we conclude that \beeq \label{eq:born} (\Hil-(\lambda^2+1\pm
i0))^{-1} = \sum_{n=0}^\infty (-1)^n
(\Hil_0-(\lambda^2+1+i0))^{-1}\Big(V(\Hil_0-(\lambda^2+1+i0))^{-1}\Big)^n
\eneq converges in the operator norm of $L^1(\R)\to L^\infty(\R)$
provided $|\lambda|\ge\lambda_1$ is sufficiently large compared to
$\|V\|_{L^1}$. Indeed, the operator norms of the $n$-th summand on
the right-hand side is bounded by
$C\lambda^{-1}(\|V\|_1\lambda^{-1})^n$. Hence\footnote{Symbols
like $1-\chi(H-I)$ and $P_s^+$ are being used in a purely formal
way --- they are defined by the $\lambda$-integrals in which they
arise.}
\begin{align}
\nn
&\big|\la e^{it\Hil} (1-\chi(H-I)) P_s^+ \phi,\psi \ra\big| \\
&\le C\sum_{n=0}^\infty \Big|\int_{0}^\infty
e^{it\lambda^2}\,\lambda(1-\chi(\lambda^2))
\big[ \big\la (\Hil_0-(\lambda^2+1+i0))^{-1}\big(V(\Hil_0-(\lambda^2+1+i0))^{-1}\big)^n \phi,\psi \big\ra \nn \\
& - \big\la
(\Hil_0-(\lambda^2+1-i0))^{-1}\big(V(\Hil_0-(\lambda^2+1-i0))^{-1}\big)^n
\phi,\psi \big\ra\big] \, d\lambda \Big|.\label{eq:born_decay}
\end{align}
The term $n=0$ represents the usual free Schr\"odinger decay and
its contribution is bounded by
$C|t|^{-\half}\|\phi\|_1\|\psi\|_1$. Indeed, the oscillatory
integral bound that arises in this case is
\begin{align*}
& \sup_{x,y\in\R} \Big|\int_{0}^\infty e^{it\lambda^2}\,(1-\chi(\lambda^2))\cos(\lambda|x-y|)\, d\lambda \Big| \\
& \le \sup_{x,y\in\R} \Big|\int_{-\infty}^\infty
[e^{it\lambda^2}]^{\vee}(u)[(1-\chi(\lambda^2))\cos(\lambda|x-y|)]^{\vee}(u)\,
du \Big| \le C |t|^{-\half}.
\end{align*}
Here the Fourier transforms are with respect to $\lambda$ and we
used that
\[ [(1-\chi(\lambda^2))\cos(\lambda a)]^{\vee}(u) \]
is a measure with  total variation norm uniformly bounded in $a$.
Next, consider the contribution by $n=1$ in~\eqref{eq:born_decay}.
 Writing $\phi=\binom{\phi_1}{\phi_2}, \psi=\binom{\psi_1}{\psi_2}$ this term becomes
(we ignore multiplicative constants and we write $dx=dx_0dx_1dx_2$
for simplicity)
\begin{align}
& \int_{\R^2}\int_{-\infty}^\infty e^{it\lambda^2}\lambda^{-1}
(1-\chi(\lambda^2)) \sin(\lambda(|x_0-x_1|+|x_1-x_2|))\,d\lambda
\, {V_1(x_1)\phi_1(x_0)
\bar{\psi}_1(x_2)}\, dx\label{eq:PU1} \\
& +\int_{\R^2}\int_{-\infty}^\infty e^{it\lambda^2}\lambda^{-1}
(1-\chi(\lambda^2)) \sin(\lambda |x_0-x_1|)
e^{-\mu |x_2-x_1|}   \,d\lambda \,  {V_2(x_1)\phi_1(x_0) \bar{\psi}_2(x_2)}\,dx \label{eq:PW1} \\
& -\int_{\R^2}\int_{-\infty}^\infty e^{it\lambda^2}\lambda^{-1}
(1-\chi(\lambda^2)) \sin(\lambda |x_2-x_1|) e^{-\mu |x_1-x_0|}
\,d\lambda \, {V_2(x_1)\phi_2(x_0) \bar{\psi}_1(x_2)}\,dx.
\label{eq:PW2}
\end{align}
To bound the oscillatory integrals, note first that
\begin{align*}
& \big\| (1+u^2)[\lambda^{-1}(1-\chi(\lambda^2))]^{\vee}(u) \big\|_\infty \\
& \le  \big\| [\lambda^{-1}(1-\chi(\lambda^2))]^{\vee}(u)
\big\|_\infty + C \big\|
\partial_\lambda^2[\lambda^{-1}(1-\chi(\lambda^2))] \big\|_1 <
\infty,
\end{align*}
since $[\lambda^{-1}]^{\vee}(u)=c\,\sign(u)$. Hence
\[ \big\| [\lambda^{-1}(1-\chi(\lambda^2))]^{\vee}(u) \big\|_1 < \infty. \]
Second, we claim that  (recall $\mu=\sqrt{2+\lambda^2}$)
\begin{align}
&\sup_{a\ge0} \big\| \int_{-\infty}^\infty
e^{-a\sqrt{2+\lambda^2}} e^{-i\lambda u}\, d\lambda \big\|_{\calM}
= \sup_{b\ge0} \big\| \int_{-\infty}^\infty
e^{-\sqrt{b+\lambda^2}} e^{-i\lambda u}\, d\lambda \big\|_{\calM}
<\infty \label{eq:L1mu}
\end{align}
where the norms refer to the total variation norms of measures. To
see this, compute
\begin{align*}
\|\partial^2_\lambda e^{-\sqrt{b+\lambda^2}}\|_1 &= \int_{-\infty}^\infty \Big|\Big(-\frac{b}{(b+\lambda^2)^{_\frac32}} +\frac{\lambda^2}{b+\lambda^2}\Big)e^{-\sqrt{b+\lambda^2}} \Big| \, d\lambda\\
&\le \int_{-\infty}^\infty  \frac{b}{(b+\lambda^2)^{_\frac32}} \,
d\lambda +
\int_{-\infty}^\infty e^{-\sqrt{b+\lambda^2}}  \, d\lambda \\
&\le \int_{-\infty}^\infty  \frac{1}{(1+\lambda^2)^{_\frac32}} \,
d\lambda + \int_{-\infty}^\infty e^{-|\lambda|}  \, d\lambda \le C
\end{align*}
uniformly  in~$b>0$. It follows that
\[ \sup_{b\ge0}\; (1+u^2)
\Big| \int_{-\infty}^\infty e^{-\sqrt{b+\lambda^2}} e^{-i\lambda
u}\, d\lambda \Big| \le C
\]
and \eqref{eq:L1mu} holds. By the same type of argument as in the
$n=0$ case, we conclude that the contribution from
\eqref{eq:PU1}-\eqref{eq:PW2} is
\[ \le C|t|^{-\half}\|V\|_1 \|\phi\|_1\|\psi\|_1,\]
as desired. Finally, the terms $n\ge2$ in \eqref{eq:born_decay}
are similar. More precisely, they lead to oscillatory integrals of
the form
\[
\sup_{b\ge0,\,a\in\R}\Big| \int_{-\infty}^\infty e^{it\lambda^2}
\lambda^{-n} (1-\chi(\lambda^2)) e^{ia\lambda}
e^{-b\sqrt{2+\lambda^2}}\, d\lambda\Big| \le C\,\lambda_1^{-n}
|t|^{-\half},
\]
which follow by the same type of arguments as before. Hence, the
entire series in~\eqref{eq:born_decay} is estimated by
\[
\sum_{n=0}^\infty C\,\lambda_1^{-n}|t|^{-\half} \|V\|_1^n
\|\phi\|_1\|\psi\|_1 \le C|t|^{-\half} \|\phi\|_1\|\psi\|_1,
\]
as claimed.

It remains to deal with those $\lambda$ that belong to the support
of $\chi$. By Lemma~\ref{lem:resolv},
\begin{align}
&\la e^{it\Hil}\chi(\Hil-I)P_s^+ \phi,\psi \ra \nn \\
&= \frac{e^{it}}{\pi i}\int_0^\infty e^{it\lambda^2}
\lambda\chi(\lambda^2)
\big\la [(\Hil-(\lambda^2+1+i0))^{-1}-(\Hil-(\lambda^2+1-i0))^{-1}]\phi, \psi\big\ra \, d\lambda \nn \\
&= -\frac{e^{it}}{\pi i}\int_{\R^2}\int_{-\infty}^\infty
e^{it\lambda^2} \lambda\chi(\lambda^2)
\big\la F_1(x,\lambda)D^{-1}(\lambda)G_2^t(y,\lambda)\sigma_3 \phi(x), \psi(y)\big\ra \, d\lambda\chi_{[x\ge y]} \,dxdy \label{eq:pot1}\\
&\quad +\frac{e^{it}}{\pi i}\int_{\R^2}\int_{-\infty}^\infty
e^{it\lambda^2} \lambda\chi(\lambda^2) \big\la
G_2(x,\lambda)D^{-1}(\lambda)F_1^t(y,\lambda)\sigma_3 \phi(x),
\psi(y)\big\ra \, d\lambda\chi_{[x< y]} \,dxdy. \label{eq:pot2}
\end{align}
The second integral~\eqref{eq:pot2} can be transformed into a
variant of the first~\eqref{eq:pot1} by means of the change of
variables $x\to-x$, $y\to-y$:
\begin{align*}
& \int_{\R^2}\int_{-\infty}^\infty e^{it\lambda^2}
\lambda\chi(\lambda^2)
\big\la G_2(x,\lambda)D^{-1}(\lambda)F_1^t(y,\lambda)\sigma_3 \phi(x), \psi(y)\big\ra \, d\lambda\chi_{[x< y]} \,dxdy \\
&= \int_{\R^2}\int_{-\infty}^\infty e^{it\lambda^2}
\lambda\chi(\lambda^2) \big\la
F_1(x,\lambda)D^{-1}(\lambda)G_2^t(y,\lambda)\sigma_3 \phi(-x),
\psi(-y)\big\ra \, d\lambda\chi_{[x> y]} \,dxdy.
\end{align*}
Hence it suffices to bound the first. To this end we need to
consider three cases: $x\ge0\ge y$, $0\ge x\ge y$, and $x\ge
y\ge0$.

\medskip
\noindent {\em Case 1:} $x\ge0\ge y$

We write the $\lambda$-integral in \eqref{eq:pot1} as the sum of
four pieces, according to the various possibilities for the
asymptotic behavior as $x\to\infty$ or $-y\to\infty$:
\begin{align}
&  \int_{-\infty}^\infty e^{it\lambda^2} \lambda\chi(\lambda^2)
 F_1(x,\lambda)D^{-1}(\lambda)G_2^t(y,\lambda)\sigma_3  \, d\lambda \nn\\
& = \int_{-\infty}^\infty e^{it\lambda^2} \lambda\chi(\lambda^2)
 [f_1(x,\lambda)\; 0]D^{-1}(\lambda)[f_1(-y,\lambda)\; 0]^t\sigma_3  \, d\lambda \label{eq:comb1}\\
& \quad + \int_{-\infty}^\infty e^{it\lambda^2}
\lambda\chi(\lambda^2)
 [f_1(x,\lambda)\; 0]D^{-1}(\lambda)[0\;f_3(-y,\lambda)]^t\sigma_3  \, d\lambda \label{eq:comb2}\\
& \quad +\int_{-\infty}^\infty e^{it\lambda^2}
\lambda\chi(\lambda^2)
 [0 \;f_3(x,\lambda)]D^{-1}(\lambda)[f_1(-y,\lambda)\; 0]^t\sigma_3  \, d\lambda \label{eq:comb3}\\
& \quad + \int_{-\infty}^\infty e^{it\lambda^2}
\lambda\chi(\lambda^2)
 [0\;f_3(x,\lambda)]D^{-1}(\lambda)[0\; f_3(-y,\lambda)]^t\sigma_3  \, d\lambda. \label{eq:comb4}
\end{align}
Now pick another cut-off function $\tilde\chi$ so that
$\tilde\chi\chi=\chi$. Then~\eqref{eq:comb1} is estimated as
follows (with the Fourier transform relative to $\lambda$):
\begin{align}
& |\eqref{eq:comb1}| \le \Big|\int_{-\infty}^\infty
e^{i[t\lambda^2+\lambda(x-y)]} \lambda\chi(\lambda^2)
 [e^{-i\lambda x}f_1(x,\lambda)\; 0]D^{-1}(\lambda)[e^{i\lambda y}f_1(-y,\lambda)\; 0]^t \sigma_3 \, d\lambda\Big| \nn \\
& \le C|t|^{-\half} \int_{-\infty}^\infty \Big|\big[
\lambda\chi(\lambda^2) D^{-1}(\lambda)\big]^{\vee}\ast
\big[\tilde\chi(\lambda^2)
 e^{-i\lambda x}f_1(x,\lambda)\big]^{\vee}\ast
\big[ \tilde\chi(\lambda^2)  e^{i\lambda y}f_1(-y,\lambda)\big]^{\vee}(\xitil)\Big| \, d\xitil  \nn \\
&\le C|t|^{-\half} \big\|\big[ \lambda\chi(\lambda^2)
D^{-1}(\lambda)\big]^{\vee}\big\|_1 \big\|
\big[\tilde\chi(\lambda^2)
 e^{-i\lambda x}f_1(x,\lambda)\big]^{\vee}\big\|_1
\big\|\big[ \tilde\chi(\lambda^2)  e^{i\lambda
y}f_1(-y,\lambda)\big]^{\vee}\big\|_1. \label{eq:drei}
\end{align}
Since $D^{-1}(\lambda)$ is smooth, the first norm on the
right-hand side of \eqref{eq:drei} is simply a constant. On the
other hand, we claim that
\[
 \sup_{x\ge0} \big\| \big[\tilde\chi(\lambda^2) e^{-i\lambda x}f_1(x,\lambda)\big]^{\vee}\big\|_1 \le C \text{\ \ and\ \ }
 \sup_{y\le0} \big\| \big[\tilde\chi(\lambda^2) e^{i\lambda y}f_1(-y,\lambda)\big]^{\vee}\big\|_1 \le C.
\]
This is an easy consequence of Lemma~\ref{lem:f1}. Indeed,
\eqref{eq:f1der} insures that
\[ |\partial_\lambda^2 [\tilde\chi(\lambda^2) e^{-i\lambda x}f_1(x,\lambda)]| \le C\, x^2 e^{-\gamma x} \]
for all $\lambda\in \R$ and all $x\ge x_0$. Since we are dealing
with a compact $\lambda$-interval it also follows that
\[ \sup_{x\ge0,\lambda\in\R} |(1-\partial_\lambda^2) [\tilde\chi(\lambda^2) e^{-i\lambda x}f_1(x,\lambda)]| \le C. \]
Hence,
\[ \sup_{\xitil\in \R, x\ge0} (1+\xitil^2) \big| [\tilde\chi(\lambda^2) e^{-i\lambda x}f_1(x,\lambda)]^{\vee}(\xitil)\big| \le C\]
and thus \beeq \label{eq:f1xlam}
 \sup_{x\ge0} \big\|[\tilde\chi(\lambda^2) e^{-i\lambda x}f_1(x,\lambda)]^{\vee}\big\|_1 \le C,
\eneq as claimed. Hence the right-hand side of~\eqref{eq:drei} is
at most $C|t|^{-\half}$, as desired.

For the second term \eqref{eq:comb2} we estimate (recall $y\le 0$)
\begin{align}
& |\eqref{eq:comb2}| \le \Big|\int_{-\infty}^\infty
e^{i[t\lambda^2+x\lambda]}e^{-\mu |y|} \lambda\chi(\lambda^2)
 [e^{-ix\lambda}f_1(x,\lambda)\; 0]D^{-1}(\lambda)[0\;e^{\mu |y|}f_3(|y|,\lambda)]^t\sigma_3  \, d\lambda\Big| \nn \\
&\le C|t|^{-\half} \big\|\big[ \lambda\chi(\lambda^2)
D^{-1}(\lambda)e^{-\mu |y|}\big]^{\vee}\big\|_1 \big\|
\big[\tilde\chi(\lambda^2)
 e^{-i\lambda x}f_1(x,\lambda)\big]^{\vee}\big\|_1
\big\|\big[ \tilde\chi(\lambda^2)  e^{\mu
|y|}f_3(|y|,\lambda)\big]^{\vee}\big\|_1. \label{eq:drei'}
\end{align}
It follows from \eqref{eq:f3est} that
\[
\sup_{y}\big\|\big[ \tilde\chi(\lambda^2)  e^{\mu
|y|}f_3(|y|,\lambda)\big]^{\vee}\big\|_1 < \infty.
\]
Furthermore, \eqref{eq:L1mu} and Lemma~\ref{lem:Dsmooth} imply
that
\[
\sup_{y}\big\|\big[ \lambda\chi(\lambda^2) D^{-1}(\lambda)e^{-\mu
|y|}\big]^{\vee}\big\|_1 \le C.
\]
Combining these estimates with \eqref{eq:f1xlam} shows that the
entire right-hand side of \eqref{eq:drei'} is uniformly bounded in
$x,y\in\R$.

The same type of arguments show that
\[ |\eqref{eq:comb3}|+|\eqref{eq:comb4}|\le C|t|^{-\half}\]
uniformly in $x\ge0\ge y$, which finishes Case~1.

\medskip
\noindent {\em Case 2:} $0\ge x\ge y$

In this case, we need to use \eqref{eq:ABdef}. Thus, we write the
$\lambda$-integral in \eqref{eq:pot1} as the sum of two integrals:
\begin{align}
& \int_{-\infty}^\infty e^{it\lambda^2} \lambda\chi(\lambda^2)
 F_1(x,\lambda)D^{-1}(\lambda)G_2^t(y,\lambda)\sigma_3  \, d\lambda
 = \int_{-\infty}^\infty e^{it\lambda^2} \lambda\chi(\lambda^2)
 G_1(x,\lambda)A(\lambda)D^{-1}(\lambda)G_2^t(y,\lambda)\sigma_3  \, d\lambda \label{eq:Aint} \\
& \quad +\int_{-\infty}^\infty e^{it\lambda^2}
\lambda\chi(\lambda^2)
 G_2(x,\lambda)B(\lambda)D^{-1}(\lambda)G_2^t(y,\lambda)\sigma_3  \, d\lambda. \label{eq:Bint},
\end{align}
each of which is itself broken up into four pieces just as in
\eqref{eq:comb1}--\eqref{eq:comb2}. Thus, starting
with~\eqref{eq:Aint},
\begin{align}
& \int_{-\infty}^\infty e^{it\lambda^2} \lambda\chi(\lambda^2)
 G_1(x,\lambda)A(\lambda)D^{-1}(\lambda)G_2^t(y,\lambda)\sigma_3  \, d\lambda \nn \\
& = \int_{-\infty}^\infty e^{it\lambda^2} \lambda\chi(\lambda^2)
 [f_2(-x,\lambda)\; 0]A(\lambda)D^{-1}(\lambda)[f_1(-y,\lambda)\; 0]^t\sigma_3  \, d\lambda \label{eq:comb1'}\\
& \quad + \int_{-\infty}^\infty e^{it\lambda^2}
\lambda\chi(\lambda^2)
 [f_2(-x,\lambda)\; 0]A(\lambda)D^{-1}(\lambda)[0\;f_3(-y,\lambda)]^t\sigma_3  \, d\lambda \label{eq:comb2'}\\
& \quad +\int_{-\infty}^\infty e^{it\lambda^2}
\lambda\chi(\lambda^2)
 [0 \;f_4(-x,\lambda)]A(\lambda)D^{-1}(\lambda)[f_1(-y,\lambda)\; 0]^t\sigma_3  \, d\lambda \label{eq:comb3'}\\
& \quad + \int_{-\infty}^\infty e^{it\lambda^2}
\lambda\chi(\lambda^2)
 [0\;f_4(-x,\lambda)]A(\lambda)D^{-1}(\lambda)[0\; f_3(-y,\lambda)]^t\sigma_3  \, d\lambda. \label{eq:comb4'}
\end{align}
We remark that in view of Corollary~\ref{cor:imbed},
\[
\lambda A(\lambda)D(\lambda)^{-1} = \bm \frac{1}{2i} & 0 \\ 0 &
\frac{-\lambda}{2\mu} \endm.
\]
In particular, the expression on the left-hand side is smooth for
all $\lambda\in\R$. Moreover, the diagonal form of the right-hand
side implies that $\eqref{eq:comb2'}=\eqref{eq:comb3'}=0$ (which
in the case of~\eqref{eq:comb3'} is crucial), as well as
\begin{align*}
\eqref{eq:comb1'} &= \frac{1}{2i}\int_{-\infty}^\infty
e^{i[t\lambda^2+(x-y)\lambda]} \chi(\lambda^2)
 [e^{-ix\lambda}f_2(-x,\lambda)\; 0] [e^{iy\lambda}f_1(-y,\lambda)\; 0]^t\sigma_3  \, d\lambda \\
\eqref{eq:comb4'} &= -\int_{-\infty}^\infty e^{it\lambda^2}
e^{-|x-y|\mu} \frac{\lambda}{2\mu} \chi(\lambda^2)
 [0\;e^{-\mu|x|}f_4(|x|,\lambda)][0\; e^{\mu|y|}f_3(|y|,\lambda)]^t\sigma_3  \, d\lambda.
\end{align*}
The same type of arguments as in Case~1 involving the Fourier
transform in $\lambda$ show that both of these expressions are
bounded by $C|t|^{-\half}$ uniformly in $0\ge x \ge y$. The only
new ingredient in this case is the estimate
\[
\sup_{x}\big\|\big[ \tilde\chi(\lambda^2)  e^{-\mu
|x|}f_4(|x|,\lambda)\big]^{\vee}\big\|_1 < \infty,
\]
which follows from Corollary~\ref{cor:f4der}. Continuing with
\eqref{eq:Bint}, we write
\begin{align*}
& \int_{-\infty}^\infty e^{it\lambda^2} \lambda\chi(\lambda^2)
 G_2(x,\lambda)B(\lambda)D^{-1}(\lambda)G_2^t(y,\lambda)\sigma_3  \, d\lambda \nn \\
& = \int_{-\infty}^\infty e^{i[t\lambda^2-(x+y)\lambda]}
\lambda\chi(\lambda^2)
 [e^{ix\lambda}f_1(-x,\lambda)\; 0]B(\lambda)D^{-1}(\lambda)[e^{iy\lambda}f_1(-y,\lambda)\; 0]^t\sigma_3  \, d\lambda \\
& \quad + \int_{-\infty}^\infty
e^{i[t\lambda^2-x\lambda]}e^{-\mu|y|} \lambda\chi(\lambda^2)
 [e^{ix\lambda}f_1(-x,\lambda)\; 0]B(\lambda)D^{-1}(\lambda)[0\;e^{\mu|y|}f_3(|y|,\lambda)]^t\sigma_3  \, d\lambda \\
& \quad +\int_{-\infty}^\infty e^{i[t\lambda^2-y\lambda ]}
e^{-\mu|x|}\lambda\chi(\lambda^2)
 [0 \;e^{\mu|x|}f_3(|x|,\lambda)]B(\lambda)D^{-1}(\lambda)[e^{iy\lambda}f_1(-y,\lambda)\; 0]^t\sigma_3  \, d\lambda \\
& \quad + \int_{-\infty}^\infty e^{it\lambda^2}
e^{-\mu(|x|+|y|)}\lambda\chi(\lambda^2)
 [0\;e^{\mu|x|}f_3(|x|,\lambda)]B(\lambda)D^{-1}(\lambda)[0\; e^{\mu|y|}f_3(|y|,\lambda)]^t\sigma_3  \, d\lambda.
\end{align*}
We remark that in view of \eqref{eq:wron2}
\[
\lambda B(\lambda)D(\lambda)^{-1} = \bm -\frac{1}{2i} & 0 \\ 0 &
\frac{\lambda}{2\mu}\endm
\calW[F_1(\cdot,\lambda),G_1(\cdot,\lambda)]^tD(\lambda)^{-1}
\]
is a smooth function in $\lambda$. Hence, the same methods as
before prove that each of these integrals are bounded by
$C|t|^{-\half}$ uniformly in $0\ge x\ge y$.

\medskip
\noindent {\em Case 3:} $x\ge y\ge 0$

This case can be reduced to the previous one. Indeed, using that
$F_1(x,\lambda)=G_2(-x,\lambda)$, as well as~\eqref{eq:Dsymm}, one
obtains
\begin{align*}
& \int_{-\infty}^\infty e^{it\lambda^2} \lambda\chi(\lambda^2)
 F_1(x,\lambda)D^{-1}(\lambda)G_2^t(y,\lambda)\sigma_3  \, d\lambda \\
& = \int_{-\infty}^\infty e^{it\lambda^2} \lambda\chi(\lambda^2)
 G_2(-x,\lambda)D^{-1}(\lambda)F_1^t(-y,\lambda)\sigma_3  \, d\lambda \\
& = \sigma_3\Big[\int_{-\infty}^\infty e^{it\lambda^2}
\lambda\chi(\lambda^2)
 F_1(-y,\lambda)D^{-1}(\lambda)G_2^t(-x,\lambda)\sigma_3  \, d\lambda \Big]^t\sigma_3.
\end{align*}
Since $0\ge -y\ge -x$, and since $\sigma_3$ and the transpose do
not affect the point-wise estimate of Case~2, we are done.
\end{proof}

The previous estimate also allows for the introduction of
derivatives. We will consider only at most two derivatives,
although more are possible.

\begin{cor}
\label{cor:inter} Let $\Hil$ be admissible. Then
\[ \big\|e^{it\Hil} P_s f\big\|_{W^{k,p'}(\R)} \le Ct^{-\frac12(\frac{1}{p}-\frac{1}{p'})} \|f\|_{W^{k,p}(\R)}\]
for  $0\le k\le 2$ and $1<p\le2$.
\end{cor}
\begin{proof}
The case $k=0$ is obtained by interpolating between
Lemmas~\ref{lem:L2stable} and~\ref{prop:disp1} and holds for the
entire range $1\le p\le2$. We need to require $p>1$ only for the
derivatives. If $a$ is sufficiently large, then
\[ (\Hil-ia)^{-1}: L^2\times L^2 \to W^{2,2}\times W^{2,2} \]
is an isomorphism. More generally,
\[ (\Hil-ia)^{-\half}: L^p\times L^p \to W^{2,p}\times W^{2,p}\]
is an isomorphism for $1<p<\infty$. This can be seen from the
resolvent identity
\[ (\Hil-ia)^{-1} = (\Hil_0-ia)^{-1}[1+V(\Hil_0-ia)^{-1}]^{-1},\]
since $\|V\|_\infty <\infty$ implies that
\[ \|V(\Hil_0-ia)^{-1}\|_{p\to p}<\half\]
if $a$ is large enough, and because
\[ (\Hil_0-ia)^{-\half}: L^p\times L^p \to W^{2,p}\times W^{2,p}\]
for any $a\ne0$ as an isomorphism. Hence,
\begin{align*}
\|\Laplace e^{it\Hil}P_s f\|_{p'} &\le C\|(\Hil-ia) e^{it\Hil} f\|_{p'} =  \|e^{it\Hil}(\Hil-ia)f\|_{p'}\\
& \le Ct^{-\frac12(\frac{1}{p}-\frac{1}{p'})} \|(\Hil-ia)f\|_p \le
Ct^{-\frac12(\frac{1}{p}-\frac{1}{p'})} \|f\|_{W^{2,p}(\R)}.
\end{align*}
This gives the case $k=2$ of the lemma, whereas $k=1$ follows by
interpolating between $k=0$ and $k=2$.
\end{proof}

Finally, as in~\cite{Sch} one can now derive Strichartz estimates
on the evolution~$e^{it\Hil}P_s$, even with derivatives.

\begin{cor}
\label{cor:strich_der} Let $\Hil$ be admissible. Then the
following Strichartz estimates hold:
\begin{align}
 \|e^{-it\Hil} P_s f\|_{L^r_t(W^{k,p}_x)} &\le C \|f\|_{W^{k,2}} \label{eq:Strich1} \\
 \Big\|\int_0^t e^{-i(t-s)\Hil} P_s F(s)\,ds \Big\|_{L^r_t(W^{k,p}_x)} &\le C \|F\|_{L^{a'}_t(W^{k,b'}_x)} \label{eq:Strich2},
\end{align}
provided $(r,p), (a,b)$ are admissible, i.e., $4<r\le\infty$ and
$\frac{2}{r}+\frac{1}{p}=\frac12$ and the same for $(a,b)$. Here
$k$ is an integer, $0\le k\le2$.
\end{cor}
\begin{proof}
We first show how to reduce matters to $k=0$. As in the previous
proof, we rely on the fact that (because of
$\|V\|_\infty<\infty$),
\[ \|\Laplace f\|_q\le C\|(\Hil-ia) f\|_q \]
for any $1<q<\infty$. Hence,
\begin{align*}
\|e^{-it\Hil} P_s f\|_{L^r_t(W^{2,p}_x)} &\le C \|(\Hil-ia)
e^{-it\Hil} P_s f\|_{L^r_t(L^p_x)}
= \| e^{-it\Hil} P_s (\Hil-ia) f\|_{L^r_t(L^p_x)}  \\
&\le C\|(\Hil-ia)f\|_2 \le C\|f\|_{W^{2,2}},
\end{align*}
which is \eqref{eq:Strich1} for $k=2$. Similarly, one
proves~\eqref{eq:Strich2} for $k=2$. The case $k=1$ is then
obtained by interpolation between $k=2$ and $k=0$. For $k=0$ we
use an argument from~\cite{RodSch}, Section~4. Let ($\calS$ for
``Strichartz'')
\[ (\calS F)(t,x)=\int_0^t (e^{-i(t-s)\Hil}P_s\,F(s,\cdot))(x)\, ds.\]
In this proof it will be understood that all times are $\ge0$.
Then by~\eqref{eq:L2stable},
\[ \|\calS F\|_{L^\infty_t(L^2_x)} \le C\|F\|_{L^1_t(L^2_x)}, \]
and more generally, by the usual fractional integration argument
based on Proposition~\ref{prop:disp1}, \beeq \label{eq:S1} \|\calS
F\|_{L^r_t(L^p_x)} \le C\|F\|_{L^{r'}_t(L^{p'}_x)} \eneq for any
admissible pair $(r,p)$. In the unitary case this
implies~\eqref{eq:Strich1} via a $TT^*$ argument, but this
reasoning does not apply here. Instead, we rely on a Kato theory
type approach as in \cite{RodSch}, Section~4. Since
$\Hil=\Hil_0+V$, Duhamel's formula yields \beeq \label{eq:du1}
 e^{-it\Hil}P_s = e^{-it\Hil_0}P_s -i\int_0^t e^{-i(t-s)\Hil_0}Ve^{-is\Hil}P_s \, ds.
\eneq Set $\rho(x)=\la x\ra^{-1-}$, say, and define \beeq
\label{eq:Mtil} \tilde{M}= \bm \rho & 0 \\ 0 & \rho \endm. \eneq
With $V=\tilde{M}\tilde{M}^{-1}V$, observe firstly that
\[
\Big\| \int_0^\infty e^{-i(t-s)\Hil_0} \tilde{M} g(s)\,
ds\Big\|_{L^r_t(L^p_x)} \les \Big\| \int_0^\infty e^{is\Hil_0}
\tilde{M} g(s)\, ds\Big\|_{L^2} \le C\|g\|_{L^2_s(L^2_x)},
\]
where the last inequality is the dual of the smoothing bound
\[ \int_0^\infty \Big\| \tilde{M}e^{-is\Hil_0^*} \psi\Big\|_2^2 \, ds \le C\|\psi\|_2^2. \]
Here ``smoothing'' is used in the sense of Kato's theory,
see~\cite{kato} and~\cite{RS4}. Now one applies the Christ-Kiselev
lemma, see~\cite{CrKi} or~\cite{RodSch}, to conclude that
\[
\Big\| \int_0^t e^{-i(t-s)\Hil_0} \tilde{M} g(s)\,
ds\Big\|_{L^r_t(L^p_x)} \le C\|g\|_{L^2_s(L^2_x)}
\]
for any admissible pair $(r,p)$. Hence, continuing
in~\eqref{eq:du1}, one obtains (using that $\|P_s f\|_2\le
C\|f\|_2$)
\[
\|e^{-it\Hil}P_s f \|_{L^r_t(L^p_x)} \le C\|f\|_2 + \Big\|
\tilde{M}^{-1}V e^{-is\Hil}P_s f\Big\|_{L^2_s(L^2_x)}.
\]
It remains to show that $\tilde{M}^{-1}V$ is $\Hil P_s$-smoothing,
i.e., \beeq \label{eq:smft} \Big\| \tilde{M}^{-1}V e^{-is\Hil}P_s
f\Big\|_{L^2_s(L^2_x)} \le C\|f\|_2. \eneq It follows from
Lemma~\ref{lem:rep} that the integrand here is the same as
\[
\tilde{M}^{-1}V e^{-is\Hil}P_s f
 = \frac{1}{2\pi i}\int_{|\lambda|\ge 1} e^{it\lambda}\,\tilde{M}^{-1}V [(\Hil-(\lambda+i0))^{-1}-(\Hil-(\lambda-i0))^{-1}]f\,d\lambda.
\]
Hence, applying Plancherel in $t$ yields
\begin{align*}
&\int_{-\infty}^\infty\|\tilde{M}^{-1}V e^{-is\Hil}P_s f\|_2^2\, ds  \\
 &= \int_{|\lambda|\ge 1}  \big\|\tilde{M}^{-1}V [(\Hil-(\lambda+i0))^{-1}-(\Hil-(\lambda-i0))^{-1}]f\big\|_2^2\,d\lambda \\
&= 2\int_{0}^\infty \lambda\big\|\tilde{M}^{-1}V [(\Hil-(\lambda^2+1+i0))^{-1}-(\Hil-(\lambda^2+1-i0))^{-1}]f\big\|_2^2\,d\lambda \\
& \quad + 2\int_{0}^\infty \lambda \big\|\tilde{M}^{-1}V
[(\Hil-(-\lambda^2+1+i0))^{-1}-(\Hil-(-\lambda^2+1-i0))^{-1}]f\big\|_2^2\,d\lambda.
\end{align*}
In view of Lemma~\ref{lem:jump},
\begin{align*}
&\int_{-\infty}^\infty\|\tilde{M}^{-1}V e^{-is\Hil}P_s f\|_2^2\, ds  \\
&= \frac{1}{2}\int_{0}^\infty \lambda^{-1}
\big\|\tilde{M}^{-1}V\calE_+(x,\lambda)\la\calE_+^*(\cdot,\lambda)\sigma_3,
f\ra\big\|_{L^2_x}^2\,d\lambda
 +\frac{1}{2}\int_{0}^\infty \lambda^{-1} \big\|\tilde{M}^{-1}V\calE_-(x,\lambda)\la\calE_-^*(\cdot,\lambda)\sigma_3, f\ra\big\|_{L^2_x}^2\,d\lambda.
\end{align*}
By the definition of $\calE$ and Lemma~\ref{lem:scat},
\[ \big\|\tilde{M}^{-1}V\calE_\pm(x,\lambda)\big\|_{L^2_x}^2 \le C|\lambda|^2(1+|\lambda|)^{-2}. \]
Hence,
\begin{align*}
&\frac{1}{2}\int_{0}^\infty \lambda^{-1}
\big\|\tilde{M}^{-1}V\calE_+(x,\lambda)\la\calE_+^*(\cdot,\lambda)\sigma_3,
f\ra\big\|_{L^2_x}^2\,d\lambda
 +\frac{1}{2}\int_{0}^\infty \lambda^{-1} \big\|\tilde{M}^{-1}V\calE_-(x,\lambda)\la\calE_-^*(\cdot,\lambda)\sigma_3, f\ra\big\|_{L^2_x}^2\,d\lambda \\
& \le C\, \int_0^\infty \big|\la\calE_+^*(\cdot,\lambda)\sigma_3,
f\ra\big|_{L^2_x}^2\,d\lambda+ C\, \int_0^\infty
\big|\la\calE_-^*(\cdot,\lambda)\sigma_3, f\ra\big|_{L^2_x}^2
\,d\lambda \le C\|f\|_2^2,
\end{align*}
where the final inequality is a Plancherel type property that
follows from Lemma~\ref{lem:scat}. Hence~\eqref{eq:smft} holds.
The conclusion is that
\[  \|e^{-it\Hil}P_s f \|_{L^r_t(L^p_x)} \le C\|f\|_2 \]
for any admissible $(r,p)$, which is~\eqref{eq:Strich1}. The proof
of \eqref{eq:Strich2} is now the usual interpolation argument.
Indeed, in view of the preceding one has the following bounds
on~$\calS$ for any admissible pair $(r,p)$:
\begin{align}
\calS: &L^1_t(L^2_x) \to L^r_t(L^p_x) \label{eq:St1}\\
\calS: &L^{r'}_t(L^{p'}_x) \to L^r_t(L^p_x) \label{eq:St2}\\
\calS: &L^{r'}_t(L^{p'}_x) \to L^\infty_t(L^2_x).  \label{eq:St3}
\end{align}
These estimates arise as follows: \eqref{eq:St2} is
exactly~\eqref{eq:S1}, whereas~\eqref{eq:St1} follows
from~\eqref{eq:Strich1} by means of Minkowski's inequality.
Finally, \eqref{eq:St3} is dual to the bound \beeq
\label{eq:again} \Big\|\int_t^\infty e^{i(t-s)\Hil^*}\tilde{P}_c
G(s)\,ds \Big\|_{ L^r_t(L^p_x)} \les \|G\|_{L^1_t(L^2_x)}. \eneq
Here $\tilde{P}_c$ corresponds to $\Hil^*$ in the same way that
$P_s$ does to $\Hil$. In particular, one has
\[ \|e^{-it\Hil^*}\tilde{P}_c\|_{1\to\infty}\le Ct^{-\frac12}\]
and therefore, \eqref{eq:again} is derived be the same methods
as~\eqref{eq:St1}. It is important to notice  that
$P_s^*=\tilde{P}_s$ which is essential for the duality argument
here. This can be seen, for example, by writing the Riesz
projections onto (generalized) eigenspaces as contour integrals
around circles surrounding the eigenvalues. Since the (complex)
eigenvalues always come in pairs, the adjoints have the desired
property. Interpolating between \eqref{eq:St1} and \eqref{eq:St2}
yields~\eqref{eq:Strich2} for the range $a'\le r'$ or $a\ge r$,
whereas interpolating between \eqref{eq:St1} and \eqref{eq:St2}
yields~\eqref{eq:Strich2} in the range $a\le r$.
\end{proof}

\section{Dispersive estimates: The weighted case}
\label{sec:weigh}

In this section we obtain the decay rate of $t^{-\frac32}$ on
$e^{it\Hil}P_s$. The latter will exploit the fact that absence of
resonances at the thresholds leads to better time-decay, albeit at
the cost of a linear weight. This property was already used by
Buslaev and Perelman~\cite{BP1}. However, they worked on $L^2$ and
with a loss of a higher power $x$, namely of $x^{3.5+\eps}$. This
would not be sufficient for our purposes. For the same estimate in
the scalar case, see~\cite{Sch2}.

\begin{prop}
\label{prop:disp2} Let $\Hil$ be admissible, see
Definition~\ref{def:admH}. Then for all $t\ne0$,
\[ \|\la x\ra^{-1}e^{it\Hil}P_s f\|_{\infty}\le C|t|^{-\frac32}\|\la x\ra f\|_1\]
with some $C=C(V)$.
\end{prop}
\begin{proof}
As in the proof of Proposition~\ref{prop:disp1}, we divide the
arguments into $\lambda$ large and not large. In the former case,
the estimates that lead to the unweighted $L^1\to L^\infty$ bound
apply here almost verbatim up to an additional integration by
parts, whereas in the latter case we will need to use the Fourier
representation from Proposition~\ref{prop:fourier}. As before, we
start from the representation formula
\begin{align*}
\la e^{it\Hil} P_s \phi,\psi \ra &= \frac{e^{it}}{\pi i}
\int_0^\infty e^{it\lambda^2}\lambda \big \la
\big((\Hil-(\lambda^2+1+i0))^{-1}-
(\Hil-(\lambda^2+1-i0))^{-1}\big)\phi,\psi \big\ra\, d\lambda \\
& + \frac{e^{-it}}{\pi i} \int_0^\infty e^{-it\lambda^2} \lambda
\big \la \big((\Hil-(-\lambda^2-1+i0))^{-1}-
(\Hil-(-\lambda^2-1-i0))^{-1}\big)\phi,\psi \big\ra\, d\lambda,
\end{align*}
which holds in the principal value sense if $\phi,\psi\in\calS$.
With the same $\chi$ as in the proof of
Proposition~\ref{prop:disp1}, we arrive at
\begin{align*}
\nn
&\big|\la e^{it\Hil} (1-\chi(H-I)) P_s^+ \phi,\psi \ra\big| \\
&\le C\sum_{n=0}^\infty |t|^{-1}\Big|\int_{0}^\infty
e^{it\lambda^2}\,\frac{d}{d\lambda}\Big\{(1-\chi(\lambda^2))
\big[ \big\la (\Hil_0-(\lambda^2+1+i0))^{-1}\big(V(\Hil_0-(\lambda^2+1+i0))^{-1}\big)^n \phi,\psi \big\ra \nn \\
& - \big\la
(\Hil_0-(\lambda^2+1-i0))^{-1}\big(V(\Hil_0-(\lambda^2+1-i0))^{-1}\big)^n
\phi,\psi \big\ra\big]\Big\} \, d\lambda \Big|.
\end{align*}
If the $\lambda$-derivative falls on $1-\chi(\lambda^2)$, then the
exact same arguments apply which we used in the unweighted case.
If the derivative falls on the resolvents, then weights
$|x_j-x_{j-1}|$ appear. However, these are bounded by
$|x_j|+|x_{j-1}|$ and can therefore be absorbed either by $V$ or
the test functions $\phi,\psi$. The conclusion is that
\[ \big|\la e^{it\Hil} (1-\chi(H-I)) P_s^+ \phi,\psi \ra\big| \le C|t|^{-\frac32} \|\la x\ra\phi\|_1\|\la x\ra\psi\|_1,\]
and therefore also
\[ \big|\la e^{it\Hil} (1-\chi(-H+I)) P_s^- \phi,\psi \ra\big| \le C|t|^{-\frac32} \|\la x\ra\phi\|_1\|\la x\ra\psi\|_1,\]
as desired. Cf.~the treatment of~\eqref{eq:PU1}--\eqref{eq:PW2}
for details.

Next, we use Proposition~\ref{prop:fourier} to write
\begin{align*}
\la e^{it\Hil} \chi(H-I) P_s^+ \phi,\psi \ra &=
\frac{e^{it}}{2\pi}\int_{-\infty}^\infty e^{it\lambda^2}
\chi(\lambda^2) \la\phi,\sigma_3 e_+(.,\lambda)\ra
\la e_+(\cdot,\lambda), \psi\ra\, d\lambda \\
\la e^{it\Hil} \chi(-H+I) P_s^- \phi,\psi \ra &=
\frac{e^{-it}}{2\pi}\int_{-\infty}^\infty e^{it\lambda^2}
\chi(\lambda^2) \la\phi,\sigma_3 e_-(.,\lambda)\ra \la
e_-(\cdot,\lambda), \psi\ra\, d\lambda.
\end{align*}
It suffices to estimate the first integral. Integrating by parts
in $\lambda$ yields
\begin{align}
& \big\la e^{it\Hil}\chi(H-I)P_s^+\phi,\psi\big\ra  \nn \\
& = \frac{1}{\pi it}\int_0^\infty e^{it\lambda^2}
\partial_\lambda\Big[\chi(\lambda^2) \lambda \la
\phi,\sigma_3\,F_1(\cdot,\lambda)D^{-1}(\lambda)\vece\ra
\la F_1(\cdot,\lambda)D^{-1}(\lambda)\vece,\psi\ra \Big] \, d\lambda \label{eq:F1teil}\\
& \quad + \frac{1}{\pi it}\int_{-\infty}^0 e^{it\lambda^2}
\partial_\lambda\Big[\chi(\lambda^2)\lambda \la
\phi,\sigma_3\,G_2(\cdot,-\lambda)D^{-1}(\lambda)\vece\ra \la
G_2(\cdot,-\lambda)D^{-1}(\lambda)\vece,\psi\ra \Big]\,
d\lambda.\nn
\end{align}
By symmetry, it will suffice to treat the integral
\eqref{eq:F1teil} involving $F_1(\cdot,\lambda)$. We distinguish
three cases, depending on where the derivative $\partial_\lambda$
falls. We start with the integral \beeq \label{eq:part1}
 \int_0^\infty e^{it\lambda^2} \omega(\lambda) \la \phi,\sigma_3\,F_1(\cdot,\lambda)D^{-1}(\lambda)\vece\ra
\la F_1(\cdot,\lambda)D^{-1}(\lambda)\vece,\psi\ra  \, d\lambda,
\eneq where we have set
$\omega(\lambda)=\partial_\lambda[\chi(\lambda^2)\lambda]$. By the
preceding, $\omega$ is a smooth function with compact support in
$[0,\infty)$. As usual, we will estimate \eqref{eq:part1} by means
of a Fourier transform in $\lambda$. Since we are working on a
half-line, this will actually be a cosine transform. Let
$\tilde\omega$ be another cut-off function satisfying
$\omega\tilde\omega=\omega$. Then
\begin{align}
& \Big| \int_0^\infty e^{it\lambda^2} \omega(\lambda) \la
\phi,\sigma_3\,F_1(\cdot,\lambda)D^{-1}(\lambda)\vece\ra
\la F_1(\cdot,\lambda)D^{-1}(\lambda)\vece,\psi\ra \, d\lambda \Big|  \nn \\
& \le C|t|^{-\half} \| \,[\omega \la
\phi,\sigma_3\,F_1(\cdot,\lambda)D^{-1}(\lambda)\vece\ra ]^{\vee}
\|_1 \|\, [\tilde\omega \la
F_1(\cdot,\lambda)D^{-1}(\lambda)\vece,\psi\ra ]^{\vee} \|_1.
\label{eq:2L1}
\end{align}
It remains to show that \beeq \nn [\omega \la
\phi,\sigma_3\,F_1(\cdot,\lambda)D^{-1}(\lambda)\vece\ra]^{\vee}(u)
:= \int_0^\infty \cos(u\lambda) \omega(\lambda) \la
\phi,\sigma_3\,F_1(\cdot,\lambda)D^{-1}(\lambda)\vece\ra \,
d\lambda \eneq satisfies \beeq \nn \int_{-\infty}^\infty\Big|
[\omega \la
\phi,\sigma_3\,F_1(\cdot,\lambda)D^{-1}(\lambda)\vece\ra]^{\vee}(u)\Big|\,
du \le C\|\la x\ra \phi\|_1. \eneq The second $L^1$-norm in
\eqref{eq:2L1} is treated the same way. This means that we need to
prove that \beeq \label{eq:omf1} \int_{-\infty}^\infty\Big|
[\omega F_1(x,\lambda)D^{-1}(\lambda)\vece]^{\vee}(u)\Big|\, du
\le C\la x\ra \eneq for all $x\in\R$. We will consider the cases
$x\ge0$ and $x\le0$ separately. In the former case,
\begin{align}
& [\omega F_1(x,\cdot)D^{-1}(\lambda)\vece]^{\vee}(u) :=
\int_0^\infty \cos(u\lambda) \omega(\lambda)
F_1(x,\lambda)D^{-1}(\lambda)\vece\, d\lambda \nn \\
&= \half\int_0^\infty e^{i(x+u)\lambda} \omega(\lambda) e^{-ix\lambda}f_1(x,\lambda)\la \vece,D^{-1}(\lambda)\vece\ra\, d\lambda \label{eq:f1one'} \\
&  + \half\int_0^\infty e^{i(x-u)\lambda} \omega(\lambda) e^{-ix\lambda}f_1(x,\lambda)\la \vece,D^{-1}(\lambda)\vece\ra\, d\lambda \label{eq:f1two}\\
& + \int_0^\infty  \cos(u\lambda) \omega(\lambda) e^{-\mu x}
e^{\mu x}f_3(x,\lambda) \la \vece',D^{-1}(\lambda)\vece\ra\,
d\lambda, \label{eq:f3}
\end{align}
where $\vece'=\binom{0}{1}$. We integrate by parts in
\eqref{eq:f1one'}--\eqref{eq:f3}:
\begin{align}
& \eqref{eq:f1one'}+ \eqref{eq:f1two} \nn \\
&= -\frac{1}{2i(x+u)} \omega(0)f_1(x,0)\la \vece,D^{-1}(0)\vece\ra - \frac{1}{2i(x-u)} \omega(0)f_1(x,0)\la \vece,D^{-1}(0)\vece\ra \label{eq:byparts1}\\
&\quad  -\frac{1}{2i(x+u)}\int_0^\infty e^{i(x+u)\lambda}
\partial_\lambda\Big[\omega(\lambda) e^{-ix\lambda}f_1(x,\lambda)
\la \vece,D^{-1}(\lambda)\vece\ra \Big]\, d\lambda \nn\\
&\quad - \frac{1}{2i(x-u)} \int_0^\infty e^{i(x-u)\lambda}
\partial_\lambda\Big[\omega(\lambda) e^{-ix\lambda}f_1(x,\lambda)
\la \vece D^{-1}(\lambda)\vece\ra \Big]\, d\lambda,\nn
\end{align}
whereas
\begin{align*}
& \int_0^\infty  \cos(u\lambda) \omega(\lambda) e^{-\mu x} e^{\mu x}f_3(x,\lambda)\la \vece',D^{-1}(\lambda)\vece\ra\, d\lambda \\
&= -u^{-1}\int_0^\infty \sin(\lambda u)
\partial_\lambda\Big[\omega(\lambda) e^{-\mu x} e^{\mu
x}f_3(x,\lambda)\la \vece',D^{-1}(\lambda)\vece\ra\Big]\,
d\lambda.
\end{align*}
Since \beeq
 \sup_{x\ge0,\,\lambda}| \partial^j_\lambda[\omega(\lambda) e^{-ix\lambda}f_1(x,\lambda)]|\le C(V),\;\;
\sup_{x\ge0,\,\lambda}| \partial^j_\lambda[\omega(\lambda)
e^{x\mu}f_3(x,\lambda)]|\le C(V), \label{eq:f1f3est} \eneq as well
as \beeq \label{eq:muxest} \sup_{x\ge0}|\partial_\lambda^j
[\omega(\lambda) e^{-\mu x}]| \le C \eneq for $j\ge0$,  it follows
via another integration by parts that
\[
|[\omega F_1(x,\cdot)D^{-1}(\lambda)\vece]^{\vee}(u)| \le C
\frac{|x|}{|x^2-u^2|} + C(u+x)^{-2} + C(u-x)^{-2} + Cu^{-2}.
\]
We will use this bound if $|\,|u|-|x|\,|>|x|$ and $|u|\ge1$. On the
other hand, if $|\,|u|-|x|\,|\le |x|$, or $|u|\le 1$,  then we simply
estimate
\[ |[\omega F_1(x,\cdot)D^{-1}(\lambda)\vece]^{\vee}(u)| \le C.\]
The conclusion is that \beeq \label{eq:xest1} \int_{\R}
|[\omega(\lambda) F_1(x,\lambda)D^{-1}(\lambda)\vece]^{\vee}(u)|\,
du \le C\la x\ra \eneq for all $x\ge0$.

If $x\le0$, then recall that
\begin{align}
& 2i\lambda F_1(x,\lambda)D^{-1}(\lambda)\vece \nn \\
&= 2i\lambda G_1(x,\lambda)A(\lambda)D^{-1}(\lambda)\vece + 2i\lambda G_2(x,\lambda)B(\lambda)D^{-1}(\lambda)\vece \nn \\
&=  G_1(x,\lambda)\bm 1 & 0\\0 & -i\frac{\lambda}{\mu} \endm\vece + G_2(x,\lambda)\vece \la 2i\lambda p B(\lambda)D^{-1}(\lambda)\vece,\vece\ra \nn \\
&\qquad + G_2(x,\lambda)\vece'\la 2i\lambda q B(\lambda)D^{-1}(\lambda)\vece,\vece'\ra \nn \\
&= f_2(-x,\lambda) + r(\lambda)f_1(-x,\lambda) + f_3(-x,\lambda) \la 2i\lambda q B(\lambda)D^{-1}(\lambda)\vece,\vece'\ra \label{eq:f2f3B}\\
&= f_1(-x,-\lambda)-f_1(-x,\lambda) +
(r(\lambda)+1)f_1(-x,\lambda) + f_3(-x,\lambda) \la 2i\lambda q
B(\lambda)D^{-1}(\lambda)\vece,\vece'\ra, \nn
\end{align}
where $r(\lambda)$ is as in Lemma~\ref{lem:scat}. Thus,
\begin{align*}
& F_1(x,\lambda)D^{-1}(\lambda)\vece \\
&= 2i\lambda G_1(x,\lambda)A(\lambda)D^{-1}(\lambda)\vece + 2i\lambda G_2(x,\lambda)B(\lambda)D^{-1}(\lambda)\vece \\
&= (2i\lambda)^{-1}[f_1(-x,-\lambda)-f_1(-x,\lambda)] +
\frac{r(\lambda)+1}{2i\lambda}f_1(-x,\lambda) + f_3(-x,\lambda)
\la q B(\lambda)D^{-1}(\lambda)\vece,\vece'\ra.
\end{align*}
By Corollary~\ref{cor:ABasymp} and Lemma~\ref{lem:Dsmooth}, $\la q
B(\lambda)D^{-1}(\lambda)\vece,\vece'\ra$ is smooth in $\lambda$.
Returning to the cosine transform (where $x=-y$ with $y\ge0$) we
conclude that
\begin{align}
&\int_0^\infty \cos(u\lambda) \omega(\lambda) F_1(-y,\lambda)D^{-1}(\lambda)\vece \, d\lambda \nn \\
& = \int_0^\infty \cos(\lambda u) \omega(\lambda)
(2i\lambda)^{-1}[f_1(y,-\lambda)-f_1(y,\lambda)] \, d\lambda \label{eq:difff1}\\
& \quad + \int_0^\infty \cos(\lambda u) e^{iy\lambda} \omega(\lambda)\frac{r(\lambda)+1}{2i\lambda}e^{-iy\lambda} f_1(y,\lambda) \, d\lambda \label{eq:r+1}\\
& \quad + \int_0^\infty \cos(\lambda u) e^{-\mu y} \la q
B(\lambda)D^{-1}(\lambda)\vece,\vece'\ra e^{\mu y}f_3(y,\lambda)
\, d\lambda \label{eq:f3lam}.
\end{align}
The easiest term to deal with is~\eqref{eq:f3lam}. Indeed,
integrating by parts in $\lambda$ twice shows that \beeq
\label{eq:newint}  \sup_{y\ge0}\Big|\int_0^\infty \cos(\lambda u)
e^{-\mu y} \la q B(\lambda)D^{-1}(\lambda)\vece,\vece'\ra e^{\mu
y}f_3(y,\lambda) \, d\lambda \Big|\le C(1+|u|)^{-2}, \eneq see
\eqref{eq:f1f3est} and \eqref{eq:muxest}. Next,
consider~\eqref{eq:r+1}. By Lemma~\ref{lem:scat}, $r(0)=-1$ and
$r(\lambda)$ is smooth. Hence, $\frac{r(\lambda)+1}{\lambda}$ is
also a smooth function. Set
$\omega_1(\lambda)=\omega(\lambda)\frac{r(\lambda)+1}{\lambda}$.
Then
\[ \eqref{eq:r+1} = \int_0^\infty \cos(u\lambda) e^{iy\lambda}\omega_1(\lambda) e^{-iy\lambda}f_1(y,\lambda)\, d\lambda. \]
By the same arguments that lead from \eqref{eq:byparts1} to
\eqref{eq:xest1} we obtain the bound \beeq \label{eq:Rint}
\int_{-\infty}^\infty \Big| \int_0^\infty \cos(u\lambda)
e^{iy\lambda}\omega_1(\lambda) e^{-iy\lambda}f_1(y,\lambda)\,
d\lambda\Big|\, du \le C\la y\ra \eneq uniformly in $y\ge0$.
Turning to~\eqref{eq:difff1}, we see that is the same as (with
$\partial_2$ being the partial derivative with respect to the
second variable of $f_1$ and ignoring constants)
\begin{align}
& \int_{-1}^1 \int_0^\infty \cos(u\lambda)\omega(\lambda) \partial_2 f_1(x,\lambda\sigma)\, d\lambda d\sigma \nn \\
&= \int_{-1}^1 \int_0^\infty \cos(u\lambda)e^{i\lambda x\sigma}\omega(\lambda) \partial_2[ e^{-i\lambda x\sigma} f_1(x,\lambda\sigma)]\, d\lambda d\sigma \label{eq:1inte}\\
& \quad + ix\int_{-1}^1 \int_0^\infty \cos(u\lambda)e^{i\lambda
x\sigma}\omega(\lambda)  e^{-i\lambda x\sigma}
f_1(x,\lambda\sigma)\, d\lambda d\sigma.  \label{eq:2inte}
\end{align}
We will focus on the second integral \eqref{eq:2inte}, since the
first one~\eqref{eq:1inte} is similar. We will integrate by parts
in $\lambda$, but only on the set $|\sigma x\pm u|\ge1$. Then
\begin{align}
& -ix\int_{-1}^1 \int_0^\infty \cos(u\lambda)e^{i\lambda
x\sigma}\omega(\lambda)  e^{-i\lambda x\sigma}
f_1(x,\lambda\sigma)\, d\lambda
 \chi_{[|\sigma x\pm u|\ge1]}\,   d\sigma \nn \\
&= \int_{-1}^1\frac{x}{2(\sigma x+u)} \omega(0)f_1(x,0)\, \chi_{[|\sigma x\pm u|\ge1]}\, d\sigma + \int_{-1}^1\frac{x}{2(\sigma x-u)} \omega(0)f_1(x,0) \, \chi_{[|\sigma x\pm u|\ge1]}\, d\sigma \nn \\
&  +\int_{-1}^1\frac{x}{2(\sigma x+u)}\int_0^\infty e^{i(\sigma x+u)\lambda} \partial_\lambda\Big[\omega(\lambda) e^{-ix\lambda\sigma}f_1(x,\lambda)\Big]\, d\lambda \, \chi_{[|\sigma x\pm u|\ge1]}\, d\sigma \label{eq:gut1}\\
& + \int_{-1}^1 \frac{x}{2(\sigma x-u)} \int_0^\infty e^{i(\sigma
x-u)\lambda} \partial_\lambda\Big[\omega(\lambda)
e^{-ix\lambda\sigma}f_1(x,\lambda)\Big]\, d\lambda\,
\chi_{[|\sigma x\pm u|\ge1]}\,  d\sigma. \label{eq:gut2}
\end{align}
The first two integrals here (which are due to the boundary
$\lambda=0$) contribute
\[
\int_{-1}^1\frac{x}{2(\sigma x+u)} \omega(0)f_1(x,0)\,
\chi_{[|\sigma x\pm u|\ge1]}\, d\sigma +
\int_{-1}^1\frac{x}{2(-\sigma x-u)} \omega(0)f_1(x,0) \,
\chi_{[|\sigma x\pm u|\ge1]}\, d\sigma =0,
\]
where we performed a change of variables $\sigma\mapsto -\sigma$
in the second one. Integrating by parts one more time in
\eqref{eq:gut1} and \eqref{eq:gut2} with respect to $\lambda$
implies
\begin{align*}
& \int_{-\infty}^\infty \Big|x\int_{-1}^1 \int_0^\infty
\cos(u\lambda)e^{i\lambda x\sigma}\omega(\lambda)  e^{-i\lambda
x\sigma} f_1(x,\lambda\sigma)\, d\lambda
 \chi_{[|\sigma x\pm u|\ge1]}\,   d\sigma\Big|\, du \nn \\
& \le C \int_{-\infty}^\infty \int_{-1}^1\frac{|x|}{(\sigma
x+u)^2} \chi_{[|\sigma x\pm u|\ge1]}\, d\sigma du + C
\int_{-\infty}^\infty \int_{-1}^1\frac{|x|}{(\sigma x-u)^2}
\chi_{[|\sigma x\pm u|\ge1]}\, d\sigma du \le C\, |x|.
\end{align*}
Finally, the cases $|\sigma x+u|\le1 $ and $|\sigma x-u|\le1$ each
contribute at most $C|x|$ to the $u$-integral. Hence
\[
 \int_{-\infty}^\infty \Big|x\int_{-1}^1 \int_0^\infty \cos(u\lambda)e^{i\lambda x\sigma}\omega(\lambda)  e^{-i\lambda x\sigma} f_1(x,\lambda\sigma)\, d\lambda d\sigma \Big|\, du
 \le C|x|.
\]
Since \eqref{eq:1inte} can be treated the same way (in fact, the
bound is $O(1)$), we obtain
\[
\int_{-\infty}^\infty\Big| \int_{-1}^1 \int_0^\infty
\cos(u\lambda)\omega(\lambda) \partial_2 f_2(x,\lambda\sigma)\,
d\lambda d\sigma \Big| \, du \le C\la x\ra.
\]
Combining this bound with \eqref{eq:newint}, \eqref{eq:Rint}, and
\eqref{eq:xest1}, we conclude that
\[
 \big\| [\omega(\lambda) F_1(x,\lambda)D^{-1}(\lambda)\vece]^{\vee} \big\|_1 \le C\la x\ra \quad \forall \;x\in\R,
\]
which in turn implies that \beeq \label{eq:uff1} \Big|
\int_0^\infty e^{it\lambda^2} \omega(\lambda) \la \phi, \sigma_3
F_1(\cdot,\lambda)D^{-1}(\lambda)\vece \ra \la
F_1(\cdot,\lambda)D^{-1}(\lambda)\vece,\psi\ra\, d\lambda \Big|
\le C\,|t|^{-\half} \|\la x\ra\phi\|_1\|\la y\ra\psi\|_1. \eneq
This is the desired estimate on \eqref{eq:F1teil}, but only for
the case when $\partial_\lambda$ falls on the factors not
involving $F_1(\cdot,\lambda)D^{-1}(\lambda)\vece$.

We now consider the case when $\partial_\lambda$ falls on
$F_1(x,\lambda)D^{-1}(\lambda)\vece$. Hence, we need to estimate
\begin{align}
&\int_0^\infty e^{it\lambda^2} \chi(\lambda^2) \lambda \la
\phi,\sigma_3\,\partial_\lambda[F_1(\cdot,\lambda)D^{-1}(\lambda)]\vece\ra
\la F_1(\cdot,\lambda)D^{-1}(\lambda)\vece,\psi\ra  \, d\lambda\nn \\
&= \int_0^\infty e^{it\lambda^2} \chi(\lambda^2) \la
\phi,\sigma_3\,\partial_\lambda[\lambda
F_1(\cdot,\lambda)D^{-1}(\lambda)]\vece\ra
\la F_1(\cdot,\lambda)D^{-1}(\lambda)\vece,\psi\ra  \, d\lambda \label{eq:F1der} \\
& - \int_0^\infty e^{it\lambda^2} \chi(\lambda^2) \la
\phi,\sigma_3\,F_1(\cdot,\lambda)D^{-1}(\lambda)\vece\ra \la
F_1(\cdot,\lambda)D^{-1}(\lambda)\vece,\psi\ra  \, d\lambda. \nn
\end{align}
The final integral we have just estimated. Hence \eqref{eq:F1der}
is the main issue. By the same reductions as before, it will
satisfy the desired bounds provided
\[ \big\|\big [\omega_1(\lambda)\partial_\lambda[\lambda F_1(x,\lambda)D^{-1}(\lambda)\vece]\,\big]^{\wedge} \big\|_1 \le  C\la x\ra \]
where $\omega_1$ is a smooth cut-off. Now
\begin{align}
& \big [\omega_1(\lambda)\partial_\lambda[\lambda F_1(x,\lambda)D^{-1}(\lambda)\vece]\,\big]^{\wedge}(u) \nn \\
& =ix \int_0^\infty \cos(u\lambda) e^{ix\lambda}\omega_1(\lambda) \lambda e^{-ix\lambda}f_1(x,\lambda)\la \vece,D^{-1}(\lambda)\vece\ra  \,d\lambda \label{eq:I}\\
&\quad + \int_0^\infty \cos(u\lambda) e^{ix\lambda}\omega_1(\lambda) \partial_\lambda[\lambda e^{-ix\lambda}f_1(x,\lambda)\la \vece,D^{-1}(\lambda)\vece\ra ] \,d\lambda \label{eq:II}\\
& \quad - \int_0^\infty \cos(u\lambda)x\mu'(\lambda)\, \omega_1(\lambda) e^{-x\mu} \lambda e^{x\mu}f_3(x,\lambda)\la \vece',D^{-1}(\lambda)\vece\ra  \,d\lambda\label{eq:III} \\
&\quad + \int_0^\infty \cos(u\lambda) \omega_1(\lambda) e^{-x\mu}
\partial_\lambda[\lambda e^{x\mu}f_3(x,\lambda)\la
\vece',D^{-1}(\lambda)\vece\ra ] \,d\lambda.\label{eq:IV}
\end{align}
We again need to distinguish $x\ge0$ from $x\le 0$. We start with
the former case. Integrating by parts in~\eqref{eq:I} leads to
\begin{align}
& ix\int_0^\infty \cos(u\lambda) e^{ix\lambda}\omega_1(\lambda)\lambda\, e^{-ix\lambda}f_1(x,\lambda)\la \vece,D^{-1}(\lambda)\vece\ra\, d\lambda \nn \\
& = -\frac{ix}{2i(x+u)}\int_0^\infty e^{i(x+u)\lambda} \partial_\lambda\Big[\omega_1(\lambda)\lambda e^{-ix\lambda}f_1(x,\lambda)\la \vece,D^{-1}(\lambda)\vece\ra\Big]\, d\lambda \label{eq:I1}\\
&\quad - \frac{ix}{2i(x-u)} \int_0^\infty e^{i(x-u)\lambda}
\partial_\lambda\Big[\omega_1(\lambda)\lambda
e^{-ix\lambda}f_1(x,\lambda)\la
\vece,D^{-1}(\lambda)\vece\ra\Big]\, d\lambda.\label{eq:I2}
\end{align}
It is important that boundary terms do not appear here (due to the
$\lambda$). On the other hand, boundary terms do arise upon
integrating~\eqref{eq:II} by parts. More precisely, for the case
of~\eqref{eq:II} we obtain
\begin{align}
& \int_0^\infty \cos(u\lambda) e^{ix\lambda}\omega_1(\lambda) \partial_\lambda[\lambda e^{-ix\lambda}f_1(x,\lambda)\la \vece,D^{-1}(\lambda)\vece\ra]\, d\lambda \nn \\
&= -\frac{1}{2i(x+u)} \omega_1(0)\partial_\lambda[ \lambda e^{-ix\lambda}f_1(x,\lambda)\la \vece,D^{-1}(\lambda)\vece\ra]\Big\vert_{\lambda=0} \label{eq:II1}\\
&  - \frac{1}{2i(x-u)} \omega_1(0) \partial_\lambda[ \lambda e^{-ix\lambda}f_1(x,\lambda)\la \vece,D^{-1}(\lambda)\vece\ra]\Big\vert_{\lambda=0} \label{eq:II2}\\
& \quad -\frac{1}{2i(x+u)}\int_0^\infty e^{i(x+u)\lambda} \partial_\lambda\Big[\omega_1(\lambda) \partial_\lambda[ \lambda e^{-ix\lambda}f_1(x,\lambda)\la \vece,D^{-1}(\lambda)\vece\ra]\Big]\, d\lambda \label{eq:II3} \\
& \quad- \frac{1}{2i(x-u)} \int_0^\infty e^{i(x-u)\lambda}
\partial_\lambda\Big[\omega_1(\lambda) \partial_\lambda[ \lambda
e^{-ix\lambda}f_1(x,\lambda)\la
\vece,D^{-1}(\lambda)\vece\ra]\Big]\, d\lambda. \label{eq:II4}
\end{align}
Integrating by parts one more time in \eqref{eq:I1} and
\eqref{eq:I2} implies
\begin{align*}
& \Big|ix\int_0^\infty \cos(u\lambda) e^{ix\lambda}\omega_1(\lambda) \lambda e^{-ix\lambda}f_1(x,\lambda)\, d\lambda \Big|  \nn \\
& \le C|x|(1+|x-u|)^{-2}+C|x|(1+|x+u|)^{-2}
\end{align*}
uniformly in $x\ge0$, whereas \eqref{eq:II1}--\eqref{eq:II4} are
treated the same way as~\eqref{eq:byparts1}. Consequently,
\eqref{eq:I} and \eqref{eq:II} each have $L^1(du)$ norm $\le C\la
x\ra$ provided~$x\ge0$. Integrating \eqref{eq:III}
and~\eqref{eq:IV} in $\lambda$ twice yields
\[ |\eqref{eq:III}|+|\eqref{eq:IV}| \le C\, (1+u^2)^{-1}\]
uniformly in $x\ge0$ since for $\ell=0,1,2$
\begin{align*}
 \sup_{x\ge0}\Big|\partial_\lambda^\ell[x\mu'(\lambda)\, \omega_1(\lambda) e^{-x\mu} \lambda e^{x\mu}f_3(x,\lambda)\la \vece',
D^{-1}(\lambda)\vece\ra]\Big| &\le C \\
 \sup_{x\ge0}\Big|\partial_\lambda^\ell[x\mu'(\lambda)\, \omega_1(\lambda) e^{-x\mu} \lambda e^{x\mu}f_3(x,\lambda)\la \vece',
D^{-1}(\lambda)\vece\ra]\Big| &\le C.
\end{align*}
Hence, \eqref{eq:III} and \eqref{eq:IV} each have $L^1(du)$ norm
$\le C$ provided~$x\ge0$.

To deal with $x\le0$, we use \eqref{eq:f2f3B}:
\[ 2i\lambda F_1(x,\lambda)D^{-1}(\lambda)\vece
= f_2(-x,\lambda) + r(\lambda)f_1(-x,\lambda) + f_3(-x,\lambda)
\la 2i\lambda q B(\lambda)D^{-1}(\lambda)\vece,\vece'\ra.
\]
This implies that
\begin{align}
& 2i\big [\omega_1(\lambda)\partial_\lambda[\lambda F_1(x,\lambda)D^{-1}(\lambda)\vece]\,\big]^{\wedge}(u) \nn \\
& = \int_0^\infty \cos(\lambda u) \omega_1(\lambda) \partial_\lambda f_2(-x,\lambda)\, d\lambda  + \int_0^\infty \cos(\lambda u) \omega_1(\lambda) \partial_\lambda [r(\lambda) f_1(-x,\lambda)]\, d\lambda \nn \\
& \quad + \int_0^\infty \cos(\lambda u) \omega_1(\lambda)
\partial_\lambda [f_3(-x,\lambda) \la 2i\lambda q
B(\lambda)D^{-1}(\lambda)\vece,\vece'\ra ]\, d\lambda \nn,
\end{align}
which further simplifies to
\begin{align}
&  = \int_0^\infty \cos(\lambda u)e^{ix\lambda} \omega_1(\lambda)
\partial_\lambda[ e^{-ix\lambda}f_2(-x,\lambda)]\, d\lambda
+ \int_0^\infty \cos(\lambda u) e^{-ix\lambda}\omega_1(\lambda) \partial_\lambda [r(\lambda)e^{ix\lambda} f_1(-x,\lambda)]\, d\lambda \label{eq:xohne} \\
& \quad +ix \int_0^\infty \cos(\lambda u)e^{ix\lambda}
\omega_1(\lambda)  e^{ix\lambda}f_2(-x,\lambda)\, d\lambda
-ix\int_0^\infty \cos(\lambda u) e^{-ix\lambda}\omega_1(\lambda) r(\lambda)e^{ix\lambda} f_1(-x,\lambda)\, d\lambda \label{eq:xmit} \\
& \quad + \int_0^\infty \cos(\lambda u)x\mu'(\lambda)e^{x\mu} \omega_1(\lambda) e^{-x\mu}f_3(-x,\lambda) \la 2i\lambda q B(\lambda)D^{-1}(\lambda)\vece,\vece'\ra\, d\lambda \label{eq:f3wieder} \\
& \quad  +\int_0^\infty \cos(\lambda u)e^{x\mu} \omega_1(\lambda)
\partial_\lambda [e^{-x\mu}f_3(-x,\lambda) \la 2i\lambda q
B(\lambda)D^{-1}(\lambda)\vece,\vece'\ra]\, d\lambda.
\label{eq:f3wieder'}
\end{align}
The two integrals in~\eqref{eq:xohne} (which are not preceded by
$ix$) are integrated by parts in the same way as
\eqref{eq:II1}--\eqref{eq:II4}. Their contribution to the
$L^1(du)$ norm is at most $C\la x\ra$ with a constant uniform
in~$x\le0$. The integrals in~\eqref{eq:xmit}, which are preceded
by~$ix$ are also integrated by parts, but we need to check in this
case that the boundary terms $\lambda=0$ cancel each other.
However, these boundary terms are
\begin{align*}
& \frac{x}{2(u-x)}\omega_1(0)r(0)f_1(-x,0)-\frac{x}{2(u+x)}\omega_1(0)r(0)f_1(-x,0) \\
& -\frac{x}{2(u+x)}\omega_1(0)f_2(-x,0)
-\frac{x}{2(x-u)}\omega_1(0)f_2(-x,0)=0,
\end{align*}
since $r(0)=-1$ and $f_1(-x,0)=f_2(-x,0)$. Hence, \eqref{eq:xmit}
can be treated as~\eqref{eq:I1}, \eqref{eq:I2}. Finally,
\eqref{eq:f3wieder} and~\eqref{eq:f3wieder'} are $\le
C(1+u^2)^{-1}$ uniformly in $x\le0$, see~\eqref{eq:III},
\eqref{eq:IV}, and we are done.
\end{proof}

Interpolating this estimate with the unweighted $L^1(\R)\to
L^\infty(\R)$ as well as the $L^2$ bound yields the following:

\begin{lemma}
\label{lem:theta} Let $\Hil$ be admissible, see
Definition~\ref{def:admH}. Then for all $0\le\theta\le1$, all
$1\le p\le2$, and all $t\ne0$,
\[ \|\la x\ra^{-\theta(\frac{1}{p}-\frac{1}{p'})}e^{it\Hil}P_s f\|_{p'}\le C|t|^{-(\half+\theta)(\frac{1}{p}-\frac{1}{p'})}\|\la x\ra^{\theta(\frac{1}{p}-\frac{1}{p'})} f\|_p.\]
Here $C$ is some absolute constant.
\end{lemma}
\begin{proof}
Interpolating between Propositions~\ref{prop:disp1}
and~\ref{prop:disp2} yields
\[ \|\la x\ra^{-\theta}e^{it\Hil}P_s f\|_{\infty}\le C|t|^{-\half-\theta}\|\la x\ra^{\theta} f\|_1. \]
The lemma now follows from a further interpolation with
Lemma~\ref{lem:L2stable}.
\end{proof}

Just as in the unweighted case, derivatives can be introduced here
as well. We restrict ourselves to two derivatives, although more
are possible.

\begin{cor}
\label{cor:weightp} Let $\Hil$ be admissible, see
Definition~\ref{def:admH}. Then for all $0\le\theta\le1$, all $1<
p\le2$, and all $t\ne0$,
\begin{equation}
 \| \la x\ra^{-\theta(\frac{1}{p}-\frac{1}{p'})}e^{it\Hil}P_s f\|_{W^{k,p'}(\R)}\le C|t|^{-(\half+\theta)(\frac{1}{p}-\frac{1}{p'})}\sum_{j=0}^k \|\la x\ra^{\theta(\frac{1}{p}-\frac{1}{p'})}
\partial_x^j f\|_{L^{p}(\R)}
\label{eq:sskm}
\end{equation}
for all $0\le k\le2$. Alternatively, \beeq \|\la
x\ra^{-\theta(\frac{1}{p}-\frac{1}{p'})}\partial_x^k\,
e^{it\Hil}P_s f\|_{L^{p'}(\R)}\le
C|t|^{-(\half+\theta)(\frac{1}{p}-\frac{1}{p'})}\sum_{j=0}^k \|\la
x\ra^{\theta(\frac{1}{p}-\frac{1}{p'})}
\partial_x^j f\|_{L^{p}(\R)}
\label{eq:sskm2} \eneq for all $0\le k\le2$. Here $C$ is a
constant that depends on $p$.
\end{cor}
\begin{proof}
The case $k=0$ is Lemma~\ref{lem:theta}. Let $a>0$ be large. As in
the proof of Corollary~\ref{cor:inter}, for $p'<\infty$,
\begin{align*}
&\Big\| \la x\ra^{-\theta(\frac{1}{p}-\frac{1}{p'})} e^{it\Hil}P_s f \Big\|_{W^{2,p'}(\R)}\\
 &\le C \Big\| (\Hil-ia) \la x\ra^{-\theta(\frac{1}{p}-\frac{1}{p'})} e^{it\Hil}P_s f \Big\|_{p'} \\
&\le C \Big\| \la x\ra^{-\theta(\frac{1}{p}-\frac{1}{p'})} e^{it\Hil}P_s  (\Hil-ia) f \Big\|_{p'} + \Big\|[\la x\ra^{-\theta(\frac{1}{p}-\frac{1}{p'})},\Hil]e^{it\Hil}P_s f\Big\|_{p'} \\
&\le C |t|^{-(\frac12+\theta)(\frac{1}{p}-\frac{1}{p'})} \Big\|\la
x\ra^{\theta(\frac{1}{p}-\frac{1}{p'})} (\Hil-ia)f\Big\|_p +
\Big\|[\la
x\ra^{-\theta(\frac{1}{p}-\frac{1}{p'})},\Hil]e^{it\Hil}P_s f
\Big\|_{p'}.
\end{align*}
The first term on the right can be controlled in terms of ordinary
derivatives, whereas the second only involves derivatives of order
zero and one. Therefore, this second term can be controlled by
means of interpolation. Hence, \eqref{eq:sskm} holds. The second
inequality follows from the first by induction in $k$.
\end{proof}

\section{The spectrum of the linearized NLS}
\label{sec:spectrum}

The linearization of the NLS
\[ i\partial_t\psi + \partial_{xx}\psi = -|\psi|^{2\sigma}\psi \]
around the ground state
\[ \phi(x) = (\sigma+1)^{\frac{1}{2\sigma}}\cosh^{-\frac{1}{\sigma}}(\sigma x) \]
of $ -\partial_{xx}\phi+\phi - \phi^{2\sigma+1}=0 $  leads to the
operator \beeq \label{eq:Hil1}
 \Hil = \bm 0 & iL_- \\ -iL_+ & 0 \endm
\eneq where
\begin{align*}
 L_- &= -\partial_{xx} + 1 - \phi^{2\sigma} = -\partial_{xx} + 1 - (\sigma+1)\cosh^{-2}(\sigma x) \\
 L_+ &= -\partial_{xx} + 1 - (2\sigma+1)\phi^{2\sigma} = -\partial_{xx} + 1 - (2\sigma+1)(\sigma+1)\cosh^{-2}(\sigma x).
\end{align*}
Equivalently, $\Hil$ can be written in the form \beeq
\label{eq:Hil2}
 \Hil= \bm -\partial_{xx}+1 - (\sigma+1)\phi^{2\sigma} & -\sigma \phi^{2\sigma} \\
 \sigma \phi^{2\sigma} & \partial_{xx} -1 + (\sigma+1)\phi^{2\sigma}\endm.
\eneq

\begin{lemma}
\label{lem:flug} If $\sigma>1$, then $L_+$ and $L_-$ have the
following properties: they have no eigenvalues in the interval
$(0,1]$ and for both $L_+$ and $L_-$ the threshold $1$ is not a
resonance.
\end{lemma}
\begin{proof}
A ``resonance'' of $L_\mp$ at energy one means  one of the
following equivalent things:
\begin{itemize}
\item There is a solution $f\in L^\infty\setminus L^2$ of $L_\mp
f=f$ \item The Wronskian of the two Jost solutions with energy one
are linearly dependent \item The transmission coefficient at
energy one does not vanish
\end{itemize}
The lemma can easily be deduced from Fl\"ugge~\cite{Flug}, for
example -- see Problem~39 on page~94. It is shown there that the
Hamiltonian $H=-\frac{d^2}{dx^2}-\alpha^2
\frac{\lambda(\lambda-1)}{\cosh^2(\alpha x)}$ with $\lambda>1$ has
a zero-energy resonance (which is the same as $T(0)=0$ for the
transmission coefficient $T(E)$) iff $\lambda$ is an integer. In
the case of $L_-$, $\alpha=\sigma$ and $\lambda=
\half+\half\sqrt{1+4\frac{\sigma+1}{\sigma^2}}$. If
$\lambda=n\ge2$, then $\frac{\sigma+1}{\sigma^2}=n(n-1)\ge2$,
which implies that $\sigma\le1$. Moreover, it is shown there that
the number of negative bound states of $H$ is the largest integer
$\le \lambda$. Here $\lambda<2$ iff $\sigma>1$ so that there is
exactly one  bound state of $L_-$ below energy one iff $\sigma>1$.
For $L_+$, we have
$\lambda=\half+\half\sqrt{1+4\frac{(2\sigma+1)(\sigma+1)}{\sigma^2}}<3$
iff $\sigma>1$. This implies that $L_+$ has exactly two bound
states below energy one, as desired. It same way, it can be
checked that $L_+$ does not have a resonance at $1$ if $\sigma>1$.
\end{proof}

The main result of this section is the following proposition.

\begin{prop}
\label{prop:spectrum} For any $\sigma\ge2$ the operator $\Hil$ on
$L^2(\R)\times L^2(\R)$ with domain $W^{2,2}(\R)\times
W^{2,2}(\R)$ satisfies:
\begin{itemize}
\item The spectrum of $\Hil$, denoted by $\spec(\Hil)$,  is
contained in $\R\cup i\R$ \item The essential spectrum equals
$(-\infty,-1]\cup[1,\infty)$. There are no embedded eigenvalues
     in the essential spectrum.
\item The only real eigenvalue in $[-1,1]$ is zero. The geometric
multiplicity of the zero eigenvalue is two, and its algebraic
multiplicity equals four or six, depending on whether $\sigma>2$
or $\sigma=2$. \item If $\sigma>2$, then there is a unique pair of
imaginary eigenvalues $\pm i\gamma$, $\gamma>0$, which are both
simple. \item The edges $\pm1$ are not resonances\footnote{see
Definition~\ref{def:res}}.
\end{itemize}
In particular, $\Hil$ is admissible\footnote{see
Definition~\ref{def:admH}} for all $\sigma\ge2$.
\end{prop}
\begin{proof}
We will rely on the techniques from~\cite{BP1} and~\cite{Pe2}. The
latter paper only deals with the critical case $\sigma=2$, which
means that we need to adapt some of Perelman's arguments to the
supercritical case $\sigma>2$.

The statement about the essential spectrum follows from Weyl's
criterion via the symmetric resolvent identity, see~\cite{Grill}.
It will be convenient to introduce the ground state
$\phi(x,\alpha)$ (which also depends on $\sigma$) of
\[ -\partial_{xx}\phi+\alpha^2\phi - \phi^{2\sigma+1}=0.\]
Then $\phi(x,\alpha)=\alpha^{\frac1\sigma}\phi(\alpha x,1)$.
Similarly,
\begin{align*}
 L_- &= -\partial_{xx} + \alpha^2 - \phi^{2\sigma} \\
 L_+ &= -\partial_{xx} + \alpha^2 - (2\sigma+1)\phi^{2\sigma}
\end{align*}
Clearly, $L_-\phi=0$, $L_+(\partial_\alpha \phi)=-2\alpha\phi$,
and $L_+(\phi')=0$. Since $\phi>0$, it follows that $L_-\ge0$.
Also, $\ker(L_+)={\rm span}\{\phi'\}$ (by the simplicity of
eigenvalues in the one-dimensional case), and $L_+$ has a unique
negative eigenvalue $E_0<0$. If $\Hil f=zf$ with
$f=\binom{f_1}{f_2}\in L^2 \setminus\{0\}$ and $z\ne0$, then
$L_-L_+f_1=z^2 f_1$ and thus $\sqrt{L_-}L_+\sqrt{L_-}g_1=z^2 g_1$
where $g_1=L_-^{-\half}f_1$ (note that $f_1\perp \phi$). It
follows that $z^2\in\R$. Moreover, this argument can be refined
(see page 1137 in~\cite{BP1}) to show that if $z\ne0$, then the
geometric and algebraic multiplicities of the eigenvalue $z$ are
equal. If $\Hil$ has imaginary eigenvalues,  then some $z^2<0$
would need to be an eigenvalue of $\sqrt{L_-}L_+\sqrt{L_-}$. This
implies that
\[ \lambda_0:=\min_{\|f\|_2=1,\,f\perp\phi}\la L_+\phi,\phi\ra <0.\]
If $g$ is a minimizer here, then by Lagrange multipliers
\[ (L_+-\lambda_0)g = c\phi,\qquad 0=\la g,\phi\ra = c\la (L_+-\lambda_0)^{-1}\phi,\phi\ra.\]
Since $c\ne0$ (otherwise, $g>0$ is the ground state of $L_+$ which
contradicts $g\perp\phi$) one has $h(\lambda_0)=0$ where
\[ h(\lambda):=\la (L_+-\lambda)^{-1}\phi,\phi\ra.\]
But $h(\lambda)$ is strictly increasing which requires
\[ 0< h(0) = \la L_+^{-1}\phi,\phi\ra = -\frac{1}{2\alpha}\la \partial_\alpha\phi,\phi\ra.\]
Since
\[ \partial_\alpha \|\phi\|_2^2 = 2\la\partial_\alpha\phi,\phi\ra < 0 \text{\ \ iff\ \ }\sigma> 2,\]
it follows that $\spec(\Hil)\subset\R$ if $\sigma\le 2$ (we will
see later that this is if and only if).

As far as the zero eigenvalue is concerned, note that (use
\eqref{eq:Hil1})
\[
\Hil\binom{0}{\phi}=0, \; \Hil\binom{\phi'}{0}=0,\;
\Hil\binom{0}{x\phi}=\binom{-2i\phi'}{0}, \;
\Hil\binom{\partial_\alpha\phi}{0}=\binom{0}{2i\alpha\phi}.
\]
It is clear that this describes $\ker(\Hil)$ and $\ker(\Hil^2)$
completely. Moreover, it is easy to check that
$\ker(\Hil^k)=\ker(\Hil^2)$ if $k\ge3$ provided $\sigma\ne2$.
Indeed, it is enough to check this for $k=3$ and that case is
settled by writing out the third power of~\eqref{eq:Hil1}
explicitly, see~\cite{Wei1} page~485. In the critical case
$\sigma=2$, one has $\partial_\alpha\phi\perp\phi$. Hence
$\partial_\alpha\phi\in\Ran(L_-)$. In fact, direct differentiation
shows that $L_-(x^2\phi)=-4\alpha\partial_\alpha\phi$. Finally,
since $x^2\phi\perp \phi'$, it follows that $L_+\rho=x^2\phi$ for
some $\rho$. In summary,
\[ \Hil\binom{0}{x^2\phi} = \binom{-4i\alpha\partial_\alpha\phi}{0}, \quad \Hil\binom{\rho}{0}=\binom{0}{-ix^2\phi}\]
for the case $\sigma=2$. This shows that $\ker(\Hil^2)\subsetneq
\ker(\Hil^3)\subsetneq\ker(\Hil^4)$ and the codimensions in each
case equal one. Again, it is easy to check that in this case
$\ker(\Hil^k)=\ker(\Hil^4)$ for $k\ge4$.

Let us now assume that the only real eigenvalue of $\Hil$ in
$[-1,1]$ is zero for all $\sigma\ge2$ (this is false if
$\sigma<2$). Then by the continuous dependence of the Riesz
projection onto the discrete spectrum of $\Hil$ on the
power~$\sigma$, it follows that the rank of this projection is
constant equal to six. Hence, if $\sigma>2$, then there are two
imaginary eigenvalues counted with multiplicity. By the
commutation properties relative to the Pauli
matrices~\eqref{eq:pauli}, it follows that this has to be a pair
$\pm i\gamma$, $\gamma>0$ of simple eigenvalues.

\noindent It remains to show the following three properties:

i) The only real eigenvalue in $[-1,1]$ is zero.

ii) The edges $\pm1$ are not resonances.

iii) There are no embedded eigenvalues
     in the essential spectrum.

\medskip
Perelman proved these statements for $\sigma=2$, cf.~Sections
2.1.2 and 2.1.3 of~\cite{Pe2}. Hence, it will suffice to consider
the case $\sigma>2$ which can be dealt with by adapting Perelman's
arguments. We will rely on Lemma~\ref{lem:flug} for that purpose.

Suppose i) fails. Then $\Hil(\alpha)^2$ has an eigenvalue $E\in
(0,\alpha^4]$. For simplicity and without loss of generality, let
us choose $\alpha=1$. Then there is $\psi\in L^2(\R)$, $\psi\ne0$,
such that
\[ L_{-}L_+\psi = E\psi \]
with $0<E\le 1$. Clearly, $\psi\perp\phi$ and $\psi\in
H^4_{loc}(\R)$ by elliptic regularity. Define $A:= PL_+P$ where
$P$ is the projection orthogonal to $\phi$. Since $\la
\phi,\partial_\alpha\phi\ra\ne0$, we conclude that
\[
\ker(A) = {\rm span} \{\phi',\,\phi\}.
\]
Moreover, let $E_0<0$ be the unique negative eigenvalue of $L_+$.
Then consider (as before) the function
\[ h(\lambda):= \la (L_+-\lambda)^{-1}\phi,\phi\ra \]
which is differentiable on the  interval $(E_0, 1)$  due to the
orthogonality of $\phi$ to the kernel of $L_+$. Moreover,
\[ h'(\lambda)=\la (L_+-\lambda)^{-2}\phi,\phi\ra>0,\; h(0)= -\half \la \phi,\partial_\alpha\phi\ra>0. \]
The final inequality here is due to the supercritical nature of
our problem. Since also $h(\lambda)\to -\infty$ as $\lambda\to
E_0$, it follows that $h(\lambda_1)=0$ for some $E_0<\lambda_1<0$.
Moreover, this is the only zero of $h(\lambda)$ with
$E_0<\lambda<1$. If we set
\[ \etatil:= (L_+-\lambda_1)^{-1} \phi,\]
then
\[ A\etatil=\lambda_1\etatil, \qquad \la \etatil,\phi\ra=0.\]
Conversely, if
\[ Af=\lambda f\]
for some $E_0<\lambda<1$, $\lambda\ne0$,  and $f\in L^2(\R)$, then
$f\perp \phi$ and
\[ (PL_+P-\lambda)f=(A-\lambda)f=0.\]
Since also
\[ E_0\la f,f\ra\le \la L_+f,f\ra = \lambda\la f,f\ra \]
it follows that $\lambda\ge E_0$. If $\lambda=E_0$, then $f$ would
necessarily have to be the ground state of $L_+$ and thus of
definite sign. But then $\la f,\phi\ra\ne0$, which is impossible.
Hence $E_0<\lambda<1$. But then $h(\lambda)=0$ implies that
$\lambda=\lambda_1$ is unique. In summary, $A$ has eigenvalues
$\lambda_1$ and $0$ in $(-\infty,1)$, with $\lambda_1$
 being a simple eigenvalue and $0$ being of multiplicity two.
Now define
\[ \calF:= {\rm span} \{ \psi,\,\etatil,\,\phi',\,\phi\}.
 \]
We claim that \beeq \label{eq:dimclaim}
 \dim(\calF)=4.
\eneq Since $\phi$ is perpendicular to the other functions, it
suffices to show that
\[ c_1\psi+c_2\etatil + c_3\phi'=0 \]
can only be the trivial linear combination. Apply $L_+$. Then
\[ c_1 L_+\psi+c_2 L_+\etatil =0\]
and therefore
\begin{align*}
c_1 \la L_+\psi,\psi\ra + c_2 \la L_+\etatil,\psi\ra &=0 \\
c_1 \la L_+\psi,\etatil\ra + c_2 \la L_+\etatil,\etatil\ra &=0.
\end{align*}
This is the same as
\begin{align*}
c_1 E\la L_{-}^{-1} \psi,\psi\ra + c_2 \lambda_1\la \etatil,\psi\ra &=0 \\
c_1 \lambda_1\la \psi,\etatil\ra + c_2 \lambda_1\la
\etatil,\etatil\ra &=0.
\end{align*}
The determinant of this system is
\[ E\lambda_1\la L_{-}^{-1} \psi,\psi\ra \la \etatil,\etatil\ra - \lambda_1^2 |\la \etatil,\psi\ra|^2 <0.\]
Hence $c_1=c_2=0$ and therefore also $c_3=0$, as desired. Thus,
\eqref{eq:dimclaim} holds. Finally, we claim that \beeq
\label{eq:qform}
 \sup_{\|f\|_2=1,\,f\in\calF} \la Af,f\ra <1.
\eneq If this is true, then by the min-max principle and
\eqref{eq:dimclaim} we would obtain that the number of eigenvalues
of $A$ in the interval $(-\infty,1)$ (counted with multiplicity)
would have to be at least four. On the other, we showed before
that this number is exactly three, leading to a contradiction.
Hence, we need to verify~\eqref{eq:qform}. Since $\la
PL_{-}^{-1}Pf,f\ra <\la f,f\ra$ for all $f\ne0$, and since $E\le
1$ by assumption,
 this in turn follows from the stronger claim that
\beeq \label{eq:qform'}
 \la Af,f\ra \le E\la PL_{-}^{-1}P f,f\ra
\eneq for all $f=a\psi+b\phi+c\phi' + d\etatil$. Clearly, we can
take $b=0$. Then the left-hand side of \eqref{eq:qform'} is equal
to
\begin{align}
&\la L_+(a\psi),a\psi+c\phi'+d\etatil\ra +\la
L_+(c\phi'+d\etatil),
a\psi+c\phi'+d\etatil\ra \nn\\
&= E\la L_-^{-1}(a\psi),a\psi+c\phi'+d\etatil\ra +E\la
c\phi'+d\etatil,
L_{-}^{-1}(a\psi)\ra + \la L_+(d\etatil),d\etatil\ra \nn\\
&= E\la L_-^{-1}(a\psi),a\psi+c\phi'+d\etatil\ra +E\la
c\phi'+d\etatil, L_{-}^{-1}(a\psi)\ra + \lambda_1\|d\etatil\|_2^2,
\label{eq:lhs1}
\end{align}
whereas the right-hand side of \eqref{eq:qform'} is \beeq = E\la
L_-^{-1}(a\psi),a\psi+c\phi'+d\etatil\ra +E\la c\phi'+d\etatil,
L_{-}^{-1}(a\psi)\ra + E\la
L_{-}^{-1}(c\phi'+d\etatil),c\phi'+d\etatil\ra. \label{eq:rhs1}
\eneq Since
\[ \lambda_1\|d\etatil\|_2^2\le 0,\quad E\la L_{-}^{-1}(c\phi'+d\etatil),c\phi'+d\etatil\ra
\ge0,
\]
we see that \eqref{eq:rhs1} does indeed dominate \eqref{eq:lhs1},
and \eqref{eq:qform'} follows.

Next, we turn to the resonances. Suppose $\Hil f=f$ where $f\in
L^\infty\setminus L^2(\R)$. By Lemma~\ref{lem:0sol}, $f = C_\pm
\binom{0}{1}+ O(e^{\mp \gamma x})$ as $x\to\pm\infty$ where
$C_+\ne0$ and $C_-\ne0$. Hence there exists $\psi\in
L^\infty\setminus L^2$ (the first component of~$f$) so that $\psi
= C_\pm+ O(e^{\mp \gamma x})$ as $x\to\pm\infty$ and such that
$L_-L_+\psi=\psi$.
This asymptotic expansion can also be differentiated. Pick a
smooth cut-off $\chi\ge0$ which is constant $=1$ around zero, and
compactly supported. Define for any $0<\eps<1$
\[ \psi^\eps := \psi \chi(\eps\cdot)+\mu(\eps)\phi, \quad \mu(\eps):= - \frac{\la \psi \chi(\eps\cdot),\phi\ra}{\la \phi,\phi\ra}.
\]
Clearly, $\la \psi^\eps,\phi\ra=0$ and $|\mu(\eps)|=o(1)$ as
$\eps\to0$ (in fact, like $e^{-C/\eps}$). It follows that
\[ \|\psi^\eps\|_2^2 = M_0(\eps)+o(1),\quad M_0(\eps):= \int_{\R}|\psi(x)|^2\chi(\eps x)^2\, dx\]
with $M_0(\eps)\to\infty$ as $\eps\to0$. We now claim that \beeq
\label{eq:L+claim}
 \la  L_+\psi^\eps,\psi^\eps \ra = \|\psi^\eps\|_2^2 + \la (L_+-1)\psi,\psi\ra +o(1)
\eneq as $\eps\to0$. From the asymptotic behavior of $\psi$ it is
clear that $ M_1:=\la (L_+-1)\psi,\psi\ra$ is a finite expression.
Write $L_{-}=-\partial_x^2+1+V_1$ and $L_+=-\partial_x^2+1+V_2$,
with Schwartz functions $V_1,V_2$ (they are of course explicitly
given in terms of $\phi$, but we are not going to use that now).
We start from the evident expression \beeq \nonumber
 \la  L_+\psi^\eps,\psi^\eps \ra = \|\psi^\eps\|_2^2 + \la (L_+-1)\psi^\eps,\psi^\eps\ra
= \|\psi^\eps\|_2^2 + \la
(-\partial_x^2+V_2)\psi^\eps,\psi^\eps\ra. \eneq By the rapid
decay of $V_2$,
\[ \la (-\partial_{x}^2+V_2)\psi^\eps,\psi^\eps\ra = \int_{\R} |\partial\psi^\eps(x)|^2\, dx
+ \int_{\R} V_2(x)|\psi(x)|^2\, dx + o(1).
\]
Since $\psi'\in L^2$, we calculate further that
\begin{align*}
\int_{\R} |\partial \psi^\eps(x)|^2\, dx &= \int_{\R} \Big| \psi'(x)\chi(\eps x)+ \eps\psi(x)\chi'(\eps x)\Big|^2\, dx \\
&= \int_{\R} |\psi'(x)|^2\, dx + \int_{\R} |\psi'(x)|^2(\chi(\eps x)^2-1)\, dx \\
& \quad + 2\eps\int_{\R} \psi(x)\chi(\eps
x)\psi'(x)\cdot\chi'(\eps x)\, dx
+ \eps^2\int_{\R}\psi(x)^2|\chi'(\eps x)|^2\, dx \\
&=  \int_{\R} |\psi'(x)|^2\, dx + o(1).
\end{align*}
To pass to the last line, use the asymptotics of $\psi$ and
$\psi'$. By these asymptotics, $f:=(L_+-1)\psi$ is rapidly
decaying. Hence, $\la (L_+-1)\psi,\psi\ra=\la f,\psi\ra$ is
well-defined as a usual scalar product. Moreover, one has
\[ L_-f=-(L_--1)\psi \text{\ \ or\ \
}\psi=-(L_--1)^{-1}L_-f=-f-(L_--1)^{-1}f.\] We conclude that \beeq
\label{eq:fpsi} \la f+\psi,f\ra = -\la (L_--1)^{-1}f,f\ra <0,
\eneq where the final inequality follows from $L_-\ge1$ on
$\{\phi\}^\perp$, as well from the property that $L_-$ has neither
an eigenvalue nor a resonance at the threshold~$\alpha^2=1$. The
inequality~\eqref{eq:fpsi} will play a crucial role in estimating
a quadratic form as in the previous paragraph dealing with the
absence of eigenvalues. To see this, let
\[ \calF_\eps:={\rm span}\{\psi^\eps,\phi',\etatil,\phi\}.
\]
As before, one shows that $\dim\calF_\eps=4$, as least if $\eps>0$
is sufficiently small (use that $\la
L_+\psi^\eps,\psi^\eps\ra\to\infty$ as $\eps\to0$). It remains to
show that for small $\eps>0$ \beeq \label{eq:tschud}
\max_{f\in\calF_\eps} \frac{\la PL_+P f,f\ra}{\la f,f\ra}<1 \eneq
where $P$ is the projection orthogonal to $\phi$. If so, then this
would imply that $A=PL_+P$ has at least four eigenvalues (with
multiplicity) in $(-\infty,1)$. However, we have shown above that
there are exactly three such eigenvalues. To
prove~\eqref{eq:tschud}, it suffices to consider the case $f\perp
\phi$.  Compute
\begin{align*}
& \frac{\la
L_+(a\psi^\eps+c\phi'+d\etatil),a\psi^\eps+c\phi'+d\etatil\ra}
{\|a\psi^\eps+c\phi'+d\etatil\|^2} \\
& =\frac{|a|^2(\|\psi^\eps\|_2^2 + M_1+o(1))+2\Re\lambda_1\la
a\psi^\eps,d\etatil\ra +
\lambda_1\|d\etatil\|_2^2}{|a|^2\|\psi^\eps\|_2^2 + 2\Re\la
a\psi^\eps, c\phi'+d\etatil\ra +
\|c\phi'+d\etatil\|_2^2} \\
& \le \max_{x\in\Compl^3} \frac{|x_1|^2(1+\delta^2
M_1+o(\delta^2))+2\delta \lambda_1\Re\la x_1\psi^\eps,x_3 e_2 \ra
+ \lambda_1|x_3|^2}{|x_1|^2+2\delta\Re\la
x_1\psi^\eps,x_2e_1+x_3e_2\ra + \|x_2e_1+x_3e_2\|_2^2}
\end{align*}
where we have set $\delta^2:=\|\psi^\eps\|_2^{-2}$ and
\[ e_1=\frac{\phi'}{\|\phi'\|_2},
 \quad e_2=\frac{\etatil}{\|\etatil\|_2}.\]
Note that $\etatil$ is an even function, since it is given by
$(L_+-\lambda_1)^{-1}\phi$ and both $\phi$ and the kernel of
$(L_+-1)^{-1}$ are even. Hence $e_1\perp e_2$. Set
\[ b_j^\eps:=  \la \psi^\eps,e_j\ra \text{\ \ for\ \ }1\le j\le 2.\]
Then $b_j^\eps\to b_j^0:= \la \psi,e_j\ra$ as $\eps\to0$ by the
exponential decay of the $e_j$. Let $B^\eps,C^\eps$ (which depend
on $\eps$) be $3\times 3$ Hermitian matrices so that
\[ C^\eps_{11}:=1+\delta^2 M_1+o(\delta^2), \;C^\eps_{13}=C^\eps_{31}:=\lambda_1
\delta b^\eps_2,\; C^\eps_{33}:= \lambda_1\] and $C^\eps_{ij}=0$
else,
\[ B^\eps_{1j}=B^\eps_{j1}:=\delta b^\eps_{j-1} \text{\ \ for\ \ } 2\le j\le 3\]
and $B^\eps_{ij}=0$ else.  In view of the preceding,
\[
\max_{f\in\calF_\eps} \frac{\la PL_+P f,f\ra}{\la f,f\ra} \le
\max_{x\in\Compl^3} \frac{\la Cx,x\ra}{\la (I+B)x,x\ra}.
\]
Clearly, the right-hand side equals the largest eigenvalue of the
Hermitian matrix \[ (I+B^\eps)^{-\half}C^\eps(I+B^\eps)^{-\half}=
C-\half(BC+CB)+\frac38(B^2C+CB^2)+\frac14BCB+O(\delta^3),
\]
where we have dropped the $\eps$ in the notation on the right-hand
side. With some patience one can check that the right-hand side
equals the matrix $D$ which is given by (dropping $\eps$ from the
notation)
\[\left[
\begin{array}{lllll} 1+\delta^2 M_2 & -\frac{\delta}{2} b_1 &  \frac{\delta}{2}(\lambda_1-1)b_2\\
-\frac{\delta}{2} b_1 & \frac{\delta^2}{4} b_1^2 &
\frac{\delta^2}{4}(1-\half\lambda_1)b_1b_2 \\
\frac{\delta}{2}(\lambda_1-1)b_2 &
\frac{\delta^2}{4}(1-\half\lambda_1)b_1b_2 & \lambda_1
+\frac{\delta^2}{4}(1-\lambda_1) b_2^2
\end{array}\right] + o(\delta^2)
\]
where $M_2:=M_1-\frac34\lambda_1 b_2^2+\frac34(b_1^2+b_2^2)$. When
$\delta=0$, this matrix has simple eigenvalues $1$, $0$,
$\lambda_1<0$. When $\delta\ne0$ but very small, the largest
eigenvalue will be close to one, of the form $1+x$ with $x$ small.
We need to see that $x<0$. Collecting powers of $x$ in $
\det(D-(1+x)I)$ we arrive at the condition \begin{align*}
 (1-\lambda_1)x &= \delta^2
[M_1(1-\lambda_1)+b_1^2(1-\lambda_1)+b_2^2(1-\lambda_1)^2]+o(\delta^2)\\
&=\delta^2
[M_1(1-\lambda_1)+(b_1^0)^2(1-\lambda_1)+(b_2^0)^2(1-\lambda_1)^2]+o(\delta^2).
\end{align*}
We have
\[ b_1^0 = \la \psi,e_1\ra = -\la (L_+-1)\psi,e_1\ra = -\la
f,e_1\ra.
\]
On the other hand,
\[ b_2^0 = \la \psi,e_2\ra = -\la f,e_2\ra + \la\psi,L_+e_2\ra =
-\la f,e_2\ra + \lambda_1 b_2^0\] and thus,
\[ b_2^0 = -(1-\lambda_1)^{-1}\la f,e_2\ra.\]
Since $\lambda_1<0$ in the supercritical case, we obtain that
\begin{align*}
 (1-\lambda_1)x &\le (1-\lambda_1)\delta^2[M_1+\sum_{j=1}^2
\la f,e_j\ra^2] + o(\delta^2)  \le (1-\lambda_1)\delta^2[M_1+
\la f,f\ra ] + o(\delta^2) \\
& = (1-\lambda_1)\delta^2 \la f+\psi,f\ra + o(\delta^2) =
-(1-\lambda_1)\delta^2\la (L_--1)^{-1}f,f\ra + o(\delta^2)
\end{align*}
which yields that $x<0$ for $\delta$ small. But $\eps>0$ small
implies that  $\delta$ is  small and we are done.

Finally, we turn to the remaining issue of embedded eigenvalues in
the essential spectrum. We will be somewhat brief, and refer the reader
to Subsection~2.1.3 of~\cite{Pe2}. Suppose that
\beeq
\label{eq:EH_def}
\Hil f = Ef \text{\ \ with\ \ } E>1.
\eneq
Then the substitution $z=\tanh(\sigma x)$ and $v(z)=f(x)$ transforms this
into the following system of differential equations with meromorphic
coefficients:
\[
\Big(-\partial_z^2+\frac{2z}{1-z^2}\partial_z + \frac{1}{\sigma^2(1-z^2)^2}\Big)
 -\frac{(\sigma+1)^2}{\sigma^2(1-z^2)}v - \frac{\sigma+1}{\sigma(1-z^2)}\sigma_1 v = \frac{E}{\sigma^2(1-z^2)^2}\sigma_3 v.
\]
We remind the reader that $\sigma>2$ is a scalar, whereas  $\sigma_1=\bm 0&1\\1&0\endm$, $\sigma_3=\bm 1&0\\0&-1\endm$ are Pauli matrices.
The singularities are exactly
$z_{\pm}=\pm1$ and $z_\infty=\infty$. In fact, these are regular singular points, see~\cite{Har} page~70.
As in~\cite{Pe2} one observes that in the vicinity of $z_{\pm}$ one can find a basis of solutions of the form
\[ (z-z_j)^{\frac{i\lambda}{2\sigma}}e_{j1}(z),\;(z-z_j)^{-\frac{i\lambda}{2\sigma}}e_{j2}(z),\; (z-z_j)^{\frac{\mu}{2\sigma}}e_{j3}(z),\]
as well as
\[
\left\{ \begin{array}{cc} (z-z_j)^{-\frac{\mu}{2\sigma}}e_{j4}(z) \text{\ \ if\ \ } \frac{\mu}{\sigma}\not\in \Z \\
\log(z-z_j)(z-z_j)^{\frac{\mu}{2\sigma}}e_{j3}(z)+(z-z_j)^{-\frac{\mu}{2\sigma}}e_{j4}(z) \text{\ \ if\ \ } \frac{\mu}{\sigma}\in \Z
\end{array}
\right.
\]
where $e_{j\ell}$, $1\le \ell\le 4$, $j=\pm$, are analytic and non-vanishing functions in some disk centered at~$z_j$, and with $E=1+\lambda^2$,
and $\mu=\sqrt{E+2}$, as before.
Since $E$ is an eigenvalue, we conclude that there would have to be a non-vanishing solution of the form
\[ (1-z^2)^{\frac{\mu}{2\sigma}}\tilde{v}\]
with an entire function $\tilde{v}$.
The remainder of the argument is identical with~\cite{Pe2}, and we skip it.
\end{proof}

We now present a simple continuity statement.

\begin{cor}
\label{cor:imag_proj} Let $\Hil(\alpha)f^{\pm}(\alpha)=\pm i\gamma
f^{\pm}(\alpha)$ where $\|f^{\pm}(\alpha)\|=1$. We can choose the
$f^{\pm}(\alpha)$ to be $\calJ$-invariant, i.e., $\calJ
f^{\pm}(\alpha)=f^{\pm}(\alpha)$. Since $\|f^{\pm}(\alpha)\|_2=1$,
they are therefore unique up to a sign. Choose this sign
consistently, i.e., so that $f^{\pm}(\alpha)$ varies continuously
with $\alpha$. In that case there is the bound \beeq
\label{eq:fal_cont}
|\gamma(\alpha_1)-\gamma(\alpha_2)|+\|f^{\pm}(\alpha_1)-f^{\pm}(\alpha_2)\|_2
\le C(\alpha_1)|\alpha_1-\alpha_2| \eneq for all
$\alpha_1,\alpha_2>0$ which are sufficiently close. Let
$\Pim^{\pm}(\alpha)$ denote the Riesz projection onto
$f^{\pm}(\alpha)$, respectively. Then one has, relative to the
operator norm on $L^2\times L^2$, \beeq \label{eq:Riesz_cont} \|
\Pim^{\pm}(\alpha_1) - \Pim^{\pm}(\alpha_2)\| \le
C(\alpha_1)|\alpha_1-\alpha_2| \eneq for all $\alpha_1,\alpha_2$
as above. Moreover, the Riesz projections admit the explicit
representation \beeq \label{eq:Riesz_rep}
 \Pim^{\pm}(\alpha) = f^{\pm}(\alpha) \la \cdot, \tilde{f}^{\pm}(\alpha)\ra,
\eneq where $\Hil(\alpha)^* \tilde{f}^{\pm}(\alpha)=\mp i\gamma
\tilde{f}^{\pm}(\alpha)$, and $\|\tilde{f}^{\pm}(\alpha)\|_2=1.$
\end{cor}
\begin{proof}
By Remark~\ref{rem:Jinv}, $\ker(\Hil(\alpha)\mp i\gamma)$ is
$\calJ$-invariant. Thus, $\calJ f^\pm(\alpha)=\lambda
f^\pm(\alpha)$ for some $\lambda\in\Compl$. It is easy to see that
this requires that $|\lambda|^2=1$. Let $e^{2i\beta}=\lambda$. It
follows that $\calJ (e^{i\beta} f^\pm(\alpha)) = e^{i\beta}
f^\pm(\alpha)$, leading to our choice of $\calJ$-invariant
eigenfunction. Using the fact that
\[ \ker[\Hil(\alpha)\mp i\gamma(\alpha)]= \ker[(\Hil(\alpha)\mp i\gamma(\alpha))^2], \]
one easily obtains (by means of the Riesz projections) that
\[
\|(\Hil(\alpha)-z)^{-1}\| \les |z\mp i\gamma(\alpha)|^{-1} \text{\
\ provided\ \ }|z\mp i\gamma(\alpha)|<r_0(\alpha).
\]
In conjunction with the resolvent identity, this yields
\[ |\gamma(\alpha_1)-\gamma(\alpha_2)| \le C(\alpha_1)|\alpha_1-\alpha_2|, \]
as well as~\eqref{eq:Riesz_cont}. However, the latter clearly
implies the remaining bound in~\eqref{eq:fal_cont}. Finally, by
the Riesz representation theorem, we necessarily have
that~\eqref{eq:Riesz_rep} holds with some choice of
$\tilde{f}^{\pm}(\alpha)\in L^2\times L^2$. Since
$\Pim^{\pm}(\alpha)^2=\Pim^{\pm}(\alpha)$, one checks that
\[ \Pim^{\pm}(\alpha)^*\tilde{f}^{\pm}(\alpha)=\tilde{f}^{\pm}(\alpha).\]
However, writing down $\Pim^{\pm}(\alpha)$ explicitly shows that
\begin{align*}
 \Pim^{+}(\alpha)^* &= \Big( \frac{1}{2\pi i}\oint_\gamma (-\Hil(\alpha)+zI)^{-1}\, dz \Big)^*
= -\frac{1}{2\pi i}\oint_\gamma (-\Hil(\alpha)^*+\bar{z}I)^{-1}\, d\bar{z} \\
&= \frac{1}{2\pi i}\oint_{-\bar{\gamma}}
(-\Hil(\alpha)^*+zI)^{-1}\, dz
\end{align*}
which is equal to the Riesz projection corresponding to the
eigenvalue $-i\gamma$ of~$\Hil(\alpha)^*$. Here $\gamma$ is a
small, positively oriented, circle around $i\gamma$. A similar
calculation applies to $\Pim^{-}(\alpha)$. Hence
$\Hil(\alpha)^*\tilde{f}^{\pm}(\alpha)=\mp i\gamma(\alpha)
\tilde{f}^{\pm}(\alpha)$, as claimed. In view
of~\eqref{eq:Riesz_rep},
\[ \|\tilde{f}^{\pm}(\alpha)\|_2^2 = \|\Pim^{+}(\alpha) \tilde{f}^{\pm}(\alpha)\|_2 \le
\|\tilde{f}^{\pm}(\alpha)\|_2.\] which implies that
$\|\tilde{f}^{\pm}(\alpha)\|_2\le1$. On the other hand,
\[ 1=\|f^{\pm}(\alpha)\|_2 = \|\Pim^{+}(\alpha) f^{\pm}(\alpha)\|_2 \le \|f^{\pm}(\alpha)\|_2
\|\tilde{f}^{\pm}(\alpha)\|_2=\|\tilde{f}^{\pm}(\alpha)\|_2,\] and
we are done.
\end{proof}

We conclude with an explicit decomposition of $L^2\times L^2$ into
the stable and unstable subspaces.

\begin{lemma}
\label{lem:L2split} Let $\calN=\ker(\Hil(\alpha)^2)$ and
$\calN^*=\ker((\Hil(\alpha)^*)^2)$ be the root spaces of $\Hil$
and $\Hil^*$, respectively, whereas $f^{\pm}(\alpha)$ and
$\tilde{f}^{\pm}(\alpha)$ are as in the previous lemma. Then there
is a direct sum decomposition
\begin{equation}
\label{eq:L2split} L^2(\R)\times L^2(\R) = \calN+{\rm
span}\{f^\pm(\alpha)\}+ \Bigl(\calN^*+  {\rm span}\{\tilde
f^\pm(\alpha)\}  \Bigr)^\perp
\end{equation}
This means that the individual summands are linearly independent,
but not necessarily orthogonal. The decomposition
\eqref{eq:L2split} is invariant under~$\Hil$. The Riesz projection
$P_s$ is precisely the projection onto the orthogonal complement
in~\eqref{eq:L2split} which is induced by the
splitting~\eqref{eq:L2split}.
\end{lemma}
\begin{proof}
This is immediate from the definition of the Riesz projections.
First,
\[
 I-P_s = \frac{1}{2\pi i} \oint_\gamma (zI-\Hil)^{-1}\, dz
\] where $\gamma$ is a simple closed curve that encloses the
entire discrete spectrum of~$\Hil$ and lies within the resolvent
set, see~\eqref{eq:rp}.  Then, on the one hand,
\[ L^2(\R)\times L^2(\R) = \ker(P_s) + \Ran(P_s) = \ker(P_s) + \ker (P_s^*)^{\perp}.\]
On the other hand,
\[ \ker(P_s) = \Ran(I-P_s) =  \calN+ {\rm span}\{f^\pm(\alpha)\} \]
as well as
\[ \ker (P_s^*) = \calN^*+  {\rm span}\{\tilde f^\pm(\alpha)\}.\]
This last equality uses that $P_s^*$ is the same as the Riesz
projection off the discrete spectrum of $\Hil^*$, as can be seen
by taking adjoints of~\eqref{eq:rp}.
\end{proof}

\begin{remark}
\label{rem:Jinv} By inspection, all root spaces in this section
are $\calJ$-invariant. This is a general fact. Indeed, one checks
easily that $J\Hil(\alpha)J=-\Hil(\alpha)$. Therefore, if
$\Hil(\alpha)f=i\sigma f$ with $\sigma\in\R$, it follows that
$\Hil(\alpha)Jf=-i\sigma J f$ where as usual $J=\bm
0&1\\1&0\endm$. Hence
\[ \calJ \ker (\Hil-i\sigma I) = \overline{J \ker (\Hil-i\sigma I)} = \ker (\Hil-i\sigma I)\]
for any $\sigma\in\R$. A similar argument shows that the root
spaces at zero are also $\calJ$-invariant. In particular, one
concludes from this that the Riesz projections $P_s, P_{\rm root},
P_{\rm im}$ preserve the space of $\calJ$-invariant functions
in~$L^2(\R)\times L^2(\R)$. This can also easily be seen directly:
Let $P$ be any Riesz projection corresponding to an eigenvalue of
$\Hil(\alpha)$ on~$i\R$, i.e.,
\[ P = \frac{1}{2\pi i}\oint_\gamma (zI-\Hil(\alpha))^{-1}\, dz \]
where $\gamma$ is a small positively oriented circle centered at
that eigenvalue. Since $J\Hil(\alpha)J=-\Hil(\alpha)$, one
concludes that
\[
JPJ = \frac{1}{2\pi i}\oint_\gamma J(zI-\Hil(\alpha))^{-1}J\, dz =
\frac{1}{2\pi i}\oint_\gamma (\Hil(\alpha)+zI)^{-1}\, dz.
\]
Thus, if $F=\binom{F_1}{F_2}$, then $-\bar{\gamma}=\gamma$ (in the
sense of oriented curves) implies that
\[
\overline{JPF} = -\frac{1}{2\pi i} \oint_\gamma
(\Hil(\alpha)+\bar{z}I)^{-1}\, d\bar{z}\; \overline{JF} =
\frac{1}{2\pi i}\oint_\gamma (zI-\Hil(\alpha))^{-1}\, dz\;
\overline{JF}= P\calJ F,
\]
so $\calJ\circ P =P\circ \calJ$, as claimed.
\end{remark}

\bibliographystyle{amsplain}

\noindent

\end{document}